\documentclass[reqno,11pt]{amsart}
\usepackage{amsthm}
\usepackage{amsmath,amsfonts,amssymb,bm}
\usepackage{xr}
\usepackage{mathrsfs}
\usepackage{mathtools}
\usepackage{hyperref}
\usepackage{mathabx}
\usepackage[initials,nobysame]{amsrefs}
\usepackage{tikz}
\usetikzlibrary{intersections,calc}
\newtheorem{theorem}{Theorem}[section]
\newtheorem{lemma}[theorem]{Lemma}
\newtheorem{proposition}[theorem]{Proposition}

\newtheorem{remark}[theorem]{Remark}
\numberwithin{figure}{section}

\def\u{\mathrm{u}}
\def\eps{\varepsilon}
\def\vv{\mathrm{v}}
\def \del{{\partial}}
\def \GRAD{\nabla_{\!\!x}}

\def \DIV{\nabla_{\!\!x} \! \cdot }

\def \DOT{{\,\cdot\,}}
\def \eps{\varepsilon}

\AtBeginDocument{
   \def\MR#1{}
}

\usepackage[top=2.0cm,bottom=2.0cm,left=3cm,right=3cm]{geometry}

\title{Low mach Number Limit of the Viscous and Heat Conductive Flow with general pressure law on torus}
\numberwithin{equation}{section}
\mathtoolsset{showonlyrefs}

\begin{document}

\author[Yuhan Chen]{Yuhan Chen}
\address[Yuhan Chen]
{\newline School of Mathematics and Statistics, Wuhan University, Wuhan 430072, P. R. China}
\email{yhchen\_math@whu.edu.cn}

\author[Guilong Gui]{Guilong Gui}
\address[Guilong Gui]
{\newline School of Mathematics and Computational Science, Xiangtan University, Xiangtan 411105, P. R. China}
\email{glgui@amss.ac.cn}

\author[Zhen Hao]{Zhen Hao}
\address[Zhen Hao]
{\newline School of Mathematics and Statistics, Wuhan University, Wuhan 430072, P. R. China}
\email{zhhao\_math@whu.edu.cn}

\author[Ning Jiang]{Ning Jiang}
\address[Ning Jiang]{\newline School of Mathematics and Statistics, Wuhan University, Wuhan 430072, P. R. China}
\email{njiang@whu.edu.cn}

\begin{abstract}
   We prove the low Mach number limit from compressible Navier-Stokes-Fourier system with the general pressure law around a constant state on the torus $\mathbb{T}^N_a$. We view this limit as a special case of the weakly nonlinear-dissipative approximation of the general hyperbolic-parabolic system with entropy. In particular, we consider the ill-prepared initial data, for which the group of fast acoustic waves is needed to be filtered. This extends the previous works, in particular Danchin [{\em Amer. J. Math.} {\bf 124} (2002), 1153-1219]  in two ways: 1. We treat the fully general non-isentropic flow, i.e. the pressure depends on the density $\rho$ and temperature $\theta$ by basic thermodynamic law. We illustrate the role played by the entropy structure of the system in the coupling of the acoustic waves and incompressible flow, and the construction of the filtering group. 2. We refine the small divisor estimate, which helps us to give the first explicit convergence rate of the filtered acoustic waves whose propogation is governed by non-local averaged system. In previous works, only convergence rate of incompressible limit was obtained.\\


\end{abstract}

\maketitle
\thispagestyle{empty}
\pagestyle{plain}
\section{Introduction}
\subsection{Compressible Navier-Stokes-Fourier system}
In this paper, we address a classic problem: the low Mach number limit from the compressible fluid model. Our starting point is the following Navier-Stokes-Fourier (NSF) system of compressible viscous and heat-conductive fluid:

\begin{equation}
    \label{NSCe}
    \left\{
       \begin{array}{l}
         \partial_{t}\rho + \nabla_x\!\cdot\!(\rho \u) = 0,\\
         \partial_{t}(\rho \u) + \nabla_{x}\!\cdot\!(\rho \u\otimes \u + pI) = \nabla_{x}\!\cdot\! \mathbb{S},\\
         \partial_{t}(\frac{1}{2}\rho| \u|^2 + \rho e) + \nabla_{x}\!\cdot\!(\frac{1}{2}\rho| \u|^2\u + \rho e \u + p \u) = \nabla_{x}\!\cdot\!(\mathbb{S}\!\cdot\! \u - q),
       \end{array}
    \right.
\end{equation}
where $\rho=\rho(t,x)>0$, $ \u=\u(t,x)\in\mathbb{R}^N$, $e=e(t,x)>0$ and $p=p(t,x)$
stand for the mass density, bulk velocity, specific internal energy and pressure respectively.
$N\geqslant2$ is the spatial dimension
and $I$ denotes the $N\times N$ identity matrix.
We are concerned with Newtonian and heat-conductive fluids
where the deviatoric strain tensor $\mathbb{S}$ and heat flux $q$ are given by the constitutive relations
\begin{equation}
    \mathbb{S}\equiv\mathbb{S}(\u, \mu, \lambda) = \mu \left[\nabla_{x} \u + (\nabla_{x} \u)^{T} - \tfrac{2}{N}(\nabla_{x}\!\cdot\! \u)I\right] + \lambda(\nabla_{x}\!\cdot\! \u)I\,,\quad\!\mbox{and}\quad\! q\equiv q(\kappa, \theta) = -\kappa\nabla_{x}\theta,
\end{equation}
where $\theta(t,x)>0$ is the temperature.
The thermodynamic independent variables are chosen to be $(\rho, \theta)$,
and the dependence of $p$ and $e$ on $(\rho, \theta)$ is given by some thermodynamic equations of state.
The shear viscosity coefficient $\mu=\mu(\rho, \theta)>0$, the bulk viscosity $\lambda=\lambda(\rho, \theta)\geqslant 0$
and the thermo-conductivity coefficient $\kappa = \kappa(\rho, \theta)>0$ come either from non-equilibrium kinetic theory or from fitting to experimental data. For the general introduction to the equations for compressible fluids, see standard books, for example \cite{Lions_book}.

The Mach number, the ratio of the characteristic velocity in the flow to the sound speed in the fluid, is an important dimensionless number in a fluid system. It measures the compressibility of the fluid. Physically, this naturally means when the Mach number is small, the system describing the fluid system should be approximated by a system of incompressible flow. The rigorous justification of this process, i.e. the Mach number approaching to zero, is a classical mathematical analytical problem in fluid dynamics, the so-called {\em low Mach number limit of a compressible system}.

In this work, we seek a mathematical justification of such an approximation
in the case of periodic boundary conditions.
Namely, the spatial domain has chosen to be the torus $\mathbb{T}_{a}^{N}$,
where $a=(a_1,...,a_N)\in\mathbb{R}_{+}^N$ denotes the aspect ratio,
meaning that for $x\in\mathbb{T}_{a}^{N}$, the $i$-th component $x_i$ varies with period $2\pi a_i$.

After suitable process of nondimensionization (see for instance \cite{Alazard, BLP-1992, Lions_book}),
the compressible Navier-Stokes system \eqref{NSCe} can be written in the following non-dimensional way:
\begin{equation}\label{CNS-scaled}
 \left\{
    \begin{aligned}
        &\partial_t \rho^{\varepsilon} + \DIV(\rho^{\varepsilon}\u^{\varepsilon}) = 0,\\
        &\partial_t (\rho^{\varepsilon}\u^{\varepsilon}) + \DIV(\rho^{\varepsilon} \u^{\varepsilon}\otimes \u^{\varepsilon})
        + \frac{1}{\varepsilon^2}\GRAD p^{\varepsilon} = \DIV \mathbb{S}^{\varepsilon},\\
        &\partial_t (\rho^{\varepsilon}e^{\varepsilon}) + \DIV(\rho^{\varepsilon}\u^{\varepsilon}e^{\varepsilon}) + p^{\varepsilon}\DIV \u^{\varepsilon}
        = \DIV(\kappa^{\varepsilon}\GRAD \theta^{\varepsilon}) + \varepsilon^2 \mathbb{S}^{\varepsilon}:\GRAD\u^{\varepsilon}\,,
    \end{aligned} \right.
\end{equation}
where $ \mathbb{S}^{\varepsilon} = \mathbb{S}(\u^\varepsilon, \mu^\varepsilon, \lambda^\varepsilon)$, and $\varepsilon$ denotes the Mach number.
$\rho^{\varepsilon}$, $\u^{\varepsilon}$, $e^{\varepsilon}$, $\theta^{\varepsilon}$,
$p^{\varepsilon}$, $\mu^{\varepsilon}$, $\lambda^{\varepsilon}$ and $\kappa^{\varepsilon}$
are the dimensionless version of the original quantities they represent.

Physically, when the Mach number $\varepsilon$ approaches to 0, the solutions of \eqref{CNS-scaled} will be asymptotically governed by an incompressible model. The compressible Navier-Stokes system includes compressible Euler, which evolutes in short time scale (Euler scaling), while the incompressible Navier-Stokes evolutes in longer time scale (Navier-Stokes scaling). In this sense, the limit $\varepsilon$ approaching to 0 in \eqref{CNS-scaled} is a singular limit. In other words, even formally, it can not expected that $(\rho^\varepsilon, \u^\varepsilon, \theta^\varepsilon)$ directly converges to some solutions of incompressible model. Among the several ways (in different scalings) that the system \eqref{CNS-scaled} asymptotically behaves like a system of incompressible model, the simplest way is fluctuating around a constant state with the size of order $O(\varepsilon)$. In this scaling, the limiting system is the incompressible Navier-Stokes system with Boussinesq relation. For the formal derivations in different regimes (Boussinesq, or dominant regimes, etc.), see Bayly-Levermore-Passot's work in early nineties \cite{BLP-1992}.  However, for general initial data, even the fluctuations around the constant state does not directly converge because of the multi-scale feature of compressible Navier-Stokes system \eqref{CNS-scaled}. More precisely, there is a linear operator (we call the acoustic operator $\mathcal{A}$) which relates both the compressible Euler and incompressible Navier-Stokes regimes. The projection on the null space of $\mathcal{A}$ converges to incompressible regime, while the orthogonal part (under some suitable inner product which determined by the entropy structure of compressible NSF) keep in the fast time scale regime. This is the common feature of all low Mach number limits. The goal of this paper is to precisely characterize both formally and rigorously this process.

\subsection{Weakly nonlinear-dissipative approximations}
Now we view the low Mach number limit of compressible NSF from a more general setting. Set $\u^\varepsilon=\frac{1}{\varepsilon}{\vv}^\varepsilon$, the rescaled compressible NSF system \eqref{CNS-scaled} has the form:
\begin{equation}\label{NSF-1}
    \left\{
       \begin{aligned}
         &\partial_{t}\rho + \tfrac{1}{\varepsilon}\nabla_x\!\cdot\!(\rho {\vv}) = 0,\\
         &\partial_{t}(\rho {\vv}) + \tfrac{1}{\varepsilon}\nabla_{x}\!\cdot\!(\rho {\vv}\otimes {\vv} + pI) = \nabla_{x}\!\cdot\! \mathbb{S},\\
         &\partial_{t}(\tfrac{1}{2}\rho| {\vv}|^2 + \rho e) + \tfrac{1}{\varepsilon}\nabla_{x}\!\cdot\!(\tfrac{1}{2}\rho| {\vv}|^2{\vv} + \rho e {\vv} + p {\vv}) = \nabla_{x}\!\cdot\!(\mathbb{S}\!\cdot\! {\vv} +\kappa\nabla_{x}\theta).
       \end{aligned}\right.
  \end{equation}

We remark that the difference between the systems \eqref{CNS-scaled} and \eqref{NSF-1} are just velocity scalings. In fact, when deriving the compressible NSF (for ideal gas) from the Boltzmann equation in the Navier-Stokes scaling \cite{BGL1991}, the system \eqref{NSF-1} will be obtained. At this scaling the so-called Knudsen number and Mach numbers will be the same order. We now explain why the format \eqref{NSF-1} is a natural way to view low Mach number limit. In fact, we can consider the general hyperbolic-parabolic (HP) system with entropy (compressible NSF \eqref{CNS-scaled} as an example):
\begin{equation}\label{HP}
  \del_t U + \DIV F(U) = \DIV \big[ D(U) \DOT \GRAD U \big] \,.
\end{equation}
 In \eqref{HP}, $U$ is a vector of densities, $F(U)$ is a twice continuously differentiable flux, and $D(U)$ is a continuously differentiable diffusion tensor. Moreover, they are assumed to possess a strictly convex, thrice
continuously differentiable, real-valued entropy density $H(U)$ such that classical solutions of \eqref{HP} also
satisfy
\begin{equation}
  \partial_t H(U)+\DIV J(U)=\DIV[H_U(U)D(U)\!\cdot\!\GRAD U]-\GRAD U\!\cdot\!H_{UU}(U)D(U)\!\cdot\!\GRAD U\,,
\end{equation}
here the entropy flux $J(U)$ satisfies $J_U(U) = H_U(U)F_U(U)$ while the tensor $H_{UU}(U)D(U)$ is symmetric and nonnegative definite. Here $H_U(U), J_U(U)$ and $F_U(U)$ denote the derivatives of $H(U), J(U)$
and $F(U)$ with respect to $U$ while $H_{UU}(U)$ denotes the Hessian of $H(U)$ with respect to $U$.

In \cite{Jiang_ARMA}, the fourth named author and Levermore studied the weakly nonlinear-dissipative (WND) approximations of \eqref{HP}. The motivation was the following: in the dimensionless form, when the time scale is $O(1)$ and the size of the diffusion term is of the order $O(\eps)$, the scaled form of \eqref{HP} reads
\begin{equation}\label{HP-short}
  \del_t U^\eps + \DIV F(U^\eps) = \eps\DIV \big[ D(U^\eps) \DOT \GRAD U^\eps \big] \,.
\end{equation}
Formally, when $\eps\rightarrow 0$, the solution of \eqref{HP-short} $U^\eps \rightarrow U$, where $U$ is the solution of the hyperbolic system:
\begin{equation}\label{Hyperbolic}
  \del_t U + \DIV F(U) = 0 \,.
\end{equation}
In the context of compressible NSF, this limit is the so-called inviscid limit. The analytic justification of this limit is easy in the context of classical solution local in time in the domains of whole space and torus. However, this inviscid limit is very hard when the boundary of the domain is nontrivial, and the solutions includes shock, contact discontinuity, etc. This is a classic hard problem in hyperbolic conservation laws, which is widely open. However, if we consider the HP in much longer time scale, say, $O(\eps)$, the rescaled form of \eqref{HP} will be
\begin{equation}\label{HP-1}
  \del_t U^\varepsilon + \tfrac{1}{\varepsilon}\DIV F(U^\varepsilon) = \DIV \big[ D(U^\varepsilon) \DOT \GRAD U^\varepsilon \big] \,.
\end{equation}
The scaled compressible NSF \eqref{NSF-1} is one of its special form.

In \cite{Jiang_ARMA}, the weakly nonlinear-dissipative approximations of \eqref{HP-1} was studied, which govern the regime in which $U$ is close to a constant
state $U^o$ and the dissipation is small. If we express the densities U in terms of any choice of dependent
variables $V$ as $U = U(V)$, so that $U^o = U(V^o)$, we then define the matrix $R^o$ by $R^o = U_V(V^o)$. The
weakly nonlinear-dissipative approximation of \eqref{HP-1} governs the perturbation $\tilde{V}$ of $V$ about $V^o$ by the
system
\begin{equation}\label{Averaged-HP}
  \partial_t \tilde{V}+\mathcal{A}\tilde{V}+\overline{\mathcal{Q}}(\tilde{V},\tilde{V})=\overline{\mathcal{D}}\tilde{V}\,,
\end{equation}
where the linear operator $\mathcal{A}$ is formally defined by
\begin{equation}
  \mathcal{A}\tilde{V}= (R^o)^{-1}F_U(U^o)R^o\!\cdot\!\GRAD \tilde{V}\,,
\end{equation}
while the averaged operators $\overline{\mathcal{Q}}$ and $\overline{\mathcal{D}}$ are formally defined by
\begin{equation}
  \begin{aligned}
    &\overline{\mathcal{Q}}(\tilde{V},\tilde{V})=\lim_{T\rightarrow\infty}\frac{1}{2T}\int^T_{-T}e^{s\mathcal{A}}
    \mathcal{Q}(e^{-s\mathcal{A}\tilde{V}},e^{-s\mathcal{A}\tilde{V}})\,\mathrm{d} s\,,\\
    &\overline{\mathcal{D}}\tilde{V}=\lim_{T\rightarrow\infty}\frac{1}{2T}\int^T_{-T}e^{s\mathcal{A}}\mathcal{D}e^{-s\mathcal{A}}\tilde{V}\,\mathrm{d}s\,,
  \end{aligned}
\end{equation}
with the operators $\mathcal{Q}$ and $\mathcal{D}$ given by
\begin{equation}
  \begin{aligned}
    &\mathcal{Q}(\tilde{V},\tilde{V})= \DIV[\tfrac{1}{2}(R^o)^{-1}F_{UU}(U^o)(R^o\tilde{V},R^o\tilde{V})]\,,\\
    &\mathcal{D}\tilde{V}=\DIV[(R^o)^{-1}D(U^o)R^o\!\cdot\!\GRAD \tilde{V}]\,.
   \end{aligned}
\end{equation}
The first two terms in \eqref{Averaged-HP} are the linearization of \eqref{HP-1} with respect to $V$ neglecting the dissipation. The
operator $\mathcal{A}$ thereby governs the fast dynamics. Nonlinearity and dissipation will modify the dynamics on
longer time scales because the perturbation $V$ is assumed to be small while the dissipation is assumed to
be weak. The operators $\mathcal{Q}$ and $\mathcal{D}$ are averages of $\mathcal{Q}$ and $\mathcal{D}$ over the fast dynamics that attempt to capture
these modifications. For detailed derivations of the averaged system \eqref{Averaged-HP}, see Section 3.2 of \cite{Jiang_ARMA}.

It was shown in Section 4 of \cite{Jiang_ARMA} that the entropy structure \eqref{HP-1} implies that the operator $\mathcal{A}$ is skew-adjoint in the Hilbert space
\begin{equation}
  \mathbb{H}=\big\{\tilde{V}\in L^2(\mathrm{d}\,;\mathbb{R}^{d+2}: \int_{\mathbb{T}^d}\tilde{V}\,\mathrm{d}x=0)\big\}\,,
\end{equation}
equipped with the natural inner product
\begin{equation}
  (\tilde{V}_1, \tilde{V}_2)_{\mathbb{H}}= \int_{\mathbb{T}^d}(R^o\tilde{V}_1)^T H_{UU}(U^o)R^o\tilde{V}_2\,\mathrm{d}x\,.
\end{equation}

Furthermore, it was shown in \cite{Jiang_AA}, applying the process of the weakly nonlinear-dissipative approximations developed in \cite{Jiang_ARMA} to the compressible Navier-Stokes-Fourier system with general pressure law, the projection of the averaged system \eqref{Averaged-HP} on the null space of $\mathcal{A}$ (which is the acoustic operator in compressible NSF) Null$(\mathcal{A})$ is exactly incompressible Navier-Stokes-Fourier system, while the project on the orthogonal space Null$(\mathcal{A})^\perp$ is a nonlocal parabolic system with 2 and 3-wave resonant terms. In this sense, the classical low Mach number limit could be understood as a special case of the WND approximation of the general hyperbolic-parabolic system with entropy. The main theme of this paper is to rigorously justify this process for compressible NSF.

\subsection{Low Mach number limit of general NSF}
In the past four decades,
the mathematical justification of the low Mach number limit
has been studied in various aspects and contexts.
The focus of interest often centers around
the absence or presence of entropy variations in the fluid
(isentropic or non-isentropic),
the types of the system (viscous or inviscid, in other words, Navier-Stokes or Euler),
the contexts of solutions being considered (weak or strong),
how general the initial data is
(“well-prepared" or “ill-prepared")
and the boundaries of the spatial domain
(whole space, periodic domain, fixed or even free boundary conditions).

Among the plenty of works, we only mention here some most representative ones. For more complete references, see the book and survey papers, such as \cite{FeireislBook2017} and \cite{Alazard2008, Schochet2007}.
More precisely, the rigorous study of the low Mach number limit
was initiated around the early eighties by Ebin \cite{Ebin},
Klainerman-Majda \cites{KM81,KM82} for local strong solutions of
compressible fluids (Euler or Navier-Stokes),
in the whole space with well-prepared data
(div$\u^{\varepsilon}_0=\mathcal{O}(\varepsilon),\ \nabla P^{\varepsilon}_{0}=\mathcal{O}(\varepsilon^2)$).
These works are later improved by Ukai \cite{Ukai} in whole space
for ill-prepared data
(div$\u^{\varepsilon}_0=\mathcal{O}(1),\ \nabla P^{\varepsilon}_{0}=\mathcal{O}(\varepsilon)$),
where acoustic waves of amplitude $\mathcal{O}(1)$
with frequency of order $\mathcal{O}(\varepsilon^{-1})$ appears.
In the context of isentropic system, the study has been fruitful
since the above pioneering works,
there have been many progress, such as Schochet \cite{Schochet_JDE}, Gallagher\cite{Gallagher1998}, Danchin \cites{Danchin_AJM, Danchin_ASENS},
Lions and Masmoudi \cites{LM_JMPA,LM_Paris}, Masmoudi \cite{Masmoudi_Poincare},
Grenier \cite{Grenier_JMPA}, Desjardins and Grenier \cite{DG},
Desjardins, Grenier, Lions and Masmoudi \cite{DGLM}, Jiang-Masmoudi \cite{Jiang-Masmoudi2015}, Hoff \cite{Hoff_CMP}, and others. One of the crucial step in these works is the relation of the fast acoustic waves and boundaries: torus, bounded domain and whole spaces. Different domains correspond to different behaviors of the acoustic waves. See the survey paper \cite{Jiang-Masmoudi2018}. In the endeavors regarding non-isentropic models,
early works include, for instance, those by Métivier-Schochet \cites{MS_ARMA,MS_JDE},
where they proved strong convergence for classical solutions of Euler equations
in the whole space or on the torus.
Their result was later extended to the boundary case by Alazard \cite{Alazard_boundary}.
Very recently, Masmoudi, Rousset and Sun \cite{Masmoudi_2022}
derived uniform higher-order regularity estimates
in bounded domains with ill-prepared initial data,
they also investigated the case of free boundary conditions in \cite{Masmoudi_2023}.

The main concerns and the novelties of the current paper are the following. First, we work the compressible NSF with the most general pressure law. In other words, we only assume the basic physical requirements: the pressure is increasing with respect to the density and the internal energy is increasing with respect to the temperature. In this sense, we do not put any technical assumption on the pressure. Second, we explicitly illustrate the role played by the entropy in this low Mach number limit. In previous works, the isentropic models are too simple to see entropy structure. In the works for non-isentropic models, most of the them treated the well-prepared initial data so that the fast waves are not generated. Third, in the previous works such as Gallagher\cite{Gallagher1998}, Danchin \cites{Danchin_AJM}, besides they worked on the isentropic cases, in the filtering process, they only obtained the convergence rate for incompressible part. In the current paper, employing some refined small divisor estimates, we analyzed the structures of 2 and 3$\mbox{-}$waves resonant sets. This allow obtain the explicit convergence rate of fast oscillating part for the first time. This is the main technical novelty of this paper.

 \subsection{Entropy structure of the system}
The entropy structure of the compressible NSF plays key roles in this work in two ways.
First, the Hessian of the entropy function is used to define a natural Hilbert space
in which the global existence and regularity of the limit system is proved.
Second, it illustrates the structure of the averaged system as projected onto
 {\em oscillating} and  {\em incompressible} modes.

No we discuss the equations of the thermodynamic state. If we select $\rho$, $\theta$ as independent thermodynamic variables,
the equations for the specific energy $e=e(\rho,\theta)$ and pressure $p=p(\rho,\theta)$
are assumed to be twice continuously differentiable over $(\rho,\theta)\in\mathbb{R}^{2}_{+}$ and to satisfy
\begin{equation}\label{es1}
    e_{\theta}(\rho,\theta)>0,\ \ p_{\rho}(\rho,\theta)>0
\end{equation}
for $\rho>0$ and $\theta>0$.
They are also assumed to satisfy the Maxwell relation
\begin{equation}\label{es2}
    \rho^2 e_{\rho} + \theta^2\partial_{\theta}(\frac{\rho}{\theta}) = 0,
\end{equation}
which implies the existence of a function $s=s(\rho,\theta)$ satisfying
\begin{equation}
    \mathrm{d}s = \frac{1}{\theta}(\mathrm{d}e + p\mathrm{d}(\frac{1}{\rho})).
\end{equation}
Here $s$ is the specific entropy of the system.
We deduce
\begin{equation}
    s_{\theta}=\frac{e_{\theta}}{\theta},\ \ s_{\rho}=\frac{1}{\theta}(e_{\rho}-\frac{p}{\rho^2}),
\end{equation}
which implies the compatibility equation
\begin{equation}
    \frac{\partial}{\partial \rho}(\frac{e_{\theta}}{\theta})
    = \frac{\partial}{\partial\theta}(\frac{1}{\theta}(e_{\rho}-\frac{p}{\rho^2}))\,.
\end{equation}
In other words,
\begin{equation}\label{eos}
    e_{\rho} = \frac{1}{\rho^2}(p-\theta p_{\theta}).
\end{equation}

We define $H(U) = -\rho s(\rho,\theta)$, where
\begin{equation}
    U \stackrel{\mathrm{def}}{=} (\rho, \rho \u_1,\dots,\rho \u_{N},\frac{1}{2}\rho|\u|^2 + \rho e(\rho,\theta))^{T}.
\end{equation}
Computations in the same manner as \cite{Jiang_AA} show that
(CNS) is a nonsingular hyperbolic-parabolic system with
a strictly convex entropy given by $H(U)$.
More precisely, we point out that
\begin{equation}\label{es3}
    R^{T}H_{UU}(U)R =
    \begin{pmatrix}
        \frac{p_{\rho}}{\rho\theta}&0&0\\
        0&\frac{\rho}{\theta} I&0\\
        0&0&\frac{\rho e_{\theta}}{\theta^2}
    \end{pmatrix},
\end{equation}
where $R$ is defined in \eqref{def_of_Rmatrix}
It is clear that the positive definiteness of $H_{UU}(U)$ is a consequence of \eqref{es1}.
Other conditions that make $H(U)$ a strictly convex entropy for the system
can also be verified using \eqref{es1} and \eqref{es2}.
Moreover, we can also show as in \cite{Jiang_ARMA} that
\begin{equation}
    \begin{aligned}
        \partial_t H(U) &+ \nabla \!\cdot\! (\frac{q}{\theta}-\rho s \u)\\
        =&-\frac{\mu}{2}|\nabla  \u + (\nabla  \u)^{T} - \frac{2}{N}(\nabla \!\cdot\! \u)I|^2
        -\lambda |\nabla \!\cdot\! \u
        -\kappa |\nabla \theta|^2.
    \end{aligned}
\end{equation}
Because $\mu>0$, $\lambda\geqslant 0$ and $\kappa>0$,
the right-hand side of the above is non-positive.
The compressible Naiver-Stokes system is therefore a hyperbolic-parabolic
system with a strictly convex entropy, which renders the system dissipative.

Next, we introduce the vector variable $V\stackrel{\mathrm{def}}{=}(\rho,\u,\theta)^{T}$, and the map
\begin{equation}\label{def_of_Rmatrix}
    R(V) \stackrel{\mathrm{def}}{=} U_{V}(V) =
    \begin{pmatrix}
        1&0&0\\
        \u&\rho I&0\\
        \frac{1}{2}|\u|^2 + e + \rho e_{\rho}& \rho \u& \rho e_{\theta}
    \end{pmatrix}.
\end{equation}
Define a Hilbert space
\begin{equation}
    \mathbb{H} = \left\{\widetilde{V}\in L^2(\mathbb{T}^{N}_a;\mathbb{C}^{N+2}):\int_{\mathbb{T}^{N}_a}\widetilde{V}\mathrm{d}x=0\right\}
\end{equation}
equipped with the inner product defined by
\begin{equation}
    \langle\widetilde{V}_1,\widetilde{V}_2\rangle_{\mathbb{H}} \stackrel{\mathrm{def}}{=} \int_{\mathbb{T}^{N}_a}(\widetilde{V}_1)^{T}(R^{\circ})^{T}H_{UU}(U^{\circ})R^{\circ}\bar{\widetilde{V}}_2\mathrm{d}x.
\end{equation}
Throughout the paper, the notation $\mathscr{U}^{\circ}$ will always be used
to denote a function $\mathscr{U}=\mathscr{U}(\rho,\theta)$ taking value at the point
$(\rho^{\circ},\theta^{\circ})$, that is $\mathscr{U}^{\circ}=\mathscr{U}(\rho^{\circ},\theta^{\circ})$.
In particular, the inner product coefficient matrix is then
\begin{equation}
    (R^{\circ})^{T}H_{UU}(U^{\circ})R^{\circ} =
    \begin{pmatrix}
        \frac{p^{\circ}_{\rho}}{\rho^{\circ}\theta^{\circ}}&0&0\\
        0&\frac{\rho^{\circ}}{\theta^{\circ}} I&0\\
        0&0&\frac{\rho^{\circ} e^{\circ}_{\theta}}{(\theta^{\circ})^2}
    \end{pmatrix}\,.
\end{equation}

\subsection{Acoustic operator $\mathcal{A}$}
Let $\mathcal{A}$ be the acoustic operator,
which is the linearization of the compressible Navier-Stokes system
neglecting dissipation about the constant state $(\rho^{\circ},0,\theta^{\circ})$
\begin{equation}\label{acoustic_operator}
    \mathcal{A}
    \left(
        \begin{array}{c}
            \tilde{\rho}\\\tilde{\u}\\\tilde{\theta}
        \end{array}
    \right)
    =\left(\begin{array}{c}
    \rho^\circ \nabla\!\cdot\!\tilde{\u} \\
    \frac{{p_\rho}^\circ}{\rho^\circ} \nabla \tilde{\rho}+\frac{p_{\theta}^\circ}{\rho^\circ} \nabla \tilde{\theta} \\
    \frac{\theta^\circ p_{\theta}^{\circ}}{\rho^\circ e_{\theta}^{\circ}} \nabla\!\cdot\! \tilde{\u}
    \end{array}\right)\,.
\end{equation}
One easily checks that $\mathcal{A}$ is a densely-defined unbounded skew-adjoint operator in the Hilbert space $\mathbb{H}$:
for $\widetilde{V}_1$, $\widetilde{V}_2$ in the domain of $\mathcal{A}$,
\begin{equation}
    \langle\widetilde{V}_1,\mathcal{A}\widetilde{V}_2\rangle_{\mathbb{H}} = -\langle\mathcal{A}\widetilde{V}_1,\widetilde{V}_2\rangle_{\mathbb{H}}.
\end{equation}
Also, from definition we see that the range and null space of $\mathcal{A}$ are given by
\begin{equation}
    \begin{gathered}
        \text{Range}(\mathcal{A}) = \left\{\begin{pmatrix}
            \rho^{\circ}\beta\\\nabla \phi\\\frac{\theta^{\circ}p^{\circ}_{\theta}}{\rho^{\circ}e^{\circ}_{\theta}}\beta
        \end{pmatrix}
        :\beta\in L^2_0(\mathbb{T}^{N}_a),\ \phi\in H^1(\mathbb{T}^{N}_a)\right\},\\
        \text{Null}(\mathcal{A}) = \left\{\begin{pmatrix}
            \gamma\\\omega\\-\frac{p^{\circ}_{\rho}}{p^{\circ}_{\theta}}\gamma
        \end{pmatrix}
        :\gamma\in L^2(\mathbb{T}^{N}_a),\ \omega\in L^2(\mathbb{T}^{N}_a;\mathbb{R}^{N}),\ \nabla \!\cdot\!\omega=0\right\},
    \end{gathered}
\end{equation}
where $L^2_0(\mathbb{T}^{N}_a)$ denotes the space of $L^2$ functions with mean zero.
Since $\mathcal{A}$ is skew-adjoint, Hilbert space theory tells us that $\text{Range}(\mathcal{A}) = \text{Null}(\mathcal{A})^{\perp}$,
the orthogonal complement of $\text{Null}(\mathcal{A})$. The Hilbert space $\mathbb{H}$ thus admits an orthogonal decomposition
\begin{equation}
    \mathbb{H} = \text{Null}(\mathcal{A}) \oplus \text{Null}(\mathcal{A})^{\perp}.
\end{equation}
We denote $\bar{f} = \frac{1}{|\mathbb{T}^N|}\int_{\mathbb{T}^N_a} f dx$ to be the mean of $f$ and $\underline{f}= f- \bar{f}$ to be the mean free part of $f$.
Moreover, we denote $\tilde{\mathcal{P}}$ and $\mathcal{P}^{\perp}$ the projections onto Null$(\mathcal{A})$ and $\text{Null}(\mathcal{A})^{\perp}$.
Algebraic computation shows that
\begin{equation}
    \tilde{\mathcal{P}}
    \begin{pmatrix}
        {\rho}\\
        {\u}\\
        {\theta}
    \end{pmatrix}
    =
    \begin{pmatrix}
        \bar{\rho}-\frac{p^{\circ}_{\theta}}{p^{\circ}_{\rho}}\left\{\frac{p^{\circ}_{\rho}}{(c^{\circ})^2}  \underline{\theta} -
        \frac{\theta^{\circ}p^{\circ}_{\theta}p^{\circ}_{\rho}}{(\rho^{\circ})^2e^{\circ}_{\theta}(c^{\circ})^2} \underline{\rho}\right\}\\
        \bar{\u}+ \Pi {\underline{\u}}\\
        \bar{\theta}+\frac{p^{\circ}_{\rho}}{(c^{\circ})^2} \underline{\theta} -
        \frac{\theta^{\circ}p^{\circ}_{\theta}p^{\circ}_{\rho}}{(\rho^{\circ})^2e^{\circ}_{\theta}(c^{\circ})^2} \underline{\rho}
    \end{pmatrix}
    ,\ \
    \mathcal{P}^{\perp}
    \begin{pmatrix}
        {\rho}\\
        {\u}\\
        {\theta}
    \end{pmatrix}
    =
    \begin{pmatrix}
        \frac{p^{\circ}_{\rho}}{(c^{\circ})^2}{\underline{\rho}} + \frac{p^{\circ}_{\theta}}{(c^{\circ})^2}{\underline{\theta}}\\
        (I-\Pi){\underline{\u}}\\
        \frac{\theta^{\circ}p^{\circ}_{\theta}}{(\rho^{\circ})^2e^{\circ}_{\theta}}
        \left\{\frac{p^{\circ}_{\rho}}{(c^{\circ})^2}{\underline{\rho}} + \frac{p^{\circ}_{\theta}}{(c^{\circ})^2}{\underline{\theta}}\right\}
    \end{pmatrix},
\end{equation}
where $c^{\circ}$ is the sound speed defined by
\begin{equation}
    c^{\circ} = \sqrt{p^{\circ}_{\rho} + \frac{\theta^{\circ}(p^{\circ}_{\theta})^2}{(\rho^{\circ})^2e^{\circ}_{\theta}}},
\end{equation}
and $\Pi$ is the Leray projection onto the space of divergence-free vector fields, which could be expressed as
\begin{equation}
    \Pi = I - \nabla \Delta^{-1}\nabla \!\cdot.
\end{equation}
For later convenience, we define the customary variables $\vartheta$, $\omega$, $\pi$ and $v$  associated to $({\rho},{\u},{\theta})^{T}$ as
\begin{equation}\label{projected variables}
    \begin{aligned}
        &\vartheta \stackrel{\mathrm{def}}{=} \frac{p^{\circ}_{\rho}}{(c^{\circ})^2}  \underline{\theta} -
        \frac{\theta^{\circ}p^{\circ}_{\theta}p^{\circ}_{\rho}}{(\rho^{\circ})^2e^{\circ}_{\theta}(c^{\circ})^2}\underline{\rho},
        \ \ \omega \stackrel{\mathrm{def}}{=} \Pi \underline{\u},\\
        &\pi \stackrel{\mathrm{def}}{=} p^{\circ}_{\rho} \underline{\rho} + p^{\circ}_{\theta} \underline{\theta},\ \ v \stackrel{\mathrm{def}}{=} (I-\Pi)\underline{\u}.
    \end{aligned}
\end{equation}
We thereby can write alternatively
\begin{equation}
    \tilde{\mathcal{P}}
    \begin{pmatrix}
        {\rho} \\
        {\u}\\
        {\theta}
    \end{pmatrix}
    =
    \begin{pmatrix}
        -\frac{p^{\circ}_{\theta}}{p^{\circ}_{\rho}}\vartheta + \bar{\rho}\\
        \omega + \bar{\u}\\
        \vartheta + \bar{\theta}
    \end{pmatrix}
    ,\ \
    \mathcal{P}^{\perp}
    \begin{pmatrix}
        {\rho}\\
        {\u}\\
        {\theta}
    \end{pmatrix}
    =
    \begin{pmatrix}
        \frac{1}{(c^{\circ})^2}\pi\\
        {v}\\
        \frac{\theta^{\circ}p^{\circ}_{\theta}}{(c^{\circ})^2(\rho^{\circ})^2e^{\circ}_{\theta}}\pi
    \end{pmatrix},
\end{equation}
and we further decompose $\tilde{\mathcal{P}}$ as $\tilde{\mathcal{P}} = \mathcal{P} + \mathcal{P}_0$ with
\begin{equation}
    \mathcal{P}
    \begin{pmatrix}
        {\rho} \\
        {\u}\\
        {\theta}
    \end{pmatrix}
    =
    \begin{pmatrix}
        -\frac{p^{\circ}_{\theta}}{p^{\circ}_{\rho}}\vartheta\\
        \omega \\
        \vartheta
    \end{pmatrix},
    \ \ \mathcal{P}_0 =
    \begin{pmatrix}
        \bar{\rho}\\
        \bar{\u} \\
        \bar{\theta}
    \end{pmatrix}.
\end{equation}

Next, we discuss the eigenvalue problem of $\mathcal{A}$.
Let $\phi_{\nu}$ be an eigenfunction of the Laplacian
\begin{equation}\label{lap_eigenvalue}
    -\Delta_{x}\phi_{\nu} = \nu^2\phi_{\nu}
\end{equation}
over $\mathbb{T}^{N}_a$ for some $\nu>0$, then a direct calculation shows
\begin{equation}
    \mathcal{A}\begin{pmatrix}
        \pm \mathrm{i}\rho^{\circ}\phi_{\nu}\\
        \frac{c^{\circ}}{\nu}\nabla \phi_{\nu}\\
        \pm \mathrm{i}\frac{\theta^{\circ}p^{\circ}_{\theta}}{\rho^{\circ}e^{\circ}_{\theta}}\phi_{\nu}
    \end{pmatrix}
    =\begin{pmatrix}
        -\rho^{\circ}c^{\circ}\nu\phi_{\nu}\\
        \pm \mathrm{i}(c^{\circ})^2\nabla \phi_{\nu}\\
        -\frac{\theta^{\circ}p^{\circ}_{\theta}}{\rho^{\circ}e^{\circ}_{\theta}}c^{\circ}\nu\phi_{\nu}
    \end{pmatrix}
    =\pm \mathrm{i}c^{\circ}\nu
    \begin{pmatrix}
        \pm \mathrm{i}\rho^{\circ}\phi_{\nu}\\
        \frac{c^{\circ}}{\nu}\nabla \phi_{\nu}\\
        \pm \mathrm{i}\frac{\theta^{\circ}p^{\circ}_{\theta}}{\rho^{\circ}e^{\circ}_{\theta}}\phi_{\nu}
    \end{pmatrix}.
\end{equation}
Thus, we see if $\nu^2$ is a positive eigenvalue of $-\Delta_{x}$ with eigenfunction $\phi_{\nu}$,
then $\pm \mathrm{i}c^{\circ}\nu$ are a conjugate pair of eigenvalues of $\mathcal{A}$ with eigenfunction
$\begin{pmatrix}
    \pm \mathrm{i}\rho^{\circ}\phi_{\nu}\\
    \frac{c^{\circ}}{\nu}\nabla \phi_{\nu}\\
    \pm \mathrm{i}\frac{\theta^{\circ}p^{\circ}_{\theta}}{\rho^{\circ}e^{\circ}_{\theta}}\phi_{\nu}
\end{pmatrix}$.

Recall that the eigenvalue problem for the Laplace equation is solved by Fourier modes.
For a spatial variable $x$ on the torus $\mathbb{T}^{N}_{a}$, the Fourier modes on this torus are given by
\begin{equation}
    \mathrm{e}^{\mathrm{i} k\cdot x},\ \ \text{with}\ \
    k\in\widetilde{\mathbb{Z}}^{N}_{a}\stackrel{\mathrm{def}}{=} \mathbb{Z}/a_1\times ... \times \mathbb{Z}/a_N
\end{equation}
In particular, the Fourier mode $\mathrm{e}^{\mathrm{i} k\cdot x}$ solves \eqref{lap_eigenvalue} with $\nu=| k|$.
Therefore, we can construct an orthonormal eigenfunction basis for $\mathcal{A}$:
the normalized eigenfunction corresponding to the eigenvalue
$\lambda_{k}^{\alpha}\stackrel{\mathrm{def}}{=}\mathrm{i}c^{\circ}\alpha \operatorname{sg}( k)| k|$ reads
\begin{equation}
    H^{\alpha}_{ k} =
    c_N
    \begin{pmatrix}
        \rho^{\circ}\\
        c^{\circ}\alpha \operatorname{sg}( k)\frac{ k}{| k|}\\
        \frac{\theta^{\circ}p^{\circ}_{\theta}}{\rho^{\circ}e^{\circ}_{\theta}}
    \end{pmatrix}
    \mathrm{e}^{\mathrm{i} k\cdot x},
\end{equation}
where $c_N$ is defined by
\begin{equation}\label{def_of_CN}
    c_N\stackrel{\mathrm{def}}{=}\frac{1}{c^{\circ}|\mathbb{T}_{a}^{N}|^{\frac{1}{2}}}
    \sqrt{\frac{\theta^{\circ}}{2\rho^{\circ}}},
\end{equation}
and $\alpha$ takes either one of $1$ or $-1$. The notation $\operatorname{sg}( k)$ stands for a generalized sign function on $\mathbb{R}^{N}\setminus\left\{0\right\}$: its value
is $1$ if and only if the first nonzero component of $ k$ is positive, $-1$ elsewhere.

The discussion above leads us to the following decomposition for any ${V}=({\rho},{\u},{\theta})^{T}\in\mathbb{H}$:
\begin{equation}
    {V} = \mathcal{P}{V} + \mathcal{P}^{\perp}{V} + \mathcal{P}_0 V=
    \begin{pmatrix}
        -\frac{p^{\circ}_{\theta}}{p^{\circ}_{\rho}}\vartheta\\
        \omega\\
        \vartheta
    \end{pmatrix}
    +
    \begin{pmatrix}
        \frac{1}{(c^{\circ})^2}\pi\\
        {\vv}\\
        \frac{\theta^{\circ}p^{\circ}_{\theta}}{(c^{\circ})^2(\rho^{\circ})^2e^{\circ}_{\theta}}\pi
    \end{pmatrix}
    +\begin{pmatrix}
        \bar{\rho}\\
        \bar{\u}\\
        \bar{\theta}
    \end{pmatrix},
\end{equation}
and the kernel orthogonal part $\mathcal{P}^{\perp}V$ can be represented by
\begin{equation}
    \mathcal{P}^{\perp}\widetilde{V} = \sum_{(\alpha, k)}V^{\alpha}_{ k} H^{\alpha}_{ k}.
\end{equation}
The coefficient in the series is defined as
\begin{equation}
    V^{\alpha}_{ k} \stackrel{\mathrm{def}}{=} \langle\widetilde{V},H^{\alpha}_{ k}\rangle_{\mathbb{H}}.
\end{equation}

\subsection{Low Mach number setting}
After suitable rescaling (see for instance \cite{Lions_book} or \cite{Alazard}),
the compressible Navier-Stokes equations,
written in a non-dimensional way \eqref{CNS-scaled}.
When $\rho^{\varepsilon}$ is bounded away from zero,
using the equations of state and \eqref{eos} to replace the equation involving $e^{\varepsilon}$ with that of $\theta^{\varepsilon}$,
the system rewrites
\begin{equation}\label{nsc_lowmach_setting}
    \left\{
    \begin{aligned}
        &\partial_t \rho^{\varepsilon} + \nabla \!\cdot\! (\rho^{\varepsilon}\u^{\varepsilon}) = 0,\\
        &\partial_t \u^{\varepsilon} + \u^{\varepsilon}\!\cdot\!\nabla \u^{\varepsilon} + \frac{1}{\varepsilon^2}\frac{1}{\rho^{\varepsilon}}\nabla p^{\varepsilon} =
        \frac{1}{\rho^{\varepsilon}}\nabla \!\cdot\!S^{\varepsilon},\\
        &\partial_t \theta^{\varepsilon} + \u^{\varepsilon}\!\cdot\!\nabla \theta^{\varepsilon} + \frac{\theta^{\varepsilon}p^{\varepsilon}_{\theta}}{\rho^{\varepsilon}e^{\varepsilon}_{\theta}}\nabla \!\cdot\! \u^{\varepsilon}
        = \frac{1}{\rho^{\varepsilon}e^{\varepsilon}_{\theta}}\nabla \!\cdot\!(\kappa^{\varepsilon}\nabla \theta^{\varepsilon}) + \frac{\varepsilon^2}{\rho^{\varepsilon}e^{\varepsilon}_{\theta}}S^{\varepsilon}:\nabla \u^{\varepsilon}.
    \end{aligned}
    \right.
\end{equation}

When the parameter $\varepsilon$ becomes small,
the asymptotic limit of \eqref{nsc_lowmach_setting} 
leads to incompressible model.
As expounded in \cite{Lions_book}, the singularity in the momentum equation
can be balanced if we assume, at least formally,
\begin{equation}
    \rho^{\varepsilon} = \rho^{\circ} + \varepsilon\tilde{\rho} + \mathcal{O}(\varepsilon^2),\ \
    \theta^{\varepsilon} = \theta^{\circ} + \varepsilon\tilde{\theta} + \mathcal{O}(\varepsilon^2),\ \
    \u^{\varepsilon} = \tilde{\u} + \mathcal{O}(\varepsilon),
\end{equation}
where $\rho^{\circ}$ and $\theta^{\circ}$ are constants.
The leading-order terms that determine the fluctuations are $\tilde{\rho}$, $\tilde{\u}$, $\tilde{\theta}$, they satisfy the Boussinesq relations
\begin{equation}
    \operatorname{div}\tilde{\u}=0,\ \ p_{\rho}(\rho^{\circ},\theta^{\circ})\nabla \tilde{\rho}+p_{\theta}(\rho^{\circ},\theta^{\circ})\nabla \tilde{\theta} = 0,
\end{equation}
and $\tilde{\u}$ is the solution to the incompressible Navier-Stokes equation
\begin{equation}
    \left\{\begin{aligned}
        &\partial_t \tilde{\u} + \mathcal{P}(\tilde{\u}\!\cdot\!\nabla \tilde{\u})-\mu\Delta \tilde{\u}=0,\ \ \operatorname{div}\tilde{\u}=0,\\
        &\tilde{\u}|_{t=0} = \tilde{\u}_{in}.
    \end{aligned}\right.
\end{equation}
One might compare the Boussinesq relation with the definition of $\text{Null}(\mathcal{A})$
to get a glimpse of how $(\tilde{\rho}^{\varepsilon},\tilde{\u}^{\varepsilon},\tilde{\theta}^{\varepsilon})^{T}$
will eventually fall into $\text{Null}(\mathcal{A})$.

In the present work, we assume that
\begin{equation}
    \rho^\varepsilon_{in} = \rho^\circ + \varepsilon \tilde{\rho}_{in},\ \ \
    \u^\varepsilon_{in} = \tilde{\u}_{in},\ \ \
    \theta^\varepsilon_{in} = \theta^\circ + \varepsilon \tilde{\theta}_{in},\ \ \
\end{equation}
and set
\begin{equation}\label{variable_fluctuations}
    \rho^\varepsilon = \rho^\circ + \varepsilon \tilde{\rho}^\varepsilon,\ \ \
    \u^\varepsilon = \tilde{\u}^\varepsilon,\ \ \
    \theta^\varepsilon = \theta^\circ + \varepsilon \tilde{\theta}^\varepsilon.\ \ \
\end{equation}
Denoting
\begin{equation}
    \widetilde{U}^{\varepsilon} = (\tilde{\rho}^{\varepsilon}, \tilde{\u}^{\varepsilon}, \tilde{\theta}^{\varepsilon})^{T}
\end{equation}
and plugging \eqref{variable_fluctuations} into \eqref{nsc_lowmach_setting},
we get the equations satisfied by the fluctuation $\widetilde{U}^{\varepsilon}$
\begin{equation}\label{nsc_epsilon}\tag{NSC$^{\hspace{0.1em}\varepsilon}$}
    \left\{
    \begin{aligned}
    \partial_t \widetilde{U}^\varepsilon + \frac{1}{\varepsilon} &\mathcal{A} \widetilde{U}^\varepsilon
    + \mathcal{Q}{(\widetilde{U}^\varepsilon,\widetilde{U}^\varepsilon)}
    = \mathcal{D}(\widetilde{U}^\varepsilon)+r^{\varepsilon},\\
    &\widetilde{U}^\varepsilon|_{t=0}=\widetilde{U}_{in} = (\tilde{\rho}_{in}, \tilde{\u}_{in},\tilde{\theta}_{in})^T,
    \end{aligned} \right.
\end{equation}
tic operator defined by \eqref{acoustic_operator}.
The quadratic term $\mathcal{Q}(\widetilde{U},\widetilde{U})$
and the dissipation term $\mathcal{D}\widetilde{U}$ are defined as follows:
\begin{equation} \notag
    \mathcal{Q}(\widetilde{U}, \widetilde{U})=\left(\begin{array}{c}
    \nabla \!\cdot\!(\tilde{\rho} \tilde{\u}) \\
    \tilde{\u} \!\cdot\! \nabla \tilde{\u}+ C_1 \tilde{\rho}\nabla \tilde{\rho} + C_2 \tilde{\theta}\nabla \tilde{\theta} + C_3 \tilde{\theta}\nabla \tilde{\rho} + C_4 \tilde{\rho}\nabla \tilde{\theta}\\
    \tilde{\u} \!\cdot\! \nabla \tilde{\theta} +\left(C_5\tilde{\theta} + C_6 \tilde{\rho}\right) \nabla \!\cdot\! \tilde{\u}
    \end{array}\right),
\end{equation}

\begin{equation} \notag
    \mathcal{D} \widetilde{U}=\left(\begin{array}{c}
    0 \\
    \frac{\mu^\circ}{\rho^\circ} \Delta \tilde{\u}+\frac{\frac{N-2}{N} \mu^\circ+\lambda^\circ}{\rho^\circ} \nabla \nabla \!\cdot\! \tilde{\u} \\
    \frac{\kappa^{\circ}}{\rho^{\circ}e_{\theta}^{\circ}}\Delta\tilde{\theta}
    \end{array}\right).
\end{equation}
Moreover, we define
\begin{equation}
    \mathcal{Q}(\widetilde{U}_1,\widetilde{U}_2) = \frac{1}{2} \left(\mathcal{Q}(\widetilde{U}_1+\widetilde{U}_2,\widetilde{U}_1+\widetilde{U}_2)-\mathcal{Q}(\widetilde{U}_1,\widetilde{U}_1)-\mathcal{Q}(\widetilde{U}_2,\widetilde{U}_2)\right).
\end{equation}
The remainder term is defined by
\begin{equation*}
    r^{\varepsilon} = \left(0,\ r^{\varepsilon}_2,\ r^{\varepsilon}_3\right)^{\text{T}},
\end{equation*}
with
\begin{equation}
    \begin{aligned}
        r^{\varepsilon}_2 &=
            -\left(
            {\frac{1}{\varepsilon^{2}}} \frac{\nabla p}{\rho^{\varepsilon}}
            - \frac{1}{\varepsilon}\frac{p^{\circ}_{\theta}}{\rho^{\circ}}\nabla \tilde{\theta}^{\varepsilon}
            -\frac{1}{\varepsilon}\frac{p^{\circ}_{\rho}}{\rho^{\circ}}\nabla \tilde{\rho}^{\varepsilon}
            -C_1 \tilde{\rho}^{\varepsilon} \nabla \tilde{\rho}^{\varepsilon}
            -C_2 \tilde{\theta}^{\varepsilon} \nabla \tilde{\theta}^{\varepsilon}
            -C_3 \tilde{\theta}^{\varepsilon} \nabla \tilde{\rho}^{\varepsilon}
            -C_4 \tilde{\rho}^{\varepsilon} \nabla \tilde{\theta}^{\varepsilon}
            \right)\\
            &+ \frac{1}{{\rho}^{\varepsilon}}\nabla \!\cdot\! S^{\varepsilon} - \left( \frac{\mu^{\circ}}{\rho^{\circ}} \Delta \tilde{\u}^{\varepsilon} + \frac{\frac{N-2}{N}\mu^{\circ}+\lambda^{\circ}}{\rho^{\circ}}\nabla \operatorname{div} \tilde{\u}^{\varepsilon} \right),\\
        r^{\varepsilon}_3 &=  -\frac{\kappa^{\circ}}{e^{\circ}_{\theta} \rho^{\circ}} \Delta  \tilde{\theta}^{\varepsilon} + \frac{1}{e_{\theta} \rho} \nabla( \kappa \nabla  \tilde{\theta}^{\varepsilon})
        + \frac{1}{\varepsilon} \left( \frac{\theta^{\circ} p^{\circ}_{\theta}}{e^{\circ}_{\theta} \rho^{\circ}} - \frac{\theta P_{\theta}}{e_{\theta} \rho} + C_5 \tilde{\theta}^{\varepsilon} + C_6 \tilde{\rho}^{\varepsilon} \right) \nabla \!\cdot\! \tilde{\u}^{\varepsilon}
        + \frac{\varepsilon}{e_{\theta} \rho} S^{\varepsilon}\!:\! \nabla  \tilde{\u}^{\varepsilon}.\\
    \end{aligned}
\end{equation}
The constants are calculated to be
\begin{align}
    &C_1 = \frac{p_{\rho\rho}^{\circ}}{\rho^{\circ}} - \frac{p_{\rho}^{\circ}}{(\rho^{\circ})^2} ,\ \ \  C_2 = \frac{p_{\theta\theta}^{\circ}}{\rho^{\circ}} ,\ \ \  C_3 = \frac{p_{\rho\theta}^{\circ}}{\rho^{\circ}} ,\ \ \  C_4 = \frac{p_{\rho\theta}^{\circ}}{\rho^{\circ}} - \frac{p_{\theta}^{\circ}}{(\rho^{\circ})^2}\\
    &C_5 = \frac{p_{\theta}^{\circ}}{\rho^{\circ}e_{\theta}^{\circ}} + \frac{\theta^{\circ}p_{\theta\theta}^{\circ}}{\rho^{\circ}e_{\theta}^{\circ}} - \frac{\theta^{\circ}p_{\theta}^{\circ}e_{\theta\theta}^{\circ}}{\rho^{\circ}(e_{\theta}^{\circ})^2} ,\ \ \
    C_6 = \frac{\theta^{\circ}p_{\theta\rho}^{\circ}}{\rho^{\circ}e_{\theta}^{\circ}} - \frac{\theta^{\circ}p_{\theta}^{\circ}}{(\rho^{\circ})^2e_{\theta}^{\circ}} - \frac{\theta^{\circ}p_{\theta}^{\circ}e_{\theta\rho}^{\circ}}{\rho^{\circ}(e_{\theta}^{\circ})^2}.
\end{align}
\begin{remark}
In the ensuing analysis, one might run into terms like $\mathcal{Q}(\widetilde{U}_1, \widetilde{U}_2)$
for non-identical $\widetilde{U}_1$ and $\widetilde{U}_2$.
The convention shall be made that the definition for $\mathcal{Q}$ is extended by polarization,
as in the case for general bilinear operators.
\end{remark}
\begin{remark}
For mathematical simplicity, we suppose that the data $(\tilde{\rho}_{in}, \tilde{\u}_{in},\tilde{\theta}_{in})$
are independent of $\varepsilon$.
Otherwise, we only need to require that $\mathcal{P}\widetilde{U}^{\varepsilon}_{in}$ and $\mathcal{P}^{\perp}\widetilde{U}^{\varepsilon}_{in}$ converge in appropriate spaces.
\end{remark}

According to the entropy structure discussed in the last section,
the system may be “decomposed" orthogonally into three parts.
First, apply $\mathcal{P}$ to both sides of \eqref{nsc_epsilon}.
\begin{equation}\label{nsc_kernel}
    \begin{aligned}
        \partial_t \mathcal{P}\widetilde{U}^{\varepsilon}
        &+ \mathcal{P}\mathcal{Q}(\mathcal{P}\widetilde{U}^{\varepsilon},\mathcal{P}\widetilde{U}^{\varepsilon})
        - \mathcal{P}\mathcal{D}\mathcal{P}\widetilde{U}^{\varepsilon}\\
        &=\mathcal{P}r^{\varepsilon} + \mathcal{P}\mathcal{D}\mathcal{P}^{\perp}\widetilde{U}^{\varepsilon}
        - \mathcal{P} \big( 2\mathcal{Q}(\mathcal{P}\widetilde{U}^{\varepsilon},\mathcal{P}^{\perp}\widetilde{U}^{\varepsilon}) + \mathcal{Q}(\mathcal{P}^{\perp}\widetilde{U}^{\varepsilon},\mathcal{P}^{\perp}\widetilde{U}^{\varepsilon})\big)\\
        &- 2\mathcal{P} \big( Q(\mathcal{P}_0 \widetilde{U}^\varepsilon , \mathcal{P}^{\perp} \widetilde{U}^{\varepsilon } ) +  Q(\mathcal{P}_0 \widetilde{U}^\varepsilon , \mathcal{P} \widetilde{U}^{\varepsilon } ) \big).
    \end{aligned}
\end{equation}

Clearly, one expects $\mathcal{P}\widetilde{U}^{\varepsilon}$
to converge to some $\mathcal{U}\in\text{Null}(\mathcal{A})$
and $\mathcal{U}$ satisfies
\begin{equation}\label{INSF}\tag{INSF}
    \left\{\begin{aligned}
        &\partial_t \mathcal{U} + \mathcal{P}\mathcal{Q}(\mathcal{U},\mathcal{U}) - \mathcal{P}\mathcal{D}\mathcal{U} = 0,\\
        &\mathcal{U}|_{t=0} = \mathcal{P}\widetilde{U}_{in}.
    \end{aligned}\right.
\end{equation}
Actually, $\mathcal{U}$ is of the form $\mathcal{U}= (-\frac{p_{\theta}^{\circ}}{p_{\rho}^{\circ}}\vartheta, \omega, \vartheta)^{T}$, and
\eqref{INSF} is nothing but the incompressible Navier-Stokes-Fourier system
\begin{equation}\nonumber
    \begin{aligned}
        \partial_t
        \begin{pmatrix}
            \omega\\
            \vartheta
        \end{pmatrix}
        +
        \begin{pmatrix}
            \Pi(\omega\!\cdot\!\nabla \omega)\\
            \omega\!\cdot\!\nabla \vartheta
        \end{pmatrix}
        =
        \begin{pmatrix}
            \frac{\mu^{\circ}}{\rho^{\circ}}\Delta\omega\\
            \frac{\kappa^{\circ}p_{\rho}^{\circ}}{(c^{\circ})^2\rho^{\circ}e_{\theta}^{\circ}}\Delta\vartheta
        \end{pmatrix}.
    \end{aligned}
\end{equation}
However, rigorous justification of such convergence will not be easy.
Even though, seemingly, the singularity appearing in the acoustic part of \eqref{nsc_epsilon} 
has been removed, terms involving
$\mathcal{P}^\perp\widetilde{U}^{\varepsilon}$
will experience high oscillations and thus
are not likely to tend to zero in any strong norm when the spatial domain is periodic.
This is in contrast to the “well-prepared" data case,
where this phenomenon does not occur, and also to the whole space case,
where oscillating waves are expected to disperse to zero,
due to Strichartz type estimates.
Thoroughly studied by various authors such as Schochet \cite{Schochet_JDE}, Gallagher \cite{Isabelle_JMKU},
Masmoudi \cite{Masmoudi_Poincare} and Danchin \cite{Danchin_AJM} among many others,
a strategy known as “filtering method"
has proved effective
in dealing with highly oscillating waves
on periodic domains.
The idea is that, we shall filter $\mathcal{P}^{\perp}\widetilde{U}^{\varepsilon}$
according to the group $e^{\frac{t}{\varepsilon}\mathcal{A}}$
generated by the acoustic operator $\mathcal{A}$.
The filtered oscillations have to be described and included into
energy estimates.
As a consequence, it will be shown that they do not affect the limit equation.
Denote $\widetilde{V}^{\varepsilon} \stackrel{\mathrm{def}}{=} \mathrm{e}^{\frac{t}{\varepsilon}\mathcal{A}}\mathcal{P}^{\perp}\widetilde{U}^{\varepsilon}$.
It satisfies
\begin{equation}\label{nsc_orthogonal}
    \begin{aligned}
        \partial_t \widetilde{V}^{\varepsilon}
        &+ \mathcal{Q}^{\varepsilon}_{2r}(\mathcal{P}\widetilde{U}^{\varepsilon},\widetilde{V}^{\varepsilon})
        + \mathcal{Q}^{\varepsilon}_{3r}(\widetilde{V}^{\varepsilon},\widetilde{V}^{\varepsilon})
        - \mathcal{D}^{\varepsilon}\widetilde{V}^{\varepsilon}\\
        &= \mathrm{e}^{\frac{t}{\varepsilon}\mathcal{A}}\mathcal{P}^{\perp}\mathcal{D}\mathcal{P}\widetilde{U}^{\varepsilon}
        + \mathrm{e}^{\frac{t}{\varepsilon}\mathcal{A}}\mathcal{P}^{\perp}r^{\varepsilon}
        - \mathrm{e}^{\frac{t}{\varepsilon}\mathcal{A}}\mathcal{P}^{\perp}\mathcal{Q}(\mathcal{P}\widetilde{U}^{\varepsilon},\mathcal{P}\widetilde{U}^{\varepsilon})\\
        &- 2\mathrm{e}^{\frac{t}{\varepsilon}\mathcal{A}}\mathcal{P}^{\perp} \big( Q(\mathcal{P}_0 \widetilde{U}^\varepsilon , \mathcal{P}^{\perp} \widetilde{U}^{\varepsilon } ) +  Q(\mathcal{P}_0 \widetilde{U}^\varepsilon , \mathcal{P} \widetilde{U}^{\varepsilon } ) \big),
    \end{aligned}
\end{equation}
where
\begin{equation}
    \begin{aligned}
        &\mathcal{Q}^{\varepsilon}_{2r}(\mathcal{P}\widetilde{U}^{\varepsilon},\widetilde{V}^{\varepsilon})  = 2\mathrm{e}^{\frac{t}{\varepsilon}\mathcal{A}}\mathcal{P}^{\perp}\mathcal{Q}(\mathcal{P}\widetilde{U}^{\varepsilon},\mathrm{e}^{-\frac{t}{\varepsilon}\mathcal{A}}\widetilde{V}^{\varepsilon}),\\
        &\mathcal{Q}^{\varepsilon}_{3r}(\widetilde{V}^{\varepsilon},\widetilde{V}^{\varepsilon}) = \mathrm{e}^{\frac{t}{\varepsilon}\mathcal{A}}\mathcal{P}^{\perp}\mathcal{Q}(\mathrm{e}^{-\frac{t}{\varepsilon}\mathcal{A}}\widetilde{V}^{\varepsilon},\mathrm{e}^{-\frac{t}{\varepsilon}\mathcal{A}}\widetilde{V}^{\varepsilon}),\\
        &\mathcal{D}^{\varepsilon}\widetilde{V}^{\varepsilon} = \mathrm{e}^{\frac{t}{\varepsilon}\mathcal{A}}\mathcal{P}^{\perp}\mathcal{D}\mathrm{e}^{-\frac{t}{\varepsilon}\mathcal{A}}\widetilde{V}^{\varepsilon}.
    \end{aligned}
\end{equation}
Formal analysis as in \cite{Jiang_AA} or \cite{Danchin_AJM} hints that, after an averaging process
(letting $\varepsilon\rightarrow\ 0$), we shall expect
\begin{equation}
    \widetilde{V}^{\varepsilon}\rightarrow\text{ some }V\in\text{Null}(\mathcal{A})^{\perp},
\end{equation}
where $V$ satisfies
\begin{equation}\label{LS}\tag{LS}
    \left\{\begin{aligned}
        &\partial_t V + \mathcal{Q}_{2r}(\mathcal{U},V) + \mathcal{Q}_{3r}(V,V)-\overline{\mathcal{D}}V = 0, 
        \\ &V|_{t=0} = \mathcal{P}^{\perp}\widetilde{U}_{in},
    \end{aligned}\right.
\end{equation}
where $\mathcal{U}$ stands for the limit already obtained in the kernel part
of the system.
The limiting operators $\mathcal{Q}_{2r}(\mathcal{U},V)$, $\mathcal{Q}_{3r}(V,V)$ and $\overline{\mathcal{D}}V$ may be computed
by means of nonstationary phase arguments in section \ref{sec_oscillating}.
\subsection{Statement of main result}
In our consideration, the evolution of the fluid
takes place in the spatial domain
$\mathbb{T}_{a}^{N}$ with $N \geqslant 2$,
where we recall $a=(a_1,...,a_N)\in\mathbb{R}_{+}^N$
is the aspect ratio of the periodic box.
Throughout the paper, the value of $a$ is fixed,
and it should be noted that the choice of $a$
is not all of $\mathbb{R}_{+}^N$, but with exceptions
that belong to a Lebesgue measure zero
set\footnote{This is due to small divisor problems
which will be made clear later}.

The functional spaces are chosen to be
the Besov spaces $B^{s}(\mathbb{T}_{a}^{N})$
(see section 3 for definition and a compendium on relevant
harmonic analysis results).
For notation simplicity,
$B^{s}$ will always stand for $B^{s}(\mathbb{T}_{a}^{N})$.
As noticed by \cite{Danchin_Inventiones}, the regularity of the compressible Navier-Stokes
equation is expected to behave differently
in “low" and “high" Fourier frequencies,
it will therefore be convenient to separate growth conditions
for different frequencies and introduce some
“hybrid Besov spaces" $B_{\eta}^{s,t}$,
where the subscript $\eta$ denotes the “anchor"
according to which we distinguish “low" and “high"
frequencies (exact definition is given in section 3).

For any Banach space $X$, $0<T\leqslant +\infty$
and $1\leqslant r\leqslant +\infty$,
we denote by $L^{r}(0,T;X)$ the set of $r$-integrable
functions on $(0,T)$ with values in $X$.
In the case $T=+\infty$, the space can be simply denoted as $L^{r}(X)$.
The space $C([0,T];X)$ (resp. $C_{b}(X)$) denotes the set
of continuous bounded functions defined on $[0,T]$ (resp. $\mathbb{R}^{+}$)
with values in $X$. We also adopt the convention that,
for a vector field $U$, “$U\in X$" means
each component of $U$ belongs to the space $X$.

We introduce the spaces
\begin{equation}
    E^{s}_{\eta} \stackrel{\mathrm{def}}{=}
        \widetilde{C}(\mathbb{R}^{+};\mathcal{E}^{s}_{\eta})\cap L^1(\mathbb{R}^{+};\dot{\mathcal{D}}^{s}_{\eta}),
\end{equation}
equipped with the norm $||\cdot||_{\widetilde{L}^{\infty}(\mathbb{R}^{+}; \mathcal{E}^{s}_{\eta}) \cap L^1(\mathbb{R}^{+};\dot{\mathcal{D}}^{s}_{\eta}) }$,
where
\begin{equation}
    \begin{aligned}
        \mathcal{E}^{s}_{\eta} \stackrel{\mathrm{def}}{=}&
        B^{s-1,s+1}_{\eta}\times B^{s-1,s}_{\eta}\times B^{s-1,s}_{\eta},\
        \mathcal{D}^{s}_{\eta} \stackrel{\mathrm{def}}{=}
        B^{s+1}\times B^{s+1,s+2}_{\eta}\times B^{s+1,s+2}_{\eta},\\
        \dot{\mathcal{E}}^{s}_{\eta} \stackrel{\mathrm{def}}{=}&
        \dot{B}^{s-1,s+1}_{\eta}\times \dot{B}^{s-1,s}_{\eta}\times \dot{B}^{s-1,s}_{\eta},\
        \dot{\mathcal{D}^{s}_{\eta}} \stackrel{\mathrm{def}}{=}
        \dot{B}^{s+1}\times \dot{B}^{s+1,s+2}_{\eta}\times \dot{B}^{s+1,s+2}_{\eta},
    \end{aligned}
\end{equation}
and
\begin{equation}
    \begin{aligned}
        F^{s}\stackrel{\mathrm{def}}{=}& \widetilde{C}(\mathbb{R}^{+};B^{s-1}) \cap L^{1}(\mathbb{R}^{+};B^{s+1}),\
        G^{s}\stackrel{\mathrm{def}}{=} \widetilde{C}(\mathbb{R}^{+};H^{s-1})\cap {L}^{2}(\mathbb{R}^{+};H^{s}).
    \end{aligned}
\end{equation}
In the above definitions, the  spaces $\widetilde{C}$ and $\widetilde{L}^2$
have the meaning that the time integration is
“localized", a specified version of which can also be found in section 3.
Equipped with these spaces,
the existence of a solution
$\widetilde{U}^{\varepsilon}=(\tilde{\rho}^{\varepsilon},
\tilde{\u}^{\varepsilon},\tilde{\theta}^{\varepsilon})$
of \eqref{nsc_epsilon} will be sought in
$E^{\frac{N}{2}}_{\varepsilon\nu^{\circ}}$,
the space for \eqref{INSF} and \eqref{LS}
is chosen to be $F^{\frac{N}{2}}$,
and convergence will be proved in a space
$G^{\frac{N}{2}-\delta}$ with a loss of regularity
$\delta>0$.

We now state our main result:
\begin{theorem}\label{MAIN}
    Given a constant state $U^{\circ} = (\rho^{\circ},0,\theta^{\circ})$, where $\rho^{\circ}$ and $\theta^{\circ}$ are two positive constants such that
    \begin{equation}
        p_{\rho}(\rho^{\circ},\theta^{\circ})>0,\ \ e_{\theta}(\rho^{\circ},\theta^{\circ})>0.
    \end{equation}
    There exist two positive constants $\alpha_0$, $C$ depending only on $U^{\circ},\  N$ and the dependence of $p,e,\mu,\lambda,\kappa$ on $\rho$ and $\theta$, 
    such that if $\widetilde{U}_{in} = (\tilde{\rho}_{in},\tilde{\u}_{in},\tilde{\theta}_{in})\in \dot{ \mathcal{E}}_{\varepsilon_0\nu^{\circ}}^{\frac{N}{2}}$,
    where $\varepsilon_{0}>0$ and
    \begin{equation}
        \|\widetilde{U}_{in}\|
        _{\mathcal{E}_{\varepsilon_0\nu^{\circ}}^{\frac{N}{2}}}
        \leqslant \alpha_0,
    \end{equation}
    then for any $0<\varepsilon<\varepsilon_0$,
    the system \eqref{nsc_epsilon} has a unique solution $U^{\varepsilon}$ such that $\widetilde{U}^{\varepsilon} \stackrel{\mathrm{def}}{=} \frac{U^{\varepsilon} - U^{\circ}}{\varepsilon}$
    is uniformly bounded in $E^{\frac{N}{2}}_{\varepsilon\nu^{\circ}}$.
    System \eqref{INSF}, \eqref{LS} has a solution $(\mathcal{U},V)\in F^{\frac{N}{2}}\times F^{\frac{N}{2}}$. They satisfy
    \begin{equation}\label{energy_ineq}
        \begin{aligned}
            \|\widetilde{U}^{\varepsilon}\|_{E^{\frac{N}{2}}_{\varepsilon\nu^{\circ}}} + \|(\mathcal{U},V)\|_{F^{\frac{N}{2}}}
            \leqslant C\|\widetilde{U}_{in}\|
            _{\mathcal{E}_{\varepsilon\nu^{\circ}}^{\frac{N}{2}}}.
        \end{aligned}
    \end{equation}
    Moreover, there exists a set $A \subset \mathbb{R}_{+}^{N}$ with Lebesgue measure zero,
    for each torus $\mathbb{T}^{N}_a$ with $a \in  A^c$, for  $\delta \in (0,1]$, $\tau \in (0,+\infty)$ and for $\sigma =\max \left\{ 2N +2\tau -1 , 5\right\}$, we have
    \begin{equation}
        \begin{aligned}
            &\|\mathcal{P}\widetilde{U}^{\varepsilon}-\mathcal{U}\|_{G^{\frac{N}{2}-\delta}} = \mathcal{O}(\varepsilon^{\delta}),\\
            &\|\mathcal{P}^{\perp}\widetilde{U}^{\varepsilon}-\mathrm{e}^{-\frac{t}{\varepsilon}\mathcal{A}}V\|_{G^{\frac{N}{2}-\delta}} =  \mathcal{O}(\varepsilon^{\frac{\delta}{2+ \sigma}}),\\
            &\| \mathcal{P}_0 \widetilde{U}^{\varepsilon} \|_{L^{\infty}(dt)} = \mathcal{O} (\varepsilon),
        \end{aligned}
    \end{equation}
    as $\varepsilon\rightarrow 0$.

\end{theorem}
\begin{remark}
    In Theorem \ref{MAIN} we require that the initial data satisfies $\widetilde{U}_{in} = (\tilde{\rho}_{in},\tilde{u}_{in},\tilde{\theta}_{in})\in \dot{ \mathcal{E}}_{\varepsilon_0\nu^{\circ}}^{\frac{N}{2}}$, which means we require that $\mathcal{P}\widetilde{U}_{in} = 0$. Similar results can also be obtained
    if we only assume $\|\widetilde{U}_{in}\|_{\mathcal{E}_{\varepsilon_0\nu^{\circ}}^{\frac{N}{2}}} \leqslant \alpha_0$. More specifically, we have
    \begin{equation}
        \begin{aligned}
            &\|\mathcal{P}\widetilde{U}^{\varepsilon}-\mathcal{U}\|_{G^{\frac{N}{2}-\delta}} = \mathcal{O}(\varepsilon^{\delta}),\\
            &\|\mathcal{P}^{\perp}\widetilde{U}^{\varepsilon}-\mathrm{e}^{-\frac{t}{\varepsilon}\mathcal{A}}V\|_{G^{\frac{N}{2}-\delta}} =  \mathcal{O}(\varepsilon^{\frac{\delta}{2+ \sigma}}),\\
            &\| \mathcal{P}_0 \widetilde{U}^{\varepsilon} - \mathcal{P}_0 \widetilde{U}_{\text{in}} \|_{L^{\infty}(dt)} = \mathcal{O} (\varepsilon),
        \end{aligned}
    \end{equation}
    as $\varepsilon \rightarrow 0$. The difference lies in that instead of \eqref{INSF} and \eqref{LS}, $\mathcal{U}$ now is the solution of \eqref{INSF2} and $V$ now is the solution of \eqref{LS2}.
    \begin{equation}\label{INSF2}
        \left\{\begin{aligned}
            &\partial_t \mathcal{U} + \mathcal{P}\mathcal{Q}(\mathcal{U}+2\mathcal{P}_0 \widetilde{U}_{\text{in}},\mathcal{U}) - \mathcal{P}\mathcal{D}\mathcal{U} = 0,\\
            &\mathcal{U}|_{t=0} = \mathcal{P}\widetilde{U}_{in}.
        \end{aligned}\right.
    \end{equation}
    \begin{equation}\label{LS2}
        \left\{\begin{aligned}
            &\partial_t V + \mathcal{Q}_{2r}(\mathcal{U}+\mathcal{P}_0 \widetilde{U}_{\text{in}},V) + \mathcal{Q}_{3r}(V,V)-\overline{\mathcal{D}}V = 0, 
            \\ &V|_{t=0} = \mathcal{P}^{\perp}\widetilde{U}_{in}.
        \end{aligned}\right.
    \end{equation}
    Note that $\mathcal{U}$ is of the form $\mathcal{U}= (-\frac{p_{\theta}^{\circ}}{p_{\rho}^{\circ}}\vartheta, \omega, \vartheta)^{T}$, and \eqref{INSF2} can be written as
    \begin{equation}
        \begin{aligned}
            \partial_t
            \begin{pmatrix}
                \omega\\
                \vartheta
            \end{pmatrix}
            +
            \begin{pmatrix}
                \Pi\left((\omega + 2 \widebar{\widetilde{u}}_{\text{in}} )\!\cdot\!\nabla \omega\right)\\
                (\omega + 2 \widebar{\widetilde{u}}_{\text{in}} )\!\cdot\!\nabla \vartheta
            \end{pmatrix}
            =
            \begin{pmatrix}
                \frac{\mu^{\circ}}{\rho^{\circ}}\Delta\omega\\
                \frac{\kappa^{\circ}p_{\rho}^{\circ}}{(c^{\circ})^2\rho^{\circ}e_{\theta}^{\circ}}\Delta\vartheta
            \end{pmatrix},
        \end{aligned}
    \end{equation}
    with initial data $\omega(0,\cdot) = \omega_{\text{in}}$ and $\vartheta(0,\cdot) = \vartheta_{\text{in}}$ defined by $\mathcal{P}\widetilde{U}_{\text{in}}$.
    A change of variable $\widetilde{\omega}=\omega + 2 \widebar{\widetilde{\u}}_{\text{in}}$ then shows that $(\widetilde{\omega},\ \vartheta)$ satisfies
    \begin{equation}
        \begin{aligned}
            \partial_t
            \begin{pmatrix}
                \widetilde{\omega}\\
                \vartheta
            \end{pmatrix}
            +
            \begin{pmatrix}
                \Pi(\widetilde{\omega} \!\cdot\!\nabla \widetilde{\omega})\\
                \widetilde{\omega} \!\cdot\!\nabla \vartheta
            \end{pmatrix}
            =
            \begin{pmatrix}
                \frac{\mu^{\circ}}{\rho^{\circ}}\Delta\widetilde{\omega}\\
                \frac{\kappa^{\circ}p_{\rho}^{\circ}}{(c^{\circ})^2\rho^{\circ}e_{\theta}^{\circ}}\Delta\vartheta
            \end{pmatrix},
        \end{aligned}
    \end{equation}
    with initial data $\widetilde{\omega}(0,\cdot) = \widetilde{\omega}_{\text{in}} + 2\widebar{\widetilde{u}}_{\text{in}}$ and $\vartheta(0,\cdot) = \vartheta_{\text{in}}$. This is again the incompressible Navier-Stokes-Fourier system.
\end{remark}
\subsection{Novelties of this paper} The main novelties of this paper can be summarized as the following  perspectives: 1, Conceptually, we view the low Mach number limit of compressible Navier-Stokes system in the much more general framework of the weakly nonlinear-dissipative hyperbolic-parabolic system with entropy. The proof in this paper will shed some light on the rigorous justification of the general case. 2, The current paper works on much more general compressible NSF system, besides that the non-isentropic fluid is considered, the pressure and the internal energy are only assumed the most basic thermodynamic requirements. 3. As mentioned above, the role played by the entropy structure in the low Mach number limit is illustrated for the first time to our best acknowledgement. The entropy was used to define the inner product and characterize the incompressible and fast oscillating regimes. 4. Technically, in \cite{Danchin_AJM}, only the convergence rate of the incompressible part was obtained. The corresponding convergence rate of fast oscillation part is much delicate since it is involves the estimate of the resonance sets, both 2-waves and 3-waves. we employ some refined small divisor estimates for the 2-wave resonance, and some geometric observation for 3-wave resonance, for the first time, obtain an explicit convergence rate for fast oscillation part.

\section{Harmonic Analysis Toolbox}
Following the lines in \cite{Danchin_AJM}.
We denote by $\mathcal{F}\u = (\widehat{\u}_k)_{k \in \widetilde{\mathbb{Z}}^{N}}$ the Fourier series of a distribution $u \in \mathcal{S}^{\prime}(\mathbb{T}^{N}_a)$ so that
\begin{equation}
    \u = \sum_{k \in \widetilde{\mathbb{Z}}^{N}} \widehat{\u}_k {\mathrm{e}^{\mathrm{i} k \cdot x}}
\end{equation}
where
\begin{equation}
    \begin{aligned}
        \widehat{\u}_k = \frac{1}{|\mathbb{T}^N_a|}u(\mathrm{e}^{-\mathrm{i} k \cdot x}),\ \widetilde{\mathbb{Z}}^{N} \stackrel{\text { def }}{=} \mathbb{Z}/ a_1 \times \cdots \times \mathbb{Z}/ a_N
        \ \text{and} \ |\mathbb{T}^N_a| \text{ stands for the measure of } \mathbb{T}^{N}_a.
    \end{aligned}
\end{equation}

We use a $C^{\infty}$ symmetric function $\varphi$ of one variable supported in
$\left\{ r \in \mathbb{R}| \frac{5}{6} \leqslant |r| \leqslant \frac{12}{5}\right\}$ such that
\begin{equation}
    \sum_{q \in \mathbb{Z}} \varphi(2^{-q}r) = 1 \text{ for } r \neq 0.
\end{equation}
One can then define the dyadic blocks as follows:
\begin{equation}
    \Delta_q \u \stackrel{\text { def }}{=} \sum_{k \in \widetilde{\mathbb{Z}}^{N}} \varphi (2^{-q}|k|) \widehat{\u}_{k} \mathrm{e}^{i k \cdot x},
\end{equation}
and the low-frequency cut-off:
\begin{equation}
    S_{q} = \widehat{\u}_{0} + \sum_{p \leqslant q-1} \Delta_{p} \u,
\end{equation}
and $\u = \widehat{\u}_{0} + \sum_{q} \Delta_{q} \u$ in $\mathcal{S}^{\prime}(\mathbb{T}^{N}_a)$.
Since we are now dealing with torus, it's obvious $\Delta_p \u = 0$ for $p$ negative enough(depending on the periodic box $\mathbb{T}^{N}_a$).

The dyadic blocks have following quasi-orthogonality:
\begin{equation}
    \begin{aligned}
    \Delta_{p} \Delta_{q} \u \equiv 0 \text{  if  } |p-q|\geqslant 2,\\
    \Delta_{p} (S_{q-1}\u \Delta_{q}\u) \equiv 0 \text{  if  } |p-q| \geqslant 4.
    \end{aligned}
\end{equation}
For $s \in \mathbb{R}$, we define inhomogeneous Besov spaces and Sobolev spaces using Littlewood-Paley decomposition as follows:
\begin{equation}
    \begin{aligned}
    & H^s\left(\mathbb{T}^N_a\right)=\left\{\u \in \mathcal{S}^{\prime}\left(\mathbb{T}^N_a\right) \ \Big|\  \|\u\|_{H^s} \stackrel{\text { def }}{=} \left( {|\widehat{\u}_0| }^2 + \sum_{q \in \mathbb{Z}} 2^{2 s q}\left\|\Delta_q \u\right\|_{L^2}^2\right)^{\frac{1}{2}}<+\infty\right\}, \\
    & B^s\left(\mathbb{T}^N_a\right)=\left\{\u \in \mathcal{S}^{\prime}\left(\mathbb{T}^N_a\right)  \ \Big|\  \|\u\|_{B^s} \stackrel{\text { def }}{=}{|\widehat{\u}_0| } + \sum_{q \in \mathbb{Z}} 2^{s q}\left\|\Delta_q \u\right\|_{L^2}<+\infty\right\} .
    \end{aligned}
\end{equation}
Similarly, we define homogeneous Besov spaces and Sobolev spaces as follows:
\begin{equation}
    \begin{aligned}
    & \dot{H}^s\left(\mathbb{T}^N_a\right)=\left\{\u \in \mathcal{S}^{\prime}\left(\mathbb{T}^N_a\right)   \ \Big|\ \widehat{\u}_0 = 0,\   \|\u\|_{\dot{H}^s} \stackrel{\text { def }}{=} \left(\sum_{q \in \mathbb{Z}} 2^{2 s q}\left\|\Delta_q \u\right\|_{L^2}^2\right)^{\frac{1}{2}}<+\infty\right\}, \\
    & \dot{B}^s\left(\mathbb{T}^N_a\right)=\left\{\u \in \mathcal{S}^{\prime}\left(\mathbb{T}^N_a\right)    \ \Big|\ \widehat{\u}_0 = 0,\  \|\u\|_{\dot{B}^s} \stackrel{\text { def }}{=} \sum_{q \in \mathbb{Z}} 2^{s q}\left\|\Delta_q \u\right\|_{L^2}<+\infty\right\} .
    \end{aligned}
\end{equation}
For $s,t \in \mathbb{R}$, $\eta > 0$, we also use this Littlewood-Paley decomposition to define hybrid Sobolev and Besov norms:
\begin{equation}
    \begin{aligned}
        \left\| \u \right\|_{H^{s,t}_{\eta}} \stackrel{\text { def }}{=}\left(  |{\hat{\u}_0}|^2 + \sum_{q < q_0 } 2^{2qs}\left\| \Delta_{q} \u \right\|^2_{L^2} +\sum_{q \geqslant q_0 } 2^{2qs}{\eta}^{2(t-s)}\left\| \Delta_{q} \u \right\|^2_{L^2}   \right)^{\frac{1}{2}},\\
        \left\| \u \right\|_{B^{s,t}_{\eta}} \stackrel{\text { def }}{=} |{\hat{\u}_0}|+  \sum_{q < q_0 } 2^{qs}\left\| \Delta_{q} \u \right\|_{L^2} +\sum_{q \geqslant q_0 } 2^{qs}{\eta}^{(t-s)}\left\| \Delta_{q} \u \right\|_{L^2}  ,
    \end{aligned}
\end{equation}
where $q_{0} = - \text{log}_2 \eta$. Note that if $s=t$, the hybrid Besov or Sobolev norms coincide with the usual Besov and Sobolev norms, and we indeed have the following equivalence:
\begin{equation}\label{norm_equivalence}
    \left\| \u \right\|_{B^{s,t}_{\eta}} \approx \left\| \u \right\|_{B^s} + {\eta}^{t-s} \left\| \u \right\|_{B^t}.
\end{equation}
Moreover the definition is irrelevant of the choice of $\varphi$. Indeed, we have following lemma:
\begin{lemma}\cite{Danchin_AJM} \label{Localising}
    Let $s \in \mathbb{R}, \u \in H_{\eta}^{s,t}$ (resp. $\u \in B_{\eta}^{s,t}$ ) and $\widetilde{\varphi} \in L^{\infty}\left(\mathbb{R}^N\right)$ be supported in the annulus $C\left(0, R_1, R_2\right)$. 
    There exists some positive constant $C$ and some nonnegative generic element $\left(c_{q,2}\right)_{q \in \mathbb{Z}}$ (resp.$ \left(c_{q,1}\right)_{q \in \mathbb{Z}}$) such that $\sum_q c_{q,2}^2 \leqslant 1$(resp.$\sum_q c_{q,1} \leqslant 1$)
    and for $q \in \mathbb{Z}$
    \begin{equation}
        \begin{aligned}
        \left\|\widetilde{\varphi}\left(2^{-q} D\right) u\right\|_{L^2} \leqslant C c_{q,2} 2^{-q \widetilde{\psi}^{s,t}} { (\widetilde{\phi}^{s,t})}^{-1}\|\u\|_{H^{s,t}_\eta}, \quad\\
        \left(\text { resp. }\left\|\widetilde{\varphi}\left(2^{-q} D\right) u\right\|_{L^2} \leqslant C c_{q,1} 2^{-q \widetilde{\psi}^{s, \ t}}{(\widetilde{\phi}^{s,\ t})}^{-1}\|\u\|_{B^{s,t}_\eta},\right)
    \end{aligned}
    \end{equation}
with
\begin{equation}
    \begin{aligned}
        \widetilde{\psi}^{s,t}(q)= s\mathbf{1}_{q \leqslant q_{0}} + t\mathbf{1}_{q > q_{0}},\ 
        \widetilde{\phi}^{s,t}(q)=\mathbf{1}_{q \leqslant q_{0}} + (\varepsilon \nu^{\circ})^{t-s} \mathbf{1}_{q > q_{0}}.
    \end{aligned}
\end{equation}
Conversely, suppose that $\u=\sum_q \u_q$ in $\mathcal{S}^{\prime}\left(\mathbb{T}^N_a\right)$ with $\operatorname{Supp} \widehat{u}_q \subset 2^q C\left(0, R_1, R_2\right)$ and that
\begin{equation}
    \begin{aligned}
        \left( \sum_{q < q_0 } 2^{2qs}\left\| \Delta_{q} \u \right\|^2_{L^2} +\sum_{q \geqslant q_0 } 2^{2qs}{\eta}^{2(t-s)}\left\| \Delta_{q} \u \right\|^2_{L^2}   \right)^{\frac{1}{2}}=K<+\infty ,\\
        \left(resp.\ \sum_{q < q_0 } 2^{qs}\left\| \Delta_{q} \u \right\|_{L^2} +\sum_{q \geqslant q_0 } 2^{qs}{\eta}^{(t-s)}\left\| \Delta_{q} \u \right\|_{L^2}=K<+\infty \right).
    \end{aligned}
\end{equation}
Then $\u \in H_\eta^{s,t}( \text{resp.} \u \in B_{\eta}^{s,t}) $ and $\|\u\|_{B_\eta^{s,t}} \lesssim K$ $( \text{resp.} \|\u\|_{B_{\eta}^{s,t}} \lesssim K) $.
\end{lemma}
We shall often use the following embedding:
\begin{lemma}\label{embedding}
    \begin{equation}
        \begin{aligned}
            &B^{s}(\mathbb{T}^N_a) \hookrightarrow  H^{s}(\mathbb{T}^N_a) \text{ for } s\in \mathbb{R},\\
            &X^{s,t_1}(\mathbb{T}^N_a) \hookrightarrow  X^{s,t_2}(\mathbb{T}^N_a)  \text{ for } s \in \mathbb{R} \text{ and } t_1 < t_2,\ (X \text{ stands for } B \text{ or } H),\\
            &B^{\frac{N}{2}}(\mathbb{T}^N_a) \hookrightarrow L^{\infty}(\mathbb{T}^N_a).
        \end{aligned}
    \end{equation}
\end{lemma}
We also need the following useful estimate which may be proved by using the definition of hybrid Besov norms
\begin{equation}\label{delta_estimates}
    \begin{aligned}
        \left\| \u \right\|_{B^{\frac{N}{2}-\delta}}
        &\leqslant   {\eta}^{\delta-1} \left\| \u \right\|_{B_{\eta}^{\frac{N}{2}-1,\frac{N}{2}}}
        \leqslant   {\eta}^{\delta-1} \left\| \u \right\|_{B_{\eta}^{\frac{N}{2}-1,\frac{N}{2}+1}}, \ \text{for } 0 \leqslant \delta  \leqslant 1.
    \end{aligned}
\end{equation}
Now we state some continuity results for the product and the composition which will be widely used in later content.
\begin{lemma}\cite{Danchin_AJM} \label{product_estimate}
    For $\eta > 0$, let $\max(s_1,\ s_2,\ t_1,\ t_2) <\frac{N}{2}$ and ${\min}(s_1 + s_2,\  t_1 + t_2 ) > 0$ then
    \begin{equation}
        \left\| \u {\vv} \right\|_{H_{\eta}^{s_1 + t_1 -\frac{N}{2},s_2 + t_2 -\frac{N}{2}}(\mathbb{T}^N_a)} \leqslant \left\| \u \right\|_{H_{\eta}^{s_1,s_2}(\mathbb{T}^N_a)} \left\| {\vv} \right\|_{H_{\eta}^{t_1,t_2}(\mathbb{T}^N_a)}.
    \end{equation}
    If $s=\frac{N}{2}$, $t < \frac{N}{2}$ $(resp.\  t=\frac{N}{2},\ s < \frac{N}{2})$, then we have
    \begin{equation}
        \begin{aligned}
            \left\| \u {\vv} \right\|_{H_{\eta}^{s + t -\frac{N}{2}}(\mathbb{T}^N_a)} \leqslant \left\| \u \right\|_{H_{\eta}^{s}(\mathbb{T}^N_a)\cap L^{\infty}(\mathbb{T}^N_a)} \left\| {\vv} \right\|_{H_{\eta}^{t}(\mathbb{T}^N_a)}
            (resp.\ \left\| {\vv} \right\|_{H_{\eta}^{s}(\mathbb{T}^N_a)} \left\| \u \right\|_{H_{\eta}^{t}(\mathbb{T}^N_a)\cap L^{\infty}(\mathbb{T}^N_a)}.).
        \end{aligned}
    \end{equation}
    If $\max(s_1,\ s_2,\ t_1,\ t_2) \leqslant \frac{N}{2}$ and ${\min}(s_1 + s_2,\  t_1 + t_2 ) > 0$, then
    \begin{equation}
        \left\| \u {\vv} \right\|_{B_{\eta}^{s_1 + t_1 -\frac{N}{2},s_2 + t_2 -\frac{N}{2}}(\mathbb{T}^N_a)} \leqslant \left\| \u \right\|_{B_{\eta}^{s_1,s_2}(\mathbb{T}^N_a)} \left\| {\vv} \right\|_{B_{\eta}^{t_1,t_2}(\mathbb{T}^N_a)}.
    \end{equation}
    For $s > 0$,
    \begin{equation}
        \left\| \u {\vv} \right\|_{B^{s}(\mathbb{T}^N_a)}
        \leqslant
        \left\| \u \right\|_{L^{\infty}(\mathbb{T}^N_a)} \left\| {\vv} \right\|_{B^{s}(\mathbb{T}^N_a)} +
        \left\| {\vv} \right\|_{L^{\infty}(\mathbb{T}^N_a)} \left\| \u \right\|_{B^{s}(\mathbb{T}^N_a)}.
    \end{equation}
\end{lemma}
\begin{lemma} \cite{Danchin_AJM}\label{composition_estimates}
    Let $s \geqslant 0$ and $\u \in X^{s}(\mathbb{T}^N_a) \cap L^{\infty}(\mathbb{T}^N_a)$ with $X = H$ or $B$. Let $F \in W^{[s]+2,\infty}_{loc}(\mathbb{R})$ be such that $F(0)=0$.
    Then $F(\u) \in X^{s} \cap L^{\infty}$ and there exists a constant $C = C(s,N,F,\left\|u\right\|_{L^{\infty}})$ such that
    \begin{equation}
        \left\| F(\u) \right\|_{X^{s}(\mathbb{T}^N_a)} \leqslant C  \left\| \u \right\|_{X^{s}(\mathbb{T}^N_a)}.
    \end{equation}
\end{lemma}
We have the following multiplier estimate:
\begin{lemma} \cite{Danchin_AJM}
    If A is an homogeneous function of degree $m \in \mathbb{R}$, then we have
    \begin{equation}
        \left\| A(D) \u \right\|_{H^{\sigma-m}(\mathbb{T}^N_a)} \lesssim \left\| \u \right\|_{H^{\sigma}(\mathbb{T}^N_a)}.
    \end{equation}
\end{lemma}
\begin{lemma} \cite{Danchin_AJM} \label{convolution_estimate}
    Let $A$ be a function on $\widetilde{\mathbb{Z}}^N \times \widetilde{\mathbb{Z}}^N $ such that for $(\alpha,\ \beta,\ \gamma) \in [-1,2]^3 $
    we have
    \begin{equation}
        |A(k,m)| \leqslant K |m|^{\alpha}|k|^{\beta}|l|^{\gamma},
    \end{equation}
    with $K$ some positive constant.

    Let
    \begin{equation}
        B(\u,{\vv}) = \sum_{(k,m)\in\widetilde{\mathbb{Z}}^N \times \widetilde{\mathbb{Z}}^N }A(k,m)  \widehat{\u}_{k} \widehat{{\vv}}_{m-k} \mathrm{e}^{i m \cdot x}.
    \end{equation}
    Then we have, for $s_1 + s_2 -\frac{N}{2}=(\alpha + s)$, $s_1,s_2\leqslant \frac{N}{2}$, $s_1+s_2 > 0$
    \begin{equation}
        \left\| B(\u,{\vv}) \right\|_{H^{s}(\mathbb{T}^N_a)} \leqslant \left\| \u \right\|_{H^{s_1+\beta}(\mathbb{T}^N_a)} \left\| {\vv} \right\|_{H^{s_2+\gamma}(\mathbb{T}^N_a)}.
    \end{equation}
\end{lemma}

Localising technique by means of operator $\Delta_q$ naturally give rise to following spaces:
\begin{equation}
    \begin{aligned}
        \widetilde{L}^{r}_{T}&(H^{s,t}_{\eta})\stackrel{\text { def }}{=} \Bigg\{u \in  \mathscr{D}^{\prime}\left((0,T)\times\mathbb{T}^N_a\right) \  \ \Big|  \\
        &\left. \left\|u \right\|_{\widetilde{L}^{r}_{T}(H^{s,t}_{\eta})} \stackrel{\text { def }}{=} \|\widehat{\u}_0\|_{L^r_T} +  \left( \sum_{q < q_0 } 2^{2qs}\left\| \Delta_{q} u \right\|^2_{L^{r}_T(L^2)} +\sum_{q \geqslant q_0 } 2^{2qs}{\eta}^{2(t-s)}\left\| \Delta_{q} \u \right\|^2_{L^{r}_{T}(L^2)}   \right)^{\frac{1}{2}} < + \infty\right\},\\
    \end{aligned}
\end{equation}
\begin{equation}
    \begin{aligned}
        \widetilde{L}^{r}_{T}&(B^{s,t}_{\eta})\stackrel{\text { def }}{=} \Bigg\{ \u \in \mathscr{D}^{\prime}\left((0,T)\times \mathbb{T}^N_a\right) \  \ \Big|  \\
        &\left. \left\|u \right\|_{\widetilde{L}^{r}_{T}(B^{s,t}_{\eta})} \stackrel{\text { def }}{=}  \|\widehat{\u}_0\|_{L^r_T} + \sum_{q < q_0 } 2^{qs}\left\| \Delta_{q} \u \right\|^2_{L^{r}_T(L^2)} +\sum_{q \geqslant q_0 } 2^{qs}{\eta}^{(t-s)}\left\| \Delta_{q} \u \right\|^2_{L^{r}_{T}(L^2)} < +\infty   \right\}.\\
    \end{aligned}
\end{equation}
We further denote $\widetilde{C}([0,T];X)= \widetilde{L}^{\infty}_T(X) \cap C([0,T];X)$,
and if $T = + \infty$ we shall simply use the notation $\widetilde{L}^{r}(B^{s,t}_{\eta})$.

Using Minkowski inequality, it's easy to deduce that
\begin{equation}\label{t_embedding}
    \begin{aligned}
        &\widetilde{L}^{r}_{T}(H_{\eta}^{s,t}) \hookrightarrow {L}^{r}_{T}(H_{\eta}^{s,t}) \ \ \text{if }\ \  r > 2,\ \
        {L}^{r}_{T}(H_{\eta}^{s,t}) \hookrightarrow \widetilde{L}^{r}_{T}(H_{\eta}^{s,t}) \ \ \text{if }\ \  r < 2,\\
        &\widetilde{L}^{r}_{T}(B_{\eta}^{s,t}) \hookrightarrow {L}^{r}_{T}(B_{\eta}^{s,t}) \ \ \text{if }\ \  r > 1,\ \
    \end{aligned}
\end{equation}
and moreover
\begin{equation}
    \widetilde{L}^{2}_{T}(H_{\eta}^{s,t}) = {L}^{2}_{T}(H_{\eta}^{s,t}), \ \widetilde{L}^{1}_{T}(B_{\eta}^{s,t}) = {L}^{1}_{T}(B_{\eta}^{s,t}).
\end{equation}
Similarly we have product estimate for above spaces as in Sobolev or Besov Space. For example
\begin{lemma} \cite{Danchin_AJM} \label{mixed_product_estimate}
    Let $s_1,\ s_2,\ t_1,\ t_2 \leqslant \frac{N}{2}$ and $\text{min}(s_1 + s_2,\  t_1 + t_2 ) > 0$ then
    \begin{equation}
        \left\| \u {\vv} \right\|_{\widetilde{L}^{r}_T(B_{\eta}^{s_1 + t_1 -\frac{N}{2},s_2 + t_2 -\frac{N}{2}})} \leqslant \left\| \u \right\|_{\widetilde{L}_T^{r_1}(B_{\eta}^{s_1,s_2})} \left\| {\vv} \right\|_{\widetilde{L}_T^{r_2}(B_{\eta}^{t_1,t_2})},
    \end{equation}
    with $\frac{1}{r} = \frac{1}{r_1}+ \frac{1}{r_2}$.
\end{lemma}
Note also that as in Lemma \ref{product_estimate} similar result holds for $\widetilde{L}^{r}(H_{\eta}^{s,t})$ type spaces.

We also state some interpolation properties:
\begin{lemma}\cite{Danchin_AJM} \label{interpolation}
    For $\theta \in [0,1]$, let $t= \theta t_1 +(1-\theta)t_2  $, $ s = \theta s_1 +(1-\theta)s_2  $, $\frac{1}{r} = \theta \frac{1}{r_1}+(1-\theta) \frac{1}{r_2}$, then
    \begin{equation}
        \begin{aligned}
            \left\| \u \right\|_{X_{\eta}^{s,t}} \lesssim  \left\| \u \right\|_{X_{\eta}^{s_1,t_1}}^{\theta}  \left\| \u \right\|_{X_{\eta}^{s_2,t_2}}^{1-\theta} ,\\
            \left\| \u \right\|_{\widetilde{L}^{r}_{T}(X_{\eta}^{s,t})} \lesssim  \left\| \u \right\|_{\widetilde{L}^{r_1}_{T}(X_{\eta}^{s_1,t_1})}^{\theta}  \left\| \u \right\|_{\widetilde{L}^{r_2}_{T}(X_{\eta}^{s_2,t_2})}^{1-\theta} ,\\
        \end{aligned}
    \end{equation}
   where $X = B \text{ or }H$ .
\end{lemma}
Dealing with the high and low frequencies separately, the following inequalities come in handy.
\begin{lemma}\cite{Danchin_AJM} \label{truncated_frequency_estimate}
For any $r \in [1,+\infty]$, $\sigma \in \mathbb{R}$ and $\delta > 0$, we have
\begin{equation}
    \begin{aligned}
        \left\| {\vv}_M \right\|_{\widetilde{L}_T^{r}(H^{\sigma})} \lesssim M^{\delta} \left\| {\vv} \right\|_{\widetilde{L}_T^{r}(H^{\sigma-\delta})},\ \forall\  \vv \in \widetilde{L}_T^{r}(H^{\sigma-\delta});
        \\
        \left\| {\vv}^M \right\|_{\widetilde{L}_T^{r}(H^{\sigma})} \lesssim M^{-\delta} \left\| {\vv} \right\|_{\widetilde{L}_T^{r}(H^{\sigma+\delta})}, \forall \ \vv \in \widetilde{L}_T^{r}(H^{\sigma+\delta}).
    \end{aligned}
\end{equation}
Where ${\vv}_M$ and ${\vv}^M$ defined by
\begin{equation}
    \begin{aligned}
        {\vv}_M \stackrel{\mathrm{def}}{=} \sum_{|k| \leqslant M} \widehat{v}_k e^{\mathrm{i} k\cdot x} \ \text{and} \ {\vv}^M \stackrel{\mathrm{def}}{=} \sum_{|k| \geqslant M} \widehat{v}_k e^{\mathrm{i} k\cdot x}.
    \end{aligned}
\end{equation}
\end{lemma}
 The operator $\mathrm{e}^{t \mathcal{A}}$ is an isomorphism on $H^s$ or $B^s$. Moreover we have
\begin{lemma}\cite{Danchin_AJM}\label{norm_conserve}
    For $\alpha,\ s \in \mathbb{R}$, the operator
    \begin{equation}
        \u \mapsto \mathrm{e}^{ \alpha t \mathcal{A}} \u
    \end{equation}
    is an isomorphism on $L^{r}((0,T);B^s(\mathbb{T}^N_a))$, $\widetilde{L}_T^{r}(B^s(\mathbb{T}^N_a))$, $L^{r}((0,T);H^s(\mathbb{T}^N_a))$ and $\widetilde{L}_T^{r}(H^s(\mathbb{T}^N_a))$.
\end{lemma}
\section{Proof of Main Theorem}
In this section, we deal with the convergence part of Theorem \ref{MAIN}. Before dive into the proof, we need to introduce some notations.
Let $s\stackrel{\mathrm{def}}{=}\frac{N}{2}-\delta$, with $\delta \in (0,1)$, we define
\begin{equation}
    \begin{aligned}
    &\begin{gathered}
        X \stackrel{\mathrm{def}}{=} \max(\|V\|_{F^{\frac{N}{2}}},\sup_{0<\varepsilon\leqslant\varepsilon_0}
        \|\mathcal{P}^{\perp}\widetilde{U}^{\varepsilon}\|_{F^{\frac{N}{2}}}),\\
        X_0 \stackrel{\mathrm{def}}{=} \|\mathcal{P}^{\perp}\widetilde{U}^{\varepsilon}(0)\|_{B^{\frac{N}{2}-1}},\\
        Y \stackrel{\mathrm{def}}{=} \max(\|\mathcal{U}\|_{F^{\frac{N}{2}}},\sup_{0<\varepsilon\leqslant\varepsilon_0}
        \|\mathcal{P}\widetilde{U}^{\varepsilon}\|_{F^{\frac{N}{2}}}),\ \
        Y_0 \stackrel{\mathrm{def}}{=} \|\mathcal{P}\widetilde{U}^{\varepsilon}(0)\|_{B^{\frac{N}{2}-1}},\\
        W^{\varepsilon} \stackrel{\mathrm{def}}{=} \|\mathcal{P}\widetilde{U}^{\varepsilon}-\mathcal{U}\|_{G^{s}},\ \
        Z^{\varepsilon} \stackrel{\mathrm{def}}{=} \|\widetilde{V}^{\varepsilon} - V\|_{G^{s}},\\
    \end{gathered}
    \\ & \begin{aligned}
        \mathcal{Y}^{\varepsilon} \stackrel{\mathrm{def}}{=}
        \max
        (\|\rho^{\circ}+\varepsilon \tilde{\rho}^{\varepsilon}\|_{L^\infty((0,+\infty)\times \mathbb{T}_a^{N})},
        \|(\rho+\varepsilon \tilde{\rho}^{\varepsilon})^{-1}\|_{L^\infty((0,+\infty)\times \mathbb{T}_a^{N})},
        \\ \|\theta^{\circ}+\varepsilon \tilde{\theta}^{\varepsilon}\|_{L^\infty((0,+\infty)\times \mathbb{T}_a^{N})},
        \|(\theta^{\circ}+\varepsilon \tilde{\theta}^{\varepsilon})^{-1}\|_{L^\infty((0,+\infty)\times \mathbb{T}_a^{N})}).
    \end{aligned}
\end{aligned}
\end{equation}
We remark that in view of Lemma \ref{norm_conserve}, the following bound is true
\begin{equation}
    \|\widetilde{V}^{\varepsilon}\|_{F^{\frac{N}{2}}} \leqslant X.
\end{equation}
The proof is divided into five steps. In the first four steps we assume that
the energy inequality \eqref{energy_ineq} is true. Then for  $\alpha_0$ small enough which will be determined, one can show that
$\mathcal{Y}^{\varepsilon}$ is bounded from both below and above and away from 0.

Our strategy for proving convergence lies in estimating the
equations satisfied by the difference between the solution
of the full problem and that of the limit problem.
Owing to the entropy structure discussed in previous sections,
the estimating process takes place separately in
$\text{Null}(\mathcal{A})$ and $\text{Null}(\mathcal{A})^{\perp}$.
Denoting $w^{\varepsilon}=\mathcal{P}\widetilde{U}^{\varepsilon}-\mathcal{U}$
and $z^\varepsilon=\widetilde{V}^\varepsilon-V$, we proceed as follows
\subsection*{First Step. Convergence of $\widetilde{V}^{\varepsilon}$ to $V$}
We have the following proposition (see the proof in Section \ref{sec_oscillating}):
\begin{proposition} \label{oscillation_convergence}
    Let $\widetilde{U}^{\varepsilon}$ be a solution of \eqref{nsc_epsilon}
    and
    $\mathcal{U}\in F^{\frac{N}{2}}$
    (resp. $V \in F^{\frac{N}{2}}$) be the solution of \eqref{INSF} (resp.\eqref{LS})
    satisfying $\| \widetilde{U}^{\varepsilon} \|_{E^{\frac{N}{2}}_{\varepsilon \nu^{\circ}}} + \left\| \mathcal{U},V \right\|_{F^{\frac{N}{2}}}
    \leqslant C \alpha_0$. 
    Then the following estimate holds:
    \begin{equation}
        \begin{aligned}
        \left\| z^{\varepsilon} \right\|&_{{L}^2(H^{\frac{N}{2}} \cap L^{\infty}) \cap \widetilde{L}^{\infty}(H^{\frac{N}{2}-1})}
        \leqslant C \mathrm{e}^{C(Y^2 + X^2)} \times \\ \bigg(
        &\varepsilon M^{1-\delta} \left( X + Y + Y^2 + XY^2 +Y^3 +( X + Y )^2 + (X + Y)^3 \right)  \\
        &+ \varepsilon C_M M^{2-\delta}(X + X^2 + XY +X^{2}Y + X Y^{2} + X^3) \\
        &+(Y+X)W^{\varepsilon} + X Z^{\varepsilon} + {(Z^{\varepsilon})}^2   + \varepsilon M^{1-\delta}\left(X_0 + Y_0 + {Y_0}^2 \right) \\
        &+ \varepsilon C_M M^{2-\delta} \left( X_0 Y_0 + {X_0}^2\right) + M^{-\delta}(Y + XY + Y^2 + X^2) \bigg),
    \end{aligned}
    \end{equation}
    where $C$ is some constant depends only on $(N,p,e,\mu,\lambda,\kappa)$.
\end{proposition}
By the energy estimate \eqref{energy_ineq} and Lemma \ref{embedding} we have
\begin{equation}\label{mutual_bound}
    \max_{0<\varepsilon\leqslant\varepsilon_0}
    (X,X_0,Y,Y_0)\leqslant C\alpha_0,
\end{equation}
then
\begin{equation}\label{z_temp}
    Z^{\varepsilon} \leqslant C \mathrm{e}^{2C\alpha^{2}_0} \left( (\alpha_0 + \alpha_0^2 + \alpha^3_0) (\varepsilon M^{1-\delta} + \varepsilon C_M M^{2-\delta} + M^{-\delta}) + \alpha_0(W^{\varepsilon}+ Z^{\varepsilon})  \right).
\end{equation}

\subsection*{Second Step. Refined small divisor estimates}
In this step we use some estimates for linear forms to give an explicit upper bound for $C_M$ (see Proposition \ref{small_divisor_estimates}), namely for $\sigma = \max{\left\{ 2N +2\tau -1 , 5\right\}}$,
\begin{equation}\label{CM_bound}
    C_M  \leqslant C(1+M)^{\sigma}.
\end{equation}
We using this result and the estimates for $W^{\varepsilon}$ to completes the proof for convergence of the oscillating part.
Combining  \eqref{z_temp}, \eqref{CM_bound} and \eqref{W_convergence} we have for $M >1$,
\begin{equation}
    Z^{\varepsilon} \leqslant C \mathrm{e}^{2C\alpha^{2}_0} \left( (\alpha_0 + \alpha_0^2 + \alpha^3_0) (\varepsilon M^{1-\delta} + \varepsilon M^{\sigma} M^{2-\delta} + M^{-\delta} + \varepsilon^{\delta}) + \alpha_0 Z^{\varepsilon}  \right).
\end{equation}
Upon requiring $\alpha_0$ small to control the linear term $\alpha_0 Z^{\varepsilon}$ and taking $\displaystyle M = \varepsilon^{-\frac{1}{2 + \sigma}}$,
we finally have
\begin{equation}
    Z^{\varepsilon} = \mathcal{O}(\varepsilon^{\frac{\delta}{2+\sigma}}).
\end{equation}

\subsection*{Third Step. Convergence of the Incompressible part}
This is guaranteed by the following proposition(see the proof in Section \ref{sec_Incompressible}).
\begin{proposition} \label{incompressible_convergence}
    Let $\widetilde{U}^{\varepsilon}$ be a solution of \eqref{nsc_epsilon}
    and
    $\mathcal{U}\in F^{\frac{N}{2}}$
    be the solution of \eqref{INSF}
    satisfying $\| \widetilde{U}^{\varepsilon} \|_{E^{\frac{N}{2}}_{\varepsilon \nu^{\circ}}} + \left\| \mathcal{U}\right\|_{F^{\frac{N}{2}}}
    \leqslant C \alpha_0$. 
    Then the following estimate holds:
    \begin{equation}
        \begin{aligned}
            W^{\varepsilon} \leqslant& \varepsilon^{\delta}C\mathrm{e}^{C(X^2+Y^2)}\times
            \bigg((X_0+X_0^2+Y_0^2)\\
            &+(X+Y+1)\left((X+Y)^2+X+Y+\varepsilon^{\delta}(X^3+X^2)\right)\\
            &+(X^3+X^2)+(X^2+XY+X)\bigg),
        \end{aligned}
    \end{equation}
    where $C$ is some constant depends only on $(N,p,e,\mu,\lambda,\kappa)$.
\end{proposition}
In view of \eqref{mutual_bound},
the conclusion of this proposition leads to
\begin{equation}\label{W_convergence}
    W^{\varepsilon}=\mathcal{O}(\varepsilon^{\delta}).
\end{equation}

\subsection*{Fourth Step. A priori estimate for $(\widetilde{\rho}^{\varepsilon},\ \widetilde{u}^{\varepsilon},\ \widetilde{\theta}^{\varepsilon})$}
In this step we show that for a suitable choice of $\alpha_0$,
$(\widetilde{\rho}^{\varepsilon}, \ \widetilde{u}^{\varepsilon},\ \widetilde{\theta}^{\varepsilon})$ has an uniform bound in $E^{\frac{N}{2}}_{\varepsilon \nu^{\circ}}$,
and that $\mathcal{Y}^{\varepsilon}$ is bounded from both above and below and boudned away from 0 (see Section \ref{sec_priori}).
\subsection*{Fifth Step. Well-posedness of the limit system}
In this step we show that the limit system is well-posed as long as \eqref{INSF} is well-posed, and that $\mathcal{U}, V $ is boudned in $F^{\frac{N}{2}}$(see section \ref{sec_limit} ),
which together with the last step completes the proof for \eqref{energy_ineq}.

\section{Convergence of the oscillating part of the system} \label{sec_oscillating}
In this section we give the proof of the convergence of the compressible part.
To do so, we first show that $\|\mathcal{P}_0 \widetilde{U}^{\varepsilon}\|_{L^{\infty}(dt)} =\mathcal{O}(\varepsilon) $, this is actually ensured by conservation law.
From \eqref{NSCe} we have
    \begin{align}
         &\partial_{t} \int_{|\mathbb{T}_a^{N}|} \rho^{\varepsilon}\,dx  = 0, \label{conservation_of_rho}\\
         &\partial_{t}\int_{|\mathbb{T}_a^{N}|}(\rho^{\varepsilon} \u^{\varepsilon})dx  = 0, \label{conservation_of_rho_u}\\
         &\partial_{t}\int_{|\mathbb{T}_a^{N}|} (\frac{1}{2}\rho^{\varepsilon} | \u^{\varepsilon}|^2 + \rho^{\varepsilon} e^{\varepsilon}) \label{conservation_of_e} dx=0.
        \end{align}
Clearly \eqref{conservation_of_rho} implies that $\widebar{\rho}^{\varepsilon} =\widebar{\rho}_{\text{in}} $ and therefore $\widebar{\widetilde{\rho}}^{\varepsilon} =\widebar{\widetilde{\rho}}_{\text{in}} $.
From \eqref{conservation_of_rho_u} we have
\begin{equation}
    \int_{|\mathbb{T}_a^{N}|} \widetilde{\u}^{\varepsilon}-{\widetilde{\u}}_{\text{in}}dx= \varepsilon \frac{1}{\rho^{\circ}}\int_{|\mathbb{T}_a^{N}|} \widetilde{\rho}^{\varepsilon} \widetilde{\u}^{\varepsilon} - {\widetilde{\rho}}_{\text{in}} {\widetilde{\u}}_{\text{in}} dx ,
\end{equation}
and therefore
\begin{equation} \label{u_bar_convergence}
    \|\widebar{\widetilde{\u}}^{\varepsilon}-\widebar{\widetilde{\u}}_{\text{in}}\|_{L^{\infty}(dt)} \lesssim \varepsilon \| \widetilde{\rho}^{\varepsilon}\|_{\widetilde{L}^{\infty}(B^{\frac{N-1}{2}})} \| \widetilde{\u}^{\varepsilon}\|_{\widetilde{L}^{\infty}(B^{\frac{N-1}{2}})}
    + \varepsilon \| \widebar{\widetilde{\rho}}_{\text{in}}\|_{B^{\frac{N-1}{2}}} \| \widebar{\widetilde{\u}}_{\text{in}}\|_{B^{\frac{N-1}{2}}}.
\end{equation}
Similarly, from \eqref{conservation_of_e} we have
\begin{equation} \label{theta_bar_convergence}
    \begin{aligned}
    \rho^{\circ} e^{\circ}_{\theta} \int_{|\mathbb{T}_a^{N}|}  \widetilde{\theta}^{\varepsilon} - \widetilde{\theta}_{\text{in}}
    = \frac{\rho^{\circ}}{2} \int_{|\mathbb{T}_a^{N}|} {\widetilde{\u}}_{\text{in}} - \widetilde{\u}^{\varepsilon} dx
    - \frac{1}{\varepsilon}\int_{|\mathbb{T}_a^{N}|} \rho^{\varepsilon}(e^{\varepsilon} - e^{\circ} - e^{\circ}_{\rho}\varepsilon \widetilde{\rho}^{\varepsilon} -e^{\circ}_{\theta}\varepsilon \widetilde{\theta}^{\varepsilon}) dx\\
    + \frac{1}{\varepsilon}\int_{|\mathbb{T}_a^{N}|} \rho_{\text{in}}(e_{\text{in}} - e^{\circ} - e^{\circ}_{\rho}\varepsilon {\widetilde{\rho}}_{\text{in}} -e^{\circ}_{\theta}\varepsilon \widetilde{\theta}_{\text{in}})dx
    - \int_{|\mathbb{T}_a^{N}|} \widetilde{\rho}^{\varepsilon} e^{\circ}_{\rho} \varepsilon \widetilde{\rho}^{\varepsilon} + \widetilde{\rho}^{\varepsilon} e^{\circ}_{\theta} \varepsilon \widetilde{\theta}^{\varepsilon} dx \\
    +  \int_{|\mathbb{T}_a^{N}|} {\widetilde{\rho}}_{\text{in}} e^{\circ}_{\rho} \varepsilon {\widetilde{\rho}}_{\text{in}} + {\widetilde{\rho}}_{\text{in}} e^{\circ}_{\theta} \varepsilon \widetilde{\theta}_{\text{in}} dx ,
    \end{aligned}
\end{equation}
Combining \eqref{u_bar_convergence} and \eqref{theta_bar_convergence} and using Lemma \ref{composition_estimates} leads to
\begin{equation}
    \begin{aligned}
    \|\widebar{\widetilde{\theta}}^{\varepsilon}-\widebar{\widetilde{\theta}}_{\text{in}}\|_{L^{\infty}(dt)}
    &\lesssim \varepsilon \| \widetilde{\rho}^{\varepsilon}\|_{\widetilde{L}^{\infty}(B^{\frac{N-1}{2}})} \| \widetilde{\u}^{\varepsilon}\|_{\widetilde{L}^{\infty}(B^{\frac{N-1}{2}})}
    +\varepsilon \| {\widetilde{\rho}}_{\text{in}}\|_{B^{\frac{N-1}{2}}} \| {\widetilde{\u}}_{\text{in}}\|_{B^{\frac{N-1}{2}}}\\
    &+\varepsilon( 1+C(\mathcal{Y}^{\varepsilon})) \| \widetilde{\rho}^{\varepsilon}\|_{\widetilde{L}^{\infty}(B^{\frac{N-1}{2}})} \| (\widetilde{\rho}^{\varepsilon},\widetilde{\theta}^{\varepsilon})\|_{\widetilde{L}^{\infty}(B^{\frac{N-1}{2}})}   \\
    &+\varepsilon( 1+C(\mathcal{Y}^{\varepsilon})) \| {\widetilde{\rho}}_{\text{in}}\|_{B^{\frac{N-1}{2}}} \| ({\widetilde{\rho}}_{\text{in}},\widetilde{\theta}_{\text{in}})\|_{B^{\frac{N-1}{2}}}.
    \end{aligned}
\end{equation}
In summary, we have
\begin{equation}
    \|\mathcal{P}_0 \widetilde{U}^{\varepsilon} - \mathcal{P}_0  \widetilde{U}_{\text{in}}\|_{L^{\infty}(dt)} = \mathcal{O}(\varepsilon),
\end{equation}
and therefore
\begin{equation}
    \|\mathcal{P}_0 \widetilde{U}^{\varepsilon}  \|_{L^{\infty}(dt)} = \mathcal{O}(\varepsilon),
\end{equation}
since we assume that $ \mathcal{P}_{0}\widetilde{U}_{\text{in}} = 0$.

We now turn to the proof of Proposition \ref{oscillation_convergence}. We first give a precise description for the operators $\mathcal{Q}^{\varepsilon}_{2r}$, $\mathcal{Q}^{\varepsilon}_{3r}$, $\mathcal{D}^{\varepsilon}$ in \eqref{nsc_orthogonal}
and averaged operators $\mathcal{Q}_{2r}$, $\mathcal{Q}_{3r}$, $\overline{\mathcal{D}}$ in \eqref{LS}.
Using the definition of  $\mathcal{Q}^{\varepsilon}_{2r}$ and  $\mathcal{Q}^{\varepsilon}_{3r}$,
 for $A=\sum_{\alpha,k}A^{\alpha}_{k}H^{\alpha}_{k}$, $B=\sum_{\alpha,k}B^{\alpha}_{k}H^{\alpha}_{k}$ and for $U  = \begin{pmatrix}
    \rho,&
    u,&
    \theta
\end{pmatrix}^{\text{T}} \in \mathrm{ker} \mathcal{A}$
we have
\begin{equation} \label{Q2r_epsilon}
    \begin{aligned}
    \mathcal{Q}^{\varepsilon}_{2r}\left(U,B\right)
    = 2 \mathrm{e}^{\frac{t}{\varepsilon}\mathcal{A}} \mathcal{P}^{\perp} \mathcal{Q}\left( U ,  \mathrm{e}^{-\frac{t}{\varepsilon} \mathcal{A}}B\right)
    = 2 \sum_{\gamma, m} \left< \mathcal{Q} \left( {U},\mathrm{e}^{-\frac{t}{\varepsilon} \mathcal{A}}B\right), H^{\gamma}_{m} \right>_{\mathbb{H}} H^{\gamma}_{m}
    \\= 2 \sum_{\gamma, \alpha}  \sum_{m,k}  \left<\mathcal{Q}\left( U, H^{\alpha}_{k}\right), H^{\gamma}_{m} \right>_{\mathbb{H}} B^{\alpha}_k H^{\gamma}_{m}
    \mathrm{e}^{-\frac{t}{\varepsilon}\left(\lambda^{\alpha}_{ {k}}-\lambda^{\gamma}_{ {m}}\right)},
    \end{aligned}
\end{equation}
\begin{equation} \label{Q3r_epsilon}
    \begin{aligned}
    \mathcal{Q}^{\varepsilon}_{3r}\left(A,B\right)
    =  \mathrm{e}^{\frac{t}{\varepsilon}\mathcal{A}} \mathcal{P}^{\perp} \mathcal{Q}\left(  \mathrm{e}^{-\frac{t}{\varepsilon} \mathcal{A}}A ,  \mathrm{e}^{-\frac{t}{\varepsilon} \mathcal{A}}B\right)
    =  \sum_{\gamma, m} \left< \mathcal{Q} \left(  \mathrm{e}^{-\frac{t}{\varepsilon} \mathcal{A}}A,\mathrm{e}^{-\frac{t}{\varepsilon} \mathcal{A}}B\right), H^{\gamma}_{m} \right>_{\mathbb{H}} H^{\gamma}_{m}
    \\=  \sum_{\gamma, m} \sum_{\alpha , k+l=m}  \left<\mathcal{Q}\left(  H^{\alpha}_{k}, H^{\beta}_{l}\right), H^{\gamma}_{m} \right>_{\mathbb{H}} A^{\alpha}_{k} B^{\beta}_{l} H^{\gamma}_{m}
    \mathrm{e}^{-\frac{t}{\varepsilon}\left(\lambda^{\alpha}_{ {k}}+ \lambda^{\beta}_{l}-\lambda^{\gamma}_{ {m}}\right)}.
    \end{aligned}
\end{equation}
Direct calculation using the definition of $\mathcal{Q}$ we have
\begin{equation}
    \begin{aligned}
        2\left<\mathcal{Q}\left(  \widehat{U}_{l} \mathrm{e}^{i l \cdot x}, H^{\alpha}_{k}\right), H^{\gamma}_{m} \right>_{\mathbb{H}}
        &=\mathrm{i} (c_N c^{\circ})^2
        \bigg[ \frac{\rho^\circ}{\theta^\circ} \left(1+\alpha \gamma \operatorname{sg}( {k})\operatorname{sg}( {m})\frac{k \cdot (l+m)}{|k| |m|}\right) \widehat{\u}_{ {l}}\!\cdot\! {m} \\
        &+ \frac{\rho^{\circ}}{\theta^{\circ}}
            \left(
                C_1 \rho^{\circ} \widehat{\rho}_l + C_2 \frac{\theta^{\circ}p^{\circ}_{\theta}}{\rho^{\circ}e^{\circ}_{\theta}} \widehat{\theta}_l
            \right) \gamma \operatorname{sg}(m)|m|
        +  \frac{\rho^{\circ}}{\theta^{\circ}}
            \left(
                \frac{p^{\circ}_{\rho}}{\theta^{\circ}} \widehat{\rho}_l \alpha \operatorname{sg}(k) \frac{k \cdot m}{|k|}
            \right) \\
        &+  \frac{\rho^{\circ}}{\theta^{\circ}}
            \left(
                C_3 \frac{\theta^{\circ}p^{\circ}_{\theta}}{\rho^{\circ}e^{\circ}_{\theta}} \widehat{\rho}_l  + C_4 \rho^{\circ} \widehat{\theta}_l
            \right) \gamma \operatorname{sg}(m) \frac{l \cdot m }{|m|} \\
        &+  \frac{p^{\circ}_{\theta}}{\theta^{\circ}} \widehat{\theta}_l \alpha \operatorname{sg}(k) \frac{k \cdot l}{|k|}
        +\frac{p^{\circ}_{\theta}}{\theta^{\circ}} \left(C_5 \widehat{\theta}_l  + C_6 \widehat{\rho}_l\right) \alpha \operatorname{sg}(k)|k|
        \bigg],
    \end{aligned}
\end{equation}

\begin{equation}
    \begin{aligned}
        \left<\mathcal{Q}\left(  H^{\alpha}_{k}, H^{\beta}_{l}\right), H^{\gamma}_{m} \right>_{\mathbb{H}}
        &=
        \frac{\mathrm{i}}{2} (c_N c^{\circ})^3
        \left[\left(\frac{\rho^\circ}{c^\circ}\right)^2 \frac{p^{\circ}_{\rho}}{ \theta^{\circ} \rho^{\circ}} \left(\alpha\operatorname{sg}({k}) \frac{k}{|{k}|}+\beta \operatorname{sg}({l})\frac{{l}}{|{l}|}\right)\!\cdot\!{m}
        \right.\\
        &+\gamma \operatorname{sg}({m})|{m}|
        \bigg( \frac{\rho^\circ}{\theta^\circ} \alpha \beta\operatorname{sg}({k})\operatorname{sg}({l}) \frac{{l}\!\cdot\!{k}}{|{k}||{l}|}
            +\frac{\rho^{\circ}}{\theta^{\circ}} \left(  \frac{\rho^\circ}{c^\circ} \right)^2 C_1
            \\ &+ \frac{\rho^{\circ}}{\theta^{\circ}}\left( \frac{\theta^\circ p^{\circ}_{\theta}}{\rho^\circ e^\circ_{\theta} c^\circ} \right)^2 C_2 + \frac{\rho^\circ p^\circ_{\theta}}{e^\circ_{\theta} (c^\circ)^2}\left( C_3 + C_4 \right)
        \bigg)
         \\
            &\left.
            + \frac{p_{\theta}^{\circ}}{ c^{\circ} {\theta}^{\circ}}
            \left( C_5\frac{\theta^\circ p^{\circ}_{\theta}}{\rho^\circ e^\circ_{\theta} c^\circ}
            + C_6\frac{\rho^\circ}{c^\circ}
            \right) \left(\beta \operatorname{sg}({l})|{l}|+\alpha\operatorname{sg}({k})|{k}|\right) \right.\\
            &\left.
            + \frac{(p_{\theta}^{\circ})^2}{\rho^{\circ} (c^{\circ})^2 e_{\theta}^{\circ}}
            \left( \alpha \operatorname{sg}({k})\frac{{k}\!\cdot\!{l}}{|{k}|}
            + \beta \operatorname{sg}({l})\frac{{k}\!\cdot\!{l}}{|{l}|}\right)
        \right],\\
    \end{aligned}
\end{equation}
where $c_N$ is defined in \eqref{def_of_CN}.

Therefore, $\mathcal{Q}^{\varepsilon}_{2r}$, $\mathcal{Q}^{\varepsilon}_{3r}$ tends (in the sense of distribution) to $\mathcal{Q}_{2r}$ and $\mathcal{Q}_{3r}$ defined by
\begin{equation}\label{Q2r}
    \begin{aligned}
        \mathcal{Q}_{2r}(U,B)& = \sum_{\substack{m,\gamma}}
        \mathrm{i} m \cdot \widecheck{C}^{\gamma}_m  H^{\gamma}_{m} = \nabla \cdot \sum_{\substack{m,\gamma}} \widecheck{C}^{\gamma}_{m}H^{\gamma}_m,
    \end{aligned}
\end{equation}
and
\begin{equation}\label{Q3r}
    \begin{aligned}
        \mathcal{Q}_{3r}(A,B)
        = \sum_{\substack{m,\gamma }} \mathrm{i} m \cdot \widetilde{C}^{\gamma}_{m} H^{\gamma}_{m}= \nabla \cdot \sum_{\substack{m,\gamma}} \widetilde{C}^{\gamma}_{m}H^{\gamma}_m,
    \end{aligned}
\end{equation}
where
\begin{equation}
    \begin{aligned}
        \widecheck{C}^{\gamma}_m
        &=\sum_{\substack{\alpha\\k+l=m\\\lambda^{\alpha}_{k} = \lambda^{\gamma}_{m}}}
        \bigg[ \frac{\rho^\circ}{\theta^\circ} \left(\frac{2k \cdot m}{|k| |m|}\right)\widehat{\u}_{l}
        +\left(
            \frac{\rho^{\circ}}{\theta^{\circ}}C_1 \rho^{\circ} + \frac{\rho^{\circ}}{\theta^{\circ}}C_3 \frac{\theta^{\circ} p^{\circ}_{\theta}}{\rho^{\circ}e^{\circ}_{\theta}}
            +\frac{p^{\circ}_{\theta}}{\theta^{\circ}}C_6
        \right) \widehat{\rho}_{l} \ \gamma \operatorname{sg}(m)\frac{m}{|m|}\\
        &+\left(
            \frac{\rho^{\circ}}{\theta^{\circ}}C_4 \rho^{\circ} + \frac{\rho^{\circ}}{\theta^{\circ}}C_2 \frac{\theta^{\circ} p^{\circ}_{\theta}}{\rho^{\circ}e^{\circ}_{\theta}}
            +\frac{p^{\circ}_{\theta}}{\theta^{\circ}}(C_5-1)
        \right) \widehat{\theta}_{l} \ \gamma \operatorname{sg}(m)\frac{m}{|m|}\\
        &+\left(
            \frac{p^{\circ}_{\rho}}{\theta^{\circ}} - \frac{p^{\circ}_{\theta}}{e^{\circ}_{\theta^{\circ}}}C_3
        \right)\widehat{\rho}_{l}\  \gamma \operatorname{sg}(m)\frac{k}{|m|}
        + \left(
            \frac{p^{\circ}_{\theta}}{\theta^{\circ}} - \frac{(\rho^{\circ})^2}{\theta^{\circ}}C_4
        \right)\widehat{\theta}_{l}\  \gamma \operatorname{sg}(m)\frac{k}{|m|}
        \bigg] (c_N c^{\circ})^2  B^{\alpha}_{ {k}},
    \end{aligned}
\end{equation}
and
\begin{equation}
    \begin{aligned}
        \widetilde{C}^{\gamma}_{m}&=
        \sum_{\substack{\alpha\\k+l=m\\ \lambda^{\alpha}_{k} + \lambda^{\beta}_{l} = \lambda^{\gamma}_{m}}}
        \left[2\left(\frac{\rho^\circ}{c^\circ}\right)^2 \frac{p^{\circ}_{\rho}}{\rho^{\circ} \theta^{\circ}}
        + \frac{(p_{\theta}^{\circ})^2}{\rho^{\circ} (c^\circ)^2 e_{\theta}^{\circ}}
        + \frac{\rho^{\circ}}{\theta^{\circ}}
        \right.\\
        &\left.
        + \frac{\rho^\circ}{\theta^\circ}
        \left(
            \left( \frac{\rho^\circ}{c^\circ} \right)^2 C_1 +\left( \frac{\theta^\circ p^{\circ}_{\theta}}{\rho^\circ e^\circ_{\theta} c^\circ} \right)^2 C_2 + \frac{\theta^\circ p^\circ_{\theta}}{e^\circ_{\theta} (c^\circ)^2}\left( C_3 + C_4 \right)
        \right)
         \right.\\
        &\left.
            +\frac{p_{\theta}^{\circ}}{c^{\circ} \theta^{\circ}}\left( C_5\left(\frac{\theta^\circ p^{\circ}_{\theta}}{\rho^\circ e^\circ_{\theta} c^\circ}\right)
            + C_6 \frac{\rho^{\circ}}{c^{\circ}}
            \right)
        \right]
        \left( c_N c^{\circ} \right)^{3} A^{\alpha}_{k} B^{\beta}_{{l}} \gamma \operatorname{sg}(m) \frac{m}{|m|} .
    \end{aligned}
\end{equation}

The operaotrs $Q_{2r}$ and $Q_{3r}$ are actually the averaged quadratic operator over the two-wave and three wave resonant sets respectively. We refer to \cites{Danchin_AJM,Masmoudi_Poincare,Jiang_AA}
for more details. Namely,
In the computation of $\mathcal{Q}_{2r}$, we need set $\lambda^{\alpha}_{k} = \lambda^{\gamma}_{m}$, i.e.,
\begin{equation}
    \alpha \operatorname{sg}(k)|k| =\gamma \operatorname{sg}(m)|m|,
\end{equation}
which implies that
\begin{equation}
    |k|  = |m|\text{, } \alpha \operatorname{sg}(k) = \gamma \operatorname{sg}(m).
\end{equation}
Similarly, in computation of $\mathcal{Q}_{3r}$, we need to set $\lambda^{\alpha}_k + \lambda^{\beta}_l = \lambda^{\gamma}_m$, i.e.,
\begin{equation}
    \alpha \operatorname{sg}(k)|k| + \beta\operatorname{sg}(l)|l|=\gamma \operatorname{sg}(m)|m|,
\end{equation}
which, combining with $k+l=m$, implies that
\begin{equation}
    \alpha = \beta = \gamma, \ k \cdot l = \operatorname{sg}(k) \operatorname{sg}(l)|k||l|,
\end{equation}
and therefore $k,\ l,\ m$ are colinear.

Moreover, we have
\begin{equation}
    \begin{aligned}
    \mathcal{D}^{\varepsilon} A &= \mathrm{e}^{\frac{t}{\varepsilon} \mathcal{A}} \mathcal{D} \mathrm{e}^{-\frac{t}{\varepsilon} \mathcal{A}} A
    = \sum_{m,\gamma}\left< \mathcal{D} \mathrm{e}^{-\frac{t}{\varepsilon} \mathcal{A}} A, H^{\gamma}_{m} \right>_{\mathbb{H}} \mathrm{e}^{\frac{t}{\varepsilon}\lambda^{\gamma}_m }
    = \sum_{m,\gamma} \sum_{\substack{\alpha \\ k = m }} \left< \mathcal{D} H^{\alpha}_{k} , H^{\gamma}_{m} \right>_{\mathbb{H}} A^{\alpha}_{k} H^{\gamma}_{m} \mathrm{e}^{\frac{t}{\varepsilon}(\lambda^{\gamma}_m - \lambda^{\alpha}_{k})},
    \end{aligned}
\end{equation}
with
\begin{equation}
    \begin{aligned}
    \left< \mathcal{D} H^{\alpha}_{k} , H^{\gamma}_{m} \right>_{\mathbb{H}}
    =-\sum_{\gamma, m}\ \sum_{\substack{\alpha\\ k=m}}\left(
        \frac{\frac{2N-2}{N} \mu^{\circ} + \lambda^{\circ}}{\theta^{\circ}}\alpha \gamma  |k|^2
        + \kappa^{\circ} \left(\frac{ {p^{\circ}_{\theta}}}{{\rho^{\circ}} {c^{\circ}} e^{\circ}_{\theta}}\right)^2 |k|^{2}
    \right)
    c_N c^{\circ}.
    \end{aligned}
\end{equation}
By setting $\lambda^{\alpha}_k = \lambda^{\gamma}_m$, we know that $\mathcal{D}^{\varepsilon}$ tends to
\begin{equation}\label{D_bar}
    \begin{aligned}
    \overline{\mathcal{D}} A &= -\sum_{\gamma, m}
    \left(
        \frac{\frac{2N-2}{N} \mu^{\circ} + \lambda^{\circ}}{\theta^{\circ}}
        + \kappa^{\circ} \left(\frac{ {p^{\circ}_{\theta}}}{{\rho^{\circ}} {c^{\circ}} e^{\circ}_{\theta}}\right)^2
    \right)
     (c_N c^{\circ})^2|m|^{2} A^{\gamma}_{m} H^{\gamma}_{m}\\
    &= \bar{\mu} \Delta A,
    \end{aligned}
\end{equation}
with \begin{equation*}
    \bar{\mu} =
    \left(
        \frac{\frac{2N-2}{N} \mu^{\circ} + \lambda^{\circ}}{\theta^{\circ}}
        + \kappa^{\circ} \left(\frac{ {p^{\circ}_{\theta}}}{{\rho^{\circ}} {c^{\circ}} e^{\circ}_{\theta}}\right)^2
    \right)
     (c_N c^{\circ})^2.
\end{equation*}

To prove proposition \ref{oscillation_convergence}, set $z^\varepsilon=\widetilde{V}^\varepsilon-V$. Subtracting \eqref{LS} from \eqref{nsc_orthogonal} yields that
\begin{equation} \label{Z_equation}
    \begin{aligned}
    \partial_t z^\varepsilon &+ \mathcal{Q}^{\varepsilon}_{3r}(z^\varepsilon,z^\varepsilon)
                            +2\mathcal{Q}^{\varepsilon}_{3r}(z^\varepsilon,V)
                            +\mathcal{Q}^{\varepsilon}_{2r}(z^\varepsilon,\mathcal{U})
                            -\overline{\mathcal{D}}z^\varepsilon \\
                            &=R^{1,\varepsilon}_c + R^{2,\varepsilon}_c + R^{3,\varepsilon}_c +R^{4,\varepsilon}_c
                            +S^\varepsilon - 2\mathrm{e}^{\frac{t}{\varepsilon}\mathcal{A}} \mathcal{P}^{\perp}\mathcal{Q}(\mathcal{P}_0 \widetilde{U}^{\varepsilon}, \widetilde{U}^{\varepsilon})
                            -\mathcal{Q}^{\varepsilon}_{2r}(w^\varepsilon,\widetilde{V}^\varepsilon)\\
                            &-\mathrm{e}^{\frac{t}{\varepsilon}\mathcal{A}} \mathcal{P}^{\perp}\mathcal{Q}(w^\varepsilon,\mathcal{P}\widetilde{U}^\varepsilon)
                            -\mathrm{e}^{\frac{t}{\varepsilon}\mathcal{A}} \mathcal{P}^{\perp}\mathcal{Q}(w^\varepsilon,\mathcal{U})
                            +\mathrm{e}^{\frac{t}{\varepsilon}\mathcal{A}}r^{\varepsilon},
    \end{aligned}
\end{equation}
where
\begin{equation}\label{definition_of_R}
    \begin{aligned}
        &R^{1,\varepsilon}_c= -\mathrm{e}^{\frac{t}{\varepsilon}\mathcal{A}} \mathcal{P}^{\perp} \mathcal{Q}(\mathcal{U},\mathcal{U}),\
        R^{2,\varepsilon}_c = (\mathcal{Q}_{2r}-\mathcal{Q}_{2r}^\varepsilon)(\mathcal{U},V), \\
        &R^{3,\varepsilon}_c = (\mathcal{Q}_{3r}-\mathcal{Q}_{3r}^\varepsilon)(V,V), \
        R^{4,\varepsilon}_c =\mathrm{e}^{\frac{t}{\varepsilon}\mathcal{A}} \mathcal{P}^{\perp} \mathcal{D} \mathcal{P} \widetilde{U}^\varepsilon ,\ \\
        &S^\varepsilon= (\mathcal{D}^\varepsilon -  \overline{\mathcal{D}}) \widetilde{V}^\varepsilon,
    \end{aligned}
\end{equation}
and
\begin{equation*}
    r^{\varepsilon} = \left(0,\ r^{\varepsilon}_2,\ r^{\varepsilon}_3\right)^{\text{T}},
\end{equation*}
with
\begin{equation}
    \begin{aligned}
        r^{\varepsilon}_2 &=
            -\left(
            {\frac{1}{\varepsilon^{2}}} \frac{\nabla p^{\varepsilon}}{\rho^{\varepsilon}}
            - \frac{1}{\varepsilon}\frac{p^{\circ}_{\theta}}{\rho^{\circ}}\nabla \tilde{\theta}^{\varepsilon}
            -\frac{1}{\varepsilon}\frac{p^{\circ}_{\rho}}{\rho^{\circ}}\nabla \tilde{\rho}^{\varepsilon}
            -C_2 \tilde{\rho}^{\varepsilon} \nabla \tilde{\rho}^{\varepsilon}
            -C_2 \tilde{\theta}^{\varepsilon} \nabla \tilde{\theta}^{\varepsilon}
            -C_3 \tilde{\theta}^{\varepsilon} \nabla \tilde{\rho}^{\varepsilon}
            -C_4 \tilde{\rho}^{\varepsilon} \nabla \tilde{\theta}^{\varepsilon}
            \right)\\
            &+\nabla \!\cdot\! S^{\varepsilon} - \left( \frac{\mu^{\circ}}{\rho^{\circ}} \Delta \tilde{\u}^{\varepsilon} + \frac{\frac{N-2}{N}\mu^{\circ}+\lambda^{\circ}}{\rho^{\circ}}\nabla \operatorname{div} \tilde{\u}^{\varepsilon} \right),
         \end{aligned}
        \end{equation}

\begin{equation}
     \begin{aligned}
        r^{\varepsilon}_3 &=  -\frac{\kappa^{\circ}}{e^{\circ}_{\theta} \rho^{\circ}} \Delta  \tilde{\theta}^{\varepsilon} + \frac{1}{e_{\theta}^{\varepsilon} \rho^{\varepsilon}} \nabla( \kappa^{\varepsilon} \nabla  \tilde{\theta}^{\varepsilon})
        + \frac{1}{\varepsilon} \left( \frac{\theta^{\circ} p^{\circ}_{\theta}}{e^{\circ}_{\theta} \rho^{\circ}} - \frac{\theta^{\varepsilon} p_{\theta}^{\varepsilon}}{e_{\theta}^{\varepsilon} \rho^{\varepsilon}} + C_5 \tilde{\theta}^{\varepsilon} + C_6 \tilde{\rho}^{\varepsilon} \right) \nabla \!\cdot\! \tilde{\u}^{\varepsilon}
        + \frac{\varepsilon}{e_{\theta}^{\varepsilon} \rho^{\varepsilon }} S^{\varepsilon}\!:\! \nabla  \tilde{\u}^{\varepsilon}.\\
    \end{aligned}
\end{equation}
Moreover, we denote $R^\varepsilon_c=R^{1,\varepsilon}_c+R^{2,\varepsilon}_c+R^{3,\varepsilon}_c+R^{4,\varepsilon}_c + S^\varepsilon$.

It's impossible to expect $R^{\varepsilon}_c$ to converge to zero in any space $\widetilde{L}^{r}(H^{\sigma})$, which actually oscillates in time. In order to get the convergence result, we need to use the so-called negative time regularity method and use the gap between $\frac{N}{2}- \delta$ and $\frac{N}{2}$.
We shall follow as in \cite{Danchin_ARMA}, and one can refer to \cite{Schochet_JDE} for more details.

Set
\begin{equation}\label{z_decomp}
    \phi^\varepsilon_M = z^\varepsilon - \varepsilon \widetilde{R}^{\varepsilon}_{c,M},
\end{equation}

Here and in the sequel we shall denote $R^{\varepsilon}_{c,M}$ the “low frequency part" of $R^{\varepsilon}_{c}$
and $R^{\varepsilon,M}_c = R^{\varepsilon}_{c} - R^{\varepsilon}_{c,M}$ the “high frequency part" of $R^{\varepsilon}_{c}$.

We choose $\widetilde{R}^{\varepsilon}_{c,M}$ such that

\begin{equation} \label{R_decomp}
    \partial_t \varepsilon \widetilde{R}^\varepsilon_{c,M} = R^{\varepsilon}_{c,M} + \varepsilon \widetilde{R}^{t,\varepsilon}_{c,M}.
\end{equation}
Note that \eqref{R_decomp} is nothing but chain rule of time derivative, and $\widetilde{R}^{\varepsilon}_{c,M}$ can be found as the time “primitive function".
For example we have
\begin{equation}
    \begin{aligned}
        S^{\varepsilon }_{M}
        = \sum_{\substack{\gamma \\|m| \leqslant M}}
         \left< \mathcal{D} H^{-\gamma}_{m} , H^{\gamma}_{m} \right>_{\mathbb{H}} A^{-\gamma}_{m} H^{\gamma}_{m} \mathrm{e}^{2\frac{t}{\varepsilon}\lambda^{\gamma}_m},
    \end{aligned}
\end{equation}
and
\begin{equation}
    \begin{aligned}
        \widetilde{S}^{\varepsilon }_{M}
        = \sum_{\substack{\gamma \\|m| \leqslant M}}
        \frac{1}{2\lambda^{\gamma}_{m}}\left< \mathcal{D} H^{-\gamma}_{m} , H^{\gamma}_{m} \right>_{\mathbb{H}} A^{-\gamma}_{m} H^{\gamma}_{m} \mathrm{e}^{2\frac{t}{\varepsilon}\lambda^{\gamma}_m},
        \\
        \widetilde{S}^{\varepsilon,t }_{M}
        = \sum_{\substack{\gamma \\|m| \leqslant M}}
        \frac{1}{2\lambda^{\gamma}_{m}}\left< \mathcal{D} H^{-\gamma}_{m} , H^{\gamma}_{m} \right>_{\mathbb{H}} \partial_t A^{-\gamma}_{m} H^{\gamma}_{m} \mathrm{e}^{2\frac{t}{\varepsilon}\lambda^{\gamma}_m}.
    \end{aligned}
\end{equation}
Similarly, we can choose proper $\widetilde{R}^{\varepsilon,1}_{c,M},\ \widetilde{R}^{\varepsilon,2},\ \widetilde{R}^{\varepsilon,3}_{c,M},\ \widetilde{R}^{\varepsilon,4}_{c,M}$ which will be given later in the proof of Lemma \ref{XR_estimates}.

Combining \eqref{Z_equation}, \eqref{z_decomp} and \eqref{R_decomp} we then have
\begin{equation}
    \begin{aligned}
        \partial_t \phi^\varepsilon_M + \mathcal{Q}^{\varepsilon}_{2r}(\phi^\varepsilon_M,\mathcal{U})
                                      &+ 2 \mathcal{Q}^{\varepsilon}_{3r}(\phi^\varepsilon_M, V)
                                      - \overline{\mathcal{D}} \phi^\varepsilon_M
                                      = R^{\varepsilon,M}_c
                                      - \mathcal{Q}^{\varepsilon}_{2r}(\widetilde{V}^\varepsilon,w^\varepsilon)
                                      - \mathcal{Q}^{\varepsilon}_{3r}(z^\varepsilon,z^\varepsilon)\\
                                      &+ \varepsilon (
                                        \widetilde{R}^{t,\varepsilon}_M
                                        - \mathcal{Q}^{\varepsilon}_{2r}(\widetilde{R}^\varepsilon_{c,M},\mathcal{U})
                                        - \mathcal{Q}^{\varepsilon}_{3r}(\widetilde{R}^\varepsilon_{c,M},V)
                                        + \overline{\mathcal{D}}\widetilde{R}^\varepsilon_{c,M}
                                      )\\
                                      &+ \mathrm{e}^{\frac{t}{\varepsilon}\mathcal{A}} \mathcal{P}^{\perp}(
                                        r^{\varepsilon}
                                        - \mathcal{Q}(w^{\varepsilon},\mathcal{\mathcal{P}} \widetilde{U}^\varepsilon)
                                        -\mathcal{Q}(w^{\varepsilon}, \mathcal{U}) - 2\mathcal{Q}(\mathcal{P}_0 \widetilde{U}^{\varepsilon}, \widetilde{U}^{\varepsilon})
                                      ).
    \end{aligned}
\end{equation}
Note that $\phi^\varepsilon_M$ indeed satisfies a nonlinear heat equation and therefore
using Lemma \ref{structure_lemma}, we have
\begin{equation} \label{phi_estimate}
    \begin{aligned}
        \left\| \phi^\varepsilon_M \right\|_{\widetilde{L}^\infty(H^{s-1})} + &\left\| \phi^\varepsilon_M \right\|_{{L}^2(H^{s})} \\
                            &\leqslant C\mathrm{ e}^{ C \int^{\infty}_0 \left\| \mathcal{U} \right\|^2_{B^\frac{N}{2}}+\left\| V \right\|^2_{B^\frac{N}{2}}}
                            \bigg(
                                \varepsilon \left\| \widetilde{R}^\varepsilon_{c,M}(0) \right\|_{H^{s-1}}
                                + \left\| R^{\varepsilon,M}_{c} \right\|_{\widetilde{L}^2 (H^{s-2})}
                                + \left\| \mathcal{Q}^\varepsilon_{3r}(z_\varepsilon,z_\varepsilon) \right\|_{\widetilde{L}^2 (H^{s-2})} \\
                            &   + \varepsilon \left\| \mathcal{Q}^\varepsilon_{3r}(\widetilde{R}^\varepsilon_{c,M},V) \right\|_{\widetilde{L}^2 (H^{s-2})}
                                + \varepsilon \left\| \mathcal{Q}^\varepsilon_{2r}(\widetilde{R}^\varepsilon_{c,M},\mathcal{U}) \right\|_{\widetilde{L}^2 (H^{s-2})}
                                + \varepsilon \left\| \overline{\mathcal{D}} \widetilde{R}^\varepsilon_{c,M} \right\|_{\widetilde{L}^2 (H^{s-2})}\\
                            &    + \varepsilon \left\|  \widetilde{R}^{t,\varepsilon}_M \right\|_{\widetilde{L}^2 (H^{s-2})}
                                + \left\| \mathcal{Q}^\varepsilon_{2r}(\widetilde{V}^\varepsilon,w^\varepsilon) \right\|_{\widetilde{L}^2 (H^{s-2})}
                                + \left\| \mathrm{e}^{\frac{t}{\varepsilon}\mathcal{A}} \mathcal{P}^{\perp}r^\varepsilon \right\|_{\widetilde{L}^2 (H^{s-2})} \\
                            &    + \left\| \mathrm{e}^{\frac{t}{\varepsilon}\mathcal{A}} \mathcal{P}^{\perp} \mathcal{Q}(w^\varepsilon,\mathcal{P}\widetilde{U}^\varepsilon) \right\|_{\widetilde{L}^2 (H^{s-2})}
                                + \left\| \mathrm{e}^{\frac{t}{\varepsilon}\mathcal{A}} \mathcal{P}^{\perp} \mathcal{Q}(w^\varepsilon,\mathcal{U}) \right\|_{\widetilde{L}^2 (H^{s-2})} \\
                             &   +2\left\|\mathrm{e}^{\frac{t}{\varepsilon}\mathcal{A}} \mathcal{P}^{\perp}\mathcal{Q}(\mathcal{P}_0 \widetilde{U}^{\varepsilon}, \widetilde{U}^{\varepsilon})\right\|_{\widetilde{L}^2 (H^{s-2})}
                            \bigg).
    \end{aligned}
\end{equation}

In view of \eqref{z_decomp}, to get the smallness of $z^\varepsilon$, we need to obtain the smallness of $\phi^{\varepsilon}_{M}$ and boundedness of $\widetilde{R}^{\varepsilon}_{c,M}$.
For the right-hand side of above inequality, we can expect the terms involving $w^{\varepsilon}$ small since we shall later prove the convergence of Incompressible part. The terms with coefficient $\varepsilon$ are expected to be bounded which coincide with
the boundedness of $\widetilde{R}^{\varepsilon}_{c,M}$. The high frequency part $R^{\varepsilon,M}_{c}$ is small using regularity gap.

Now we give a more detailed description.
Using Lemma \ref{Qe_estimate} and the Lemma \ref{embedding} we have
\begin{equation}\label{source_estimate}
    \begin{aligned}
        &\left\| \mathrm{e}^{\frac{t}{\varepsilon}\mathcal{A}} \mathcal{P}^{\perp} \mathcal{Q}(w^\varepsilon,\mathcal{P}\widetilde{U}^\varepsilon) \right\|_{{L}^{2}(H^{s-2})}
            \lesssim \left \| \mathcal{P} \widetilde{U}^\varepsilon \right \|_{\widetilde{L}^2(B^{\frac{N}{2}})}
                  \left \| w^\varepsilon \right \|_{\widetilde{L}^{\infty}(H^{s-1})} ,\\
        &\left\| \mathrm{e}^{\frac{t}{\varepsilon}\mathcal{A}} \mathcal{P}^{\perp} \mathcal{Q}(w^\varepsilon,\mathcal{U}) \right\|_{{L}^{2}(H^{s-2})}
            \lesssim \left \|  w^\varepsilon \right \|_{{L}^2(H^{s})}
                  \left \| \mathcal{U} \right \|_{\widetilde{L}^{\infty}(H^{\frac{N}{2}-1})} ,\\
        &\left\| \mathrm{e}^{\frac{t}{\varepsilon}\mathcal{A}} \mathcal{P}^{\perp} \mathcal{Q}(\mathcal{P}_0 \widetilde{U}^{\varepsilon},\widetilde{U}^{\varepsilon}) \right\|_{L^2({H}^{s-2})}
            \lesssim \left\| \mathcal{P}_0 \widetilde{U}^{\varepsilon} \right\|_{L^{\infty}(dt)} \left\| \widetilde{U} \right\|_{L^{2}(\dot{H}^{\frac{N}{2}})} \\
        &\left\| \mathcal{Q}^\varepsilon_{2r}(\widetilde{V}^\varepsilon,w^\varepsilon) \right\|_{{L}^{2}(H^{s-2})}
            \lesssim \left \| w^\varepsilon \right \|_{{L}^{2}(H^{s})}
                  \left \| \widetilde{V}^\varepsilon \right \|_{\widetilde{L}^2(B^{\frac{N}{2}})} ,\\
        &\left\| \mathcal{Q}^\varepsilon_{3r}(z^\varepsilon,z^\varepsilon) \right\|_{{L}^{2}(H^{s-2})}
            \lesssim (\left \| \widetilde{V}^\varepsilon \right \|_{\widetilde{L}^2(B^{\frac{N}{2}})}+\left \| V \right \|_{\widetilde{L}^2(B^{\frac{N}{2}})})
                  \left \| z^\varepsilon \right \|_{\widetilde{L}^{\infty}(H^{s-1})}, \\
        &\left\|\mathcal{Q}^\varepsilon_{3r}(\widetilde{R}^\varepsilon_{c,M},V)\right\|_{{L}^{2}(H^{s-2})}
            \lesssim \left\|\widetilde{R}^\varepsilon_{c,M}\right\|_{\widetilde{L}^{\infty}(H^{s-1})}
                  \left\|V\right\|_{\widetilde{L}^2(B^{\frac{N}{2}})},\\
        &\left\|\mathcal{Q}^\varepsilon_{2r}(\widetilde{R}^\varepsilon_{c,M},\mathcal{U})\right\|_{{L}^{2}(H^{s-2})}
            \lesssim \left\|\widetilde{R}^\varepsilon_{c,M}\right\|_{\widetilde{L}^{\infty}(H^{s-1})}
                   \left\|\mathcal{U}\right\|_{\widetilde{L}^2(B^{\frac{N}{2}})}.
    \end{aligned}
\end{equation}
Note that $\overline{D}$ is just Laplacian, so we have
\begin{equation}\label{Laplacian_estimate}
    \left\|\overline{\mathcal{D}}\widetilde{R}^\varepsilon_{c,M}\right\|_{{L}^{2}(H^{s-2})}
            \lesssim \left\|\widetilde{R}^\varepsilon_{c,M}\right\|_{{L}^{2}(H^{s})}.
\end{equation}
Now combining \eqref{phi_estimate}, \eqref{source_estimate} and \eqref{Laplacian_estimate} we have
\begin{equation} \label{phi_full_estimate}
    \begin{aligned}
        \| \phi^\varepsilon_M& \|_{\widetilde{L}^\infty(H^{s-1})} + \left\| \phi^\varepsilon_M \right\|_{{L}^2(H^{s})} \\
        &\leqslant C\mathrm{ e}^{C\int^t_0 \left\| \mathcal{U} \right\|^2_{B^\frac{N}{2}}+\left\| V \right\|^2_{B^\frac{N}{2}}}
        \bigg(
            \varepsilon \left\| \widetilde{R}^\varepsilon_{c,M}(0) \right\|_{H^{s-1}}
            + \left\| R^{\varepsilon,M}_{c} \right\|_{\widetilde{L}^2 (H^{s-2})}
            + \varepsilon \left\| R^{t,\varepsilon}_{c,M} \right\|_{\widetilde{L}^2 (H^{s-2})}
        \\ &    + \varepsilon (1+ X +Y)\left\| \widetilde{R}^{\varepsilon}_{c,M} \right\|_{{L}^{2}(H^s) \cap \widetilde{L}^{\infty}(H^{s-1})}
            + XW + YW + \varepsilon(X+Y) + \left\| \mathcal{P}^{\perp} r^{\varepsilon} \right\|_{\widetilde{L}^{2}({H}^{s-2})} \bigg).
    \end{aligned}
\end{equation}
What left now is the estimates for $r^{\varepsilon}$,  $\widetilde{R}^\varepsilon_{c,M}$, $\widetilde{R}^{\varepsilon,M}_{c}$ and $\widetilde{R}^{t,\varepsilon}_{c,M}$.
\begin{lemma} \label{r_estimate}
    We have the following estimate for $r_\varepsilon$:
    \begin{equation}
        \begin{aligned}
            \left\| r_{\varepsilon} \right\|_{\widetilde{L}^2_{T}(\dot{H}^{s-2})}
            \lesssim
            C(\mathcal{Y}^{\varepsilon}) \varepsilon^{\delta}
            \left((X+Y)^2 + (X + Y)^3\right),
        \end{aligned}
    \end{equation}
    where $C(\mathcal{Y}^{\varepsilon}) $ is a polynomial in $\mathcal{Y}^{\varepsilon}$.
\end{lemma}
\begin{lemma} \label{XR_estimates}
    We have the following estimate for $\widetilde{R}^{\varepsilon}_{c,M}$, $\widetilde{R}^{\varepsilon,t}_{c,M}$, and ${R}^{\varepsilon,M}_{c}$
    \begin{equation}
        \begin{aligned}
            &\left\| \widetilde{R}^{\varepsilon}_{c,M} \right\|_{\widetilde{L}^{\infty}(H^{s-1})\cap {L}^2(H^{s})}
            \lesssim
            M^{1-\delta}\left( Y^{2} + X + Y \right) + C^1_{M} M^{2-\delta} X Y + C^2_{M} M^{2-\delta} X^{2}, \\
            & \left\| {R}^{\varepsilon,M}_{c} \right\|_{{L}^{2}(H^{s-2})}
            \lesssim
            M^{-\delta}\left(X^2 + Y^2 + X Y +Y\right),
        \end{aligned}
    \end{equation}
    \begin{equation}
        \begin{aligned}
            \| &\widetilde{R}^{\varepsilon,t}_{c,M} \|_{{L}^2(H^{s-2})}
            \lesssim \left( M^{1-\delta} Y + C_{M}^{1}M^{2-\delta} X \right) \left\| \partial_t \mathcal{U} \right\|_{{L}^{2}(H^{\frac{N}{2}-2})}\\
                &+ \left( C^{1}_{M} M^{2-\delta} Y + C^{2}_{M} M^{2-\delta} \right) \left\| \partial_t V \right\|_{{L}^{2}(H^{\frac{N}{2}-2})}
                + M^{1-\delta} \left( \left\| \partial_t \mathcal{P} \widetilde{U}^{\varepsilon} \right\|_{\widetilde{L}^{2}(H^{\frac{N}{2}})} + \left\|  \partial_t \widetilde{V}^{\varepsilon} \right\|_{\widetilde{L}^{2}(H^{\frac{N}{2}})} \right),
           \end{aligned}
    \end{equation}
    with $C_{M}^1, \ C_{M}^2$ polynomials in $M$,
    and in addition
    \begin{equation}
        \begin{aligned}
            \left\| \widetilde{{R}}^{\varepsilon}_{c,M}(0) \right\|_{H^{s-2}}
           &\lesssim
           M^{1-\delta}\left( Y_{0}^{2} + X_{0} + Y_{0} \right) + C^1_{M} M^{2-\delta} X_{0} Y_{0} + C^2_{M} M^{2-\delta} X_{0}^{2}, \\
        \end{aligned}
    \end{equation}

\end{lemma}
Therefore, we further need the estimates for time derivative of $\mathcal{U}$, $V$ and $\widetilde{U}^{\varepsilon}$. Indeed, we have the following estimates:
\begin{lemma}\label{time_derivative_estimate}
\begin{equation}
    \begin{aligned}
        &\left\| \partial_t \mathcal{U} \right\|_{{L}^2(H^{\frac{N}{2}-2})}
        \lesssim
        Y^2+Y,
        \\
        &\left\| \partial_t V \right\|_{{L}^2(H^{\frac{N}{2}-2})}
        \lesssim
        X(1+Y+X),
        \\
        &\left\| \partial_t \mathcal{P}\widetilde{U}^{\varepsilon} \right\|_{{L}^2(H^{\frac{N}{2}-2})}
        \lesssim
        X + XY + X^2 + Y^2 +C(\mathcal{Y}^\varepsilon)\left( (X+Y)^2+(X+Y)^3\right),\\
        \\
        &\left\| \partial_t \widetilde{V}^{\varepsilon} \right\|_{{L}^2(H^{\frac{N}{2}-2})}
        \lesssim
        X+Y+(X+Y)^2+C(\mathcal{Y}^\varepsilon)\left( (X+Y)^2+(X+Y)^3\right).
    \end{aligned}
\end{equation}
\end{lemma}

We can now complete the proof of Proposition \ref{oscillation_convergence}.
\begin{proof}[Proof of Proposition \ref{oscillation_convergence}]

Note that
\begin{equation}
    \left\| z^{\varepsilon} \right\|_{{L}^2(H^{\frac{N}{2}} \cap L^{\infty}) \cap \widetilde{L}^{\infty}(H^{\frac{N}{2}-1})}
    \lesssim
    \varepsilon \left\| \widetilde{R}^{\varepsilon}_{c,M} \right\|_{{L}^2(H^{\frac{N}{2}} \cap L^{\infty}) \cap \widetilde{L}^{\infty}(H^{\frac{N}{2}-1})}
    + \left\| \phi^{\varepsilon}_{M} \right\|_{{L}^2(H^{\frac{N}{2}} \cap L^{\infty}) \cap \widetilde{L}^{\infty}(H^{\frac{N}{2}-1})},
\end{equation}
combining with \eqref{phi_full_estimate}  Lemma \ref{r_estimate}, Lemma \ref{time_derivative_estimate} and using the fact $\mathcal{Y}^{\varepsilon} $ is uniformly bounded, we have
\begin{equation}
    \begin{aligned}
    \left\| z^{\varepsilon} \right\|&_{{L}^2(H^{\frac{N}{2}} \cap L^{\infty}) \cap \widetilde{L}^{\infty}(H^{\frac{N}{2}-1})}
    \leqslant C \mathrm{e}^{C(Y^2 + X^2)} \times \\
    &\bigg(\varepsilon M^{1-\delta} \left( X + Y + Y^2 + XY^2 +Y^3 +( X + Y )^2 + (X + Y)^3 \right)  \\
    &+ \varepsilon C_M M^{2-\delta}(X + X^2 + XY +X^{2}Y + X Y^{2} + X^3) \\
    &+(Y+X)W^{\varepsilon} + X Z^{\varepsilon} + {(Z^{\varepsilon})}^2   + \varepsilon M^{1-\delta}\left(X_0 + Y_0 + {Y_0}^2 \right) \\
    &+ \varepsilon C_M M^{2-\delta} \left( X_0 Y_0 + {X_0}^2\right) + M^{-\delta}(Y + XY + Y^2 + X^2)\bigg),
\end{aligned}
\end{equation}
where
\begin{equation}
    C_M = \max \left( C_{M}^1 ,\ C_{M}^{2}\right).
\end{equation}
\end{proof}
\begin{remark}
    We can actually give an explicit form of $C_M$ using some refined small divisor estimates, which plays an essential role in finding the convergence rate of the oscillating part. We shall postpone this to next section.
\end{remark}
\begin{proof}[Proof of Lemma \ref{r_estimate}]

    For $r^{\varepsilon}_2$ note that $\| \mathcal{P}^{\perp} r_\varepsilon\|_{L^2({H^{s-2}})} \lesssim \| r_\varepsilon \|_{L^2(\dot{H}^{s-2})}$ and note that
\begin{equation}
    \begin{aligned}
         \frac{1}{{\rho}^{\varepsilon}}\nabla \!\cdot\! S^{\varepsilon} - \left( \frac{\mu^{\circ}}{\rho^{\circ}} \Delta \tilde{\u}^{\varepsilon} + \frac{\frac{N-2}{N}\mu^{\circ}+\lambda^{\circ}}{\rho^{\circ}}\nabla \operatorname{div} \tilde{\u}^{\varepsilon} \right)
        = \frac{1}{{\rho}^{\varepsilon}} \nabla\left[ \left( \mu^{\varepsilon} - \mu^{\circ}\right) \nabla \widetilde{\u}^{\varepsilon} \right] + \left(\frac{1}{\rho^{\varepsilon}} - \frac{1}{\rho^{\circ}}\right) \mu^{\circ}\Delta \widetilde{\u}^{\varepsilon}\\
        + \left(\frac{1}{\rho^{\varepsilon}} - \frac{1}{\rho^{\circ}}\right)\left({\frac{N-2}{N}\mu^{\circ}+\lambda^{\circ}}\right)\nabla \operatorname{div} \tilde{\u}^{\varepsilon}
        + \frac{1}{{\rho}^{\varepsilon}} \nabla\left(\left({\frac{N-2}{N} (\mu^{\varepsilon} - \mu^{\circ})+ \lambda^{\varepsilon} - \lambda^{\circ}}\right) \operatorname{div} \tilde{\u}^{\varepsilon} \right).
    \end{aligned}
\end{equation}
Using Lemma \ref{composition_estimates} and \eqref{delta_estimates} we have
\begin{equation}
    \begin{aligned}
        \left\| \frac{1}{\rho^{\circ}} - \frac{1}{{\rho}^{\varepsilon}} \right\|_{H^{s}} + \left\| \mu^{\varepsilon} - \mu^{\circ} \right\|_{H^{s}} + \left\| \lambda^{\varepsilon} - \lambda^{\circ} \right\|_{H^{s}}
        \leqslant C(\mathcal{Y}^{\varepsilon}) \left\| \left( \varepsilon  \widetilde{\rho}^{\varepsilon} ,  \varepsilon \widetilde{\theta}^{\varepsilon} \right) \right\|_{H^{s}}
        \lesssim
        C(\mathcal{Y}^{\varepsilon}) \varepsilon^{\delta} \left\| \left(   \widetilde{{\rho}}^{\varepsilon} ,  \ \widetilde{\theta}^{\varepsilon} \right) \right\|_{B^{\frac{N}{2}-1, \frac{N}{2}}_{\varepsilon \nu^{\circ} }}.
    \end{aligned}
\end{equation}
Note also that $\frac{1}{\rho^{\varepsilon}}$ is bounded by $ \mathcal{Y}^{\varepsilon} $ , then using Lemma \ref{embedding}, Lemma \ref{product_estimate} and Lemma \ref{composition_estimates}, it's easy to get
\begin{equation}
    \begin{aligned}
    \left\| \frac{1}{{\rho}^{\varepsilon}} \nabla\left[ \left( \mu^{\varepsilon} - \mu^{\circ}\right) \nabla \widetilde{u}^{\varepsilon} \right] \right\|_{{L}^2(H^{s-2})}
    &= \left\|\frac{1}{{\rho}^{\varepsilon}} \nabla\left[ \left( \mu^{\varepsilon} - \mu^{\circ}\right) \nabla \widetilde{u}^{\varepsilon} \right] \right\|_{L^2(H^{s-2})}
    \lesssim \left\| \frac{1}{\rho^{\varepsilon}} \right\|_{L^{\infty}(L^{\infty})}
    \left\| \left( \mu^{\varepsilon} - \mu^{\circ}\right) \nabla \widetilde{u}^{\varepsilon}  \right\|_{L^2(H^{s-1})} \\
    &\lesssim
    \left\| \frac{1}{\rho^{\varepsilon}} \right\|_{L^{\infty}(L^{\infty})}
     \left\|  \mu^{\varepsilon} - \mu^{\circ} \right\|_{L^{\infty}(H^{s})}
    \left\| \widetilde{\u}^{\varepsilon} \right\|_{L^{2}(\dot{H}^{\frac{N}{2}} \cap L^{\infty})}
    \\ &\lesssim
    C(\mathcal{Y}^{\varepsilon}) \varepsilon^{\delta} \left\| \left(   \widetilde{{\rho}}^{\varepsilon} ,  \ \widetilde{\theta}^{\varepsilon} \right) \right\|_{\widetilde{L}^{\infty}(B^{\frac{N}{2}-1, \frac{N}{2}}_{\varepsilon \nu^{\circ} })}
    \left\| \widetilde{\u}^{\varepsilon} \right\|_{L^{2}(\dot{B}^{\frac{N}{2}})},
    \end{aligned}
\end{equation}
and
\begin{equation}
    \left\|  \left(\frac{1}{\rho^{\varepsilon}} - \frac{1}{\rho^{\circ}}\right) \mu^{\circ}\Delta \widetilde{\u}^{\varepsilon} \right\|_{{L}^2(H^{s-2})}
    \lesssim
    C(\mathcal{Y}^{\varepsilon}) \varepsilon^{\delta}
    \left\| \widetilde{\u}^{\varepsilon} \right\|_{\widetilde{L}^2(\dot{B}^{\frac{N}{2}})}
    \left\| \left(   \widetilde{\rho}^{\varepsilon} ,  \ \widetilde{\theta}^{\varepsilon} \right) \right\|_{\widetilde{L}^{\infty}(B^{\frac{N}{2}-1, \frac{N}{2}}_{\varepsilon \nu^{\circ} })}.
\end{equation}
Quite similar, we have
\begin{equation}
    \begin{aligned}
    \big\|\left(\frac{1}{\rho^{\varepsilon}} - \frac{1}{\rho^{\circ}}\right)\left({\frac{N-2}{N}\mu^{\circ}+\lambda^{\circ}}\right)\nabla \operatorname{div} \tilde{\u}^{\varepsilon} &+ \frac{1}{{\rho}^{\varepsilon}} \nabla\left(\left({\frac{N-2}{N} (\mu^{\varepsilon} - \mu^{\circ})+ \lambda^{\varepsilon} - \lambda^{\circ}}\right) \operatorname{div} \tilde{\u}^{\varepsilon} \right) \big\|_{{L}^2(H^{s-2})}
    \\ \lesssim
    &C(\mathcal{Y}^{\varepsilon}) \varepsilon^{\delta} \left\| \widetilde{\u}^{\varepsilon} \right\|_{\widetilde{L}^2(\dot{B}^{\frac{N}{2}})}
    \left\| \left(   \widetilde{{\rho}}^{\varepsilon} ,  \ \widetilde{\theta}^{\varepsilon} \right) \right\|_{\widetilde{L}^{\infty}(B^{\frac{N}{2}-1, \frac{N}{2}}_{\varepsilon \nu^{\circ} })}.
    \end{aligned}
\end{equation}
For the rest term in $r^{\varepsilon}_2$ we actually need a second order expansion to rule out the factor $\frac{1}{\varepsilon}$ and gain smallness. Indeed, we have
\begin{equation}
    \begin{aligned}
    {\frac{1}{\varepsilon^{2}}} \frac{\nabla p^{\varepsilon}}{\rho^{\varepsilon}}
    - \frac{1}{\varepsilon}\frac{p^{\circ}_{\theta}}{\rho^{\circ}}\nabla \tilde{\theta}^{\varepsilon}
    -\frac{1}{\varepsilon}\frac{p^{\circ}_{\rho}}{\rho^{\circ}}\nabla \tilde{\rho}^{\varepsilon}
    -C_1 \tilde{\rho}^{\varepsilon} \nabla \tilde{\rho}^{\varepsilon}
    -C_2 \tilde{\theta}^{\varepsilon} \nabla \tilde{\theta}^{\varepsilon}
    -C_3 \tilde{\theta}^{\varepsilon} \nabla \tilde{\rho}^{\varepsilon}
    -C_4 \tilde{\rho}^{\varepsilon} \nabla \tilde{\theta}^{\varepsilon} \\
    =\frac{1}{\varepsilon} \left( \frac{p_{\rho}^{\varepsilon}}{\rho^{\varepsilon}} - \frac{p^{\circ}_{\rho}}{\rho^{\circ}} - \left(\frac{p_{\rho}}{\rho}\right)^{\circ}_{\rho}(\varepsilon \widetilde{\rho}^{\varepsilon})-\left(\frac{p_{\rho}}{\rho}\right)^{\circ}_{\theta}(\varepsilon \widetilde{\theta}^{\varepsilon})\right) \nabla \widetilde{\rho}^{\varepsilon}
    \\+\frac{1}{\varepsilon} \left( \frac{p_{\theta}^{\varepsilon}}{\rho^{\varepsilon}} - \frac{p^{\circ}_{\theta}}{\rho^{\circ}} - \left(\frac{p_{\theta}}{\rho}\right)^{\circ}_{\rho}(\varepsilon \widetilde{\rho}^{\varepsilon} )-\left(\frac{p_{\theta}}{\rho}\right)^{\circ}_{\theta}(\varepsilon \widetilde{\theta}^{\varepsilon})\right) \nabla \widetilde{\theta}^{\varepsilon}   .
    \end{aligned}
\end{equation}
Again using Lemma \ref{delta_estimates}, Lemma \ref{product_estimate} and Lemma \ref{composition_estimates} we have
\begin{equation}\label{r2_1}
    \begin{aligned}
        \left\| \frac{1}{\varepsilon} \left( \frac{p_{\rho}^{\varepsilon}}{\rho^{\varepsilon}} - \frac{p^{\circ}_{\rho}}{\rho^{\circ}} - \left(\frac{p_{\rho}}{\rho}\right)^{\circ}_{\rho}(\varepsilon \widetilde{\rho}^{\varepsilon})-\left(\frac{p_{\rho}}{\rho}\right)^{\circ}_{\theta}(\varepsilon \widetilde{\theta}^{\varepsilon})\right) \nabla \widetilde{\rho}^{\varepsilon} \right\|_{{L}^{2}(H^{s-2})}
        \\=\left\| \frac{1}{\varepsilon} \left( \frac{p_{\rho}^{\varepsilon}}{\rho^{\varepsilon}} - \frac{p^{\circ}_{\rho}}{\rho^{\circ}} - \left(\frac{p_{\rho}}{\rho}\right)^{\circ}_{\rho}(\varepsilon \widetilde{\rho}^{\varepsilon})-\left(\frac{p_{\rho}}{\rho}\right)^{\circ}_{\theta}(\varepsilon \widetilde{\theta}^{\varepsilon})\right) \nabla \widetilde{\rho}^{\varepsilon} \right\|_{{L}^{2}(H^{s-2})}
        \\ \lesssim
        \left\| \left( \frac{p_{\rho \rho}\rho - p_{\rho}}{\rho^2} \bigg|_{(\rho^{\circ}+ \tau \varepsilon \widetilde{\rho}^{\varepsilon}, \theta^{\circ} + \tau \varepsilon \widetilde{\theta}^{\varepsilon} )} - \frac{p_{\rho \rho}\rho - p_{\rho}}{\rho^2} \bigg|_{(\rho^{\circ},\theta^{\circ})} \right)\widetilde{\rho}^{\varepsilon} \nabla \widetilde{\rho}^{\varepsilon}\right\|_{L^{2}(H^{s-2})}
        \\+\left\| \left(\frac{p_{\rho \theta} }{\rho} \bigg|_{(\rho^{\circ}+ \tau \varepsilon \widetilde{\rho}^{\varepsilon}, \theta^{\circ} + \tau \varepsilon \widetilde{\theta}^{\varepsilon} )} - \frac{p_{\rho \theta}}{\rho} \bigg|_{(\rho^{\circ},\theta^{\circ})} \right)  \widetilde{\theta}^{\varepsilon} \nabla \widetilde{\rho}^{\varepsilon} \right\|_{L^{2}(H^{s-2})}
        \\ \lesssim
        C(\mathcal{Y}^{\varepsilon}) \left\| (\varepsilon \widetilde{\rho}^{\varepsilon} , \varepsilon \widetilde{\theta}^{\varepsilon}) \right\|_{L^{\infty}(H^{s})}
        \left(
        \left\| \widetilde{\rho}^{\varepsilon} \nabla \widetilde{\rho}^{\varepsilon}\right\|_{L^2(H^{\frac{N}{2}-2})}
        +\left\| \widetilde{\theta}^{\varepsilon} \nabla \widetilde{\rho}^{\varepsilon}\right\|_{L^2(H^{\frac{N}{2}-2})}
        \right)
        \\ \lesssim
        C(\mathcal{Y}^{\varepsilon}) \varepsilon^{\delta}\left\| ( \widetilde{\rho}^{\varepsilon} ,  \widetilde{\theta}^{\varepsilon}) \right\|_{\widetilde{L}^{\infty}(B^{\frac{N}{2}-1, \frac{N}{2}}_{\varepsilon \nu^{\circ}})}
            \left\| (\widetilde{\rho}^{\varepsilon},\ \widetilde{\theta}^{\varepsilon})\right\|_{\widetilde{L}^{\infty}(B^{\frac{N}{2}-1,\frac{N}{2}}_{\varepsilon \nu^{\circ}})} \left\| \widetilde{\rho}^{\varepsilon} \right\|_{\widetilde{L}^2(\dot{B}^{\frac{N}{2}})}.
    \end{aligned}
\end{equation}
Same arguments lead to
 \begin{equation}\label{r2_2}
    \begin{aligned}
        \left\| \frac{1}{\varepsilon} \left( \frac{p_{\theta}^{\varepsilon}}{\rho^{\varepsilon}} - \frac{p^{\circ}_{\theta}}{\rho^{\circ}} - \left(\frac{p_{\theta}}{\rho}\right)^{\circ}_{\rho}(\varepsilon \widetilde{\rho} )-\left(\frac{p_{\theta}}{\rho}\right)^{\circ}_{\theta}(\varepsilon \widetilde{\theta}^{\varepsilon})\right) \nabla \widetilde{\theta}^{\varepsilon}  \right\|_{{L}^2({H}^{s-2})}
        \\ \lesssim
        C(\mathcal{Y}^{\varepsilon}) \varepsilon^{\delta}
            \left\| (\widetilde{\rho}^{\varepsilon},\ \widetilde{\theta}^{\varepsilon})\right\|_{\widetilde{L}^{\infty}(B^{\frac{N}{2}-1,\frac{N}{2}}_{\varepsilon \nu^{\circ}})}^2
            \| \widetilde{\theta}^{\varepsilon} \|_{\widetilde{L}^2(\dot{B}^{\frac{N}{2}})}.
    \end{aligned}
 \end{equation}
 Combining \eqref{r2_1} and \eqref{r2_2} we have
\begin{equation}
    \begin{aligned}
        \left\|
            {\frac{1}{\varepsilon^{2}}} \frac{\nabla p^{\varepsilon}}{\rho^{\varepsilon}}
            \right.&\left.- \frac{1}{\varepsilon}\frac{p^{\circ}_{\theta}}{\rho^{\circ}}\nabla \tilde{\theta}^{\varepsilon}
            -\frac{1}{\varepsilon}\frac{p^{\circ}_{\rho}}{\rho^{\circ}}\nabla \tilde{\rho}^{\varepsilon}
            -C_1 \tilde{\rho}^{\varepsilon} \nabla \tilde{\rho}^{\varepsilon}
            -C_2 \tilde{\theta}^{\varepsilon} \nabla \tilde{\theta}^{\varepsilon}
            -C_3 \tilde{\theta}^{\varepsilon} \nabla \tilde{\rho}^{\varepsilon}
            -C_4 \tilde{\rho}^{\varepsilon} \nabla \tilde{\theta}^{\varepsilon}
            \right\|_{{L}^{2}(H^{s-2})}\\
            &\lesssim
            C(\mathcal{Y}^\varepsilon) \varepsilon^{\delta} 
            \left\| (\tilde{\rho}^{\varepsilon},\ \tilde{\theta}^{\varepsilon}) \right\|_{\widetilde{L}^{\infty}(B^{\frac{N}{2}-1,\frac{N}{2}}_{\varepsilon \nu^{\circ}})}^2
            \left\| (\tilde{\rho}^{\varepsilon},\ \tilde{\theta}^{\varepsilon}) \right\|_{\widetilde{L}^{2}(\dot{B}^{\frac{N}{2},\frac{N}{2}+1}_{\varepsilon \nu^{\circ}})}.
    \end{aligned}
\end{equation}
In summary
\begin{equation}
    \left\| r^{\varepsilon}_2 \right\|_{\widetilde{L}^2_{T}(\dot{H}^{s-2})}  \lesssim
    C(\mathcal{Y}^\varepsilon) \varepsilon^{\delta}
    \left(
    \left\| \widetilde{U}^{\varepsilon} \right\|_{E^{\frac{N}{2}}_{\varepsilon \nu^{\circ}}}^2
    +\left\| \widetilde{U}^{\varepsilon} \right\|_{E^{\frac{N}{2}}_{\varepsilon \nu^{\circ}}}^3
    \right)
    \lesssim
    C(\mathcal{Y}^\varepsilon) \varepsilon^{\delta}
    \left((X+Y)^2 + (X+Y)^3\right).
\end{equation}

We now focus on the estimates for $r^{\varepsilon}_3$. Note that
\begin{equation}
    \frac{\kappa^{\circ}}{e^{\circ}_{\theta} \rho^{\circ}} \Delta  \tilde{\theta}^{\varepsilon}
     - \frac{1}{e_{\theta}^{\varepsilon} \rho^{\varepsilon}} \nabla( \kappa^{\varepsilon} \nabla  \tilde{\theta}^{\varepsilon})
    =
    \left(\frac{1}{e_{\theta}^{\circ} \rho^{\circ}} - \frac{1}{{e_{\theta}^{\varepsilon}} \rho^{\varepsilon}}\right) \kappa^{\circ} \Delta \widetilde{\theta}^{\varepsilon}
    +\frac{1}{{e_{\theta}^{\varepsilon}} \rho^{\varepsilon}} \nabla [\left(  \kappa^{\circ} - \kappa^{\varepsilon}\right) \nabla {\theta}^{\varepsilon}],
\end{equation}
and that
\begin{equation}
    \begin{aligned}
    &\frac{1}{\varepsilon} \left( \frac{\theta^{\circ} p^{\circ}_{\theta}}{e^{\circ}_{\theta} \rho^{\circ}} - \frac{\theta^{\varepsilon} p_{\theta}^{\varepsilon}}{e_{\theta}^{\varepsilon} \rho^{\varepsilon}} + C_5 \varepsilon\tilde{\theta}^{\varepsilon} + C_6 \varepsilon \tilde{\rho}^{\varepsilon} \right) \nabla \!\cdot\! \tilde{\u}^{\varepsilon}
    \\ &=\frac{1}{\varepsilon} \left(
        \frac{\theta^{\circ} p^{\circ}_{\theta}}{e^{\circ}_{\theta} \rho^{\circ}} - \frac{\theta^{\varepsilon} p_{\theta}^{\varepsilon}}{e_{\theta}^{\varepsilon} \rho^{\varepsilon}} + \left(\frac{\theta p_{\theta}}{e_{\theta} \rho}\right)_{\theta} \bigg|_{(\rho^{\circ}, \theta^{\circ})} \varepsilon \widetilde{\theta}^{\varepsilon} +
        \left(\frac{\theta p_{\theta}}{e_{\theta} \rho}\right)_{\rho} \bigg|_{(\rho^{\circ}, \theta^{\circ})} \varepsilon \widetilde{\rho}^{\varepsilon}
    \right) \nabla \!\cdot\! \widetilde{\u}^{\varepsilon}
    \\ & = \frac{1}{\varepsilon}
        \left(
            \left(\frac{\theta p_{\theta}}{e_{\theta} \rho}\right)_{\theta} \bigg|_{(\rho^{\circ}, \theta^{\circ})}  - \left(\frac{\theta p_{\theta}}{e_{\theta} \rho}\right)_{\theta} \bigg|_{(\rho^{\circ}+ \tau \varepsilon \widetilde{\rho}^{\varepsilon}, \theta^{\circ}+ \tau \varepsilon \widetilde{\theta}^{\varepsilon})}
        \right)\varepsilon \widetilde{\theta}^{\varepsilon}\nabla \!\cdot\! \widetilde{\u}^{\varepsilon}
    \\ &+\frac{1}{\varepsilon}
        \left(
            \left(\frac{\theta p_{\theta}}{e_{\theta} \rho}\right)_{\rho} \bigg|_{(\rho^{\circ}, \theta^{\circ})}  - \left(\frac{\theta p_{\theta}}{e_{\theta} \rho}\right)_{\rho} \bigg|_{(\rho^{\circ}+ \tau \varepsilon \widetilde{\rho}^{\varepsilon}, \theta^{\circ}+ \tau \varepsilon \widetilde{\theta}^{\varepsilon})}
        \right)\varepsilon \widetilde{\rho}^{\varepsilon}\nabla \!\cdot\! \widetilde{\u}^{\varepsilon}.
    \end{aligned}
\end{equation}
Using again Lemma , we have
\begin{equation}
    \begin{aligned}
        \left\|\frac{\kappa^{\circ}}{e^{\circ}_{\theta} \rho^{\circ}} \Delta  \tilde{\theta}^{\varepsilon} - \frac{1}{e_{\theta}^{\varepsilon} \rho^{\varepsilon}} \nabla( \kappa^{\varepsilon} \nabla  \tilde{\theta}^{\varepsilon})\right\|_{{L}^{2}(H^{s-2})}
        &\lesssim
        C(\mathcal{Y}^{\varepsilon}) \varepsilon^{\delta} \left\| (\tilde{\rho}^{\varepsilon},\ \tilde{\theta}^{\varepsilon}) \right\|_{\widetilde{L}^{\infty}(B^{\frac{N}{2}-1,\frac{N}{2}}_{\varepsilon \nu^{\circ}})}
        \| \tilde{\theta}^{\varepsilon} \|_{\widetilde{L}^{2}(\dot{B}^{\frac{N}{2},\frac{N}{2}+1}_{\varepsilon \nu^{\circ}})},\\
        \left\|
        \frac{1}{\varepsilon} \left( \frac{\theta^{\circ} p^{\circ}_{\theta}}{e^{\circ}_{\theta} \rho^{\circ}} - \frac{\theta^{\varepsilon} p_{\theta}^{\varepsilon}}{e_{\theta}^{\varepsilon} \rho^{\varepsilon}} + C_5 \tilde{\theta}^{\varepsilon} + C_6 \tilde{\rho}^{\varepsilon} \right) \nabla \!\cdot\! \tilde{\u}^{\varepsilon}
        \right\|_{{L}^{2}(H^{s-2})}
        &\lesssim
        C(\mathcal{Y}^\varepsilon) \varepsilon^{\delta}
        \left\| (\tilde{\rho}^{\varepsilon},\ \tilde{\theta}^{\varepsilon}) \right\|_{\widetilde{L}^{\infty}(B^{\frac{N}{2}-1,\frac{N}{2}}_{\varepsilon \nu^{\circ}})}^{2}
        \left\| \tilde{\u}^{\varepsilon} \right\|_{\widetilde{L}^{2}(\dot{B}^{\frac{N}{2},\frac{N}{2}+1}_{\varepsilon \nu^{\circ}})}.\\
    \end{aligned}
\end{equation}
For the last term in $r^{\varepsilon}_3$ we have
\begin{equation*}
    \begin{aligned}
        \frac{\varepsilon}{e_{\theta}^{\varepsilon} \rho^{\varepsilon}} S^{\varepsilon}\!:\! \nabla  \tilde{\u}^{\varepsilon}
        = \varepsilon \frac{\mu^{\varepsilon}}{e_{\theta}^{\varepsilon} \rho^{\varepsilon}} \left(\nabla \widetilde{\u}^\varepsilon + {\nabla \widetilde{\u}^{\varepsilon}}^{\text{T}}\right): \nabla \widetilde{\u}^{\varepsilon}
        +\varepsilon \left( \lambda^{\varepsilon} - \frac{N}{2} \right) \frac{1}{e_{\theta}^{\varepsilon} \rho^{\varepsilon}} (\nabla \!\cdot\! \widetilde{\u}^{\varepsilon})^2,
    \end{aligned}
\end{equation*}
since $\mu$, $\lambda$ and $e_{\theta}$ are continuous in $\rho^{\varepsilon}$ and $\theta^{\varepsilon}$, they and their reciprocal are bounded by a continuous modulus $C(\mathcal{Y}^{\varepsilon})$.
Therefore, it's easy to get
\begin{equation}
    \begin{aligned}
        \left\|
        \frac{\varepsilon}{e_{\theta} \rho^{\varepsilon}} S^{\varepsilon}\!:\! \nabla  \tilde{\u}^{\varepsilon}
        \right\|_{{L}^{2}(H^{s-2})}
        &\lesssim
        C(\mathcal{Y}^\varepsilon) \varepsilon^{\delta} \left\| \tilde{\u}^{\varepsilon} \right\|_{\widetilde{L}^{2}(\dot{B}^{\frac{N}{2},\frac{N}{2}+1}_{\varepsilon \nu^{\circ}})}
        \left\| \tilde{\u}^{\varepsilon} \right\|_{\widetilde{L}^{\infty}(B^{\frac{N}{2}-1,\frac{N}{2}}_{\varepsilon \nu^{\circ}})}.\\
    \end{aligned}
\end{equation}
Combining the estimate for each term of $r^{\varepsilon}_3$, we have
\begin{equation}
    \begin{aligned}
        \left\| r^{\varepsilon}_3 \right\|_{\widetilde{L}^2_{T}(\dot{H}^{s-2})}  \lesssim  C(\mathcal{Y}^\varepsilon) \varepsilon^{\delta} \left((X+Y)^2 + (X+Y)^3\right).
    \end{aligned}
\end{equation}
\end{proof}
\begin{proof}[Proof of Lemma \ref{XR_estimates}]

First observe that $\mathcal{U}$ is of the form $\mathcal{U}= (-\frac{p_{\theta}^{\circ}}{p_{\rho}^{\circ}}\vartheta, \omega, \vartheta)^{T}$, then we have
\begin{equation}
    \begin{aligned}
        R^{1,\varepsilon}_c&= -\mathrm{e}^{\frac{t}{\varepsilon}\mathcal{A}} \mathcal{P}^{\perp} \mathcal{Q}(\mathcal{U},\mathcal{U}),\  \\
        &=-\mathrm{e}^{\frac{t}{\varepsilon}\mathcal{A}} \mathcal{P}^{\perp}
        \begin{pmatrix}
            p_{\theta}^{\circ} \omega \cdot \nabla \vartheta \\
            \omega \cdot \nabla \omega + \left( C_1(p_\theta^{\circ})^2
            + C_2(p_{\rho}^{\circ})^2 - (C_3 + C_4)p_{\theta}^{\circ} p_{\rho}^{\circ}\right) \vartheta \nabla \vartheta\\
            - p_{\rho}^{\circ} \omega \cdot \nabla \vartheta \\
        \end{pmatrix}\\
        &= -\sum_{\gamma, m}\left< \begin{pmatrix}
            p_{\theta}^{\circ} \omega \cdot \nabla \vartheta \\
            \omega \cdot \nabla \omega + \left( C_1(p_\theta^{\circ})^2
            + C_2(p_{\rho}^{\circ})^2 - (C_3 + C_4)p_{\theta}^{\circ} p_{\rho}^{\circ}\right)
             \vartheta \nabla \vartheta\\
            - p_{\rho}^{\circ} \omega \cdot \nabla \vartheta \\
        \end{pmatrix},\ H^{\gamma}_m\right>_{\mathbb{H}} \mathrm{e}^{\frac{t}{\varepsilon}\mathcal{A}} H^{\gamma}_m\\
        &= - \mathrm{i} c_N c^{\circ} \sum_{\gamma, m} \sum_{k+l =m}\left( \frac{\widehat{\omega}_l \!\cdot\! m \ \widehat{\omega}_k \!\cdot\! m}{|m|}
        + \left( C_1(p_\theta^{\circ})^2 + C_2(p_{\rho}^{\circ})^2 - (C_3 + C_4)p_{\theta}^{\circ} p_{\rho}^{\circ}\right) \widehat{\vartheta}_{l} \widehat{\vartheta}_{k} \frac{l \!\cdot\! m}{|m|}
        \right) \mathrm{e}^{\frac{t}{\varepsilon} \lambda^{\gamma}_{m}} H^{\gamma}_m ,
    \end{aligned}
\end{equation}
and
\begin{equation}
    \begin{aligned}
        \widetilde{R}^{1,\varepsilon}_{c,M}
        = - \mathrm{i} c_N c^{\circ} \sum_{\gamma, |m| \leqslant M} \sum_{k+l =m}\left( \frac{\widehat{\omega}_l \!\cdot\! m \ \widehat{\omega}_k \!\cdot\! m}{|m|}
        + \left( C_1(p_\theta^{\circ})^2 + C_2(p_{\rho}^{\circ})^2 - (C_3 + C_4)p_{\theta}^{\circ} p_{\rho}^{\circ}\right) \widehat{\vartheta}_{l} \widehat{\vartheta}_{k} \frac{l \!\cdot\! m}{|m|}
        \right)
        \\ \times \frac{1}{\lambda^{\gamma}_{m}}\mathrm{e}^{\frac{t}{\varepsilon} \lambda^{\gamma}_{m}} H^{\gamma}_m,
    \end{aligned}
\end{equation}
\begin{equation}
    \begin{aligned}
        \widetilde{R}^{1,\varepsilon,t}_{c,M}
        = - \mathrm{i} c_N c^{\circ}
        \sum_{\gamma, |m| \leqslant M} \sum_{k+l =m}
        \partial_t
        \left( \frac{\widehat{\omega}_l \!\cdot\! m \ \widehat{\omega}_k \!\cdot\! m}{|m|}
        +  \left( C_1(p_\theta^{\circ})^2 + C_2(p_{\rho}^{\circ})^2 - (C_3 + C_4)p_{\theta}^{\circ} p_{\rho}^{\circ}\right) \widehat{\vartheta}_{l} \widehat{\vartheta}_{k} \frac{l \!\cdot\! m}{|m|}
        \right)
        \\ \times \frac{1}{\lambda^{\gamma}_{m}}\mathrm{e}^{\frac{t}{\varepsilon} \lambda^{\gamma}_{m}} H^{\gamma}_m,
    \end{aligned}
\end{equation}
\begin{equation}
    \begin{aligned}
        {R}^{1,\varepsilon,M}_{c}
        = - \mathrm{i} c_N c^{\circ} \sum_{\gamma, m > |M|} \sum_{k+l =m}\left( \frac{\widehat{\omega}_l \!\cdot\! m \ \widehat{\omega}_k \!\cdot\! m}{|m|}
        + \left( C_1(p_\theta^{\circ})^2 + C_2(p_{\rho}^{\circ})^2 - (C_3 + C_4)p_{\theta}^{\circ} p_{\rho}^{\circ}\right) \widehat{\vartheta}_{l} \widehat{\vartheta}_{k} \frac{l \!\cdot\! m}{|m|}
        \right)
        \\ \times \mathrm{e}^{\frac{t}{\varepsilon} \lambda^{\gamma}_{m}} H^{\gamma}_m.
    \end{aligned}
\end{equation}
Note that the first and the last components of $H^{\gamma}_{m}$ are just constants, and that $|\frac{m}{|m|}| =1$, so we can just treat the term $H^{\gamma}_m$ as $\mathrm{e}^{\mathrm{i}m \cdot x}$ in the sequel. 

Since $\displaystyle \left|\frac{l \cdot m}{|m|} \right| \leqslant |l|$ and $\displaystyle \left| \frac{m_i m_j}{|m|} \right|\leqslant |m|$, we then have by Lemma \ref{convolution_estimate} and Lemma \ref{truncated_frequency_estimate}
\begin{equation}
    \begin{aligned}
        \left\| \sum_{|m| > M}\frac{l \!\cdot\! m}{|m|} \widehat{\vartheta}_{l} \widehat{\vartheta}_{k} \mathrm{e}^{\mathrm{i} m \cdot x}\right\|_{{L}^2(H^{s-2})}
        \lesssim
        M^{-\delta}\left\| \sum_{m}\frac{l \!\cdot\! m}{|m|} \widehat{\vartheta}_{l} \widehat{\vartheta}_{k} \mathrm{e}^{\mathrm{i} m \cdot x}\right\|_{{L}^2(H^{\frac{N}{2}-2})}
        \\\lesssim
        M^{-\delta}\left\| \vartheta \right\|_{{L}^{2}(H^{\frac{N}{2}}\cap L^{\infty})}  \left\| \vartheta \right\|_{\widetilde{L}^{\infty}(H^{\frac{N}{2}-1})},
        \\
        \left\| \sum_{|m| > M}
        \frac{\widehat{\omega}^i_l m_i \ \widehat{\omega}^j_k  m_j}{|m|} \mathrm{e}^{\mathrm{i} m \cdot x}
        \right\|_{{L}^2(H^{s-2})}
        \lesssim
        M^{-\delta}\left\| \sum_{m}
        \frac{\widehat{\omega}^i_l m_i \ \widehat{\omega}^j_k  m_j}{|m|} \mathrm{e}^{\mathrm{i} m \cdot x}
        \right\|_{{L}^2(H^{\frac{N}{2}-2})}
        \\\lesssim
        M^{-\delta}\left\| \omega^{i} \right\|_{{L}^{2}(H^{\frac{N}{2}}\cap L^{\infty})}  \left\| \omega^{j} \right\|_{\widetilde{L}^{\infty}(H^{\frac{N}{2}-1})},
    \end{aligned}
\end{equation}
and therefore
\begin{equation}
    \begin{aligned}
        \left\| {R}^{1,\varepsilon,M}_{c} \right\|_{{L}^{2}(H^{s-2})}
        \lesssim
        M^{-\delta} \| \mathcal{U} \|_{F^{\frac{N}{2}}}^2
        \lesssim M^{-\delta} X^2.
    \end{aligned}
\end{equation}
Noticing also that  $\displaystyle \left|\frac{l \cdot m}{|m| \lambda^{\gamma}_{m}} \right| \leqslant \frac{|l|}{|m|}$ and $\displaystyle \left| \frac{m_i m_j}{|m|\lambda^{\gamma}_{m}} \right|\leqslant 1$,
we have
\begin{equation}
    \begin{aligned}
        \left\| \sum_{|m| \leqslant M}\frac{l \!\cdot\! m}{|m| \lambda^{\gamma}_{m}} \widehat{\vartheta}_{l} \widehat{\vartheta}_{k} \mathrm{e}^{\mathrm{i} m \cdot x}\right\|_{{L}^2(H^{s})\cap \widetilde{L}^{\infty}(H^{s-1})}
        &\lesssim
        M^{1-\delta}\left\| \sum_{m}\frac{l \!\cdot\! m}{|m|\lambda^{\gamma}_{m}} \widehat{\vartheta}_{l} \widehat{\vartheta}_{k} \mathrm{e}^{\mathrm{i} m \cdot x}\right\|_{{L}^2(H^{\frac{N}{2}-1}) \cap \widetilde{L}^{\infty} (H^{\frac{N}{2}}\cap L^{\infty})  }
        \\ &\lesssim
        M^{1-\delta}\left\| \vartheta \right\|_{{L}^{2}(H^{\frac{N}{2}}\cap L^{\infty}) \cap \widetilde{L}^{\infty}(H^{\frac{N}{2}-1})}  \left\| \vartheta \right\|_{\widetilde{L}^{\infty}(H^{\frac{N}{2}-1})},
    \end{aligned}
\end{equation}
\begin{equation}
    \begin{aligned}
         \left\| \sum_{|m| \leqslant M}
        \frac{\widehat{\omega}^i_l m_i \ \widehat{\omega}^j_k  m_j}{|m|\lambda^{\gamma}_{m}} \mathrm{e}^{\mathrm{i} m \cdot x}
        \right\|_{{L}^2(H^{s}) \cap\widetilde{L}^{\infty}(H^{s-1}) }
        &\lesssim
        M^{1-\delta}\left\| \sum_{m}
        \frac{\widehat{\omega}^i_l m_i \ \widehat{\omega}^j_k  m_j}{|m|\lambda^{\gamma}_{m}} \mathrm{e}^{\mathrm{i} m \cdot x}
        \right\|_{{L}^2(H^{\frac{N}{2}-1})}
        \\ &\lesssim
        M^{1-\delta}
        \left\| \omega^{i} \right\|_{{L}^{2}(H^{\frac{N}{2}}\cap L^{\infty}) \cap \widetilde{L}^{\infty}(H^{\frac{N}{2}-1})}  \left\| \omega^{j} \right\|_{\widetilde{L}^{\infty}(H^{\frac{N}{2}-1})},
    \end{aligned}
\end{equation}
and therefore
\begin{equation}
    \begin{aligned}
        \left\| \widetilde{R}^{1,\varepsilon}_{c,M} \right\|_{{{L}^2(H^{s}) \cap\widetilde{L}^{\infty}(H^{s-1}) }}
        \lesssim
        M^{1-\delta} \left\| \mathcal{U} \right\|_{F^{\frac{N}{2}}}^2
        \lesssim
        M^{1-\delta} X^2.
    \end{aligned}
\end{equation}
The estimate for $\widetilde{R}^{1,\varepsilon,t}_{c,M}$ is quite similar using chain rule, then we have
\begin{equation}
    \begin{aligned}
        \left\| \widetilde{R}^{1,\varepsilon,t}_{c,M} \right\|_{{L}^2(H^{s-2})}
        \lesssim
        M^{1-\delta} \left\| \mathcal{U} \right\|_{\widetilde{L}^{\infty}(H^{\frac{N}{2}-1})} \left\| \partial_t \mathcal{U} \right\|_{{L}^2(H^{\frac{N}{2}-2})}
        \lesssim
        M^{1-\delta} X\left\| \partial_t \mathcal{U} \right\|_{{L}^2(H^{\frac{N}{2}-2})}.
    \end{aligned}
\end{equation}
We currently skip the estimates for $R^{2,\varepsilon}_{c}$, $R^{3,\varepsilon}_{c}$ and their corresponding terms, which requires a “small divisor estimate” we will deal with later.

Note that

\begin{equation}
    \begin{aligned}
        S^{\varepsilon}=- \sum_{\alpha,{k}}
        \left(\frac{2\mu^\circ-\frac{2}{N}\mu^\circ
        +\lambda^\circ}{\rho^\circ}\frac{\rho^\circ}{\theta^\circ}\alpha \gamma \operatorname{sg}({k})\operatorname{sg}({m})+\frac{\kappa^\circ (p^\circ_\theta)^2}{(\rho^\circ)^2(c^\circ)^2 e^\circ_{\theta}}
        \right)
        |{k}|^2 \mathrm{ e}^{2\frac{t}{\varepsilon} \lambda^{\alpha}_{k}} V^{\alpha,\varepsilon}_{k} H^{-\alpha}_{k},
    \end{aligned}
\end{equation}
\begin{equation}
    \begin{aligned}
        \widetilde{S}^{\varepsilon}_M=- \sum_{\alpha,|{k}|\leqslant M}
        \left(\frac{2\mu^\circ-\frac{2}{N}\mu^\circ
        +\lambda^\circ}{\rho^\circ}\frac{\rho^\circ}{\theta^\circ}\alpha \gamma \operatorname{sg}({k})\operatorname{sg}({m})+\frac{\kappa^\circ (p^\circ_\theta)^2}{(\rho^\circ)^2(c^\circ)^2 e^\circ_{\theta}}
        \right)
        \frac{1}{2\lambda^\alpha_{k}}|{k}|^2 \mathrm{ e}^{2\frac{t}{\varepsilon} \lambda^{\alpha}_{k}} V^{\alpha,\varepsilon}_{k} H^{-\alpha}_{k},
    \end{aligned}
\end{equation}
\begin{equation}
    \begin{aligned}
        \widetilde{S}^{\varepsilon,t}_M=- \sum_{\alpha,|{k}|\leqslant M}
        \left(\frac{2\mu^\circ-\frac{2}{N}\mu^\circ
        +\lambda^\circ}{\rho^\circ}\frac{\rho^\circ}{\theta^\circ}\alpha \gamma \operatorname{sg}({k})\operatorname{sg}({m})+\frac{\kappa^\circ (p^\circ_\theta)^2}{(\rho^\circ)^2(c^\circ)^2 e^\circ_{\theta}}
        \right)
        \frac{1}{2\lambda^\alpha_{k}}|{k}|^2 \mathrm{ e}^{2\frac{t}{\varepsilon} \lambda^{\alpha}_{k}} \partial_t V^{\alpha,\varepsilon}_{k} H^{-\alpha}_{k},
    \end{aligned}
\end{equation}
\begin{equation}
    \begin{aligned}
        {S}^{\varepsilon,M}=- \sum_{\alpha,|{k}|\geq M}
        \left(\frac{2\mu^\circ-\frac{2}{N}\mu^\circ
        +\lambda^\circ}{\rho^\circ}\frac{\rho^\circ}{\theta^\circ}\alpha \gamma \operatorname{sg}({k})\operatorname{sg}({m})+\frac{\kappa^\circ (p^\circ_\theta)^2}{(\rho^\circ)^2(c^\circ)^2 e^\circ_{\theta}}
        \right)
       |{k}|^2 \mathrm{ e}^{2\frac{t}{\varepsilon} \lambda^{\alpha}_{k}} V^{\alpha,\varepsilon}_{k} H^{-\alpha}_{k},
    \end{aligned}
\end{equation}
using again Lemma \ref{convolution_estimate}, Lemma \ref{embedding} as in estimates of $R^{1,\varepsilon}_c$ we have
\begin{equation}\label{S_estimates}
    \begin{aligned}
        &\left\| \widetilde{S}^{\varepsilon}_{M} \right\|_{L^2(H^s) \cap \widetilde{L}^{\infty}(H^{s-1})} \lesssim M^{1-\delta}\left\| \widetilde{V}^\varepsilon   \right\|_{{{L}}^{2}(H^{\frac{N}{2}} \cap L^{\infty})\cap \widetilde{L}^{\infty}(H^{\frac{N}{2}-1})}
        \lesssim M^{1-\delta}X, \\
        &\left\| \widetilde{S}^{\varepsilon,t}_{M} \right\|_{{{L}}^{2}(H^{s-2})} \lesssim M^{1-\delta}\left\| \partial_t\widetilde{V}^\varepsilon \right\|_{{{L}}^{2}(H^{\frac{N}{2}-2})},\\
        &\left\| {S}^{\varepsilon,M}  \right\|_{{{L}}^{2}(H^{s-2})}  \lesssim M^{-\delta} \left\|\widetilde{V}^\varepsilon \right\|_{{{L}}^{2}(H^{\frac{N}{2}} \cap L^{\infty})}
        \lesssim M^{-\delta} X,\\
    \end{aligned}
\end{equation}
\begin{equation} \label{R4_0}
    \begin{aligned}
        R^{4,\varepsilon}_{c}
        =\mathrm{e}^{\frac{t}{\varepsilon}\mathcal{A}} \mathcal{P}^{\perp} \mathcal{D} \mathcal{P} \widetilde{U}^\varepsilon
        =
        \sum_{\alpha,{k}}\frac{-\kappa^\circ p^{\circ}_{\rho}\theta^\circ p^\circ_{\theta}}{(e^{\circ}_{\theta})^2 (\rho^\circ)^2 c^{\circ}} \widehat{\vartheta}_{k}|{k}|^2  \mathrm{e}^{\frac{t}{\varepsilon}\lambda^{\alpha}_{k}} H^{\alpha}_{k}.
    \end{aligned}
\end{equation}
For $R^{4,\varepsilon}_c$ we have
\begin{equation}\label{R4_1}
    \begin{aligned}
        \widetilde{R}^{4,\varepsilon}_{c,M}=\sum_{\alpha,|{k}|\leqslant M}-\frac{1}{\lambda^{\alpha}_{k}}\frac{\kappa^\circ p^{\circ}_{\rho}\theta^\circ p^\circ_{\theta}}{(e^{\circ}_{\theta})^2 (\rho^\circ)^2 c^{\circ}} \widehat{\vartheta}_{k}|{k}|^2 \mathrm{e}^{\frac{t}{\varepsilon}\lambda^{\alpha}_{k}} H^{\alpha}_{k},
    \end{aligned}
\end{equation}
\begin{equation}\label{R4_2}
    \begin{aligned}
        \widetilde{R}^{4,\varepsilon,t}_{c,M}=\sum_{\alpha,|{k}|\leqslant M}-\frac{1}{\lambda^{\alpha}_{k}}\frac{\kappa^\circ p^{\circ}_{\rho}\theta^\circ p^\circ_{\theta}}{(e^{\circ}_{\theta})^2 (\rho^\circ)^2 c^{\circ}} \partial_t\widehat{\vartheta}_{k}|{k}|^2 \mathrm{e}^{\frac{t}{\varepsilon}\lambda^{\alpha}_{k}} H^{\alpha}_{k},
    \end{aligned}
\end{equation}
\begin{equation}\label{R4_3}
    \begin{aligned}
        R^{4,\varepsilon,M}_{c}=\sum_{\alpha,|{k}|>M}\frac{-\kappa^\circ p^{\circ}_{\rho}\theta^\circ p^\circ_{\theta}}{(e^{\circ}_{\theta})^2 (\rho^\circ)^2 c^{\circ}} \widehat{\vartheta}_{k}|{k}|^2 \mathrm{e}^{\frac{t}{\varepsilon}\lambda^{\alpha}_{k}} H^{\alpha}_{k}.
    \end{aligned}
\end{equation}
Same method applying to \eqref{R4_0}, \eqref{R4_1}, \eqref{R4_2} and \eqref{R4_3} yields
\begin{equation}\label{R4_estimates}
    \begin{aligned}
        &\left\| \widetilde{R}^{4,\varepsilon}_{M} \right\|_{L^2(H^s)\cap {\widetilde{L}}^{\infty}(H^{s-1}) } \lesssim M^{1-\delta}\left\|\mathcal{P} \widetilde{U}^\varepsilon   \right\|_{{\widetilde{L}}^{2}(B^{\frac{N}{2}})\cap {\widetilde{L}}^{\infty}(H^{\frac{N}{2}-1})}
        \lesssim  M^{1-\delta}X,\\
        &\left\| \widetilde{R}^{4,\varepsilon,t}_{M} \right\|_{{{L}}^{2}(H^{s-2})} \lesssim M^{1-\delta}\left\| \partial_t \mathcal{P} \widetilde{U}^\varepsilon \right\|_{{{L}}^{2}(H^{\frac{N}{2}-2})},\\
        &\left\| {R}^{4,\varepsilon,M}  \right\|_{{{L}}^{2}(H^{s-2})}  \lesssim M^{-\delta} \left\| \mathcal{P} \widetilde{U}^\varepsilon \right\|_{{{L}}^{2}(H^{\frac{N}{2}} \cap L^{\infty})}
        \lesssim M^{-\delta} X.\\
    \end{aligned}
\end{equation}
Finally, we turn to the estimates for $R^{2,\varepsilon}_{c}$ and $R^{3,\varepsilon}_{c}$.
For $\gamma \in \{-1,1\}$ and ${m} \in \widetilde{\mathbb{Z}}^{N} \backslash \{0\}$, denote similarly as in Danchin's work the following index set
\begin{equation*}
    \begin{aligned}
          \mathcal{B}^{1,\gamma}_m=&\{(\alpha,{k},{l})\in\{-1,1\}\times\widetilde{\mathbb{Z}}^{N} \backslash \{0\} \times \widetilde{\mathbb{Z}}^{N} \backslash \{0\} \mid \alpha \operatorname{sg} ({k})|{k}| \neq \gamma \operatorname{sg}({m})|{m}|,\ {k}+{l}={m}\},\\
          \mathcal{B}^{2,\gamma}_m=&\{(\alpha,\beta,{k},{l})\in\{-1,1\}^2\times\widetilde{\mathbb{Z}}^{N} \backslash \{0\} \times \widetilde{\mathbb{Z}}^{N} \backslash \{0\} \mid \alpha \operatorname{sg} ({k})|{k}| +\beta \operatorname{sg}({l})|{l}|\neq \gamma \operatorname{sg}({m})|{m}|,\ \\
          &{k}+{l}={m}\},
    \end{aligned}
\end{equation*}

\begin{equation*}
    \begin{aligned}
          \mathcal{B}^{1,\gamma}_{m,M}=&\{(\alpha,{k},{l})\in\{-1,1\}\times\widetilde{\mathbb{Z}}^{N} \backslash \{0\} \times \widetilde{\mathbb{Z}}^{N} \backslash \{0\} \mid \alpha \operatorname{sg} (k)|{k}| \neq \gamma  \operatorname{sg} (m)|{m}||,\\
           &{k}+{l}={m}, \ |l|,|k|\leqslant M\},\\
          \mathcal{B}^{2,\gamma}_{m,M}=&\{(\alpha,\beta,{k},{l})\in\{-1,1\}^2\times\widetilde{\mathbb{Z}}^{N} \backslash \{0\} \times \widetilde{\mathbb{Z}}^{N} \backslash \{0\} \mid \alpha \operatorname{sg} (k)|{k}| +\beta  \operatorname{sg} (l)|{l}|\neq \gamma  \operatorname{sg} (m)|{m}||,\ \\
          &{k}+{l}={m},\ |l|,|k|\leqslant M \},
    \end{aligned}
\end{equation*}
and
\begin{equation*}
    \begin{aligned}
        \mathcal{B}^{1,\gamma,M}_{m}=&\{(\alpha,{k},{l})\in\{-1,1\}\times\widetilde{\mathbb{Z}}^{N} \backslash \{0\} \times \widetilde{\mathbb{Z}}^{N} \backslash \{0\} \mid \alpha \operatorname{sg} (k)|{k}| \neq \gamma  \operatorname{sg} (m)|{m}||,\\
        & {k}+{l}={m},
        |k|\leqslant M,\ |l| > M \}
        \\ &\cup
        \{(\alpha,{k},{l})\in\{-1,1\}\times\widetilde{\mathbb{Z}}^{N} \backslash \{0\} \times \widetilde{\mathbb{Z}}^{N} \backslash \{0\} \mid \alpha \operatorname{sg} (k)|{k}| \neq \gamma  \operatorname{sg} (m)|{m}||,\\
        & {k}+{l}={m},
        |k| >M\},
        \\
          \mathcal{B}^{2,\gamma,M}_{m}=&\{(\alpha,\beta,{k},{l})\in\{-1,1\}^2\times\widetilde{\mathbb{Z}}^{N} \backslash \{0\} \times \widetilde{\mathbb{Z}}^{N} \backslash \{0\} \mid \alpha \operatorname{sg} (k)|{k}| +\beta  \operatorname{sg} (l)|{l}|\neq \gamma  \operatorname{sg} (m)|{m}||,\ \\
          &{k}+{l}={m},\ |k|\leqslant M, \ |l| > M \}
        \\ & \cup
        \{(\alpha,\beta,{k},{l})\in\{-1,1\}^2\times\widetilde{\mathbb{Z}}^{N} \backslash \{0\} \times \widetilde{\mathbb{Z}}^{N} \backslash \{0\} \mid \alpha \operatorname{sg} (k)|{k}| +\beta  \operatorname{sg} (l)|{l}|\neq \gamma  \operatorname{sg} (m)|{m}||,\ \\
          &{k}+{l}={m},\ |k| > M \}.
    \end{aligned}
\end{equation*}
Also, we denote
\begin{equation}
    \begin{aligned}
        C^1_M  \stackrel{\text { def }}{=}
        \max_{\substack{m,\gamma \\ (\alpha,k,l)\in \mathcal{B}_{m,M}^{1,\gamma}}}
        \big|\alpha \operatorname{sg} (k)|{k}| - \gamma  \operatorname{sg} (m)|{m}| \big|^{-1},
        \\
        C_M^2 \stackrel{\text { def }}{=}
        \max_{\substack{m, \gamma  \\(\alpha, \beta, k,l) \in \mathcal{B}_{m, M}^{2,\gamma}}}
        \big|\alpha {\operatorname{sg}} (k)| k|+\beta {\operatorname{sg}} (l)|l|-\gamma {\operatorname{sg}} m |m| \big|^{-1}.
    \end{aligned}
\end{equation}
We will show in next section that $C^1_M,\ C^2_M$  can actually be achieved and are polynomials in $M$.
Using \eqref{Q2r_epsilon} and \eqref{Q2r} we have
\begin{equation}
    \begin{aligned}
        R^{2,\varepsilon}_{c}=\sum_{\substack{m,\gamma\\(\alpha, {k} ,{l} ) \in \mathcal{B}^{1,\gamma}_m}} -\mathrm{i}
        &\left[\frac{\rho^\circ}{\theta^\circ} \left(1+\alpha \gamma \operatorname{sg}({k})\operatorname{sg}({m})\right)\widehat{\omega}_{{l}}\!\cdot\!{m}
        +\left(
            p^\circ_\rho C_5 - p^\circ_\theta C_6 - \frac{p^\circ_\rho p^\circ_\theta}{c^{\circ} \theta^\circ}
        \right)\alpha \operatorname{sg}({k}) \frac{{k} \! \cdot\! {l} }{|{k}|} \widehat{\theta}_{{l}} \right.\\
        &\left.
        +\left(
            \tilde{c}_1\frac{{k} \!\cdot \! {m}}{|{m}|} + \tilde{c}_2 \frac{{l} \!\cdot \! {m}}{|{m}|}
        \right)
            \gamma \operatorname{sg}({m}) \widehat{\theta}_{{l}}
        \right]
        (c_N c^{\circ})^2
        V^{\alpha}_{k} H^{\gamma}_{{m}} \mathrm{e}^{-\frac{t}{\varepsilon}\left(\lambda^{\alpha}_{k}-\lambda^{\gamma}_{{m}}\right)},
    \end{aligned}
\end{equation}
\begin{equation}
    \begin{aligned}
        \widetilde{R}^{2,\varepsilon}_{c,M}=
        \sum_{\substack{m,\gamma\\(\alpha, {k} ,{l} ) \in \mathcal{B}^{1,\gamma}_{m,M}}} \mathrm{i}
        &\frac{1}{\lambda^{\alpha}_{k}-\lambda^{\gamma}_{{m}}}\left[\frac{\rho^\circ}{\theta^\circ} \left(1+\alpha \gamma \operatorname{sg}({k})\operatorname{sg}({m})\right)  \widehat{\omega}_{{l}} \!\cdot\!{m} \right. \\
        &\left.
        +\left(
            p^\circ_\rho C_5 - p^\circ_\theta C_6 - \frac{p^\circ_\rho p^\circ_\theta}{c^{\circ} \theta^\circ}
        \right)\alpha \operatorname{sg}({k}) \frac{{k} \! \cdot\! {l} }{|{k}|}   \widehat{\theta}_{{l}}  \right.\\
        &\left.
        +\left(
            \tilde{c}_1\frac{{k} \!\cdot \! {m}}{|{m}|} + \tilde{c}_2 \frac{{l} \!\cdot \! {m}}{|{m}|}
        \right)\gamma \operatorname{sg}({m})   \widehat{\theta}_{{l}}
        \right]
        (c_N c^{\circ})^2
        V^{\alpha}_{k} H^{\gamma}_{{m}} \mathrm{e}^{-\frac{t}{\varepsilon}\left(\lambda^{\alpha}_{k}-\lambda^{\gamma}_{{m}}\right)},
    \end{aligned}
\end{equation}
Note that,  up to a constant  $C(\rho^{\circ}, \theta^{\circ}, e, p)$, $V^{\alpha}_k$ is $\displaystyle \pi + \alpha\operatorname{sg}(k)\frac{k\cdot v}{|k|} (\pi,v \text{ defined in \eqref{projected_variables} })$.
It's then easy to see
\begin{equation*}
    \begin{aligned}
        \big| \frac{1}{\lambda^{\alpha}_k - \lambda^{\gamma}_m}\big| \big|(1+\alpha \gamma \operatorname{sg}(k) \operatorname{sg}(m)) m_i\left(1 + \alpha \operatorname{sg}(k) \frac{k_j}{|k|}\right) \big|
        \lesssim C^1_M |m|,
    \end{aligned}
\end{equation*}
therefore, using again Lemma \ref{convolution_estimate} and noticing that $|m| \leqslant 2M$ we have
\begin{equation}
    \begin{aligned}
    \big\|
    \sum_{\substack{m,\gamma\\(\alpha, {k} ,{l} ) \in \mathcal{B}^{1,\gamma}_{m,M}}} \mathrm{i}
    \frac{1}{\lambda^{\alpha}_{k}-\lambda^{\gamma}_{{m}}}
    &\left(1+\alpha \gamma \operatorname{sg}({k})\operatorname{sg}({m})\right)  \widehat{\omega}_{{l}}^i \!\cdot\!{m}_i
    (\pi + \alpha \operatorname{sg}(k)\frac{{\vv}_j k_j}{|k|}) H^{\gamma}_{{m}} \mathrm{e}^{-\frac{t}{\varepsilon}\left(\lambda^{\alpha}_{k}-\lambda^{\gamma}_{{m}}\right)}
    \big\|_{{L}^2(H^s)\cap\widetilde{L}^{\infty}(H^{s-1})}
    \\
    &\begin{aligned}
    &\lesssim
    C^1_M M^{2-\delta} \left\| \omega^i \right\|_{{L}^2(H^{\frac{N}{2}}\cap L^{\infty})\cap\widetilde{L}^{\infty}(H^{\frac{N}{2}-1})}
    \left\| (\pi,\ {\vv}^j) \right\|_{\widetilde{L}^{\infty}(H^{\frac{N}{2}-1})}
    \\ &\lesssim
    \left\| \mathcal{U} \right\|_{{L}^2(H^{\frac{N}{2}}\cap L^{\infty})\cap\widetilde{L}^{\infty}(H^{\frac{N}{2}-1})}
    \left\| V \right\|_{\widetilde{L}^{\infty}(H^{\frac{N}{2}-1})}
    \\ &\lesssim C^1_M M^{2-\delta}XY.
    \end{aligned}
    \end{aligned}\\
\end{equation}

Now for the term
\begin{equation}
    \begin{aligned}
        \widetilde{R}^{2,\varepsilon,t}_{c,M}=\sum_{\substack{m,\gamma\\(\alpha, {k} ,{l} ) \in \mathcal{B}^{1,\gamma}_{m,M}}} \mathrm{i}
        &\frac{1}{\lambda^{\alpha}_{k}-\lambda^{\gamma}_{{m}}}\left[\frac{\rho^\circ}{\theta^\circ} \left(1+\alpha \gamma \operatorname{sg}({k})\operatorname{sg}({m})\right)\partial_t \left(\widehat{\omega}_{{l}}V^{\alpha}_{k}\right) \!\cdot\!{m} \right. \\
        &\left.
        +\left(
            p^\circ_\rho C_5 - p^\circ_\theta C_6 - \frac{p^\circ_\rho p^\circ_\theta}{c_0 \theta^\circ}
        \right)\alpha \operatorname{sg}({k}) \frac{{k} \! \cdot\! {l} }{|{k}|} \partial_t \left(\widehat{\theta}_{{l}} V^{\alpha}_{k}\right) \right.\\
        &\left.
        +\left(
            \tilde{c}_1\frac{{k} \!\cdot \! {m}}{|{m}|} + \tilde{c}_2 \frac{{l} \!\cdot \! {m}}{|{m}|}
        \right)\gamma \operatorname{sg}({m}) \partial_t \left(\widehat{\theta}_{{l}} V^{\alpha}_{k} \right)
        \right] (c_N c^{\circ})^2  H^{\gamma}_{{m}} \mathrm{e}^{-\frac{t}{\varepsilon}\left(\lambda^{\alpha}_{k}-\lambda^{\gamma}_{{m}}\right)} ,
    \end{aligned}
\end{equation}
same methods yields
\begin{equation}
    \begin{aligned}
        &\left\| \widetilde{R}^{2,\varepsilon,t}_{M} \right\|_{{{L}}^{2}(H^{s-2})} \lesssim C^{1}_{M} M^{2-\delta} \left(\left\| V \right\|_{{\widetilde{L}}^{\infty}(H^{\frac{N}{2}-1})} \left\| \partial_t \mathcal{U} \right\|_{{{L}}^{2}(H^{\frac{N}{2}-2})}   +     \left\| \mathcal{U} \right\|_{{\widetilde{L}}^{\infty}(H^{\frac{N}{2}-1})} \left\| \partial_t {V} \right\|_{{{L}}^{2}(H^{\frac{N}{2}-2})}\right).\\
    \end{aligned}
\end{equation}
Treatment for $R^{2,\varepsilon,M}_c$ is a little bit different, we first give its expression
\begin{equation}
    \begin{aligned}
        R^{2,\varepsilon,M}_{c}=\sum_{\substack{m,\gamma\\(\alpha, {k} ,{l} ) \in \mathcal{B}^{1,\gamma,M}_m}} -\mathrm{i}
        &\left[\frac{\rho^\circ}{\theta^\circ} \left(1+\alpha \gamma \operatorname{sg}({k})\operatorname{sg}({m})\right)\widehat{\omega}_{{l}}\!\cdot\!{m}
        +\left(
            p^\circ_\rho C_5 - p^\circ_\theta C_6 - \frac{p^\circ_\rho p^\circ_\theta}{c^{\circ} \theta^\circ}
        \right)\alpha \operatorname{sg}({k}) \frac{{k} \! \cdot\! {l} }{|{k}|} \widehat{\theta}_{{l}} \right.\\
        &\left.
        +\left(
            \tilde{c}_1\frac{{k} \!\cdot \! {m}}{|{m}|} + \tilde{c}_2 \frac{{l} \!\cdot \! {m}}{|{m}|}
        \right)
            \gamma \operatorname{sg}({m}) \widehat{\theta}_{{l}}
        \right] (c_N c^{\circ})^2 V^{\alpha}_{k} H^{\gamma}_{{m}} \mathrm{e}^{-\frac{t}{\varepsilon}\left(\lambda^{\alpha}_{k}-\lambda^{\gamma}_{{m}}\right)}.
    \end{aligned}
\end{equation}
note that
\begin{equation}
    \begin{aligned}
        &\sum_{\substack{m,\gamma\\(\alpha, {k} ,{l} ) \in \mathcal{B}^{1,\gamma,M}_m}} -\mathrm{i}
        \left(1+\alpha \gamma \operatorname{sg}({k})\operatorname{sg}({m})\right)\widehat{\omega}_{{l}}\!\cdot\!{m}
         (\pi + \alpha \operatorname{sg}(k)\frac{{\vv}_j k_j}{|k|}) H^{\gamma}_{{m}} \mathrm{e}^{-\frac{t}{\varepsilon}\left(\lambda^{\alpha}_{k}-\lambda^{\gamma}_{{m}}\right)}
         \\ &=\sum_{\substack{m,\gamma\\(\alpha, {k} ,{l} ) \in \mathcal{B}^{1,\gamma}_m}}-\mathrm{i}
        \left(1+\alpha \gamma \operatorname{sg}({k})\operatorname{sg}({m})\right)\widehat{\omega}_{{l}}\!\cdot\!{m}
         (\pi^M + \alpha \operatorname{sg}(k)\frac{{\vv}^M_j k_j}{|k|}) H^{\gamma}_{{m}} \mathrm{e}^{-\frac{t}{\varepsilon}\left(\lambda^{\alpha}_{k}-\lambda^{\gamma}_{{m}}\right)}
         \\ &+\sum_{\substack{m,\gamma\\(\alpha, {k} ,{l} ) \in \mathcal{B}^{1,\gamma}_m}}-\mathrm{i}
        \left(1+\alpha \gamma \operatorname{sg}({k})\operatorname{sg}({m})\right)\widehat{\omega^M}_{{l}}\!\cdot\!{m}
         (\pi_M + \alpha \operatorname{sg}(k)\frac{({\vv}_M)_j k_j}{|k|}) H^{\gamma}_{{m}} \mathrm{e}^{-\frac{t}{\varepsilon}\left(\lambda^{\alpha}_{k}-\lambda^{\gamma}_{{m}}\right)}.
    \end{aligned}
\end{equation}
Same as in estimates of $\widetilde{R}^{2,\varepsilon}_{c,M}$, by Lemma we have
\begin{equation*}
    \begin{aligned}
        \big\|
        \sum_{\substack{m,\gamma\\(\alpha, {k} ,{l} ) \in \mathcal{B}^{1,\gamma}_m}}-\mathrm{i}
        \left(1+\alpha \gamma \operatorname{sg}({k})\operatorname{sg}({m})\right)\widehat{\omega}_{{l}}\!\cdot\!{m}
         (\pi^M + \alpha \operatorname{sg}(k)\frac{{\vv}^M_j k_j}{|k|}) H^{\gamma}_{{m}} \mathrm{e}^{-\frac{t}{\varepsilon}\left(\lambda^{\alpha}_{k}-\lambda^{\gamma}_{{m}}\right)}
        \big\|_{{L}^2(H^{s-2})}
        \\
         \begin{aligned}
         &
        \lesssim
        \left\| \mathcal{U} \right\|_{{L}^2(H^{\frac{N}{2}}\cap L^{\infty})}
        \left\| V^M \right\|_{\widetilde{L}^{\infty}(H^{s-1})}
        \lesssim M^{-\delta}
        \left\| \mathcal{U} \right\|_{{L}^2(H^{\frac{N}{2}}\cap L^{\infty})}
        \left\| V \right\|_{\widetilde{L}^{\infty}(H^{\frac{N}{2}-1})}
        \\ &
        \lesssim
        M^{-\delta} XY,
         \end{aligned}
    \end{aligned}
\end{equation*}
and
\begin{equation}
    \begin{aligned}
       \big\|
        \sum_{\substack{m,\gamma\\(\alpha, {k} ,{l} ) \in \mathcal{B}^{1,\gamma}_m}}-\mathrm{i}
        \left(1+\alpha \gamma \operatorname{sg}({k})\operatorname{sg}({m})\right)\widehat{\omega^M}_{{l}}\!\cdot\!{m}
         (\pi_M + \alpha \operatorname{sg}(k)\frac{({\vv}_M)_j k_j}{|k|}) H^{\gamma}_{{m}} \mathrm{e}^{-\frac{t}{\varepsilon}\left(\lambda^{\alpha}_{k}-\lambda^{\gamma}_{{m}}\right)}
         \big\|_{{L}^2(H^{s-2})}
        \\ \begin{aligned}
            &\lesssim
        \left\| \mathcal{U}^M \right\|_{{L}^2(H^{s})}
        \left\| V_M \right\|_{\widetilde{L}^{\infty}(H^{\frac{N}{2}-1})}
        \lesssim M^{-\delta}
        \left\| \mathcal{U} \right\|_{{L}^2(H^{\frac{N}{2}}\cap L^{\infty})}
        \left\| V \right\|_{\widetilde{L}^{\infty}(H^{\frac{N}{2}-1})}
        \\ &
        \lesssim
        M^{-\delta} XY.
        \end{aligned}
    \end{aligned}
\end{equation}
The rest terms in $R^{2,\varepsilon,M}$ can be dealt with similarly, and we have
\begin{equation}
    \lVert R^{2,\varepsilon,M} \rVert_{{L}^2(H^{s-2})}
    \lesssim M^{-\delta} XY.
\end{equation}
For terms corresponding to $R^{3,\varepsilon}_c$ using \eqref{Q3r_epsilon} and \eqref{Q3r}
\begin{equation}
    \begin{aligned}
        R^{3,\varepsilon}_{c}=
        &\sum_{\substack{m,\gamma\\(\alpha, {k} ,{l} ) \in \mathcal{B}^{2,\gamma}_m}} -\mathrm{i}
        \left[\left(\frac{\rho^\circ}{\theta^\circ}\right)^2\frac{1}{\theta^\circ} \frac{p^{\circ}_{\rho}}{\rho^{\circ}} \left(\alpha({k})\operatorname{sg} \frac{k}{|{k}|}+\beta \operatorname{sg}({l})\frac{{l}}{|{l}|}\right)\!\cdot\!{m}
        \right.\\
        &\left.+\gamma \operatorname{sg}({m})|{m}|
        \left( \frac{\rho^\circ}{\theta^\circ} \alpha \beta\operatorname{sg}({k})\operatorname{sg}({l}) \frac{{l}\!\cdot\!{k}}{|{k}||{l}|}
            +\left( \frac{\rho^\circ}{\theta^\circ} \right)^2 C_1 +\left( \frac{\theta^\circ p^{\circ}_{\theta}}{\rho^\circ e^\circ_{\theta} c^\circ} \right)^2 C_2 + \frac{\theta^\circ p^\circ_{\theta}}{e^\circ_{\theta} (c^\circ)^2}\left( C_3 + C_4 \right)
        \right)
         \right.\\
        &\left.
            +\frac{\rho^\circ}{(\theta^\circ)^2 e^\circ_\theta }\left(\frac{\theta^\circ p^{\circ}_{\theta}}{\rho^\circ e^\circ_{\theta} c^\circ}\right)^2
            \left( \alpha \operatorname{sg}({k})\frac{{k}\!\cdot\!{l}}{|{k}|}
            + \beta \operatorname{sg}({l})\frac{{k}\!\cdot\!{l}}{|{l}|}\right) \right.\\
            &\left.
            +\left( C_5\left(\frac{\theta^\circ p^{\circ}_{\theta}}{\rho^\circ e^\circ_{\theta} c^\circ}\right)^2
            + C_6\frac{p^\circ}{c^\circ}\frac{\theta^\circ p^{\circ}_{\theta}}{\rho^\circ e^\circ_{\theta} c^\circ}
            \right)\left(\beta \operatorname{sg}({l})|{l}|+\alpha\operatorname{sg}({k})|{k}|\right)
        \right] \\
        &
        (c_N c^{\circ})^3
        V^{\alpha}_{k} V^{\beta}_{{l}} H^{\gamma}_{{m}} \mathrm{e}^{-\frac{t}{\varepsilon}\left(\lambda^{\alpha}_{k}+\lambda^{\beta}_{{l}}-\lambda^{\gamma}_{{m}}\right)},
    \end{aligned}
\end{equation}
\begin{equation}
    \begin{aligned}
        \widetilde{R}^{3,\varepsilon}_{c,M}=
        &\sum_{\substack{m,\gamma\\(\alpha, {k} ,{l} ) \in \mathcal{B}^{2,\gamma}_{m,M}}} \mathrm{i}
        \left[\left(\frac{\rho^\circ}{\theta^\circ}\right)^2\frac{1}{\theta^\circ} \frac{p^{\circ}_{\rho}}{\rho^{\circ}} \left(\alpha\operatorname{sg} ({k})\frac{k}{|{k}|}+\beta \operatorname{sg}({l})\frac{{l}}{|{l}|}\right)\!\cdot\!{m}
        \right.\\
        &\left.+\gamma \operatorname{sg}({m})|{m}|
        \left( \frac{\rho^\circ}{\theta^\circ} \alpha \beta\operatorname{sg}({k})\operatorname{sg}({l}) \frac{{l}\!\cdot\!{k}}{|{k}||{l}|}
            +\left( \frac{\rho^\circ}{\theta^\circ} \right)^2 C_1 +\left( \frac{\theta^\circ p^{\circ}_{\theta}}{\rho^\circ e^\circ_{\theta} c^\circ} \right)^2 C_2 + \frac{\theta^\circ p^\circ_{\theta}}{e^\circ_{\theta} (c^\circ)^2}\left( C_3 + C_4 \right)
        \right)
         \right.\\
        &\left.
            +\frac{\rho^\circ}{(\theta^\circ)^2 e^\circ_\theta }\left(\frac{\theta^\circ p^{\circ}_{\theta}}{\rho^\circ e^\circ_{\theta} c^\circ}\right)^2
            \left( \alpha \operatorname{sg}({k})\frac{{k}\!\cdot\!{l}}{|{k}|}
            + \beta \operatorname{sg}({l})\frac{{k}\!\cdot\!{l}}{|{l}|}\right) \right.\\
            &\left.
            +\left( C_5\left(\frac{\theta^\circ p^{\circ}_{\theta}}{\rho^\circ e^\circ_{\theta} c^\circ}\right)^2
            + C_6\frac{p^\circ}{c^\circ}\frac{\theta^\circ p^{\circ}_{\theta}}{\rho^\circ e^\circ_{\theta} c^\circ}
            \right)\left(\beta \operatorname{sg}({l})|{l}|+\alpha\operatorname{sg}({k})|{k}|\right)
        \right] \\
        &(c_N c^{\circ})^3
         \frac{1}{\lambda^{\alpha}_{k}+\lambda^{\beta}_{{l}}-\lambda^{\gamma}_{{m}}} V^{\alpha}_{k} V^{\beta}_{{l}} H^{\gamma}_{{m}} \mathrm{e}^{-\frac{t}{\varepsilon}\left(\lambda^{\alpha}_{k}+\lambda^{\beta}_{{l}}-\lambda^{\gamma}_{{m}}\right)},
    \end{aligned}
\end{equation}
\begin{equation}
    \begin{aligned}
        \widetilde{R}^{3,\varepsilon,t}_{c,M}=
        &\sum_{\substack{m,\gamma\\(\alpha, {k} ,{l} ) \in \mathcal{B}^{2,\gamma}_{m,M}}} \mathrm{i}
        \left[\left(\frac{\rho^\circ}{\theta^\circ}\right)^2\frac{1}{\theta^\circ} \frac{p^{\circ}_{\rho}}{\rho^{\circ}} \left(\alpha({k})\operatorname{sg} \frac{k}{|{k}|}+\beta \operatorname{sg}({l})\frac{{l}}{|{l}|}\right)\!\cdot\!{m}
        \right.\\
        &\left.+\gamma \operatorname{sg}({m})|{m}|
        \left( \frac{\rho^\circ}{\theta^\circ} \alpha \beta\operatorname{sg}({k})\operatorname{sg}({l}) \frac{{l}\!\cdot\!{k}}{|{k}||{l}|}
            +\left( \frac{\rho^\circ}{\theta^\circ} \right)^2 C_1 +\left( \frac{\theta^\circ p^{\circ}_{\theta}}{\rho^\circ e^\circ_{\theta} c^\circ} \right)^2 C_2 + \frac{\theta^\circ p^\circ_{\theta}}{e^\circ_{\theta} (c^\circ)^2}\left( C_3 + C_4 \right)
        \right)
         \right.\\
        &\left.
            +\frac{\rho^\circ}{(\theta^\circ)^2 e^\circ_\theta }\left(\frac{\theta^\circ p^{\circ}_{\theta}}{\rho^\circ e^\circ_{\theta} c^\circ}\right)^2
            \left( \alpha \operatorname{sg}({k})\frac{{k}\!\cdot\!{l}}{|{k}|}
            + \beta \operatorname{sg}({l})\frac{{k}\!\cdot\!{l}}{|{l}|}\right) \right.\\
            &\left.
            +\left( C_5\left(\frac{\theta^\circ p^{\circ}_{\theta}}{\rho^\circ e^\circ_{\theta} c^\circ}\right)^2
            + C_6\frac{p^\circ}{c^\circ}\frac{\theta^\circ p^{\circ}_{\theta}}{\rho^\circ e^\circ_{\theta} c^\circ}
            \right)\left(\beta \operatorname{sg}({l})|{l}|+\alpha\operatorname{sg}({k})|{k}|\right)
        \right] \\
        &(c_N c^{\circ})^3
        \frac{1}{\lambda^{\alpha}_{k}+\lambda^{\beta}_{{l}}-\lambda^{\gamma}_{{m}}}
        \partial_t \left({V^{\alpha}_{k} V^{\beta}_{{l}}} \right) H^{\gamma}_{{m}} \mathrm{e}^{-\frac{t}{\varepsilon}\left(\lambda^{\alpha}_{k}+\lambda^{\beta}_{{l}}-\lambda^{\gamma}_{{m}}\right)},
    \end{aligned}
\end{equation}

\begin{equation}
    \begin{aligned}
        R^{3,\varepsilon,M}_{c}=
        &\sum_{\substack{m,\gamma\\(\alpha, {k} ,{l} ) \in \mathcal{B}^{2,\gamma,M}_m}} -\mathrm{i}
        \left[\left(\frac{\rho^\circ}{\theta^\circ}\right)^2\frac{1}{\theta^\circ} \frac{p^{\circ}_{\rho}}{\rho^{\circ}} \left(\alpha({k})\operatorname{sg} \frac{k}{|{k}|}+\beta \operatorname{sg}({l})\frac{{l}}{|{l}|}\right)\!\cdot\!{m}
        \right.\\
        &\left.+\gamma \operatorname{sg}({m})|{m}|
        \left( \frac{\rho^\circ}{\theta^\circ} \alpha \beta\operatorname{sg}({k})\operatorname{sg}({l}) \frac{{l}\!\cdot\!{k}}{|{k}||{l}|}
            +\left( \frac{\rho^\circ}{\theta^\circ} \right)^2 C_1 +\left( \frac{\theta^\circ p^{\circ}_{\theta}}{\rho^\circ e^\circ_{\theta} c^\circ} \right)^2 C_2 + \frac{\theta^\circ p^\circ_{\theta}}{e^\circ_{\theta} (c^\circ)^2}\left( C_3 + C_4 \right)
        \right)
         \right.\\
        &\left.
            +\frac{\rho^\circ}{(\theta^\circ)^2 e^\circ_\theta }\left(\frac{\theta^\circ p^{\circ}_{\theta}}{\rho^\circ e^\circ_{\theta} c^\circ}\right)^2
            \left( \alpha \operatorname{sg}({k})\frac{{k}\!\cdot\!{l}}{|{k}|}
            + \beta \operatorname{sg}({l})\frac{{k}\!\cdot\!{l}}{|{l}|}\right) \right.\\
            &\left.
            +\left( C_5\left(\frac{\theta^\circ p^{\circ}_{\theta}}{\rho^\circ e^\circ_{\theta} c^\circ}\right)^2
            + C_6\frac{p^\circ}{c^\circ}\frac{\theta^\circ p^{\circ}_{\theta}}{\rho^\circ e^\circ_{\theta} c^\circ}
            \right)\left(\beta \operatorname{sg}({l})|{l}|+\alpha\operatorname{sg}({k})|{k}|\right)
        \right] \\
        &(c_N c^{\circ})^3
        V^{\alpha}_{k} V^{\beta}_{{l}} H^{\gamma}_{{m}} \mathrm{e}^{-\frac{t}{\varepsilon}\left(\lambda^{\alpha}_{k}+\lambda^{\beta}_{{l}}-\lambda^{\gamma}_{{m}}\right)}.
    \end{aligned}
\end{equation}
Using same method as in estimates of $\widetilde{R}^{2,\varepsilon}_{c,M}$, $\widetilde{R}^{2,\varepsilon,t}_{c,M}$ and $R^{2,\varepsilon,M}$ we have,
\begin{equation}
    \begin{aligned}
        &\left\| \widetilde{R}^{3,\varepsilon}_{M} \right\|_{{\widetilde{L}}^{\infty}(H^{s-1})\cap {{L}}^{2}(H^{s})} \lesssim C^{2}_{M} M^{2-\delta}X^2,\\
        &\left\| \widetilde{R}^{3,\varepsilon,t}_{M} \right\|_{{{L}}^{2}(H^{s-2})} \lesssim C^{2}_{M} M^{2-\delta}  X \left\| \partial_t V \right\|_{{{L}}^{2}(H^{\frac{N}{2}-2})},\\
        &\left\| {R}^{3,\varepsilon,M} \right\|_{{{L}}^{2}(H^{s-2})} \lesssim M^{- \delta} X^2.
    \end{aligned}
\end{equation}
\end{proof}

\begin{proof}[Proof of Lemma \ref{time_derivative_estimate}]

From \eqref{nsc_orthogonal},
we have
\begin{equation}
    \begin{aligned}
    \left\| \partial_t \widetilde{V}^{\varepsilon} \right\|_{{L}^2(H^{\frac{N}{2}-2})}
    \lesssim
     \left\| \mathcal{Q}^{\varepsilon}_{2r}(\mathcal{P}\widetilde{U}^{\varepsilon},\widetilde{V}^{\varepsilon}) \right\|_{{L}^2(H^{\frac{N}{2}-2})}
    +\left\| \mathcal{Q}^{\varepsilon}_{3r}(\widetilde{V}^{\varepsilon},\widetilde{V}^{\varepsilon})\right\|_{{L}^2(H^{\frac{N}{2}-2})}
    +\left\| \mathcal{D}^{\varepsilon}\widetilde{V}^{\varepsilon}\right\|_{{L}^2(H^{\frac{N}{2}-2})}\\
    +\left\| \mathrm{e}^{\frac{t}{\varepsilon}\mathcal{A}}\mathcal{P}^{\perp}\mathcal{D}\mathcal{P}\widetilde{U}^{\varepsilon} \right\|_{{L}^2(H^{\frac{N}{2}-2})}
    +\left\| \mathrm{e}^{\frac{t}{\varepsilon}\mathcal{A}}\mathcal{P}^{\perp}r^{\varepsilon} \right\|_{{L}^2(H^{\frac{N}{2}-2})}
    +\left\| \mathrm{e}^{\frac{t}{\varepsilon}\mathcal{A}}\mathcal{P}^{\perp}\mathcal{Q}(\mathcal{P}\widetilde{U}^{\varepsilon},\mathcal{P}\widetilde{U}^{\varepsilon})\right\|_{{L}^2(H^{\frac{N}{2}-2})}.
    \end{aligned}
\end{equation}
Using Lemma \ref{Qe_estimate} and the definition \eqref{definition_of_R} we get
\begin{equation}
    \begin{aligned}
        \left\| \partial_t \widetilde{V}^{\varepsilon} \right\|_{{L}^2(H^{\frac{N}{2}-2})}
        &\lesssim \left\| {R}^{4,\varepsilon}_{c} \right\|_{{L}^2(H^{\frac{N}{2}-2})}
        +\left\| r^{\varepsilon} \right\|_{{L}^2(H^{\frac{N}{2}-2})}
        +\left\| S^{\varepsilon} \right\|_{{L}^2(H^{\frac{N}{2}-2})}
        +\left\| \overline{\mathcal{D}} \widetilde{V}^{\varepsilon} \right\|_{{L}^2(H^{\frac{N}{2}-2})}
        \\ &\begin{aligned}
        &+\left\| \widetilde{V}^{\varepsilon} \right\|_{\widetilde{L}^{\infty}(H^{\frac{N}{2}-1})}
        \left\| \mathcal{P}\widetilde{\mathcal{U}}^{\varepsilon} \right\|_{\widetilde{L}^2(B^{\frac{N}{2}})}
        +\left\| \widetilde{V}^{\varepsilon} \right\|_{\widetilde{L}^{\infty}(H^{\frac{N}{2}-1})}
        \left\| \widetilde{V}^{\varepsilon} \right\|_{\widetilde{L}^2(B^{\frac{N}{2}} )}
        \\ &+\left\| \mathcal{P}\widetilde{\mathcal{U}}^{\varepsilon} \right\|_{\widetilde{L}^{\infty}(H^{\frac{N}{2}-1})}
        \left\|  \mathcal{P}\widetilde{\mathcal{U}}^{\varepsilon} \right\|_{\widetilde{L}^2(B^{\frac{N}{2}})},
        \end{aligned}
    \end{aligned}
\end{equation}
using Lemma \ref{r_estimate} and combining with \eqref{S_estimates} and \eqref{R4_estimates} then we have
\begin{equation}
    \begin{aligned}
            \|\partial_t\widetilde{V}^{\varepsilon}\|_{L^2_{T}(H^{\frac{N}{2}-2})}
            \lesssim X+Y+(X+Y)^2+C(\mathcal{Y}^\varepsilon)\left( (X+Y)^2+(X+Y)^3\right).\\
    \end{aligned}
\end{equation}
From \eqref{nsc_kernel} we have
\begin{equation}
    \begin{aligned}
        \left\| \partial_t \mathcal{P}\widetilde{U}^{\varepsilon} \right\|_{{L}^2(H^{\frac{N}{2}-2})}
        &\lesssim
        \left\| \mathcal{P} \mathcal{Q} (\mathcal{P}\widetilde{U}^{\varepsilon},\mathcal{P}\widetilde{U}^{\varepsilon}) \right\|_{{L}^2(H^{\frac{N}{2}-2})}
        + \left\| \mathcal{P}\mathcal{D}\mathcal{P}\widetilde{U}^{\varepsilon} \right\|_{{L}^2(H^{\frac{N}{2}-2})}
        \\ &+\left\| \mathcal{P}\mathcal{D}\mathcal{P}^{\perp}\widetilde{U}^{\varepsilon}  \right\|_{{L}^2(H^{\frac{N}{2}-2})}
        +\left\| \mathcal{P}\mathcal{Q}(\mathcal{P}\widetilde{U}^{\varepsilon},\mathcal{P}^{\perp}\widetilde{U}^{\varepsilon}) \right\|_{{L}^2(H^{\frac{N}{2}-2})}
        \\ &+\left\| \mathcal{Q}(\mathcal{P}^{\perp}\widetilde{U}^{\varepsilon},\mathcal{P}^{\perp}\widetilde{U}^{\varepsilon}) \right\|_{{L}^2(H^{\frac{N}{2}-2})}
         +\left\| r^{\varepsilon} \right\|_{{L}^2(H^{\frac{N}{2}-2})}.
    \end{aligned}
\end{equation}
Direct use of Lemma \ref{r_estimate} and Lemma \ref{Qe_estimate} yields that
\begin{equation}
    \begin{aligned}
        \left\| \partial_t \mathcal{P}\widetilde{U}^{\varepsilon} \right\|_{{L}^2(H^{\frac{N}{2}-2})}
        &\lesssim
        X + XY + X^2 + Y^2 +C(\mathcal{Y}^\varepsilon)\left( (X+Y)^2+(X+Y)^3\right).\\
    \end{aligned}
\end{equation}
For the term $\partial_t \mathcal{P} \mathcal{U}$,  by \eqref{INSF} we have
\begin{equation}
    \begin{aligned}
        \left\| \partial_t \mathcal{P}\mathcal{U} \right\|_{{L}^2(H^{\frac{N}{2}-2})}
        \lesssim  \left\| \mathcal{P} \mathcal{Q} (\mathcal{U},\mathcal{U}) \right\|_{{L}^2(H^{\frac{N}{2}-2})} +
        \left\| \mathcal{P} \mathcal{D} \mathcal{U} \right\|_{{L}^2(H^{\frac{N}{2}-2})}
         \lesssim
        Y^2 + Y.
    \end{aligned}
\end{equation}
For the term $\partial_t V$, from \eqref{LS} we know that
\begin{equation}
    \begin{aligned}
        \left\| \partial_t V \right\|_{{L}^2(H^{\frac{N}{2}-2})}
        &\lesssim  \left\|  \mathcal{Q}_{2r} (\mathcal{U},V) \right\|_{{L}^2(H^{\frac{N}{2}-2})}
        +  \left\|  \mathcal{Q}_{3r} (V,V) \right\|_{{L}^2(H^{\frac{N}{2}-2})}
        +  \left\| \overline{\mathcal{D}} V \right\|_{{L}^2(H^{\frac{N}{2}-2})}.
    \end{aligned}
\end{equation}
Using the expression \eqref{Q2r}, \eqref{Q3r} and using Lemma \ref{convolution_estimate} we have
\begin{equation}
    \begin{aligned}
    \left\|  \mathcal{Q}_{2r} (\mathcal{U},V) \right\|_{{L}^2(H^{\frac{N}{2}-2})}
        \lesssim \left\| \mathcal{U} \right\|_{\widetilde{L}^2(B^{\frac{N}{2}} )}
                 \left\| V \right\|_{\widetilde{L}^{\infty}(H^{\frac{N}{2}-1})},
    \\
    \left\|  \mathcal{Q}_{3r} (V,V) \right\|_{{L}^2(H^{\frac{N}{2}-2} )}
        \lesssim \left\| V \right\|_{\widetilde{L}^2(B^{\frac{N}{2}} )}
                 \left\| V \right\|_{\widetilde{L}^{\infty}(H^{\frac{N}{2}-1})}.
    \end{aligned}
\end{equation}
Again using Lemma \ref{embedding} and Lemma \ref{interpolation} we have
\begin{equation}
    \begin{aligned}
        \left\| \partial_t V \right\|_{{L}^2(H^{\frac{N}{2}-2})}
        \lesssim XY + X^2 + X.
    \end{aligned}
\end{equation}
\end{proof}

\section{Refined small divisor estimates}
    In this section we deal with the problem concerning the terms $ C^1_M $ and $ C^2_M $.
    Note that
    \begin{equation}
        \begin{aligned}
            C^1_M =
            \max_{\substack{m,\gamma \\ (\alpha,k,l)\in \mathcal{B}_{m,M}^{1,\gamma}}}
            \big|\alpha \operatorname{sg} (k)|{k}| - \gamma  \operatorname{sg} (m)|{m}| \big|^{-1},
            \\
            C_M^2 =
            \max_{\substack{m, \gamma  \\(\alpha, \beta, k, l) \in \mathcal{B}_{m, M}^{2,\gamma}}}
            \big|\alpha \operatorname{sg} (k)| k|+\beta \operatorname{sg} (l)|l|-\gamma \operatorname{sg} (m)| m |\big|^{-1}.
        \end{aligned}
    \end{equation}

    The terms $\left|\alpha\operatorname{sg}(k)|k|-\beta\operatorname{sg}(l)|l| \right|$  and $\left| \alpha\operatorname{sg}(k)|k|+\beta\operatorname{sg}(l)|l|-\gamma\operatorname{sg}(m)|m|\right|$
    can become very small when $M$ becomes large due to the arbitrary choice of the aspect ratios $a_2/a_1, \cdots , a_N/a_1$, which might lead to some singularity. Luckily, as has been illustrated in \cite{Danchin_AJM} and \cite{Isabelle_JMKU},
    $C^1_M$ and $C^2_M$ can be bounded by some polynomial in $M$ for almost all aspect ratios, which makes the convergence of the oscillating part possible. However this is an existence result and it's impossible to get an explicit convergence rate for the oscillating part of the system.
    To overcome this diffculty we shall use some number theory involving linear forms. Indeed, we have the following refined small divisor estimates.

    \begin{proposition}[Estimate of $C^1_M$ ] \label{CM1_estimate}
        For any $\tau>0$, $\alpha,\beta\in\{1,-1\}$, for almost all choices of
        $a=(a_1,\dots,a_N)\in\mathbb{R}_{+}^{N}$ in the sense of Lebesgue measure,
        the following inequalities are true
        for $k,l\in\widetilde{\mathbb{Z}}^{N}_{a}\setminus\{0\}$ with $|k|,|l|\leqslant M$.

        \begin{equation}
            \begin{aligned}
                \alpha\operatorname{sg}(k)|k|-&\beta\operatorname{sg}(l)|l|\neq 0\\
                \Longrightarrow\
                &\frac{1}{\left|\alpha\operatorname{sg}(k)|k|-\beta\operatorname{sg}(l)|l| \right|}
                \leqslant C(1+M)^{2N+2\tau-1},
            \end{aligned}
        \end{equation}
        where $C$ is a constant depending only on $a$, $N$ and $\tau$.
    \end{proposition}
    \begin{proposition}[Estimate of $C^2_M$ ] \label{CM2_estimate}
        \label{small_divisor_estimates}
        For any $a=(a_1,\dots,a_N)\in\mathbb{R}_{+}^{N}$, $\alpha,\beta,\gamma \in\{1,-1\}$,
        and for $k,l,m\in\widetilde{\mathbb{Z}}^{N}_{a}\setminus\{0\}$ with $|k|,|l|,|m|\leqslant M$ and $k+l=m$, the following inequality holds
        \begin{equation}
            \begin{aligned}
                \alpha\operatorname{sg}(k)|k|+&\beta\operatorname{sg}(l)|l|-\gamma\operatorname{sg}(m)|m|\neq 0\\
                \Longrightarrow\
                &\frac{1}{\left|\alpha\operatorname{sg}(k)|k|+\beta\operatorname{sg}(l)|l|-\gamma\operatorname{sg}(m)|m| \right|}
                \leqslant C(1+M)^{5},
            \end{aligned}
        \end{equation}
    \end{proposition}

    The proof of proposition \ref{CM1_estimate} relies heavily upon
    the following lemma on the lower-bound estimate of linear forms
    proven by V. G. Sprindzhuk in 1977 \cite{Sprindzhuk}:
    \begin{lemma}[Sprindzhuk]\label{Sprindzhuk_Lemma}
        For any $\tau>0$, for almost all points
        $\lambda=(\lambda_1,\dots,\lambda_N)\in\mathbb{R}^{N}$ in the sense of Lebesgue measure,
        and for any $x=(x_1,\dots,x_N)\in\mathbb{Z}^{N}$, $x\neq 0$, we have
        \begin{equation}
            |\lambda_1x_1+\dots+\lambda_Nx_N|\geqslant
            C\frac{1}{|x|^{N-1+\tau}},
        \end{equation}
        where $C$ is a constant depending only on $\lambda$ and $\tau$.
    \end{lemma}

    \begin{proof}[Proof of proposition \ref{CM1_estimate}]

    Without loss of generality, we only need to prove estimates for
    divisors of the form
    \begin{equation}
        \frac{1}{\left| |k|-|l| \right|}.
    \end{equation}
    Simple computation shows that
    \begin{equation}\label{RSD_temp_0}
            \frac{1}{\left||k|-|l| \right|} = \frac{|k|+|l|}
            {\left| \frac{ \check{k}_1^2-\check{l}_1^2}{a_1^2}
            +\dots+\frac{\check{k}_N^2-\check{l}_N^2}{a_N^2}\right|}.
    \end{equation}
    Set $x_i = \check{k}_i^2 - \check{l}_i^2 $, by lemma \ref{Sprindzhuk_Lemma} $\forall \tau >0$ there exists a set $B \subset \mathbb{R}^{n}_{+}$ such that $\text{Meas}(\mathbb{R}^{n}_{+}\backslash B) = 0$ and that
    $\forall \lambda \in B$
    \begin{equation}
        |\lambda_1x_1+\dots+\lambda_Nx_N|\geqslant
            C\frac{1}{|x|^{N-1+\tau}}.
    \end{equation}
    Now we only need to show that there exists a set $A \subset \mathbb{R}^{n}_{+}$ with $\text{Meas}(\mathbb{R}^{n}_{+}\backslash A) = 0$ such that $\forall a \in A$, $(\frac{1}{a_1^2},\frac{1}{a_2^2},\cdots, \frac{1}{a_N^2}) \in B$.
    We consider the following map
    \begin{equation}
        \begin{aligned}
            T :\  \mathbb{R}^N_{+} &\rightarrow \mathbb{R}^N_{+}\\
                (y_1,y_2,\cdots,y_{\scriptscriptstyle N}) &\mapsto (\frac{1}{\sqrt{y_1}},\frac{1}{\sqrt{y_2}},\cdots,\frac{1}{\sqrt{y_{\scriptscriptstyle N}}}),
        \end{aligned}
    \end{equation}
    which is clearly differentiable on $\mathbb{R}^N_{+}$ and therefore maps measure zero set to measure zero set.
    We can thus take $A = T(B)$ since $T$ is bijection, and then $\forall a \in A$ we have
    \begin{equation} \label{RSD_temp_1}
        \left| \sum_{1\leqslant i \leqslant N}\frac{ \check{k}_i^2-\check{l}_i^2}{a_i^2} \right| \geqslant C \frac{1}{|x|^{N-1+\tau}}.
    \end{equation}
    On the other hand we have
    \begin{equation}\label{RSD_temp_2}
        \begin{aligned}
        |x| = \left( \sum_{1\leqslant i \leqslant N} (\check{k}_i^2-\check{l}_i^2)^2 \right)^{\frac{1}{2}} = \left( \sum_{1\leqslant i \leqslant N} ({(a_ik_i)}^2-(a_il_i)^2)^2 \right)^{\frac{1}{2}}
        \leqslant (\max_{i} a_i)^2 \left(  |k|^2 + |l|^2 \right).
        \end{aligned}
    \end{equation}
    Combining \eqref{RSD_temp_0}, \eqref{RSD_temp_1} and \eqref{RSD_temp_2} we get
    \begin{align}
        \frac{1}{\left| |k| - |l| \right|  }&\leqslant C \left(|k| +|l|\right) \left((\max_{i} a_i)^2 \left(  |k|^2 + |l|^2 \right)  \right)^{N-1+\tau}\\
        &\leqslant C M^{2N+2\tau-1},
    \end{align}
    where $C$ depends only on $a,\ \tau,\  N$.
    \end{proof}
    The estimate of $C^2_M$  relies heavily on the geometric property $k+l=m$ and is indeed simpler then the estimate of $C^1_M$.
    \begin{proof}[Proof of Proposition \ref{CM2_estimate}]
        We first define $\displaystyle {C_a \stackrel{\mathrm{def}}{=}   \inf_{1\leqslant i \leqslant N} \frac{1}{a_i}}$. Since $k,\ l,\ m \in \widetilde{\mathbb{Z}}^{N}_{a}\setminus\{0\} $,
        we readily have $|k|,\ |l|,\ |m| \geqslant C_a$.
        Without loss of generality, we only need to prove estimates for divisors of the form
        \begin{equation}
          \frac{1} { \left| |k| + |l| -|m| \right|}.
        \end{equation}
        It is clear that $|k|+|l|-|m| >0 $ since $k+l=m$ and $|k|+|l| \neq |m|$.
        Note that
        \begin{align}\label{CM2_temp_1}
            |k|+|l|-|m| &= \frac{1}{|k|+|l| + |m|} \left( |k|^2 + |l|^2 -|m|^2 + 2|k||l|  \right)\\
            &=\frac{2|k||l|}{|k|+|l| + |m|}(1 - \cos \theta),
        \end{align}
    where $\theta $ is the angle between $k$ and $l$.
    If $k$ and $l$ are collinear then $\cos \theta = \pm 1$. Since $|k|+|l| \neq |m|$ $\cos \theta$ can only be $-1$, and then
    \begin{equation} \label{CM2_parital_result_1}
        \frac{1}{|k|+|l|-|m|} = \frac{|k|+|l| + |m|}{4|k||l|} \leqslant \frac{3M}{4(C_a)^2}.
    \end{equation}
    Now if $k$ and $l$ are not collinear, we only need to consider the angle between the line where $k$ la y s on and the line where $l$ lies on.
    Given $k = (k_1,k_2,\cdots,k_N) \stackrel{\mathrm{def}}{=} (\check{k}_1,\check{k}_2,\cdots,\check{k}_N) $, we consider the line $p = t k$,
    where $p=(p_1,p_2,\cdots,p_N) \in \mathbb{R}^N$ and $t \in \mathbb{R}$.

    We claim that for all the lines of the form $p= tk + r$ which passes through  points $p \in \widetilde{\mathbb{Z}}^{N}_{a}\setminus\{0\}$ with $r \in \mathbb{R}^N\setminus\{0\}$,
    there is a lower bound $d_\text{min}$ for the distance between the two parallel lines $p= tk + r$ and $p= tk$, and that
    \begin{equation}
        d_{\text{min}} = \frac{1}{|k|} \sqrt{\sum_{1\leqslant i < j \leqslant N}\left(\frac{1}{a_i a_j} \text{gcd}(\hat{k}_i,\hat{k}_j)\right)} < M.
    \end{equation}
    This claim immediately implies that
    the region
    \begin{equation}
        \{\ q \in \mathbb{R}^N | \ 0< dist(q, p) < d_{\text{min}},\  p=tk \ \}
    \end{equation}
    does not contain any points in $\widetilde{\mathbb{Z}}^{N}_{a}\setminus\{0\}$.
    The figure \ref{2d_demonstration} is a demonstration of the case $N=2$.
    \begin{figure}[h]
        \centering
    \begin{tikzpicture}[scale = 0.5]
         \draw[->,thick] (-1,0) -- (10,0) node[right] {$x$};
         \draw[->,thick] (0,-1) -- (0,10) node[above] {$y$};
         \draw (0,0) node[below left] {$O$};
        \foreach \x in {0,2,4,6,8}{
          \foreach \y in {0,2,4}{
            \fill (\x,\y*2) circle (2pt);
          }
        }
        \def\linefunction#1{1*#1 }
        \draw[ name path=line1] (0, {\linefunction{0}}) -- (8, {\linefunction{8}});
        \draw[red] (0, {\linefunction{0}}) -- (4, {\linefunction{4}});
        \def\linefunction_2#1{1*#1 +2 }
        \draw[ name path=line2] (0, {\linefunction_2{0}}) -- (6, {\linefunction_2{6}});
        \def\linefunction_3#1{1*#1 -2 }
        \draw[ name path=line3] (0, {\linefunction_3{0}}) -- (9, {\linefunction_3{9}});
        \node[below right, red] at (4,4) {k};
        \node[below right, red] at (2,8) {l};
        \draw[red] (0,0) -- (2,8);
        \path[name path=circle] (0,0) -- ++(0:7.5) arc (0:90:7.5);
        \node[below right,font=\tiny] at (9,0) {(M,0)};
        \node[below,font=\tiny] at (2,0) {($\frac{1}{a_1}$,0)};
        \node[below,font=\tiny] at (4,0) {($\frac{2}{a_1}$,0)};
        \node[below,font=\tiny] at (6,0) {($\frac{3}{a_1}$,0)};
        \node[below,font=\tiny] at (8,0) {($\frac{4}{a_1}$,0)};
        \node[left,font=\tiny] at (0,4) {(0,$\frac{1}{a_2}$)};
        \node[left,font=\tiny] at (0,8) {(0,$\frac{2}{a_2}$)};
        \node[left,font=\tiny] at (0,9) {(0,M)};
        \draw[name path = circle ] (0,0) -- ++(0:9.2) arc (0:90:9.2);
        \path[name intersections={of=line2 and circle, by=intersection}];
        \coordinate (A) at (intersection);
        \draw[dashed, line width=0.8pt] (0, 0) -- (intersection) node[midway,above=24pt,right=18pt] {M};
        \coordinate (P) at (intersection-1);
        \path[name intersections={of=line1 and circle, by=intersection}];
        \coordinate (Q) at (intersection-2);
        \path[name intersections={of=line3 and circle, by=intersection}];
        \coordinate (R) at (intersection-2);
        \path[name path=perpline] (P) -- ($(P) + (P) - (1, 1)$);
        \draw[dashed,<->, line width=0.8pt] (P) -- ($(P) - 0.5*(P) + 0.5*(R)$) node [midway,above right] {$\tiny d_{\text{min}}$};
        \end{tikzpicture}
        \caption{Case N = 2.}
        \label{2d_demonstration}
    \end{figure}

        Since $l \in \widetilde{\mathbb{Z}}^{N}_{a}\setminus\{0\} $, then we have
        \begin{equation}
            \cos \theta \leqslant \frac{\sqrt{M^2- d_{\text{min}}^2}}{M},
        \end{equation}
        and therefore
        \begin{equation} \label{CM2_temp_2}
            1 - \cos \theta \geqslant 1 - \frac{\sqrt{M^2- d_{\text{min}}^2}}{M}
            =  \frac{d_{\text{min}}^2}{M\left(\sqrt{M^2- d_{\text{min}}^2} + M\right)}
            \geqslant \frac{N(N-1)(C_a)^2}{4M^4}>0.
        \end{equation}
    Combining \eqref{CM2_temp_1} and \eqref{CM2_temp_2} we have
    \begin{equation} \label{CM2_parital_result_2}
        \frac{1} { \left| |k| + |l| -|m| \right|} \leqslant C M^5,
    \end{equation}
    where $C$ is a constant depends only on $a$, and  $N$.
    In summary, whether $k$ and $l$ are collinear or not, \eqref{CM2_parital_result_1} and \eqref{CM2_parital_result_2} leads to
    \begin{equation}
        \frac{1} { \left| |k| + |l| -|m| \right|} \leqslant C (1+M)^5,
    \end{equation}
    where $C$ is a constant depends only on $a$, and  $N$.

    We now turn to the proof of our claim. The distance $d$ between the two parallel lines $p=tk$ and $p = tk + r$ is calculated as
    \begin{equation}
        d = \frac{1}{|k|} \sqrt{r^2k^2 - (r \cdot k)^2} = \frac{1}{|k|} \sqrt{ \sum_{1 \leqslant i < j \leqslant N} \left( r_i k_j -r_j k_i\right)^2}.
    \end{equation}\label{calculation_of_d}
    Note that we require $p = tk +r $ passes through points in $\widetilde{\mathbb{Z}}^{N}_{a}\setminus\{0\}$, which means there exists $p \in \widetilde{\mathbb{Z}}^{N}_{a}\setminus\{0\} $ and $t \in \mathbb{R}$ such that
    \begin{equation}
        \left\{
            \begin{matrix}
                p_1 = t k_1 + r_1,\\
                p_2 = t k_2 + r_2,\\
                \vdots\\
                p_N = t k_N + r_N,
            \end{matrix}
            \right.
    \end{equation}
    which implies that
    \begin{equation} \label{pre_Bezout}
        p_i k_j - p_j k_i = r_i k_j - r_j k_i, \ \ 1\leqslant i < j \leqslant N.
    \end{equation}
    Multiplying $a_ia_j$ to both sides of \eqref{pre_Bezout} we have
    \begin{equation}
        a_ip_i \check{k}_j - a_jp_j \check{k}_i = (r_i k_j - r_j k_i)a_ia_j.
    \end{equation}
    Since $p \in  \widetilde{\mathbb{Z}}^{N}_{a}\setminus\{0\} $, we know that for $ 1\leqslant i \leqslant N$, $a_i p_i $ are integers. Classic Bézout's Lemma then implies that
    \begin{equation}
        |(r_i k_j - r_j k_i)a_ia_j| \geqslant |\text{gcd}(\check{k}_i, \check{k}_j)|,
    \end{equation}
    where $\text{gcd}(\check{k}_i, \check{k}_j)$ means the greatest common divisor of integers $\check{k}_i$ and $\check{k}_j$.
    And therefore
    \begin{equation}
        d \geqslant  \frac{1}{|k|} \sqrt{ \sum_{1 \leqslant i < j \leqslant N} \left( \frac{\text{gcd}(\check{k}_i ,\check{k}_j)}{a_i a_j}\right)^2} = d_{\text{min}}.
    \end{equation}
    Next, we show that $d_{\text{min}} < M$. Note that $\frac{\text{gcd}(\check{k}_i ,\check{k}_j)}{a_i a_j} \leqslant k_i k_j$, we then have
    \begin{equation}
        \frac{1}{|k|} \sqrt{ \sum_{1 \leqslant i < j \leqslant N} \left( \frac{\text{gcd}(\check{k}_i ,\check{k}_j)}{a_i a_j}\right)^2}
        \leqslant \frac{1}{|k|} \sqrt{ \sum_{1 \leqslant i < j \leqslant N} \left( k_i k_j\right)^2 } < \frac{1}{|k|}\sqrt{|k|^4} = |k| \leqslant M.
    \end{equation}
\end{proof}
    \begin{remark}
        Using Proposition \ref{CM1_estimate} and \ref{CM2_estimate} we readily get that: for $M>1$,
        \begin{equation}
            C_M = \max\{C^1_M,C^2_M\} \leqslant CM^\sigma,
        \end{equation}
        where $\sigma = \max \left\{ 2N +2\tau -1 , 5\right\}$
    \end{remark}

\section{Convergence of the incompressible part of the system}\label{sec_Incompressible}
This section is devoted to the proof of Proposition \ref{incompressible_convergence}.
Although the equation is projected onto $\text{Null}(\mathcal{A})$,
terms that are liable to oscillate in time do not disappear
due to the nonlinear nature of the problem.
The strategy of using norms with negative regularity index in time
is again employed to deal with these terms.
We shall also point out that the projection $\mathcal{P}$
plays an important role in carring out some of the computations.

\begin{proof}[Proof of proposition \ref{incompressible_convergence}]

Subtracting \eqref{INSF} from \eqref{nsc_kernel}, we get

\begin{equation}\label{difference_equation_kernel}
    \begin{aligned}
        &\partial_t w^{\varepsilon} + 2\mathcal{P}\mathcal{Q}(A^{\varepsilon},w^{\varepsilon}) - \mathcal{P}\mathcal{D}w^{\varepsilon}\\
        &=
        - \mathcal{P}R^{1,\varepsilon}
        - \mathcal{P}R^{2,\varepsilon}
        + \mathcal{P}R^{3,\varepsilon}
        + \mathcal{P}r^{\varepsilon}
        - \mathcal{P}\mathcal{Q}(w^{\varepsilon},w^{\varepsilon})
        - 2\mathcal{P}\mathcal{Q}(\mathcal{P}_0 \widetilde{U}^{\varepsilon}, \widetilde{U}^{\varepsilon}),
    \end{aligned}
\end{equation}
with
\begin{equation}
    A^{\varepsilon} = \mathcal{P}^{\perp}\widetilde{U}^{\varepsilon} + \mathcal{U}.
\end{equation}
Here, the oscillating terms are
\begin{equation}
    R^{1,\varepsilon} = \mathcal{Q}(\mathcal{P}^{\perp}\widetilde{U}^{\varepsilon},\mathcal{P}^{\perp}\widetilde{U}^{\varepsilon}),\quad
    R^{2,\varepsilon} = 2\mathcal{Q}(\mathcal{P}^{\perp}\widetilde{U}^{\varepsilon},\mathcal{U}),\quad
    \text{and}\quad R^{3,\varepsilon} = \mathcal{D}\mathcal{P}^{\perp}\widetilde{U}^{\varepsilon}.
\end{equation}
Their low-frequency-truncated parts
$R^{1,\varepsilon}_{M}$, $R^{2,\varepsilon}_{M}$,
and $R^{3,\varepsilon}_{M}$ are defined similarly as in the previous section.
Also, we define
$R^{\varepsilon} \stackrel{\mathrm{def}}{=} R^{1,\varepsilon}+R^{2,\varepsilon}-R^{3,\varepsilon}$,
$R^{\varepsilon}_{M} \stackrel{\mathrm{def}}{=} R^{1,\varepsilon}_{M}+R^{2,\varepsilon}_{M}-R^{3,\varepsilon}_{M}$,
and $R^{\varepsilon,M}\stackrel{\mathrm{def}}{=}R^{\varepsilon}-R^{\varepsilon}_{M}$.
The change of function to be made in this section is
\begin{equation}
    \psi^{\varepsilon}_{M} \stackrel{\mathrm{def}}{=} w^{\varepsilon} + \varepsilon\mathcal{P}\widetilde{R}^{\varepsilon}_{M},
\end{equation}
where the fluctuation
\begin{equation}\label{negative_time_relation}
    \varepsilon\partial_t\mathcal{P}\widetilde{R}^{\varepsilon}_{M} = \mathcal{P}R^{\varepsilon}_{M} + \varepsilon\mathcal{P}\widetilde{R}^{t,\varepsilon}_{M}.
\end{equation}
The difference equation in terms of $\psi^{\varepsilon}_{M}$ becomes
\begin{equation}\label{difference_equation_psi}
    \begin{aligned}
        &\partial_t\psi^{\varepsilon}_{M} + 2\mathcal{P}\mathcal{Q}(A^{\varepsilon},\psi^{\varepsilon}_{M}) -
        \mathcal{P}\mathcal{D}\psi^{\varepsilon}_{M}\\
        &= \varepsilon\mathcal{P}\widetilde{R}^{t,\varepsilon}_{M} +
        2\varepsilon\mathcal{P}\mathcal{Q}(A^{\varepsilon},\mathcal{P}\widetilde{R}^{\varepsilon}_{M})
        - \varepsilon\mathcal{P}\mathcal{D}\mathcal{P}\widetilde{R}^{\varepsilon}_{M}
        + \mathcal{P}r^{\varepsilon}
        - \mathcal{P}R^{\varepsilon,M} - \mathcal{P}\mathcal{Q}(w^{\varepsilon},w^{\varepsilon})
        +2\mathcal{P}\mathcal{Q}(\mathcal{P}_0 \widetilde{U}^{\varepsilon}, \widetilde{U}^{\varepsilon}).
    \end{aligned}
\end{equation}
By Lemma \ref{structure_lemma}, one gets
\begin{equation}\label{estimate_of_psi}
    \begin{aligned}
        \|&\psi^{\varepsilon}_{M}\|_{\widetilde{L}^{\infty}(H^{s-1})} + c\mu^{\prime}\|\psi^{\varepsilon}_{M}\|_{L^{2}(H^{s})}\\
        \leqslant& \frac{C}{\mu^\prime}\text{exp}\Bigg(\frac{C}{\mu^{\prime}}\int_{0}^{\infty}\|A^{\varepsilon}\|^2_{B^{\frac{N}{2}}}\Bigg)\times
        \Big(\varepsilon\|\mathcal{P}\widetilde{R}^{\varepsilon}_{M}(0)\|_{H^{s-1}}\\
        &+ \varepsilon(\|\mathcal{P}\widetilde{R}^{t,\varepsilon}_{M}\|_{{L}^{2}(H^{s-2})} + \|\mathcal{P}\mathcal{Q}(A^{\varepsilon},\mathcal{P}\widetilde{R}^{\varepsilon}_{M})\|_{{L}^{2}(H^{s-2})}
        + \|\mathcal{P}\mathcal{D}\widetilde{R}^{\varepsilon}_{M}\|_{L^{2}(H^{s-2})})\\
        &+ \|\mathcal{P}r^{\varepsilon}\|_{L^{2}(H^{s-2})} + \|\mathcal{P}R^{\varepsilon,M}\|_{L^{2}(H^{s-2})} \\
        &+  \|\mathcal{P}\mathcal{Q}(w^{\varepsilon},w^{\varepsilon})\|_{{L}^{2}(H^{s-2})}
        + 2\|\mathcal{P}\mathcal{Q}(\mathcal{P}_0 \widetilde{U}^{\varepsilon}, \widetilde{U}^{\varepsilon})\|_{{L}^{2}(H^{s-2})}
        \Big),
    \end{aligned}
\end{equation}
where we recall $\mu^{\prime} = \text{min}\big\{\frac{\mu^{\circ}}{\rho^{\circ}},\frac{\kappa^{\circ}p^{\circ}_{\rho}}{(c^{\circ})^2\rho^{\circ}e^{\circ}_{\theta}}\big\}$,
and $C$ is a constant.
According to Lemma \ref{product_estimate} and Lemma \ref{convolution_estimate} , we have
\begin{equation}\label{PQWW_estimate}
    \|\mathcal{P}\mathcal{Q}(w^{\varepsilon},w^{\varepsilon})\|_{{L}^{2}(H^{s-2})}
    \leqslant \|w^{\varepsilon}\|_{\widetilde{L}^2(B^{\frac{N}{2}})}\|\nabla w^{\varepsilon}\|_{\widetilde{L}^\infty(H^{s-2})}
    \leqslant PW^{\varepsilon},
\end{equation}
and
\begin{equation}
    \|\mathcal{P}\mathcal{Q}(\mathcal{P}_0 \widetilde{U}^{\varepsilon}, \widetilde{U}^{\varepsilon})\|_{{L}^{2}(H^{s-2})}
    \lesssim \| \mathcal{P}_0 \widetilde{U}^{\varepsilon} \|_{L^{\infty}(dt)} \| \widetilde{U}^{\varepsilon} \|_{L^{2}(\dot{H}^{\frac{N}{2}})} \lesssim \varepsilon (X+Y) .
\end{equation}
Lemma \ref{r_estimate} guarantees
\begin{equation}\label{remainder_kernel}
    \|r^{\varepsilon}\|_{L^2(H^{s-2})} \leqslant C(\mathcal{Y}^{\varepsilon})\varepsilon^{\delta}\left((X+Y)^2+(X+Y)^3\right).
\end{equation}
For the oscillating terms, first we have by Lemma \ref{product_estimate} that
\begin{equation}\label{PQ(A,PR)&DPR}
    \begin{aligned}
        \|\mathcal{P}\mathcal{Q}(A^{\varepsilon},\mathcal{P}\widetilde{R}^{\varepsilon}_{M})\|_{{L}^{2}(H^{s-2})}
        &\leqslant \|A^{\varepsilon}\|_{\widetilde{L}^2(B^{\frac{N}{2}})}\|\mathcal{P}\widetilde{R}^{\varepsilon}_{M}\|_{\widetilde{L}^\infty(H^{s-1})}.\\
    \end{aligned}
\end{equation}
And it's easily seen
\begin{equation}\label{DPR}
    \|\mathcal{D}\mathcal{P}\widetilde{R}^{\varepsilon}_{M}\|_{L^2 (H^{s-2})}
    \leqslant \|\mathcal{P}\widetilde{R}^{\varepsilon}_{M}\|_{L^2 (H^{s})}.
\end{equation}
As in the preceding section, proving smallness of the right-hand side of \eqref{estimate_of_psi}
is reduced to proving suitable bounds on $\mathcal{P}\widetilde{R}^\varepsilon_{M}$,
$\mathcal{P}\widetilde{R}^{t,\varepsilon}_{M}$, and $\mathcal{P}\widetilde{R}^{\varepsilon,M}$.
In fact, we have
\begin{lemma}\label{estimates_of_R^epsilon}
    \begin{align}
        &\|\mathcal{P}\widetilde{R}^{\varepsilon}_{M}\|_{\widetilde{L}^{\infty}(H^{s-1})\cap L^{2}(H^{s})}
        \lesssim M^{1-\delta}(X+X^2+Y^2),\label{PR}\\
        &\|\mathcal{P}\widetilde{R}^{\varepsilon}_{M}(0)\|_{H^{s-1}} \lesssim M^{1-\delta}(X_0+X_0^2+Y_0^2),\label{PR0}\\
        \begin{split}
            &\|\mathcal{P}\widetilde{R}^{t,\varepsilon}_{M}\|_{L^2(H^{s-2})}
            \lesssim M^{1-\delta}(X+Y+1)\bigg((X+Y)^2+X+Y\\
            &\quad\quad\quad\quad\quad\quad\quad +C(\mathcal{Y}^{\varepsilon})\varepsilon^{\delta}\left((X+Y)^2+(X+Y)^3\right)\bigg),\\
        \end{split}\label{PRt}\\
        &\|\mathcal{P}R^{\varepsilon,M}\|_{L^2(H^{s-2})}\lesssim M^{-\delta}(X^2+XY+X).\label{PR^M}
    \end{align}
\end{lemma}

Collecting  \eqref{PQWW_estimate}, \eqref{remainder_kernel},
\eqref{PQ(A,PR)&DPR}, \eqref{DPR}, and Lemma \ref{estimates_of_R^epsilon}
we arrive at
\begin{equation}
    \begin{aligned}
        W^{\varepsilon} \leqslant& C(\mathcal{Y}^{\varepsilon})\mathrm{e}^{(X^2+Y^2)}\times
        \bigg(\varepsilon M^{1-\delta}(X_0+X_0^2+Y_0^2)\\
        &+\varepsilon M^{1-\delta}((X+Y+1)(X+X^2+Y^2))\\
        &+\varepsilon M^{1-\delta}(X+Y+1)\big((X+Y)^2+X+Y+\varepsilon^{\delta}(X^3+X^2)\big)\\
        &+\varepsilon^{\delta}(X^3+X^2)+M^{-\delta}(X^2+XY+X)+YW^{\varepsilon}\bigg).
    \end{aligned}
\end{equation}
The term $YW^{\varepsilon}$ can be cancelled out using the smallness of the estimate \eqref{mutual_bound}.
The conclusion of the proposition follows by choosing $M=\varepsilon^{-1}$.
\end{proof}

Now we turn to the proof of Lemma \ref{estimates_of_R^epsilon}.
\begin{proof}[Proof of Lemma \ref{estimates_of_R^epsilon}]

    By definition, the spectral representation of the oscillating terms
    and the low-frequency-truncated parts thereof read
    \begin{equation}\nonumber
        \begin{aligned}
            R^{1,\varepsilon} &= \mathcal{Q}(\mathrm{e}^{-\frac{t}{\varepsilon}\mathcal{A}}\widetilde{V}^{\varepsilon},\mathrm{e}^{-\frac{t}{\varepsilon}\mathcal{A}}\widetilde{V}^{\varepsilon})
            =c_N^2\sum_{m}\sum_{\substack{\alpha,\beta\\ k+ l= m}}
            \mathrm{e}^{-\frac{t}{\varepsilon}(\lambda_{k}^{\alpha}+\lambda_{l}^{\beta})}
            V^{\alpha,\varepsilon}_{ k}V^{\beta,\varepsilon}_{ l}
            \mathcal{Q}(H^{\alpha}_{ k},H^{\beta}_{ l}),\\
            R^{2,\varepsilon} &= 2\mathcal{Q}(\mathrm{e}^{-\frac{t}{\varepsilon}\mathcal{A}}\widetilde{V}^{\varepsilon},\mathcal{U})
            = 2c_N\sum_{m}\sum_{\substack{\alpha\\ k+ l= m}}
            \mathrm{e}^{-\frac{t}{\varepsilon}\lambda_{k}^{\alpha}}V^{\alpha,\varepsilon}_{ k}
            \mathcal{Q}(H^{\alpha}_{ k},\widehat{\mathcal{U}}_{l}\mathrm{e}^{\mathrm{i} l\cdot x}),\\
            R^{3,\varepsilon} &= \mathcal{D}\mathrm{e}^{-\frac{t}{\varepsilon}\mathcal{A}}\widetilde{V}^{\varepsilon}
            = c_N\sum_{\alpha, k}\mathrm{e}^{-\frac{t}{\varepsilon}\lambda_{k}^{\alpha}}
            V^{\alpha,\varepsilon}_{ k}
            \mathcal{D}H^{\alpha}_{ k},\\
            R^{1,\varepsilon}_{M} &= c_N^2\sum_{m}\sum_{\substack{\alpha,\beta\\ k+ l= m\\|k|+|l|\leqslant M}}
            \mathrm{e}^{-\frac{t}{\varepsilon}(\lambda_{k}^{\alpha}+\lambda_{l}^{\beta})}
            V^{\alpha,\varepsilon}_{ k}V^{\beta,\varepsilon}_{ l}
            \mathcal{Q}(H^{\alpha}_{ k},H^{\beta}_{ l}),\\
            R^{2,\varepsilon}_{M} &= 2c_N\sum_{m}\sum_{\substack{\alpha\\ k+ l= m\\|k|+|l|\leqslant M}}
            \mathrm{e}^{-\frac{t}{\varepsilon}\lambda_{k}^{\alpha}}V^{\alpha,\varepsilon}_{ k}
            \mathcal{Q}(H^{\alpha}_{ k},\widehat{\mathcal{U}}_{l}\mathrm{e}^{\mathrm{i} l\cdot x}),\\
            R^{3,\varepsilon}_{M} &= c_N\sum_{\substack{\alpha, k\\|k|\leqslant M}}\mathrm{e}^{-\frac{t}{\varepsilon}\lambda_{k}^{\alpha}}
            V^{\alpha,\varepsilon}_{ k}
            \mathcal{D}H^{\alpha}_{ k}.
        \end{aligned}
    \end{equation}
    Direct calculations give
    \begin{equation}\nonumber
        \begin{aligned}
        \mathcal{Q}&(H^{\alpha}_{ k},H^{\beta}_{ l}) \\
        &=\frac{1}{2}\left(\begin{array}{c}
            \mathrm{i}\rho^{\circ} c^{\circ} m\!\cdot\!(\alpha \operatorname{sg}( k)\frac{ k}{| k|} + \beta \operatorname{sg}( l)\frac{ l}{| l|})\\

            \mathrm{i} m(c^{\circ})^2\alpha\beta \operatorname{sg}( k)\operatorname{sg}( l)\frac{ k\cdot l}{| k|| l|}
            +\mathrm{i} m(C_1(\rho^{\circ})^2 + C_2\frac{(\theta^{\circ})^2(p^{\circ}_{\theta})^2}{(\rho^{\circ})^2(e^{\circ}_{\theta})^2} +
            (C_3+C_4)\frac{\theta^{\circ}p^{\circ}_{\theta}}{e^{\circ}_{\theta}})\\

            \mathrm{i} k\!\cdot\! lc^{\circ}\frac{\theta^{\circ}p^{\circ}_{\theta}}{\rho^{\circ}e^{\circ}_{\theta}}(\alpha \operatorname{sg}( k)\frac{1}{| k|} + \beta \operatorname{sg}( l)\frac{1}{| l|})
            +\mathrm{i}c^{\circ}(\rho^{\circ}C_6 + \frac{\theta^{\circ}p^{\circ}_{\theta}}{\rho^{\circ}e^{\circ}_{\theta}}C_5)(\alpha \operatorname{sg}( k)| k| + \beta \operatorname{sg}( l)| l|)
        \end{array}\right)
        \mathrm{e}^{\mathrm{i} m\cdot x},\\
        \end{aligned}
    \end{equation}
    \begin{equation}\nonumber
        \begin{gathered}
            \mathcal{Q}(H^{\alpha}_{ k},\widehat{\mathcal{U}}_{ l}\mathrm{e}^{\mathrm{i} l\cdot x}) =
            \frac{1}{2}\left(\begin{array}{c}
                \mathrm{i} m\!\cdot\!(\rho^{\circ}\hat{\omega}_{ l}-c^{\circ}\alpha \operatorname{sg}( k)\frac{ k}{| k|}\frac{p^{\circ}_{\theta}}{p^{\circ}_{\rho}}\hat{\vartheta}_{ l})\\

                \mathrm{i}c^{\circ}\alpha \operatorname{sg}( k)\frac{( k\cdot l)\hat{\omega}_{ l} + ( k\cdot\hat{\omega}_{ l}) k}{| k|}
                +\mathrm{i}C_7 l\hat{\vartheta}_{ l} + \mathrm{i}C_8 k\hat{\vartheta}_{ l} + \mathrm{i}C_9 m\hat{\vartheta}_{ l}\\

                \mathrm{i}c^{\circ}\alpha \operatorname{sg}( k)\frac{ k\cdot l}{| k|}\hat{\vartheta}_{l}
                +\mathrm{i}\hat{\omega}_{ l}\!\cdot\! k\frac{\theta^{\circ}p^{\circ}_{\theta}}{\rho^{\circ}e^{\circ}_{\theta}}
                \mathrm{i}C_{10} l\!\cdot\!\hat{\omega}_{ l} + \mathrm{i}C_{11} c^{\circ}\alpha \operatorname{sg}( k)| k|\hat{\vartheta}_{l}
            \end{array}\right)
            \mathrm{e}^{i m\cdot x},\\
        \end{gathered}
    \end{equation}
    and
    \begin{equation}
        \mathcal{D}H^{\alpha}_{ k}=
        \left(\begin{array}{c}
            0\\
            -\frac{2\mu^\circ-\frac{N}{2}\mu^\circ+\lambda^\circ}{\rho^\circ}c^{\circ}\alpha \operatorname{sg}( k) k| k|\\
            -\frac{\kappa^{\circ}\theta^{\circ}p^{\circ}_{\theta}}{(\rho^{\circ})^2(e^{\circ}_{\theta})^2}| k|^2
        \end{array}\right)
        \mathrm{e}^{\mathrm{i} m\cdot x}.
    \end{equation}
    In the above expressions
    \begin{equation}
        \begin{gathered}
            C_7 = C_3\rho^{\circ} - C_2\frac{\theta^{\circ}(p^{\circ}_{\theta})^2}{\rho^{\circ}e^{\circ}_{\theta}p^{\circ}_{\rho}},\ \ \
            C_8 = C_2\rho^{\circ} - C_3\frac{\theta^{\circ}(p^{\circ}_{\theta})^2}{\rho^{\circ}e^{\circ}_{\theta}p^{\circ}_{\rho}},\ \ \
            C_9 = \frac{\theta^{\circ}p^{\circ}_{\theta}}{\rho^{\circ}e^{\circ}_{\theta}}C_4 - \rho^{\circ}\frac{p^{\circ}_{\theta}}{p^{\circ}_{\rho}}C_1\\
            C_{10} = \rho^{\circ}C_5 + \frac{\theta^{\circ}p^{\circ}_{\theta}}{\rho^{\circ}e^{\circ}_{\theta}}C_6,\ \ \
            C_{11} = C_6 - \frac{p^{\circ}_{\theta}}{p^{\circ}_{\rho}}C_1.
        \end{gathered}
    \end{equation}
    The terms $\widetilde{R}^{2,\varepsilon}_{M}$ and $\widetilde{R}^{3,\varepsilon}_{M}$ are defined as
    goes as frictionless as
    \begin{equation}\nonumber
        \begin{aligned}
            &\widetilde{R}^{2,\varepsilon}_{M} = -2c_N\sum_{m}\sum_{\substack{\alpha\\ k+ l= m\\|k|+|l|\leqslant M}}
            \frac{1}{\lambda_{k}^{\alpha}}
            \mathrm{e}^{-\frac{t}{\varepsilon}\lambda_{k}^{\alpha}}V^{\alpha,\varepsilon}_{ k}
            \mathcal{Q}(H^{\alpha}_{ k},\widehat{\mathcal{U}}_{l}\mathrm{e}^{\mathrm{i} l\cdot x}),\\
            &\widetilde{R}^{3,\varepsilon}_{M} = -c_N\sum_{\substack{\alpha, k\\|k|\leqslant M}}
            \frac{1}{\lambda_{k}^{\alpha}}
            \mathrm{e}^{-\frac{t}{\varepsilon}\lambda_{k}^{\alpha}}
            V^{\alpha,\varepsilon}_{ k}
            \mathcal{D}H^{\alpha}_{ k}.
        \end{aligned}
    \end{equation}
    However, the attempt to define $\widetilde{R}^{1,\varepsilon}_{M}$ in the same way would fail
    since $\lambda^{\alpha}_{k}+\lambda^{\beta}_{l}$ could vanish for certain aligned eigendirections.
    Fortunately, the projection $\mathcal{P}$ annihilates such bad cases.
    In fact, direct computations give
    \begin{equation}\nonumber
    \mathcal{P}R^{1,\varepsilon}
    =c_N^2\sum_{m}\sum_{\substack{\alpha,\beta\\ k+ l= m}}
        \mathrm{e}^{-\frac{t}{\varepsilon}(\lambda_{k}^{\alpha}+\lambda_{l}^{\beta})}
        V^{\alpha,\varepsilon}_{ k}V^{\beta,\varepsilon}_{ l}
        \mathcal{P}\mathcal{Q}(H^{\alpha}_{ k},H^{\beta}_{ l}),
\end{equation}
where

\begin{equation}\nonumber
    \begin{aligned}
    \mathcal{P}\mathcal{Q}(H^{\alpha}_{ k},H^{\beta}_{ l}) =&
    \frac{1}{2}
    \frac{\theta^{\circ}p^{\circ}_{\rho}}{(c^{\circ})^2 e^{\circ}_{\theta}}\left(p^{\circ}_{\rho\theta} - \frac{p^{\circ}_{\theta}e^{\circ}_{\rho\theta}}{e^{\circ}_{\theta}}
    + \frac{(p^{\circ}_{\theta})^2}{(\rho^{\circ})^2e^{\circ}_{\theta}} + \frac{\theta^{\circ}p^{\circ}_{\theta}p^{\circ}_{\theta\theta}}{(\rho^{\circ})^2e^{\circ}_{\theta\theta}}
    - \frac{\theta^{\circ}(p^{\circ}_{\theta})^2e^{\circ}_{\theta\theta}}{(\rho^{\circ})^2(e^{\circ}_{\theta})}\right)\\
    &\times
    (\lambda_{k}^{\alpha}+\lambda_{l}^{\beta})
    \left(\begin{array}{c}
        -i\frac{p^{\circ}_{\theta}}{p^{\circ}_{\rho}}\\

        0\\

        1
    \end{array}\right)
    \mathrm{e}^{i m\cdot x}.
    \end{aligned}
\end{equation}
Therefore, we may only define the $\text{Null}(\mathcal{A})$ part of “$\widetilde{R}^{1,\varepsilon}_{M}$".
With a slight abuse of notation, we write
\begin{equation}\nonumber
    \mathcal{P}\widetilde{R}^{1,\varepsilon}_{M} = -c_N^2\sum_{m}\sum_{\substack{\alpha,\beta\\ k+ l= m\\|k|+|l|\leqslant M}}
    \frac{\mathrm{e}^{-\frac{t}{\varepsilon}(\lambda_{k}^{\alpha}+\lambda_{l}^{\beta})}}
    {\lambda_{k}^{\alpha}+\lambda_{l}^{\beta}}
    V^{\alpha,\varepsilon}_{ k}V^{\beta,\varepsilon}_{ l}
    \mathcal{P}\mathcal{Q}(H^{\alpha}_{ k},H^{\beta}_{ l}),
\end{equation}

Next, we have by Lemma \ref{embedding}, Lemma \ref{convolution_estimate}
and Lemma \ref{truncated_frequency_estimate}
\begin{equation}
    \begin{aligned}
        &\|\mathcal{P}\widetilde{R}^{1,\varepsilon}_{M}\|_{\widetilde{L}^{\infty}(H^{s-1})}\leqslant M^{1-\delta}\|\widetilde{V}^{\varepsilon}\|^{2}_{\widetilde{L}^{\infty}(H^{\frac{N}{2}-1})},\
        \|\mathcal{P}\widetilde{R}^{1,\varepsilon}_{M}\|_{L^{2}(H^{s})}\leqslant M^{1-\delta}\|\widetilde{V}^{\varepsilon}\|_{\widetilde{L}^{\infty}(H^{\frac{N}{2}-1})}\|\widetilde{V}^{\varepsilon}\|_{\widetilde{L}^2(B^{\frac{N}{2}})},\\
        &\|\widetilde{R}^{2,\varepsilon}_{M}\|_{\widetilde{L}^{\infty}(H^{s-1})}\leqslant M^{1-\delta}\|\widetilde{V}^{\varepsilon}\|_{\widetilde{L}^{\infty}(H^{\frac{N}{2}-1})}\|\mathcal{U}\|_{\widetilde{L}^{\infty}(H^{\frac{N}{2}-1})},\
        \|\widetilde{R}^{2,\varepsilon}_{M}\|_{L^{2}(H^{s})}\leqslant M^{1-\delta}\|\widetilde{V}^{\varepsilon}\|_{\widetilde{L}^{\infty}(H^{\frac{N}{2}-1})}\|\mathcal{U}\|_{\widetilde{L}^2(B^{\frac{N}{2}})},\\
        &\|\widetilde{R}^{3,\varepsilon}_{M}\|_{\widetilde{L}^{\infty}(H^{s-1})\cap L^{2}(H^{s})}\leqslant M^{1-\delta}\|\widetilde{V}^{\varepsilon}\|_{\widetilde{L}^{\infty}(H^{\frac{N}{2}-1})\cap \widetilde{L}^2(B^{\frac{N}{2}})},
    \end{aligned}
\end{equation}
and therefore, \eqref{PR} is proved.
We can also treat \eqref{PR0} similarly.

The “high frequenct" terms can be expressed as
\begin{equation}
    \begin{aligned}
            &\mathcal{P}R^{1,\varepsilon,M} = \mathcal{P}\mathcal{Q}(\mathrm{e}^{-\frac{t}{\varepsilon}\mathcal{A}}\widetilde{V}^{\varepsilon},\mathrm{e}^{-\frac{t}{\varepsilon}\mathcal{A}}\widetilde{V}^{\varepsilon,M})
            + \mathcal{P}\mathcal{Q}(\mathrm{e}^{-\frac{t}{\varepsilon}\mathcal{A}}\widetilde{V}^{\varepsilon,M},\mathrm{e}^{-\frac{t}{\varepsilon}\mathcal{A}}\widetilde{V}^{\varepsilon}_{M}),\\
            &R^{2,\varepsilon,M} = 2\mathcal{Q}(\mathrm{e}^{-\frac{t}{\varepsilon}\mathcal{A}}\widetilde{V}^{\varepsilon},\mathcal{U}^{M})
            + 2\mathcal{Q}(\mathrm{e}^{-\frac{t}{\varepsilon}\mathcal{A}}\widetilde{V}^{\varepsilon,M},\mathcal{U}_{M}),\\
            &R^{3,\varepsilon,M} = \mathcal{D}\mathrm{e}^{-\frac{t}{\varepsilon}\mathcal{A}}\widetilde{V}^{\varepsilon,M},
    \end{aligned}
\end{equation}
whence also by Lemma \ref{convolution_estimate} and Lemma
\ref{truncated_frequency_estimate} we have
\begin{equation}
    \begin{aligned}
        \|\mathcal{P}R^{1,\varepsilon,M}\|_{L^2(H^{s-2})}
        &\lesssim \|\widetilde{V}^{\varepsilon}\|_{\widetilde{L}^{\infty}(H^{\frac{N}{2}-1})}\|\widetilde{V}^{\varepsilon,M}\|_{L^2(H^{s})}
        + \|\widetilde{V}^{\varepsilon,M}\|_{\widetilde{L}^{\infty}(H^{s-1})}\|\widetilde{V}^{\varepsilon}_{M}\|_{\widetilde{L}^2(B^{\frac{N}{2}})}\\
        &\lesssim M^{-\delta} (\|\widetilde{V}^{\varepsilon}\|_{\widetilde{L}^{\infty}(H^{\frac{N}{2}-1})}\|\widetilde{V}^{\varepsilon}\|_{\widetilde{L}^2(B^{\frac{N}{2}})}
        + \|\widetilde{V}^{\varepsilon}\|_{\widetilde{L}^{\infty}(H^{\frac{N}{2}-1})}\|\widetilde{V}^{\varepsilon}_{M}\|_{\widetilde{L}^2(B^{\frac{N}{2}})})\\
        &\lesssim M^{-\delta}X^2,\\
        \|R^{2,\varepsilon,M}\|_{L^2(H^{s-2})}
        &\lesssim\|\widetilde{V}^{\varepsilon}\|_{\widetilde{L}^{\infty}(H^{\frac{N}{2}-1})}\|\mathcal{U}^{M}\|_{L^2(H^{s})}
        + \|\widetilde{V}^{\varepsilon,M}\|_{\widetilde{L}^{\infty}(H^{s-1})}\|\mathcal{U}_{M}\|_{\widetilde{L}^2(B^{\frac{N}{2}})}\\
        &\lesssim M^{-\delta}XY,\\
        \|R^{3,\varepsilon,M}\|_{L^2(H^{s-2})}
        &\lesssim \|\widetilde{V}^{\varepsilon}\|_{L^2(H^{s})}
        \lesssim M^{-\delta}X.
    \end{aligned}
\end{equation}
Hence we know \eqref{PR^M} is true.

From \eqref{negative_time_relation} we compute
\begin{equation}\nonumber
    \begin{aligned}
        &\mathcal{P}R^{1,t,\varepsilon}_{M} =
         -c_N^2\sum_{m}\sum_{\substack{\alpha,\beta\\ k+ l= m\\|k|+|l|\leqslant M}}
        \frac{\mathrm{e}^{-\frac{t}{\varepsilon}(\lambda_{k}^{\alpha}+\lambda_{l}^{\beta})}}
        {\lambda_{k}^{\alpha}+\lambda_{l}^{\beta}}
        \partial_t (V^{\alpha,\varepsilon}_{ k}V^{\beta,\varepsilon}_{ l})
        \mathcal{P}\mathcal{Q}(H^{\alpha}_{ k},H^{\beta}_{ l}),\\
        &R^{2,t,\varepsilon}_{M}  = -2c_N\sum_{m}\sum_{\substack{\alpha\\ k+ l= m\\|k|+|l|\leqslant M}}
        \frac{1}{\lambda_{k}^{\alpha}}
        \mathrm{e}^{-\frac{t}{\varepsilon}\lambda_{k}^{\alpha}}
        \partial_t (V^{\alpha,\varepsilon}_{ k}
        \mathcal{Q}(H^{\alpha}_{ k},\widehat{\mathcal{U}}_{l}\mathrm{e}^{\mathrm{i} l\cdot x})),\\
        &R^{3,t,\varepsilon}_{M} = -c_N\sum_{\substack{\alpha, k\\|k|\leqslant M}}
        \frac{1}{\lambda_{k}^{\alpha}}
        \mathrm{e}^{-\frac{t}{\varepsilon}\lambda_{k}^{\alpha}}
        \partial_t (V^{\alpha,\varepsilon}_{ k})
        \mathcal{D}H^{\alpha}_{ k}.
    \end{aligned}
\end{equation}
Then, Lemma \ref{time_derivative_estimate} gives
\begin{equation}
    \begin{aligned}
        \|\mathcal{P}\widetilde{R}^{1,t,\varepsilon}_{M}\|_{L^2(H^{s-2})}
        &\leqslant M^{1-\delta}\|\widetilde{V}^{\varepsilon}\|_{\widetilde{L}^{\infty}(B^{\frac{N}{2}-1})}\|\partial_t \widetilde{V}^{\varepsilon}\|_{L^2_{T}(H^{\frac{N}{2}-2})}\\
        &\leqslant M^{1-\delta}(X+Y)\bigg((X+Y)^2+X+Y
        +C (\mathcal{Y}^{\varepsilon})\varepsilon^{\delta}\left((X+Y)^2+(X+Y)^3\right)\bigg),\\
    \end{aligned}
\end{equation}
\begin{equation}
    \begin{aligned}
        \|\widetilde{R}^{2,t,\varepsilon}_{M}\|_{L^2(H^{s-2})}
        &\leqslant M^{1-\delta}\|\mathcal{U}\|_{\widetilde{L}^{\infty}(H^{\frac{N}{2}-1})}\|\partial_t\widetilde{V}^{\varepsilon}\|_{L^2_{T}(H^{\frac{N}{2}-2})}
        + M^{1-\delta}\|\widetilde{V}^{\varepsilon}\|_{\widetilde{L}^{\infty}(B^{\frac{N}{2}-1})}\|\partial_t\mathcal{U}\|_{L^2_{T}(H^{\frac{N}{2}-2})}\\
        &\leqslant M^{1-\delta}Y\bigg((X+Y)^2+X+Y
         +C (\mathcal{Y}^{\varepsilon})\varepsilon^{\delta}\left((X+Y)^2+(X+Y)^3\right)\bigg)\\
        & +M^{1-\delta}(X+Y)(Y^2+Y),\\
    \end{aligned}
\end{equation}
\begin{equation}
    \begin{aligned}
        \|\widetilde{R}^{3,t,\varepsilon}_{M}\|_{L^2(H^{s-2})}
        &\leqslant M^{1-\delta}\|\partial_t\widetilde{V}^{\varepsilon}\|_{\widetilde{L}^2(B^{\frac{N}{2}})}
        \leqslant M^{1-\delta}\bigg((X+Y)^2+X+Y\\
        & +C (\mathcal{Y}^{\varepsilon})\varepsilon^{\delta}\big((X+Y)^2+(X+Y)^3\big)\bigg).
    \end{aligned}
\end{equation}
Collecting all terms, \eqref{PRt} is thus proved to hold.

\end{proof}
\section{A priori estimate for $\text{NSC}^\varepsilon$}\label{sec_priori}
In this section we give uniform in $\varepsilon$ a priori estimates for \eqref{nsc_epsilon}. For the well-posedness of \eqref{nsc_epsilon}, we refer to \cites{Chen-Gui-Jiang,Danchin_ARMA}.

As in  \cite{Chen-Gui-Jiang} and \cite{Danchin_ARMA}, we set
\begin{equation}
    \begin{aligned}
        &\rho_j = \Delta_{j} \tilde{\rho}^{\varepsilon},\
        u_{j}=\Delta_{j} \tilde{\u}^{\varepsilon},\
        \theta_{j} =\Delta_{j} \tilde{\theta}^{\varepsilon},\ \\
        &d_j =\Delta_{j} (\Lambda^{-1} \text{div } \tilde{\u}^{\varepsilon}), \
        \Omega_{j} = \Delta_{j} (\Lambda^{-1} \text{curl } \tilde{\u}^{\varepsilon}).\ \\
    \end{aligned}
\end{equation}
We then perform a localization to \eqref{nsc_epsilon} and rewrite it as
\begin{equation}
    \left\{ \begin{array}{lr}
        \partial_t \rho_{j} + \frac{1}{\varepsilon} \rho^{\circ} \Lambda d_{j} + \left(\widetilde{\u}^{\varepsilon} \!\cdot\! \nabla \widetilde{\rho}^{\varepsilon} \right)_{j}=F_{j}, \\
        \partial_t d_{j} + \frac{\nu^{\circ}}{\rho^{\circ}} \Lambda^{2} d_{j} - \frac{1}{\varepsilon} \frac{p^{\circ}_{\rho}}{\rho^{\circ}} \Lambda \rho_{j} -\frac{1}{\varepsilon}\frac{P^{\circ}_{\theta}}{\rho^{\circ}} \Lambda \theta_{j} +\left(\widetilde{\u}^{\varepsilon} \!\cdot\! \nabla d \right)_{j} = G_j,\\
        \partial_t \theta_{j} + \frac{\kappa^{\circ} }{\rho^{\circ} e^{\circ}_{\theta}} \Lambda^{2} \theta_{j} +\frac{1}{\varepsilon} \frac{\theta^{\circ} P^{\circ}_{\theta}}{e^{\circ}_{\theta} \rho^{\circ}}  \Lambda d_{j} +\left(\widetilde{\u}^{\varepsilon} \!\cdot\! \nabla \widetilde{\theta}^{\varepsilon} \right)_{j} = H_j,\\
        \partial_t \Omega_{j} - \frac{\mu^{\circ}}{\rho^{\circ}}\Delta \Omega_{j} = \widetilde{G}_{j},\\
        \u_{j} = -\Lambda^{-1} \nabla d_{j} - \Lambda^{-1} \text{div } \Omega_{j},
    \end{array}
    \right.
\end{equation}

where
\begin{equation}
    \begin{aligned}
        F = & - \tilde{\rho}^{\varepsilon} \nabla \!\cdot\!  \tilde{\u}^{\varepsilon},\\
        G = & \tilde{\u}^{\varepsilon} \!\cdot\! \nabla d - \Lambda^{-1} \operatorname{div} ( \tilde{\u}^{\varepsilon} \!\cdot\! \nabla  \tilde{\u}^{\varepsilon})
        + \frac{\nu^{\circ}}{\rho ^ {\circ}} \Lambda^{2} d + \Lambda^{-1} \operatorname{div} \left( \frac{\nabla \!\cdot\! S^{\varepsilon}}{\rho^{\varepsilon}} \right) 
        -\frac{1}{\varepsilon}\Lambda^{-1} \mathrm{div}  \left(\frac{p_{\rho}^{\varepsilon}}{\rho^{\varepsilon}} - \frac{p^{\circ}_{\rho}}{\rho^{\circ}} \right)\nabla  \tilde{\rho}^{\varepsilon}\\
        &-\frac{1}{\varepsilon}\Lambda^{-1} \mathrm{div} \left(\frac{p_{\theta}^{\varepsilon}}{\rho^{\varepsilon}} - \frac{p^{\circ}_{\theta}}{\rho^{\circ}} \right)\nabla  \tilde{\theta}^{\varepsilon},\\
        H = & -\frac{\kappa^{\circ}}{e^{\circ}_{\theta} \rho^{\circ}} \Delta  \tilde{\theta}^{\varepsilon} + \frac{1}{e_{\theta} \rho^{\varepsilon}} \nabla( \kappa^{\varepsilon} \nabla  \tilde{\theta}^{\varepsilon})
          + \frac{1}{\varepsilon} \Lambda^{-1} \mathrm{div}\left( \frac{\theta^{\circ} p^{\circ}_{\theta}}{e^{\circ}_{\theta} \rho^{\circ}} - \frac{\theta^{\varepsilon} p_{\theta}^{\varepsilon}}{e_{\theta}^{\varepsilon} \rho^{\varepsilon}}\right)\Lambda d
          + \frac{\varepsilon}{e_{\theta}^{\varepsilon} \rho^{\varepsilon}} S^{\varepsilon}: \nabla  \tilde{\u}^{\varepsilon},\\
        \widetilde{G} = & \Lambda^{-1} \text{curl } (  \tilde{\u}^{\varepsilon} \!\cdot\! \nabla  \tilde{\u}^{\varepsilon})
                        +\frac{\mu^{\circ}}{\rho^{\circ}} \Lambda^2 \Omega_{j} + \Lambda^{-1} \text{curl } \left( \nabla \!\cdot\! S^{\varepsilon} \right). 
    \end{aligned}
\end{equation}
\begin{proposition}\label{A_priori_estimates}
    Given a constant state $U^{\circ}= \left(\rho^{\circ},0,\theta^{\circ}\right)^{\mathrm{T}}$, where $\rho^{\circ}$ and $\theta^{\circ}$ are positive constants such that
    \begin{equation}
        p_{\rho}(\rho^{\circ},\theta^{\circ}) > 0,\ e_{\theta}(\rho^{\circ},\theta^{\circ}) >0.
    \end{equation}
    There exists two positive constants $\alpha_0 = \alpha_0(U^{\circ},N,p,e,\mu,\lambda,\kappa)$ and $M= M(U^{\circ},N,p,e,\mu,\lambda,\kappa)$ such that for some $\varepsilon_0 > 0$
    \begin{equation}
        \begin{aligned}
            \left\| \widetilde{U}^{in} \right\|_{\mathcal{E}^{\frac{N}{2}}_{\varepsilon_0 \nu^{\circ}} } \leqslant \alpha_0 ,
        \end{aligned}
    \end{equation}
    then the solution of system \eqref{nsc_epsilon} is uniformly bounded in $E^{\frac{N}{2}}_{\varepsilon\nu^{\circ}}$, whenever $\varepsilon \leqslant \varepsilon_0$. Moreover,
    we have
    \begin{equation}
        \left\| (\tilde{\rho}^{\varepsilon},\ \tilde{u}^{\varepsilon},\ \tilde{\theta}^{\varepsilon}) \right\|_{E_{\varepsilon \nu^{\circ}}^{\frac{N}{2}}}
        \leqslant M \alpha^{\circ}.
    \end{equation}
\end{proposition}
Throughout this section we denote $\| \cdot \|$ the $L^2(\mathbb{T}^{N}_a)$ norm and $< \cdot\  , \ \cdot >$ the $L^2(\mathbb{T}^{N}_a)$ inner product.

\subsection{Decay effect at low frequency}

We focus on the first three equations, namely
\begin{align}
    \partial_t \rho_{j} + \frac{1}{\varepsilon} \rho^{\circ} \Lambda d_{j} + \left(\widetilde{\u}^{\varepsilon} \!\cdot\! \nabla \widetilde{\rho}^{\varepsilon} \right)_{j}=F_{j},\label{rho_equation}   \\
    \partial_t d_{j} + \frac{\nu^{\circ}}{\rho^{\circ}} \Lambda^{2} d_{j} - \frac{1}{\varepsilon} \frac{p^{\circ}_{\rho}}{\rho^{\circ}} \Lambda \rho_{j} -\frac{1}{\varepsilon}\frac{p^{\circ}_{\theta}}{\rho^{\circ}} \Lambda \theta_{j} +\left(\widetilde{\u}^{\varepsilon} \!\cdot\! \nabla d \right)_{j} = G_j,\label{d_equation}\\
    \partial_t \theta_{j} + \frac{\kappa^{\circ} }{\rho^{\circ} e^{\circ}_{\theta}} \Lambda^{2} \theta_{j} +\frac{1}{\varepsilon} \frac{\theta^{\circ} p^{\circ}_{\theta}}{e^{\circ}_{\theta} \rho^{\circ}}  \Lambda d_{j} +\left(\widetilde{\u}^{\varepsilon} \!\cdot\! \nabla \widetilde{\theta}^{\varepsilon} \right)_{j} = H_j\label{theta_equation}.
\end{align}
Above three equations taking entropy inner product with $\left( \rho_j, d_j, \theta_j\right)^{\text{T}}$, more precisely 
\begin{equation}
    \begin{aligned}
        &\frac{1}{2} \frac{p^{\circ}_{\rho} }{\theta^{\circ} \rho^{\circ}} \frac{d}{dt} \left\| \rho_j \right\|^{2}
             + \frac{1}{\varepsilon} \frac{p^{\circ}_{\rho} }{\theta^{\circ}}  \left<\Lambda d_j, \rho_j\right>
             + \frac{p^{\circ}_{\rho} }{\theta^{\circ} \rho^{\circ}} \left< \left(\widetilde{\u}^{\varepsilon}\cdot \nabla \widetilde{\rho}^{\varepsilon}\right)_j ,\rho_j \right>
             = \frac{p^{\circ}_{\rho} }{\theta^{\circ} \rho^{\circ}} \left< F_j, \rho_j \right>,\\
        &\frac{1}{2} \frac{\rho^{\circ}}{\theta^{\circ}}\frac{d}{dt}\left\| d_{j} \right\|^{2} + \frac{\nu^{\circ}}{\theta^{\circ}} \left\| \Lambda d_j \right\|^2
        - \frac{1}{\varepsilon} \frac{p_{\rho}^{\circ}}{\theta^{\circ}}\left< \Lambda \rho_j, d_j \right>
        -\frac{1}{\varepsilon} \frac{p_{\theta}^{\circ}}{\theta^{\circ}}\left< \Lambda \theta_j, d_j \right>
        + \frac{\rho^{\circ}}{\theta^{\circ}}  \left< \left(\widetilde{\u}^{\varepsilon}\cdot \nabla d \right)_j ,d_j \right>
        = \frac{\rho^{\circ}}{\theta^{\circ}} \left< G_j , d_j \right>,\\
        &\frac{1}{2}\frac{ \rho^{\circ} e^{\circ}_{\theta}}{(\theta^{\circ})^2} \frac{d}{dt} \left\| \theta_{j} \right\|^{2}
        + \frac{\kappa^{\circ}}{(\theta^{\circ})^2} \left\| \Lambda \theta_j \right\|^2 + \frac{1}{\varepsilon} \frac{p_{\theta}^{\circ}}{\theta^{\circ}} \left< \Lambda d_j , \theta_j \right>
        + \frac{\rho^{\circ} e_{\theta}^{\circ}}{(\theta^{\circ})^2}   \left< \left(\widetilde{\u}^{\varepsilon}\cdot \nabla \widetilde{\theta}^{\varepsilon}\right)_j ,\theta_j \right>
        = \frac{\rho^{\circ} e_{\theta}^{\circ}}{(\theta^{\circ})^2} \left< H_j , \theta_j \right>,
    \end{aligned}
\end{equation}
and therefore
\begin{equation}\label{Energy1}
    \begin{aligned}
        &\frac{1}{2} \frac{d}{dt}
        \left(
            \frac{p^{\circ}_{\rho} }{\theta^{\circ} \rho^{\circ}}
            \left\| \rho_{j} \right \|^{2} +
            \frac{\rho^{\circ}}{\theta^{\circ}}\left\| d_{j} \right\|^{2}
            + \frac{ \rho^{\circ} e^{\circ}_{\theta}}{(\theta^{\circ})^2} \left\| \theta_{j} \right\|^{2}
            \right) + \frac{\nu^{\circ}}{\theta^{\circ}}
         \left\| \Lambda d_{j} \right\|^{2}
          + \frac{\kappa^{\circ}}{(\theta^{\circ})^2} \left\| \Lambda \theta_{j} \right\|^{2} \\
        &= -\left< \frac{p^{\circ}_{\rho}}{\theta^{\circ} \rho^{\circ}}
                \left( \widetilde{\u}^{\varepsilon} \!\cdot\! \nabla \widetilde{\rho}^{\varepsilon} \right)_{j} ,  \rho_{j}
          \right>
          + \frac{p^{\circ}_{\rho}}{\theta^{\circ} \rho^{\circ} } \left< F_{j}, \rho_{j} \right>
          -\left<
            \frac{\rho^{\circ}}{\theta^{\circ}} \left( \widetilde{\u}^{\varepsilon} \!\cdot\! \nabla d \right)_{j} ,  d_{j}
          \right> +\frac{\rho^{\circ}}{\theta^{\circ}}  \left< G_{j}, d_{j} \right> \\
        & -\left<
            \frac{\rho^{\circ} e^{\circ}_{\theta}}{(\theta^{\circ})^2} \left( \widetilde{\u}^{\varepsilon} \!\cdot\! \nabla \widetilde{\theta}^{\varepsilon} \right)_{j} ,  \theta_{j}
          \right> +\frac{\rho^{\circ} e^{\circ}_{\theta}}{(\theta^{\circ})^2}  \left< H_{j}, \theta_{j} \right>.
    \end{aligned}
\end{equation}
In order to get the decay estimate for $\rho$ at low frequency, $\varepsilon \nu^{\circ} \Lambda$ \eqref{rho_equation} $-(\rho^{\circ})^2$ \eqref{d_equation} we have
\begin{equation}
    \begin{aligned}
        \partial_t \left( \varepsilon \nu^{\circ} \Lambda \rho_j - (\rho^{\circ})^{2} d_j \right)
        + \frac{1}{\varepsilon}p_{\rho}^{\circ} \rho^{\circ} \Lambda \rho_j
        + \frac{1}{\varepsilon}p_{\theta}^{\circ} \rho^{\circ} \Lambda \theta_j
        +\varepsilon \nu^{\circ} \Lambda (\widetilde{\u}^{\varepsilon} \cdot \nabla \widetilde{\rho}^{\varepsilon})_j - (\rho^{\circ})^{2} \left(\widetilde{\u}^{\varepsilon}\cdot \nabla d\right)_j
        \\ = \varepsilon \nu^{\circ} \Lambda F_j - (\rho^{\circ})^{2} G_j.
    \end{aligned}
\end{equation}
Basic $L^2$ estimate then leads to
\begin{equation}
    \begin{aligned}
        &\frac{1}{2} \frac{d}{dt} \left\| \varepsilon \nu^{\circ} \Lambda \rho_{j}
        - ({\rho^{\circ}})^{2} d_{j} \right\|^{2}
        + \rho^{\circ} \nu^{\circ} p^{\circ}_{\rho} \left\| \Lambda \rho_{j} \right\|^2
        + \rho^{\circ} \nu^{\circ} p^{\circ}_{\theta} \left< \Lambda \theta_{j}, \Lambda \rho_{j} \right>
        - \frac{1}{\varepsilon}(\rho^{\circ})^{3} p_{\rho}^{\circ} \left< \Lambda \rho_j , d_j \right>
        \\ &- \frac{1}{\varepsilon}(\rho^{\circ})^{3} p_{\theta}^{\circ} \left< \Lambda \theta_j , d_j \right>
        =  \varepsilon \nu^{\circ} (\rho^{\circ})^{2}
        \left(
            \left< \Lambda \left( \widetilde{\u}^{\varepsilon} \!\cdot\! \nabla \widetilde{\rho}^{\varepsilon} \right)_{j} , d_{j} \right>
           +\left<  \left( \widetilde{\u}^{\varepsilon} \!\cdot\! \nabla d \right)_{j} , \Lambda \rho_{j} \right>
         \right)
        -\varepsilon \nu^{\circ} (\rho^{\circ})^{2} \left< \Lambda F_j ,  d_{j} \right> \\
        &-\varepsilon \nu^{\circ} (\rho^{\circ})^{2}   \left< G_{j} , \Lambda \rho_{j} \right>
        -\varepsilon^{2} (\nu^{\circ})^2 \left< \Lambda \left( \widetilde{\u}^{\varepsilon} \!\cdot\! \nabla \widetilde{\rho}^{\varepsilon} \right)_{j} , \Lambda \rho_{j} \right>
        +\varepsilon^{2} (\nu^{\circ})^2 \left< \Lambda F_{j}, \Lambda \rho_{j} \right>
        - (\rho^{\circ})^{4} \left< \left( \widetilde{\u}^{\varepsilon} \!\cdot\! \nabla d \right)_{j} , d_{j} \right>\\
        &+ (\rho^{\circ})^{4} \left< G_{j}, d_{j} \right>.
    \end{aligned}
\end{equation}
Moreover, using again \eqref{rho_equation} and \eqref{theta_equation} to subtract the terms $\left< \Lambda \rho_j , d_j\right> $ and $\left< \Lambda \theta_j , d_j\right> $ from above equation then we have
\begin{equation}\label{Energy2}
    \begin{aligned}
        &\frac{1}{2} \frac{d}{dt}
        \left(
         \left\| \varepsilon \nu^{\circ} \Lambda \rho_{j}- ({\rho^{\circ}})^{2} d_{j} \right\|^{2}
        +p^{\circ}_{\rho} (\rho^{\circ})^2  \left\| \rho_{j} \right\|^{2}
        +\frac{(\rho^{\circ})^4 e^{\circ}_{\theta}}{\theta^{\circ}} \left\| \theta_{j} \right\|^{2}
        \right)
        + \nu^{\circ} p^{\circ}_{\rho} \rho^{\circ} \left\| \Lambda \rho_{j} \right\|^{2}\\
        &+ \nu^{\circ} p^{\circ}_{\theta} \rho^{\circ} \left< \Lambda \rho_{j} ,\Lambda \theta_{j} \right>
        + \frac{(\rho^{\circ})^{3} \kappa^{\circ}}{\theta^{\circ}} \left\| \Lambda \theta_{j} \right\|^{2}
        =
        - \varepsilon \nu^{\circ} (\rho^{\circ})^{2}
        \left(
            \left< \Lambda \left( \widetilde{\u}^{\varepsilon} \!\cdot\! \nabla \widetilde{\rho}^{\varepsilon} \right)_{j} , d_{j} \right>
           +\left<  \left( \widetilde{\u}^{\varepsilon} \!\cdot\! \nabla d \right)_{j} , \Lambda \rho_{j} \right>
        \right)\\
        &+\varepsilon \nu^{\circ} (\rho^{\circ})^{2} \left< \Lambda F_j ,  d_{j} \right>
        +\varepsilon \nu^{\circ} (\rho^{\circ})^{2}   \left< G_{j} , \Lambda \rho_{j} \right>
        -\varepsilon^{2} (\nu^{\circ})^2 \left< \Lambda \left( \widetilde{\u}^{\varepsilon} \!\cdot\! \nabla \widetilde{\rho}^{\varepsilon} \right)_{j} , \Lambda \rho_{j} \right>
        +\varepsilon^{2} (\nu^{\circ})^2 \left< \Lambda F_{j}, \Lambda \rho_{j} \right> \\
        &- (\rho^{\circ})^{4} \left< \left( \widetilde{\u}^{\varepsilon} \!\cdot\! \nabla d \right)_{j} , d_{j} \right>
        + (\rho^{\circ})^{4} \left< G_{j}, d_{j} \right>
        + P^{\circ}_{\rho} (\rho)^{2} \left< F_{j} , \rho_{j} \right>
        - P^{\circ}_{\rho} (\rho)^{2} \left< \left(\widetilde{\u}^{\varepsilon} \!\cdot\! \nabla \widetilde{\rho}^{\varepsilon} \right)_j , \rho_{j} \right> \\
        &+ \frac{(\rho^{\circ})^{4}e^{\circ}_{\theta}}{\theta^{\circ}} \left< H_j, \theta_{j} \right>
        -\frac{(\rho^{\circ})^4 e^{\circ}_{\theta}}{\theta^{\circ}}  \left< \left(\widetilde{\u}^{\varepsilon} \!\cdot\! \nabla \widetilde{\theta}^{\varepsilon} \right)_j , \theta_{j} \right>.
    \end{aligned}
\end{equation}
Denote the right-hand side of \eqref{Energy1} $R_A$ and the right-hand side of \eqref{Energy2} $R_B$,
multiply \eqref{Energy1} by $K (\rho^\circ)^3 \theta^{\circ}$ and plus \eqref{Energy2}, then we obtain
\begin{equation}
    \begin{aligned}
        &\frac{1}{2} \frac{d}{dt} \left(
            \left\| \varepsilon \nu^{\circ} \Lambda \rho_{j} - ({\rho^{\circ}})^{2} d_{j} \right\|^{2}
           + (1 + K) p^{\circ}_{\rho} (\rho^{\circ})^2  \left\| \rho_{j} \right\|^{2}
           + (1 + K) \frac{(\rho^{\circ})^4 e^{\circ}_{\theta}}{\theta^{\circ}} \left\| \theta_{j} \right\|^{2}
           + K (\rho^{\circ})^4 \left\| d_{j} \right\|^{2}
           \right) \\
        &  + \nu^{\circ} p^{\circ}_{\rho} \rho^{\circ} \left\| \Lambda \rho_{j} \right\|^{2}
        + \nu^{\circ} p^{\circ}_{\theta} \rho^{\circ} \left< \Lambda \rho_{j} ,\Lambda \theta_{j} \right>
        + (1 + K)\frac{(\rho^{\circ})^{3} \kappa^{\circ}}{\theta^{\circ}} \left\| \Lambda \theta_{j} \right\|^{2}
        + K (\rho^{\circ})^{3} \nu^{\circ} \left\| \Lambda d_{j} \right\|^{2}\\
        &= R_A + K (\rho^\circ)^3 \theta^{\circ}R_B.
    \end{aligned}
\end{equation}
For convenience, we set
\begin{equation*}
    \begin{aligned}
        &f^{2}_{j} \equiv
            \left\| \varepsilon \nu^{\circ} \Lambda \rho_{j}- ({\rho^{\circ}})^{2} d_{j} \right\|^{2}
           +( 1 + K )\left(p^{\circ}_{\rho} (\rho^{\circ})^2  \left\| \rho_{j} \right\|^{2}
           +\frac{(\rho)^4 e^{\circ}_{\theta}}{\theta^{\circ}} \left\| \theta_{j} \right\|^{2} \right)
           + K (\rho^{\circ})^{4}\left\| \Lambda d_{j} \right\|^{2},\\
        &g^{2}_{j} \equiv
        K (\rho^{\circ})^{3} \nu^{\circ} \left\| \Lambda d_{j} \right\|^{2}
        + \frac{1}{2} \nu^{\circ} p^{\circ}_{\rho} \rho^{\circ}  \left\| \Lambda \rho_{j} \right\|^{2}
        + \rho^{\circ} \left[ ( 1 + K ) \frac{\kappa^{\circ}}{\theta^{\circ}} (\rho^{\circ})^{2} - \frac{1}{2} \nu^{\circ} \frac{ (p^{\circ}_{\theta})^2}{p^{\circ}_{\rho}} \right]\left\| \Lambda \theta_{j} \right\|^{2} .
    \end{aligned}
\end{equation*}
By Hölder's and Young's inequality, we have
\begin{equation}
    |2\varepsilon \nu^{\circ} (\rho^{\circ})^{2} \left< \Lambda \rho_{j},d_j \right>|
    \leqslant
     \frac{1}{2} (\nu^{\circ})^{2} \left\| \varepsilon \Lambda \rho_{j} \right\|^{2}
     + 2(\rho^{\circ})^{4} \left\|d_j \right\|^{2},
\end{equation}
then it's easy to see that
\begin{equation}
    \begin{aligned}
        &f^{2}_{j} \geq \frac{1}{2} \left\| \varepsilon \nu^{\circ} \Lambda \rho_{j}\right\|^{2}
        + (K-1) (\rho^{\circ})^{4} \left\| d_j \right\|^{2}
        + ( 1 + K )
        \left(p^{\circ}_{\rho} (\rho^{\circ})^2  \left\| \rho_{j} \right\|^{2}
        + \frac{(\rho^{\circ})^4 e^{\circ}_{\theta}}{\theta^{\circ}} \left\| \theta_{j} \right\|^{2}
        \right),\\
        &f^{2}_{j} \leqslant \frac{3}{2} \left\| \varepsilon \nu^{\circ} \Lambda \rho_{j}\right\|^{2}
        + ( K + 3 ) (\rho^{\circ})^{4} \left\| d_j \right\|^{2}
        + ( 1 + K )
        \left(p^{\circ}_{\rho} (\rho^{\circ})^2  \left\| \rho_{j} \right\|^{2}
        + \frac{(\rho^{\circ})^4 e^{\circ}_{\theta}}{\theta^{\circ}} \left\| \theta_{j} \right\|^{2}
        \right),\\
    \end{aligned}
\end{equation}
also we need to choose $K$ such that $K > 1 $, and $( 1 + K ) \frac{\kappa^{\circ}}{\theta^{\circ}} (\rho^{\circ})^{2} - \frac{1}{2} \nu^{\circ} \frac{ (p^{\circ}_{\theta})^2}{p^{\circ}_{\rho}}>0$,
to ensure that 
\begin{equation}
    \begin{aligned}
        f_j^2 \approx \left\| \varepsilon \nu^{\circ} \Lambda \rho_j \right\|^2 + \left\| \rho_j \right\|^2 + \left\| d_j \right\|^2 + \left\| \theta_j \right\|^2,\\
        g_j^2
        \approx (\left\| \Lambda \rho_j \right\|_{L^2} + \left\| \Lambda \theta_j \right\|_{L^2}+\left\| \Lambda d_j \right\|_{L^2}).
    \end{aligned}
\end{equation}

Note that by Bernstein's Lemma we have
\begin{equation}
    \left\| \varepsilon \nu^{\circ} \Lambda \rho_{j}\right\| \approx \varepsilon \nu^{\circ}2^{j}\left\| \rho_{j}\right\|,
\end{equation}
therefore for low frequency $j \leqslant j_0 := - \log_2(\varepsilon \nu^{\circ})$ we have
\begin{equation}
    f^2_j \approx \left\| \rho_j \right\|^2 + \left\| d_j \right\|^2 +\left\| \theta_j \right\|^2 , \ \ g^2_{j} \gtrsim 2^{j} f^2_{j},
\end{equation}
and for high frequency $j > j_0 $ we have
\begin{equation}
    f^2_j \approx \left\| \varepsilon \nu^{\circ} \Lambda \rho_j \right\|^2 + \left\| d_j \right\|^2 +\left\| \theta_j \right\|^2,
    \ \   g^2_{j} \gtrsim \left( \frac{1}{\varepsilon \nu^{\circ}}\right)^{2} f^2_{j}.
\end{equation}

Now we turn to the estimate of right-hand side, using Lemma \ref{convection_estimates} we have
\begin{equation}
    \begin{aligned}
        (\nu^{\circ})^2 \varepsilon^2 (\widetilde{\rho}^{\circ})^{2}
        \left|
            \left< \Lambda \left( \widetilde{\u}^{\varepsilon} \!\cdot\! \nabla \widetilde{\rho}^{\varepsilon} \right)_{j} , \Lambda \rho_{j} \right>
        \right|
        \lesssim \
        &\varepsilon^{2} (\nu^{\circ})^2 C_j 2^{-j(\widetilde{\psi}^{\frac{N}{2}-1,\frac{N}{2}+1}(j)-1)} (\widetilde{\phi}^{\frac{N}{2}-1,\frac{N}{2}+1}(j))^{-1}\\
        &\left\| \nabla \widetilde{\u}^{\varepsilon} \right\|_{B^{\frac{N}{2}}} \left\| \widetilde{\rho}^{\varepsilon} \right\|_{B^{\frac{N}{2}-1,\frac{N}{2}+1}_{\varepsilon \nu^{\circ}}} \left\| \Lambda \rho_{j} \right\|_{L^2},
    \end{aligned}
\end{equation}
and
\begin{equation}\label{convection_terms}
    \begin{aligned}
        &\begin{aligned}
        \left|
            \left<   \left( \widetilde{\u}^{\varepsilon} \!\cdot\! \nabla \widetilde{\rho}^{\varepsilon} \right)_{j} ,   \rho_{j} \right>
        \right|
        \lesssim \
        &C_j 2^{-j(\widetilde{\psi}^{\frac{N}{2}-1,\frac{N}{2}+1}(j))} (\widetilde{\phi}^{\frac{N}{2}-1,\frac{N}{2}+1}(j))^{-1}
        \left\| \nabla \widetilde{\u}^{\varepsilon} \right\|_{B^{\frac{N}{2}}} \left\| \widetilde{\rho}^{\varepsilon} \right\|_{B^{\frac{N}{2}-1,\frac{N}{2}+1}_{\varepsilon \nu^{\circ}}} \left\|   \rho_{j} \right\|_{L^2},\\
        \left|
            \left<   \left( \widetilde{\u}^{\varepsilon} \!\cdot\! \nabla d \right)_{j} ,   d_{j} \right>
        \right|
        \lesssim \
        &C_j 2^{-j(\widetilde{\psi}^{\frac{N}{2}-1,\frac{N}{2}}(j))} (\widetilde{\phi}^{\frac{N}{2}-1,\frac{N}{2}}(j))^{-1}
        \left\| \nabla \widetilde{\u}^{\varepsilon} \right\|_{B^{\frac{N}{2}}} \left\| d \right\|_{B^{\frac{N}{2}-1,\frac{N}{2}}_{\varepsilon \nu^{\circ}}} \left\|   d_{j} \right\|_{L^2},\\
        \left|
            \left<   \left( \widetilde{\u}^{\varepsilon} \!\cdot\! \nabla \widetilde{\theta}^{\varepsilon} \right)_{j} ,   \theta_{j} \right>
        \right|
        \lesssim \
        &C_j 2^{-j(\widetilde{\psi}^{\frac{N}{2}-1,\frac{N}{2}}(j))} (\widetilde{\phi}^{\frac{N}{2}-1,\frac{N}{2}}(j))^{-1}
        \left\| \nabla \widetilde{\u}^{\varepsilon} \right\|_{B^{\frac{N}{2}}} \left\| \widetilde{\theta}^{\varepsilon} \right\|_{B^{\frac{N}{2}-1,\frac{N}{2}}_{\varepsilon \nu^{\circ}}} \left\|   \theta_{j} \right\|_{L^2},\\
    \end{aligned}\\
    &\begin{aligned}
        \nu^{\circ} \varepsilon (\rho^{\circ})^{2} \bigg|
            \left< \Lambda \left( \widetilde{\u}^{\varepsilon} \!\cdot\! \nabla \widetilde{\rho}^{\varepsilon} \right)_{j} , d_{j} \right>
           &+\left<  \left( \widetilde{\u}^{\varepsilon} \!\cdot\! \nabla d \right)_{j} , \Lambda \rho_{j} \right>
         \bigg|
         \lesssim \nu^{\circ} \varepsilon (\rho^{\circ})^{2} C_{j}\left\| \nabla \widetilde{\u}^{\varepsilon} \right\|_{B^{\frac{N}{2}}}
         \\ &
         \times\left(2^{-j(\widetilde{\psi}^{\frac{N}{2}-1,\frac{N}{2}}(j))} (\widetilde{\phi}^{\frac{N}{2}-1,\frac{N}{2}}(j))^{-1}
        \left\| d \right\|_{B^{\frac{N}{2}-1,\frac{N}{2}}_{\varepsilon \nu^{\circ}}}
        \left\| \Lambda \rho_{j} \right\|_{L^2}\right.
        \\ &
        \left.+2^{-j(\widetilde{\psi}^{\frac{N}{2}-1,\frac{N}{2}+1}(j)-1)} (\widetilde{\phi}^{\frac{N}{2}-1,\frac{N}{2}+1}(j))^{-1}
        \left\| \widetilde{\rho}^{\varepsilon} \right\|_{B^{\frac{N}{2}-1,\frac{N}{2}+1}_{\varepsilon \nu^{\circ}}}
        \left\|  d_{j} \right\|_{L^2}\right).
    \end{aligned}
\end{aligned}
\end{equation}
For the estimate involving $F$, $G$, $H$,
using the relation
\begin{equation}
    \begin{aligned}
        2^j \varepsilon \nu^{\circ} \leqslant 1,\ \ j\leqslant j_0,\\
        2^{-j} \leqslant \varepsilon \nu^{\circ},\ \ j > j_0,
    \end{aligned}
\end{equation} and using Hölder's inequality we have
\begin{equation} \label{source_terms}
    \begin{aligned}
         &|  \left< F_{j}, \rho_{j} \right>|
         \leqslant
         \left\| F_{j} \right\| \left\| \rho_{j} \right\|, \ \ j \leqslant j_0;\
         |  \left< F_{j}, \rho_{j} \right>|
         \leqslant
         \left\|\varepsilon \nu^{\circ} \Lambda F_{j} \right\| \left\| \varepsilon \nu^{\circ} \Lambda \rho_{j} \right\|,\  j > j_0;\ \\
         &|\varepsilon \nu^{\circ} \left< \Lambda F_{j},  d_{j} \right> |
         \leqslant
         \left\|  F_{j} \right\| \left\|  d_{j} \right\|, \ \ j \leqslant j_0;\
         |\varepsilon \nu^{\circ} \left< \Lambda F_{j},  d_{j} \right> |
         \leqslant
         \left\| \varepsilon \nu^{\circ} F_{j} \right\| \left\|  d_{j} \right\|, \ \ j > j_0;\\
         &| \varepsilon \nu^{\circ} (\rho^{\circ})^{2}   \left< G_{j} , \Lambda \rho_{j} \right>|
         \leqslant \left\| G_{j} \right\| \left\|  \rho_{j} \right\|, \ \ j \leqslant j_0;\
         | \varepsilon \nu^{\circ} (\rho^{\circ})^{2}   \left< G_{j} , \Lambda \rho_{j} \right> |
         \leqslant \left\| G_{j} \right\| \left\|\varepsilon \nu^{\circ}  \rho_{j} \right\|, \ \ j >j_0;\\
         &|\varepsilon^{2} (\nu^{\circ})^2 \left< \Lambda F_{j},\Lambda \rho_{j} \right>|
         \leqslant
         \left\|  F_{j} \right\| \left\|  \rho_{j} \right\|, \ \ j \leqslant j_0;\
         |\varepsilon^{2} (\nu^{\circ})^2  \left< \Lambda F_{j},\Lambda \rho_{j} \right>|
         \leqslant
         \left\| \varepsilon \nu^{\circ} \Lambda F_{j} \right\| \left\| \varepsilon \nu^{\circ} \Lambda \rho_{j} \right\|, j > j_0;\\
         &|  \left< G_{j}, d_{j} \right>|
         \leqslant
         \left\| G_{j} \right\| \left\| d_{j} \right\|, \ \
         |  \left< H_{j}, \theta_{j} \right>|
         \leqslant
         \left\| H_{j} \right\| \left\| \theta_{j} \right\|.\ \\\
    \end{aligned}
\end{equation}
In summary, we have for low frequency
\begin{equation}\label{low}
    \begin{aligned}
        \frac{1}{2} \frac{d}{dt} f^{2}_{j} + 2^{2j} f^{2}_{j}
        \lesssim
        &f_{j} \bigg(
            \left\| \left(F_j,G_j , H_j\right)\right\|
            + C_j\left\|
             \nabla \widetilde{\u}^{\varepsilon}
            \right\|_{B^{\frac{N}{2}}} 2^{-j(\frac{N}{2}-1)}\\
            &\left. \times \left(
                \left\| \widetilde{\rho}^{\varepsilon} \right\|_{B^{\frac{N}{2}-1,\frac{N}{2}+1}_{\varepsilon \nu^{\circ}}}
            +\left\|
            \left( d,\ \widetilde{\theta}^{\varepsilon} \right)
            \right\|_{B^{\frac{N}{2}-1,\frac{N}{2}}_{\varepsilon \nu^{\circ}}}
            \right)
            \right),
    \end{aligned}
\end{equation}
and for high frequency
\begin{equation}\label{high}
    \begin{aligned}
        \frac{1}{2} \frac{d}{dt} f^{2}_{j} + \left(\frac{1}{\varepsilon \nu^{\circ}}\right)^{2} f^{2}_{j}
        \lesssim
        &f_{j} \bigg(
            \left\| \left(\varepsilon \nu^{\circ} \Lambda F_j,G_j ,H_j\right)\right\|
            + C_j 2^{-j{\frac{N}{2}}}\frac{1}{\varepsilon \nu^{\circ}}
            \left\| \nabla \widetilde{\u}^{\varepsilon} \right\|_{B^{\frac{N}{2}}} \\
            &\left.\times \left(
                \left\| \widetilde{\rho}^{\varepsilon} \right\|_{B^{\frac{N}{2}-1,\frac{N}{2}+1}_{\varepsilon \nu^{\circ}}}
            +\left\|
            \left(
                d,\ \widetilde{\theta}^{\varepsilon}
            \right)
            \right\|_{B^{\frac{N}{2}-1,\frac{N}{2}}_{\varepsilon \nu^{\circ}}}
            \right)
            \right).
    \end{aligned}
\end{equation}
\subsection{Smoothing effect at high frequency}
Now we state the estimate for $\Omega$ and smoothing effect at high frequency.
For $j> j_0$ we have
\begin{equation}
    \left\{ \begin{array}{lr}
        \partial_t d_{j} + \frac{\nu^{\circ}}{\rho^{\circ}} \Lambda^{2} d_{j}  -\frac{1}{\varepsilon}\frac{p^{\circ}_{\theta}}{\rho^{\circ}} \Lambda \theta_{j} +\left(\widetilde{\u}^{\varepsilon} \!\cdot\! \nabla d \right)_{j} = G_j + \frac{1}{\varepsilon} \frac{p^{\circ}_{\rho}}{\rho^{\circ}} \Lambda \rho_{j},\\
        \partial_t \theta_{j} + \frac{\kappa^{\circ} }{\rho^{\circ} e^{\circ}_{\theta}} \Lambda^{2} \theta_{j} +\frac{1}{\varepsilon} \frac{\theta^{\circ} p^{\circ}_{\theta}}{e^{\circ}_{\theta} \rho^{\circ}}  \Lambda d_{j} +(\widetilde{\u}^{\varepsilon} \!\cdot\! \nabla \widetilde{\theta}^{\varepsilon} )_{j} = H_j,\\
    \end{array}
    \right.
\end{equation}
hence, we obtain
\begin{equation}
    \begin{aligned}
    & \frac{1}{2} \frac{d}{d t}
    \left(
        \left\|d_j\right\|_{L^2}^2
        +\frac{{e^{ \circ }_{\theta}}}{  { \theta^{\circ} } }
        \left\|\theta_{j}\right\|_{L^2}^2
    \right)
        + \frac{\nu^{\circ}}{{ \rho^{\circ} } }\left\|\Lambda d_j\right\|_{L^2}^2
        +\frac{{e^{ \circ }_{\theta}} \kappa^{\circ}}{{ \rho^{\circ} } { \theta^{\circ} } }\left\|\Lambda \theta_{j}\right\|_{L^2}^2
         \\
        & =
        \frac{1}{\varepsilon} \frac{{p^{\circ}_\rho}}{{ \rho^{\circ} } } \left<\Lambda \rho_{j} , d_j\right>
        -\left<(\widetilde{\u}^{\varepsilon} \!\cdot\! \nabla d)_{j} , d_j\right>
        -\frac{{e^{ \circ }_{ {\theta} }}}{  { \theta^{\circ} } }
        \left<( \widetilde{\u}^{\varepsilon} \!\cdot\! \nabla \widetilde{\theta}^{\varepsilon} )_j , \theta_{j}\right>
        +\left<G_j, d_j\right>
        +\frac{{e^{ \circ }_{ \theta }}}{ { \theta^{ \circ } } }
        \left<H_j , \theta_{j}\right>.
    \end{aligned}
\end{equation}
Note that
\begin{equation}
    \begin{aligned}
        |\frac{1}{\varepsilon} \frac{{p^{\circ}_\rho}}{{ \rho^{\circ} } } \left<\Lambda \rho_{j} , d_j\right>|
        \lesssim  \frac{1}{\varepsilon \nu^{\circ}} 2^{j}  \left\| \rho_j \right\|
        \left\| d_j \right\|
        \lesssim \frac{1}{\varepsilon \nu^{\circ}}2^{-js} \left\| \widetilde{\rho}^{\varepsilon} \right\|_{B^{\frac{N}{2}+1}}\left\| d_j \right\|,
    \end{aligned}
\end{equation}
and the estimates for last four terms has been achieved in \eqref{convection_terms} and \eqref{source_terms}.
Then we have for high frequency
\begin{equation}
    \begin{aligned}
        & \frac{1}{2} \frac{d}{d t}
    \left(
        \left\|d_j\right\|_{L^2}^2
        +\frac{{e^{ \circ }_{\theta}}}{{ \rho^{\circ} }  { \theta^{\circ} } }
        \left\|\theta_{j}\right\|_{L^2}^2
    \right)
        + \frac{\nu}{{ \rho^{\circ} } }\left\|\Lambda d_j\right\|_{L^2}^2
        +\frac{{e^{ \circ }_{\theta}} \kappa^{\circ}}{{ \rho^{\circ} } ^2 { \theta^{\circ} } }\left\|\Lambda \theta_{j}\right\|_{L^2}^2
         \\
        & \lesssim
        C_j 2^{-js}
            \frac{1}{\varepsilon \nu^{\circ}}
        \left\|
             \nabla \widetilde{\u}^{\varepsilon}
        \right\|_{B^{\frac{N}{2}}}
        \left(
            \left\|
            \left(
                d, \widetilde{\theta}^{\varepsilon}
            \right)
            \right\|_{B^{\frac{N}{2}-1,\frac{N}{2}}_{\varepsilon \nu^{\circ}}}
            +\left\| \widetilde{\rho}^{\varepsilon} \right\|_{B^{\frac{N}{2}+1}}
        \right).
    \end{aligned}
\end{equation}
Note also that the equation satisfied by $\Omega$ is, up to non-linear terms, a heat equation,
and thus standard heat equation estimate yields that
\begin{equation}\label{heat}
    \begin{aligned}
        \frac{1}{2} \left\| \Omega_{j} \right\| + \frac{\mu^{\circ}}{\rho^{\circ}} 2^{2j}  \left\| \Omega_{j} \right\|
        \lesssim  \big\| \widetilde{G}_{j} \big\|.
    \end{aligned}
\end{equation}
\subsection{Conclusion}
Combining the estimate \eqref{high}, \eqref{low}, and \eqref{heat} we get
\begin{equation}\label{homogeneous_energy_estimates}
    \begin{aligned}
        \left\|( \tilde{\rho}^{\varepsilon},\  \tilde{\u}^{\varepsilon},\  \tilde{\theta}^{\varepsilon} )\right\|_{\widetilde{L}_{t}^{\infty}(\dot{\mathcal{E}}_{\varepsilon \nu^{\circ}}^{\frac{N}{2}})}
        + \left\| (\tilde{\rho}^{\varepsilon},\  \tilde{\u}^{\varepsilon},\  \tilde{\theta}^{\varepsilon} )\right\|_{\widetilde{L}_{t}^{1}(\dot{\mathcal{D}}_{\varepsilon \nu^{\circ}}^{\frac{N}{2}})}
        \lesssim  \left\| (F,\  G,\ \widetilde{G},\  H) \right\|_{\widetilde{L}^{1}(\dot{\mathcal{E}}_{\varepsilon \nu^{\circ}}^{\frac{N}{2}})} \\
         +\left\|( \tilde{\rho}_{\text{in}},\  \tilde{\u}_{\text{in}},\  \tilde{\theta}_{\text{in}} )\right\|_{(\dot{\mathcal{E}}_{\varepsilon \nu^{\circ}}^{\frac{N}{2}})}
        +\int^t_0 \left\| \nabla \tilde{\u}^{\varepsilon} (\tau)\right\|_{B^{\frac{N}{2}}}  \left\| \tilde{\rho}^{\varepsilon},\  \tilde{\u}^{\varepsilon},\  \tilde{\theta}^{\varepsilon} \right\|_{\widetilde{L}_{\tau}^{\infty}(\dot{\mathcal{E}}_{\varepsilon \nu^{\circ}}^{\frac{N}{2}})} d\tau.
    \end{aligned}
\end{equation}

We have the following estimates for non-linear terms $F$, $G$, $\widetilde{G}$, $H$.
\begin{lemma}\label{RightHandSide_estimates}
\begin{equation}
    \begin{aligned}
        \left\| (F,\ G,\ \widetilde{G},\ H) \right\|_{\widetilde{L}^{1}(\dot{\mathcal{E}}_{\varepsilon \nu^{\circ}}^{\frac{N}{2}})}
        \leqslant
        C(\mathcal{Y}^{\varepsilon})\left\| (\tilde{\rho}^{\varepsilon},\ \tilde{\u}^{\varepsilon},\ \tilde{\theta}^{\varepsilon}) \right\|^{2}_{E_{\varepsilon \nu^{\circ}}^{\frac{N}{2}}}.
    \end{aligned}
\end{equation}
\end{lemma}
\begin{proof}

For $F$, using \eqref{norm_equivalence}, Lemma \ref{product_estimate} and Lemma \ref{embedding} we have
\begin{equation}
    \begin{aligned}
        \left\| F  \right\|_ {{L}^{1}(B_{\varepsilon \nu^{\circ}}^{\frac{N}{2}-1,\frac{N}{2}+1})} =
        \left\| \tilde{\rho}^{\varepsilon} \text{ div } \tilde{\u}^{\varepsilon} \right\|_ {{L}^{1}(B_{\varepsilon \nu^{\circ}}^{\frac{N}{2}-1,\frac{N}{2}+1})}
        &\leqslant \left\| \tilde{\rho}^{\varepsilon}\right\|_{\widetilde{L}^{\infty}(B_{\varepsilon \nu^{\circ}}^{\frac{N}{2}-1,\frac{N}{2}+1})}
        \left\| \tilde{\u}^{\varepsilon} \right\|_{L^{1}(\dot{B}^{\frac{N}{2}+1})}. \\
    \end{aligned}
\end{equation}
For $G$, by definition
\begin{equation}
    \begin{aligned}
        \frac{\nu^{\circ}}{\rho^{\circ}} \Lambda^{2} d + \Lambda^{-1} \operatorname{div} \left( \nabla \!\cdot\! \frac{S^{\varepsilon}}{\rho^{\varepsilon}}  \right)
        &=\frac{\mu^{\circ}}{\rho^{\circ}} \Lambda^2 d + \Lambda^{-1} \mathrm{div} \left[\frac{1}{\rho^{\varepsilon}}\nabla \!\cdot\! ({\mu^{\varepsilon}} \nabla \widetilde{\u}^{\varepsilon}) \right]
        \\ &
        +\left(\frac{2}{N} \mu^{\circ }+ \lambda^{\circ} \right)\Lambda^2 d  +
        \Lambda^{-1} \mathrm{div} \frac{1}{\rho^{\varepsilon}} \left[ \nabla \!\cdot\! \left(\mu^{\varepsilon} \frac{2}{N} \nabla \!\cdot\! \widetilde{\u}^{\varepsilon}\mathrm{I}\right)\right]\\
        &+ \Lambda^{-1} \mathrm{div} \frac{1}{\rho^{\varepsilon}} \nabla \left(\lambda^{\varepsilon} \nabla \!\cdot\! \widetilde{\u}^{\varepsilon} \right)+
        +\frac{\mu^{\circ}}{\rho^{\circ}} \Lambda^2 d + \Lambda^{-1}  \mathrm{div} \left[\frac{1}{\rho^{\varepsilon}}\nabla \!\cdot\! ({\mu^{\varepsilon}} {(\nabla \widetilde{\u}^{\varepsilon})}^{\mathrm{T}}) \right].
    \end{aligned}
\end{equation}
Note that
\begin{equation}
    \begin{aligned}
        \Lambda^{-1} \mathrm{div} \nabla \!\cdot\! (\nabla \!\cdot\! \widetilde{\u}^{\varepsilon}\mathrm{I}) = - \Lambda^{2} d =  \Lambda^{-1} \mathrm{div} \nabla \!\cdot\! {(\nabla\widetilde{\u}^{\varepsilon})}^{\mathrm{T}},
    \end{aligned}
\end{equation}
then we have
\begin{equation}
    \begin{aligned}
    &\begin{aligned}
        \frac{\mu^{\circ}}{\rho^{\circ}} \Lambda^2 d + \Lambda^{-1} \mathrm{div} &\left[\frac{1}{\rho^{\varepsilon}}\nabla \!\cdot\! ({\mu^{\varepsilon}} \nabla \widetilde{\u}^{\varepsilon}) \right]
        \\ &=
        \Lambda^{-1} \mathrm{div} \frac{1}{\rho^{\varepsilon}}\nabla \!\cdot\! \left[ \left( \mu^{\varepsilon} - {\mu^{\circ}}\right) \nabla \widetilde{\u}^{\varepsilon}\right]
        + \mu^{\circ}\Lambda^{-1} \mathrm{div}  \left(\frac{1}{\rho^{\varepsilon}} - \frac{1}{\rho^{\circ}} \right)\Delta \widetilde{\u}^{\varepsilon},
    \end{aligned}\\
    &\begin{aligned}
        \frac{\mu^{\circ}}{\rho^{\circ}} \Lambda^2 d + \Lambda^{-1} & \mathrm{div} \left[\frac{1}{\rho^{\varepsilon}}\nabla \!\cdot\! ({\mu^{\varepsilon}} {(\nabla \widetilde{\u}^{\varepsilon})}^{\mathrm{T}}) \right] \\
        &=
        \Lambda^{-1} \mathrm{div} \frac{1}{\rho^{\varepsilon}} \nabla \!\cdot\! \left[ \left( \mu^{\varepsilon} -\mu^{\circ}\right) {(\nabla \widetilde{\u}^{\varepsilon})}^{\mathrm{T}}\right]
        + \mu^{\circ} \Lambda^{-1} \mathrm{div} \left[\left(\frac{1}{\rho^{\varepsilon}} - \frac{1}{\rho^{\circ}}\right)\nabla \!\cdot\! {( \nabla \widetilde{\u}^{\varepsilon})}^{\mathrm{T}} \right],
    \end{aligned}\\
    &\begin{aligned}
        \left(\frac{2}{N} \mu^{\circ }+ \lambda^{\circ} \right)& \Lambda^2 d  +
        \Lambda^{-1} \mathrm{div} \frac{1}{\rho^{\varepsilon}} \left[ \nabla \!\cdot\! \left(\mu^{\varepsilon} \frac{2}{N} \nabla \!\cdot\! \widetilde{\u}^{\varepsilon}\mathrm{I}\right)\right]
        + \Lambda^{-1} \mathrm{div} \frac{1}{\rho^{\varepsilon}} \nabla \left(\lambda^{\varepsilon} \nabla \!\cdot\! \widetilde{\u}^{\varepsilon} \right)
        \\ &=
        \frac{2}{N}\Lambda^{-1} \mathrm{div} \frac{1}{\rho^{\varepsilon}} \nabla \!\cdot\! \left[ \left(\mu^{\varepsilon} - \mu^{\circ}\right) \nabla \!\cdot\! \widetilde{\u}^{\varepsilon} \mathrm{I} \right]
        +\mu^{\circ}\frac{2}{N}\Lambda^{-1} \mathrm{div} \left[\left(\frac{1}{\rho^{\varepsilon}}-\frac{1}{\rho^{\circ}}\right) \nabla \!\cdot\! \left(\nabla \!\cdot\! \widetilde{\u}^{\varepsilon} \mathrm{I} \right) \right]
        \\ &+ \Lambda^{-1} \mathrm{div} \frac{1}{\rho^{\varepsilon}}\nabla\left[\left(\lambda^{\varepsilon} - \lambda^{\circ}\right)\nabla \!\cdot\! \widetilde{\u}^{\varepsilon} \right]
        + \lambda^{\circ}\Lambda^{-1} \mathrm{div} \left[\left(\frac{1}{\rho^{\varepsilon}}-\frac{1}{\rho^{\circ}}\right) \nabla(\nabla \!\cdot\! \widetilde{\u}^{\varepsilon})\right].
    \end{aligned}
\end{aligned}
\end{equation}
Using Lemma \ref{product_estimate}, Lemma \ref{composition_estimates} and \eqref{delta_estimates} we have
\begin{equation}
    \begin{aligned}
        \left\|\frac{\nu^{\circ}}{\rho^{\circ}} \Lambda^{2} d + \Lambda^{-1} \operatorname{div} \left( \nabla \!\cdot\! \frac{S^{\varepsilon}}{\rho^{\varepsilon}}  \right) \right\|_ {L^{1}(B_{\varepsilon \nu^{\circ}}^{\frac{N}{2}-1,\frac{N}{2}})}
        &\leqslant
        C(\mathcal{Y}^{\varepsilon})\left\| \left(\varepsilon \tilde{\rho}^{\varepsilon},\ \varepsilon \tilde{\theta}^{\varepsilon}\right)\right\|_{\widetilde{L}^{\infty}(B^{\frac{N}{2}}) } \left\| \widetilde{\u}^{\varepsilon} \right\|_{L^{1}(\dot{B}_{\varepsilon \nu^{\circ}}^{\frac{N}{2}+1,\frac{N}{2}+2})}\\
        &\leqslant
        C(\mathcal{Y}^{\varepsilon})\left\| \left( \tilde{\rho}^{\varepsilon},\  \tilde{\theta}^{\varepsilon}\right)  \right\|_{\widetilde{L}^{\infty}(B_{\varepsilon \nu^{\circ}}^{\frac{N}{2}-1,\frac{N}{2}})}
        \left\| \widetilde{\u}^{\varepsilon} \right\|_{L^{1}(\dot{B}_{\varepsilon \nu^{\circ}}^{\frac{N}{2}+1,\frac{N}{2}+2})},\\
    \end{aligned}
\end{equation}
and similarly we have
\begin{equation}
    \begin{aligned}
        \left\|\frac{1}{\varepsilon} \left(\frac{p_{\rho}^{\varepsilon}}{\rho^{\varepsilon}} - \frac{p^{\circ}_{\rho}}{\rho^{\circ}} \right)\nabla \tilde{\rho}\right\|_ {L^{1}(\dot{B}_{\varepsilon \nu^{\circ}}^{\frac{N}{2}-1,\frac{N}{2}})}
        &\leqslant C(\mathcal{Y}^{\varepsilon}) \frac{1}{\varepsilon}\left\| \left( \varepsilon\tilde{\rho}^{\varepsilon},\ \varepsilon\tilde{\theta}^{\varepsilon}\right)  \right\|_{\widetilde{L}^{2}(\dot{B}^{\frac{N}{2}})}
        \left\| \tilde{\rho}^\varepsilon \right\|_{\widetilde{L}^{2}(\dot{B}_{\varepsilon \nu^{\circ}}^{\frac{N}{2},\frac{N}{2}+1})}\\
        &\leqslant  C(\mathcal{Y}^{\varepsilon}) \left\| \left( \tilde{\rho}^{\varepsilon},\ \tilde{\theta}^{\varepsilon}\right)  \right\|_{\widetilde{L}^{2}(\dot{B}_{\varepsilon \nu^{\circ}}^{\frac{N}{2},\frac{N}{2}+1})}
        \left\| \tilde{\rho}^\varepsilon \right\|_{\widetilde{L}^{2}(\dot{B}_{\varepsilon \nu^{\circ}}^{\frac{N}{2},\frac{N}{2}+1})},\\
        \left\|\frac{1}{\varepsilon} \left(\frac{p_{\theta}^{\varepsilon}}{\rho^{\varepsilon}} - \frac{p^{\circ}_{\theta}}{\rho^{\circ}} \right)\nabla \tilde{\theta}\right\|_ {L^{1}(\dot{B}_{\varepsilon \nu^{\circ}}^{\frac{N}{2}-1,\frac{N}{2}})}
        &\leqslant C(\mathcal{Y}^{\varepsilon}) \frac{1}{\varepsilon}\left\| \left( \varepsilon\tilde{\rho}^{\varepsilon},\ \varepsilon\tilde{\theta}^{\varepsilon}\right)  \right\|_{\widetilde{L}^{2}(\dot{B}^{\frac{N}{2}})}
        \left\| \tilde{\theta}^\varepsilon \right\|_{\widetilde{L}^{2}(\dot{B}_{\varepsilon \nu^{\circ}}^{\frac{N}{2},\frac{N}{2}+1})}\\
        &\leqslant  C(\mathcal{Y}^{\varepsilon}) \left\| \left( \tilde{\rho}^{\varepsilon},\  \tilde{\theta}^{\varepsilon}\right)  \right\|_{\widetilde{L}^{2}(\dot{B}_{\varepsilon \nu^{\circ}}^{\frac{N}{2},\frac{N}{2}+1})}
        \left\| \tilde{\theta}^\varepsilon \right\|_{\widetilde{L}^{2}(\dot{B}_{\varepsilon \nu^{\circ}}^{\frac{N}{2},\frac{N}{2}+1})},\\
        \left\| \tilde{\u}^{\varepsilon} \!\cdot\! \nabla d - \Lambda^{-1} \operatorname{div} ( \tilde{\u}^{\varepsilon} \!\cdot\! \nabla  \tilde{\u}^{\varepsilon})\right\|_{L^{1}(\dot{B}_{\varepsilon \nu^{\circ}}^{\frac{N}{2}-1,\frac{N}{2}})}
        &\leqslant
        C\left\| \tilde{\u}^{\varepsilon} \right\|_{\widetilde{L}^{\infty}(\dot{B}_{\varepsilon \nu^{\circ}}^{\frac{N}{2}-1, \frac{N}{2}})}
        \left\| \tilde{\u}^{\varepsilon} \right\|_{\widetilde{L}^{1}(\dot{B}_{\varepsilon \nu^{\circ}}^{\frac{N}{2}+1, \frac{N}{2}+2})}.\\
    \end{aligned}
\end{equation}
For $H$, note that $ {e_{\theta} \rho^{\varepsilon}} $ is bounded away from $0$ under the assumption $\mathcal{Y}^{\varepsilon}$ is small and $e_{\theta}^{\circ} > 0$,
and note that
\begin{equation}
    \begin{aligned}
        -\frac{\kappa^{\circ}}{e^{\circ}_{\theta} \rho^{\circ}} \Delta  \tilde{\theta}^{\varepsilon}
        + \frac{1}{e_{\theta}^{\varepsilon} \rho^{\varepsilon}} \nabla( \kappa^{\varepsilon} \nabla  \tilde{\theta}^{\varepsilon})
        = \frac{1}{e_{\theta}^{\varepsilon} \rho^{\varepsilon}} \nabla \left[\left( \kappa^{\varepsilon} - \kappa^{\circ} \right)\nabla  \tilde{\theta}^{\varepsilon}\right]
        + \kappa^{\circ} \left( \frac{1}{e_{\theta}^{\varepsilon} \rho^{\varepsilon}} - \frac{1}{e_{\theta}^{\circ} \rho^{\circ}} \right) \Delta \widetilde{\theta}^{\varepsilon}.
    \end{aligned}
\end{equation}
As in estimate of $G$, Lemma \ref{product_estimate}, Lemma \ref{composition_estimates} and \eqref{delta_estimates} leads to
\begin{equation}
    \begin{aligned}
        \left\| -\frac{\kappa^{\circ}}{e^{\circ}_{\theta} \rho^{\circ}} \Delta  \tilde{\theta}^{\varepsilon} + \frac{1}{e_{\theta}^{\varepsilon} \rho^{\varepsilon}} \nabla( \kappa^{\varepsilon} \nabla  \tilde{\theta}^{\varepsilon}) \right\|_{L^{1}(\dot{B}_{\varepsilon \nu^{\circ}}^{\frac{N}{2}-1,\frac{N}{2}})}
        &\leqslant
        C(\mathcal{Y}^{\varepsilon}) \left\| \left( \tilde{\rho}^{\varepsilon},\  \tilde{\theta}^{\varepsilon}\right)  \right\|_{\widetilde{L}^{\infty}(\dot{B}_{\varepsilon \nu^{\circ}}^{\frac{N}{2}-1,\frac{N}{2}})}
        \left\| \tilde{\theta}^\varepsilon \right\|_{\widetilde{L}^{1}(\dot{B}_{\varepsilon \nu^{\circ}}^{\frac{N}{2}+1,\frac{N}{2}+2})},\\
        \left\|  \frac{1}{\varepsilon} \left( \frac{\theta^{\circ} p^{\circ}_{\theta}}{e^{\circ}_{\theta} \rho^{\circ}} - \frac{\theta^{\varepsilon} p_{\theta}^{\varepsilon}}{e_{\theta}^{\varepsilon} \rho^{\varepsilon}}\right)\Lambda d\right\|_{L^{1}(\dot{B}_{\varepsilon \nu^{\circ}}^{\frac{N}{2}-1,\frac{N}{2}})}
        &\leqslant
        C(\mathcal{Y}^{\varepsilon}) \left\| \left( \tilde{\rho}^{\varepsilon},\  \tilde{\theta}^{\varepsilon}\right)  \right\|_{\widetilde{L}^{2}(\dot{B}_{\varepsilon \nu^{\circ}}^{\frac{N}{2},\frac{N}{2}+1})}
        \left\| d \right\|_{\widetilde{L}^{2}(\dot{B}_{\varepsilon \nu^{\circ}}^{\frac{N}{2},\frac{N}{2}+1})},\\
        \left\| \frac{\varepsilon}{e_{\theta}^{\varepsilon} \rho^{\varepsilon}} S^{\varepsilon}: \nabla  \tilde{\u}^{\varepsilon} \right\|_{\widetilde{L}^{1}(\dot{B}_{\varepsilon \nu^{\circ}}^{\frac{N}{2}-1,\frac{N}{2}})}
        &\leqslant
        \varepsilon C(\mathcal{Y}^{\varepsilon}) \left\| \tilde{\u}^{\varepsilon} \right\|^{2}_{\widetilde{L}^{2}(\dot{B}_{\varepsilon \nu^{\circ}}^{\frac{N}{2}, \frac{N}{2}+1})},
    \end{aligned}
\end{equation}
the estimate for $\widetilde{G}$ is just as $G$.
Now using Proposition \ref{interpolation} we have
\begin{equation}
    \begin{aligned}
        \left\| (\widetilde{\rho}^{\varepsilon},\ \widetilde{\u}^{\varepsilon}, \ \widetilde{\theta}^{\varepsilon}) \right\|_{\widetilde{L}^{2}(\dot{B}^{\frac{N}{2},\frac{N}{2}+1}_{\varepsilon \nu^{\circ}})}^2
        &\lesssim
        \left\| (\widetilde{\rho}^{\varepsilon},\ \widetilde{\u}^{\varepsilon}, \ \widetilde{\theta}^{\varepsilon}) \right\|_{\widetilde{L}^{\infty}(\mathcal{E}^{\frac{N}{2}})}
        \left\| (\widetilde{\rho}^{\varepsilon},\ \widetilde{\u}^{\varepsilon}, \ \widetilde{\theta}^{\varepsilon}) \right\|_{\widetilde{L}^{1}(\dot{\mathcal{D}}_{\varepsilon \nu^{\circ}}^{\frac{N}{2}})}
        \\ &\lesssim
        \| (\widetilde{\rho}^{\varepsilon},\ \widetilde{\u}^{\varepsilon}, \ \widetilde{\theta}^{\varepsilon}) \|_{E^{\frac{N}{2}}_{\varepsilon \nu^{\circ}}}^2,
    \end{aligned}
\end{equation}
which completes the proof.
\end{proof}
Now we deal with the mean part of $\widetilde{\u}^{\varepsilon}$. In fact $\widebar{{\widetilde{\rho}}}^{\varepsilon}$ is conserved and therefore  $\widebar{{\widetilde{\rho}}}^{\varepsilon} =  \widebar{{\widetilde{\rho}}}_{\text{in}}$, we only need to consider
$\widebar{{\widetilde{\u}}}^{\varepsilon}$ and  $\widebar{{\widetilde{\theta}}}^{\varepsilon}$.
\begin{equation}
    \begin{aligned}
        \partial_t  \widebar{{\widetilde{\u}}}^{\varepsilon}
        &= \frac{1}{|\mathbb{T}^N_a|}\int_{\mathbb{T}^N_a}
         - \widetilde{\u}^{\varepsilon} \cdot \nabla \widetilde{\u}^{\varepsilon}
         - \frac{1}{\varepsilon} \frac{p^{\varepsilon}_{\rho}}{\rho^{\varepsilon}}\nabla \widetilde{\rho}^{\varepsilon}
         - \frac{1}{\varepsilon}  \frac{p^{\varepsilon}_{\theta}}{\theta^{\varepsilon}}\nabla \widetilde{\theta}^{\varepsilon}
         + \frac{1}{\rho^{\varepsilon}} \nabla S^{\varepsilon} dx,\\
        \partial_t  \widebar{{\widetilde{\theta}}}^{\varepsilon}
        &= \frac{1}{|\mathbb{T}^N_a|}\int_{\mathbb{T}^N_a}
         \frac{1}{\varepsilon}\left( \frac{\theta^{\circ} p^{\circ}_{\theta}}{e^{\circ}_{\theta} \rho^{\circ}} - \frac{\theta p_{\theta}}{e_{\theta} \rho}\right) \nabla \widetilde{\u}^{\varepsilon}
          - \widetilde{\u}^{\varepsilon} \cdot \nabla \widetilde{\theta}^{\varepsilon}
        +\frac{1}{e_{\theta}\rho} \nabla(\kappa^{\varepsilon} \nabla \widetilde{\theta}^{\varepsilon}) 
         + \frac{\varepsilon}{e_{\theta}\rho} \nabla \widetilde{\u}^{\varepsilon} : S^{\varepsilon} dx
    \end{aligned}
\end{equation}
Using Lemma \ref{composition_estimates} we have
\begin{equation}
   \begin{aligned}
   \| \frac{1}{|\mathbb{T}^N_a|}\int_{\mathbb{T}^N_a} \widetilde{\u}^{\varepsilon} \cdot \nabla \widetilde{\theta}^{\varepsilon} dx \|_{L^1{dt}}
   &\lesssim
   \|\underline{\widetilde{\u}^{\varepsilon}} \cdot \nabla \underline{\widetilde{\theta}^{\varepsilon}}\|_{L^1(dtdx)} \\
   &\lesssim \| \widetilde{\u}^{\varepsilon} \|_{\widetilde{L}^{\infty}(\dot{B}^{\frac{N-1}{2}})} \| \widetilde{\theta}^{\varepsilon}\|_{L^{1}(\dot{B}^{\frac{N+1}{2}})}
   \lesssim \| (\widetilde{\rho}^{\varepsilon},\ \widetilde{\u}^{\varepsilon}, \ \widetilde{\theta}^{\varepsilon}) \|_{E^{\frac{N}{2}}_{\varepsilon \nu^{\circ}}}^2,\\
   \end{aligned}
\end{equation}
\begin{equation}
    \begin{aligned}
   \|\frac{1}{|\mathbb{T}^N_a|}\int_{\mathbb{T}^N_a}
         \frac{1}{\varepsilon}\left( \frac{\theta^{\circ} p^{\circ}_{\theta}}{e^{\circ}_{\theta} \rho^{\circ}} - \frac{\theta p_{\theta}}{e_{\theta} \rho}\right) \nabla \widetilde{\u}^{\varepsilon}\|_{L^1(dt)}
    &\lesssim \| \frac{1}{\varepsilon}{\left( \frac{\theta^{\circ} p^{\circ}_{\theta}}{e^{\circ}_{\theta} \rho^{\circ}} - \frac{\theta p_{\theta}}{e_{\theta} \rho}\right)} \|_{L^{\infty}(\dot{B}^{\frac{N-1}{2}})} \| \nabla \widetilde{\u}^{\varepsilon}\|_{L^1(B^{\frac{N}{2}})} \\
    &\lesssim \| (\widetilde{\rho}^{\varepsilon},\ \widetilde{\u}^{\varepsilon}, \ \widetilde{\theta}^{\varepsilon}) \|_{E^{\frac{N}{2}}_{\varepsilon \nu^{\circ}}}^2,\\
   \end{aligned}
\end{equation}

\begin{equation}
    \begin{aligned}
    \left\|\frac{1}{|\mathbb{T}^N_a|}\int_{\mathbb{T}^N_a} \frac{1}{e_{\theta}\rho} \nabla(\kappa^{\varepsilon} \nabla \widetilde{\theta}^{\varepsilon}) -\frac{\kappa^{\circ}}{e^{\circ}_{\theta}\rho^{\circ}} \Delta \widetilde{\theta}^{\varepsilon} dx\right\|_{L^1(dt)}
    &\lesssim \left\|\frac{1}{|\mathbb{T}^N_a|}\int_{\mathbb{T}^N_a} \frac{1}{e_{\theta}\rho} [\nabla(\kappa^{\varepsilon}-\kappa^{\circ}) \nabla \widetilde{\theta}^{\varepsilon}] dx\right\|_{L^1(dt)} \\
    &+ \left\|\frac{1}{|\mathbb{T}^N_a|}\int_{\mathbb{T}^N_a}\kappa^{\circ} ( \frac{1}{e_{\theta}\rho}- \frac{1}{e_{\theta}^{\circ}\rho^{\circ}} ) \Delta\widetilde{\theta}^{\varepsilon}) dx\right\|_{L^1(dt)} \\
    &\lesssim \varepsilon C(\mathcal{Y}^{\varepsilon})\| (\rho^{\varepsilon}, \theta^{\varepsilon})\|_{L^{\infty}(\dot{B}^{\frac{N}{2}-1})} \| \theta^{\varepsilon}\|_{L^1(\dot{B}^{\frac{N}{2}+1})} ,
    \end{aligned}
\end{equation}
\begin{equation}
    \begin{aligned}
    \|\frac{1}{|\mathbb{T}^N_a|}\int_{\mathbb{T}^N_a}  \frac{\varepsilon}{e_{\theta}\rho} \nabla \widetilde{\u}^{\varepsilon} : S^{\varepsilon} dx \|_{L^{1}(dt)}
    &\lesssim \varepsilon C(\mathcal{Y})  \left\| \|\nabla \widetilde{\u}^{\varepsilon} \|_{L^2(dx)} \| S^{\varepsilon} \|_{L^2(dx)} \right\|_{L^1(dt)}\\
    &\lesssim \varepsilon C(\mathcal{Y})^2 \left\| \|\nabla \widetilde{\u}^{\varepsilon} \|_{L^2(dx)} \| \nabla \widetilde{\u}^{\varepsilon} \|_{L^2(dx)} \right\|_{L^1(dt)}\\
    &\lesssim \varepsilon C(\mathcal{Y})^2   \| \widetilde{\u}^{\varepsilon} \|^2_{L^2(\dot{B}^{\frac{N}{2}})},
    \end{aligned}
\end{equation}
and therefore we have
\begin{equation}\label{theta_mean_estimate}
    \| \widebar{\widetilde{\theta}}^{\varepsilon} \|_{L^{\infty}(dt)}
    \lesssim |\widebar{\widetilde{\theta}}_{\text{in}} |
    + C(\mathcal{Y}^{\varepsilon})
    \| (\widetilde{\rho}^{\varepsilon},\ \widetilde{\u}^{\varepsilon}, \ \widetilde{\theta}^{\varepsilon}) \|_{E^{\frac{N}{2}}_{\varepsilon \nu^{\circ}}}^2.
\end{equation}
The term $ \| \widebar{\widetilde{\u}}^{\varepsilon} \|_{L^{\infty}(dt)}  $ is estimated in a similar way, and we omit it here. We then have
\begin{equation}\label{mean_value_estimates}
    \| (\widebar{\widetilde{\rho}}^{\varepsilon},\widebar{\widetilde{\u}}^{\varepsilon},\widebar{\widetilde{\theta}}^{\varepsilon}) \|_{L^{\infty}(dt)}
    \lesssim |(\widebar{\rho}_{\text{in}},\widebar{\widetilde{\u}}_{\text{in}},\widebar{\widetilde{\theta}}_{\text{in}}) |
     + C(\mathcal{Y}^{\varepsilon})
    \| (\widetilde{\rho}^{\varepsilon},\ \widetilde{\u}^{\varepsilon}, \ \widetilde{\theta}^{\varepsilon}) \|_{E^{\frac{N}{2}}_{\varepsilon \nu^{\circ}}}^2.
\end{equation}
Combining \eqref{homogeneous_energy_estimates} and \eqref{mean_value_estimates} leads to
\begin{equation}\label{inhomogeneous_energy_estimates}
    \begin{aligned}
        \left\|( \tilde{\rho}^{\varepsilon},\  \tilde{\u}^{\varepsilon},\  \tilde{\theta}^{\varepsilon} )\right\|_{\widetilde{L}_{t}^{\infty}(\mathcal{E}^{\frac{N}{2}})}
        + \left\| (\tilde{\rho}^{\varepsilon},\  \tilde{\u}^{\varepsilon},\  \tilde{\theta}^{\varepsilon} )\right\|_{\widetilde{L}_{t}^{1}(\dot{\mathcal{D}}_{\varepsilon \nu^{\circ}}^{\frac{N}{2}})}
        \lesssim  \left\| (F,\  G,\ \widetilde{G},\  H) \right\|_{\widetilde{L}^{1}(\dot{\mathcal{E}}_{\varepsilon \nu^{\circ}}^{\frac{N}{2}})}
        +
        \left\|( \tilde{\rho}_{\text{in}},\  \tilde{\u}_{\text{in}},\  \tilde{\theta}_{\text{in}} )\right\|_{(\mathcal{E}^{\frac{N}{2}})}
        \\
             +C(\mathcal{Y}^{\varepsilon})
             \| (\widetilde{\rho}^{\varepsilon},\ \widetilde{\u}^{\varepsilon}, \ \widetilde{\theta}^{\varepsilon}) \|_{E^{\frac{N}{2}}_{\varepsilon \nu^{\circ}}}^2
        +\int^t_0 \left\| \nabla \tilde{\u}^{\varepsilon} (\tau)\right\|_{B^{\frac{N}{2}}}  \left\| (\tilde{\rho}^{\varepsilon},\  \tilde{\u}^{\varepsilon},\  \tilde{\theta}^{\varepsilon}) \right\|_{\widetilde{L}_{\tau}^{\infty}(\mathcal{E}^{\frac{N}{2}})} d\tau.
    \end{aligned}
\end{equation}
Lemma \ref{RightHandSide_estimates} and  simple application of Gronwall's inequality ensures that
\begin{equation}\label{LPH_estimates}
    \begin{aligned}
        &\left\|( \tilde{\rho}^{\varepsilon},\  \tilde{\u}^{\varepsilon},\  \tilde{\theta}^{\varepsilon} )\right\|_{\widetilde{L}_{t}^{\infty}({\mathcal{E}}_{\varepsilon \nu^{\circ}}^{\frac{N}{2}})}
        + \left\| (\tilde{\rho}^{\varepsilon},\  \tilde{\u}^{\varepsilon},\  \tilde{\theta}^{\varepsilon} )\right\|_{\widetilde{L}_{t}^{1}(\dot{\mathcal{D}}_{\varepsilon \nu^{\circ}}^{\frac{N}{2}})}
        \leqslant
        C\mathrm{e}^{C\left\|\nabla \tilde{\u}^{\varepsilon}\right\|_{\widetilde{L}_{t}^{1}(B^{\frac{N}{2}})}}\\
        &\times\left(
            \left\|( \tilde{\rho}^{\varepsilon},\  \tilde{\u}^{\varepsilon},\  \tilde{\theta}^{\varepsilon} )(0)\right\|_{{\mathcal{E}}_{\varepsilon \nu^{\circ}}^{\frac{N}{2}}}
         + + C(\mathcal{Y}^{\varepsilon})
         \| (\widetilde{\rho}^{\varepsilon},\ \widetilde{\u}^{\varepsilon}, \ \widetilde{\theta}^{\varepsilon}) \|_{E^{\frac{N}{2}}_{\varepsilon \nu^{\circ}}}^2
        \right),
    \end{aligned}
\end{equation}
where $C$ is a constant independent of $\varepsilon$.
We now use standard bootstrap argument to close the estimate.
\begin{proof}[Proof of Proposition \ref{A_priori_estimates}]

Given that
\begin{equation}
    \begin{aligned}
        \left\| \left(\tilde{\rho}_{in},\ \tilde{\u}_{in},\ \tilde{\theta}_{in}\right) \right\|_{\mathcal{E}_{\varepsilon \nu^{\circ}}^{{\frac{N}{2}}}}
        \leqslant \alpha_{0},
    \end{aligned}
\end{equation}
with $\alpha_{0}$ to be determined later. For a fixed $M > 1$ define following set
\begin{equation}
    \begin{aligned}
        \mathcal{H}=\left\{  t \in [0,+\infty)|\left\| \left(\tilde{\rho}^{\varepsilon},\ \tilde{\u}^{\varepsilon},\ \tilde{\theta}^{\varepsilon} \right) \right\|_{\widetilde{L}^{\infty}_t(\mathcal{E}_{\varepsilon \nu^{\circ}}^{{\frac{N}{2}}}) \cap \widetilde{L}^{1}_t(\dot{\mathcal{D}}_{\varepsilon \nu^{\circ}}^{\frac{N}{2}})}
        \leqslant M \alpha_{0} \right\}.
    \end{aligned}
\end{equation}
It's clearly a closed set by continuity, and certainly not empty since the initial condition ensures $t=0$ in this set.

We now show that it's also an open set.
Again, by continuity, $\forall t \in \mathcal{H}$, there exists $\delta_0 > 0$ such that for every $t^{*}$ satisfying $|t-t^{*}| < \delta_0$ we have
\begin{equation}\label{bootstrap_conclusion}
\left\| \left(\tilde{\rho}^{\varepsilon},\ \tilde{\u}^{\varepsilon},\ \tilde{\theta}^{\varepsilon}\right) \right\|_{\widetilde{L}^{\infty}_{t^*}(\mathcal{E}_{\varepsilon \nu^{\circ}}^{{\frac{N}{2}}}) \cap \widetilde{L}^{1}_{t^*}(\dot{\mathcal{D}}_{\varepsilon \nu^{\circ}}^{\frac{N}{2}})}
 \leqslant 2M \alpha_{0}.
\end{equation}
To avoid $\left\|\rho^{\varepsilon}\right\|_{L^{\infty}(dtdx)} = 0 $ and $\left\|\theta^{\varepsilon}\right\|_{L^{\infty}(dtdx)} = 0 $,
we use
\begin{equation}
    \begin{aligned}
        \left\|\rho^{\varepsilon}\right\|_{L^{\infty}(dtdx)} = \left\|\rho^{\circ}+\varepsilon\tilde{\rho}^{\varepsilon}\right\|_{L^{\infty}(dtdx)} \geqslant \| \rho^{\circ} - \mathcal{C}_B \left\| \varepsilon \tilde{\rho}^{\varepsilon} \right\|_{B^{\frac{N}{2}}}\|_{L^{\infty}(dt)},
    \end{aligned}
\end{equation}
where $\mathcal{C}_B$ is the continuity modulus of $B^{\frac{N}{2}} \hookrightarrow L^{\infty}$.
Note that we have
\begin{equation}
    \left\| \varepsilon \tilde{\rho}^{\varepsilon} \right\|_{B^{\frac{N}{2}}} \leqslant \frac{1}{\nu^{\circ}}  \left\|  \tilde{\rho}^{\varepsilon} \right\|_{B^{\frac{N}{2}-1, \frac{N}{2}}_{\varepsilon \nu^{\circ}}} \leqslant \frac{1}{\nu^{\circ}}2M \alpha_{0},
\end{equation}
and similarly
\begin{equation}
    \begin{aligned}
        \left\|\theta^{\varepsilon} \right\|_{L^{\infty}} \geqslant \theta^{\circ} - \frac{\mathcal{C}_B}{\nu^{\circ}}2M \alpha_0.
    \end{aligned}
\end{equation}
Therefore, we require that
\begin{equation}
    \alpha_0 \leqslant \min \left\{\frac{1}{4} \frac{\nu^{\circ} \rho^{\circ}}{\mathcal{C}_B M},\ \frac{1}{4} \frac{\nu^{\circ} \theta^{\circ}}{\mathcal{C}_B M} \right\},
\end{equation}
which implies that
\begin{equation} \label{Y_estimates}
    \max \left\{ \frac{1}{2} \rho^{\circ},\ \frac{1}{2} \theta^{\circ},
    \frac{2}{3} \frac{1}{\rho^{\circ}},\ \frac{2}{3} \frac{1}{\theta^{\circ}} \right\}
    \leqslant
    \left\| \mathcal{Y}^{\varepsilon} \right\|_{L^{\infty}([0,+\infty) \times \mathbb{R}^{n})}
    \leqslant \min \left\{ \frac{3}{2} \rho^{\circ},\ \frac{3}{2} \theta^{\circ},
    2 \frac{1}{\rho^{\circ}},\ 2 \frac{1}{\theta^{\circ}} \right\} .
\end{equation}
Now using Lemma \ref{RightHandSide_estimates}, and combining \eqref{LPH_estimates}, \eqref{bootstrap_conclusion} and \eqref{Y_estimates}, we get
\begin{equation}
    \begin{aligned}
        \left\| (\tilde{\rho}^{\varepsilon},\ \tilde{\u}^{\varepsilon},\ \tilde{\theta}^{\varepsilon}) \right\|_{\widetilde{L}^{\infty}_{t^*}(\mathcal{E}_{\varepsilon \nu^{\circ}}^{{\frac{N}{2}}}) \cap {L}^{1}_{t^*}(\dot{\mathcal{D}}_{\varepsilon \nu^{\circ}}^{\frac{N}{2}})}
        \leqslant C \mathrm{e}^{C 2M\alpha_{0}}
        \left(
        \alpha_{0}
        +M^{2} \alpha_{0}^{2}
        \right).
    \end{aligned}
\end{equation}
Choosing $M=4C$ and choosing $\alpha_{0}$ small such that
\begin{equation}
    M^2 \alpha_{0} \leqslant 1, \ \mathrm{e}^{2CM\alpha_{0}}\leqslant 2,\  \text{ and } \alpha_{0} \leqslant  \frac{\nu^{\circ}}{\mathcal{C}_B 16C},
\end{equation}
and therefore
\begin{equation}
    \begin{aligned}
        \left\| (\tilde{\rho}^{\varepsilon},\ \tilde{\u}^{\varepsilon},\ \tilde{\theta}^{\varepsilon}) \right\|_{\widetilde{L}^{\infty}_{t^*}(\mathcal{E}_{\varepsilon \nu^{\circ}}^{{\frac{N}{2}}}) \cap \widetilde{L}^{1}_{t^*}(\dot{\mathcal{D}}_{\varepsilon \nu^{\circ}}^{\frac{N}{2}})}
        \leqslant M \alpha^{\circ},
    \end{aligned}
\end{equation}
which implies that $\mathcal{H}$ is an open set and completes the proof.
\end{proof}
\section{About the limit system}\label{sec_limit}
In this section, we prove that the limit system \eqref{LS}
\begin{equation}
    \left\{\begin{aligned}
        &\partial_t V + \mathcal{Q}_{2r}(\mathcal{U},V) + \mathcal{Q}_{3r}(V,V)-\bar{\mu} \Delta V =0,\\
        &V|_{t=0} = \mathcal{P}^{\perp}\widetilde{U}_{in},
    \end{aligned}\right.
\end{equation}
is globally well-posed in $F^{\frac{N}{2}}$
as far as \eqref{INSF} has a global solution in
$F^{\frac{N}{2}}$.
Since \eqref{INSF} rewrites
\begin{equation}\nonumber
    \left\{
    \begin{aligned}
        &\partial_t\omega + \Pi(\omega\!\cdot\!\nabla \omega) = \frac{\mu^{\circ}}{\rho^{\circ}}\Delta\omega\\
        &\partial_t\vartheta + \omega\!\cdot\!\nabla \vartheta = \frac{\kappa^{\circ}p_{\rho}^{\circ}}{(c^{\circ})^2\rho^{\circ}e_{\theta}^{\circ}}\Delta\vartheta,\\
        &\omega|_{t=0} = \omega_{in},\ \ \vartheta|_{t=0} = \vartheta_{in},
    \end{aligned}
    \right.
\end{equation}
which is just the incompressible Naiver-Stokes equation
coupled with a passive scalar transport-diffusion equation.
Theory on well-posedness of these types of equations in
Besov spaces has been thorough, see for example \cite{Chemin_Danchin_book}.
We state the following proposition without proof
\begin{proposition}
    Suppose $N\geqslant 2$. Let $\mathcal{U}_{in}=
    (-\frac{p_{\theta}^{\circ}}{p_{\rho}^{\circ}}\vartheta_{in}, \omega_{in}, \vartheta_{in})^{T}\in B^{\frac{N}{2}-1}$,
    where $\omega_{in}$ is a divergence free vector field.
    Then there exist a constant $C$ depending only on $N$ such that,
    under the condition
    \begin{equation}
        \|\mathcal{U}_{in}\|_{B^{\frac{N}{2}-1}} \leqslant C\mu^{\prime},
    \end{equation}
    the incompressible Navier-Stokes-Fourier system
    \begin{equation}\nonumber
        \left\{\begin{aligned}
            &\partial_t \mathcal{U} + \mathcal{P}\mathcal{Q}(\mathcal{U},\mathcal{U}) - \mathcal{P}\mathcal{D}\mathcal{U} = 0,\\
            &\mathcal{U}|_{t=0} = \mathcal{U}_{in}.
        \end{aligned}\right.
    \end{equation}
    has a unique solution $\mathcal{U}
    = (-\frac{p_{\theta}^{\circ}}{p_{\rho}^{\circ}}\vartheta, \omega, \vartheta)^{T}$
    in $F^{\frac{N}{2}}$.
    Moreover, it satisfies the estimate
    \begin{equation}
        \|\mathcal{U}\|_{F^{\frac{N}{2}}}\leqslant
        C\|\mathcal{U}_{in}\|_{B^{\frac{N}{2}-1}}
    \end{equation}
    where $\mu^{\prime} = \text{min}\big\{\frac{\mu^{\circ}}{\rho^{\circ}},
    \frac{\kappa^{\circ}p^{\circ}_{\rho}}{(c^{\circ})^2\rho^{\circ}e^{\circ}_{\theta}}\big\}$.
\end{proposition}

The main theorem in this section reads
\begin{theorem}[A priori estimates for the limit system]
    \label{apriori_limit_system}
    Let $s\geqslant1$, $V_{in}\in B^{s-1}\cap \text{Null}(\mathcal{A})^{\perp}$
    and $\mathcal{U}\in F^{s}$ is a
    fixed solution to the incompressible Navier-Stokes-Fourier system.
    Then the limit system \eqref{LS} has a unique solution $V\in F^{s}$
    which remains in $\text{Null} (\mathcal{A})^{\perp}$ for all time.
    Moreover, for $t\in(0,\infty]$, the solution $V$ satisfies the energy estimates
    \begin{equation}\label{l2_estimate_for_limit_system}
        \|V(t)\|^2_{L^2} + 2\bar{\mu}\int_{0}^{t}\|\nabla V(\tau)\|^2_{L^2}\mathrm{d}\tau \leqslant\|V_{in}\|^2_{L^2},
    \end{equation}
    and
    \begin{equation}
        \|V(t)\|_{B^{s-1}} + c\bar{\mu}\|V\|_{L^1_t (B^{s+1})} \leqslant \mathrm{e}^{\frac{C}{\bar{\mu}}\|V_{in}\|^2_{L^2}}\|V_{in}\|_{B^{s-1}}.
    \end{equation}
    Where the constants $c$ and $C$ depend only on $s$ and $N$.
\end{theorem}

Owing to the non-local nature of the problem,
proving well-posedness using energy arguments
may require a more in-depth inspection of
the structure of the resonance operators.
First, an alternative expression for the three-wave resonant term
stemming from previous computations reads
\begin{equation}
    \begin{aligned}
        \mathcal{Q}_{3r}(V,V) =
        \sum_{(\gamma, m)}
        \sum_{\substack{ k+ l= m\\ \operatorname{sg}( k)| k| + \operatorname{sg}( l)| l| = \operatorname{sg}( m)| m|}}
        \chi^{\gamma}_{m}
        V^{\gamma}_{ k}V^{\gamma}_{ l}H^{\gamma}_{ m},
    \end{aligned}
\end{equation}
where
\begin{equation}
    \chi^{\gamma}_{m} = \mathrm{i} \widetilde{C}_1\gamma \operatorname{sg}(m)|m|.
\end{equation}
Here, for simplicity's sake, we use the notation $\widetilde{C}_1$
to represent some generic constant independent of the summation.
Following \cite{Masmoudi_Poincare}, we introduce the set of
“prime vector representatives" for frequencies:
\begin{equation}
    \mathscr{P} \stackrel{\mathrm{def}}{=}
    \{p\in\widetilde{\mathbb{Z}}^{N}\;|\;
    a_1p_1\wedge\dots\wedge a_Np_N = 1\}.
\end{equation}
This is equivalent to saying that,
for any $p\in\mathscr{P}$, there does not
exist any $(n,q)\in\mathbb{N}\times\widetilde{\mathbb{Z}}^{N}$
such that $p=nq$ and $n\geqslant2$.
For a vector $V(x)=\sum_{\alpha,k}V^{\alpha}_{k}H^{\alpha}_{k}(x)\in\mathrm{Null}(\mathcal{A})^{\perp}$,
we define
\begin{equation}
    V^p (x) \stackrel{\mathrm{def}}{=}
    \sum_{n\in\mathbb{Z}^{\ast}}
    V^{\operatorname{sg}(p)}_{np}H^{\operatorname{sg}(p)}_{np}(x).
\end{equation}
It's easily seen that for all $s\in\mathbb{R}$
\begin{equation}\label{energy_equivalence_1}
    V(x) = \sum_{p\in\mathscr{P}}W^p(x),\quad
    \|V\|_{H^s} = \sum_{p\in\mathscr{P}}\|V^p\|_{H^s},\quad\mathrm{and}\quad
    \|V\|_{B^s} = \sum_{p\in\mathscr{P}}\|V^p\|_{B^s}.
\end{equation}
And also, for $p,q\in\mathscr{P}$, $p\neq q$,
we have
\begin{equation}
    \mathcal{Q}_{3r}(V^p,V^q)=0.
\end{equation}
Thus, $\mathcal{Q}_{3r}(V,V)$
is formed by resonances between modes
oscillating in the same direction, namely
\begin{equation}\label{primely_reduced_Q3r}
    \mathcal{Q}_{3r}(V,V) = \sum_{p\in\mathcal{P}}
    \mathcal{Q}_{3r}(V^p,V^p).
\end{equation}

Next, we move to see how $\mathcal{Q}_{2r}(\mathcal{U},V)$
can be reduced using the prime vector representation.
The spectral decomposition rewrites
\begin{equation}
    \mathcal{Q}_{2r}(\mathcal{U},V)
        =\sum_{(\gamma, m)}\sum_{\substack{ k+ l=m\\\alpha \operatorname{sg}( k)=\gamma \operatorname{sg}( m)\\| k|=| m|}}
        V^{\alpha}_{ k}H^{\gamma}_{ m}\sigma^{\gamma}_{klm},
\end{equation}
where
\begin{equation}
    \sigma^{\gamma}_{klm}
    =\widetilde{C}_{2}\frac{ k\!\cdot\! m}{|k|| m|}(\hat{\omega}_{l}\!\cdot\! k)
    + \widetilde{C}_{3} \gamma \operatorname{sg}( m)| m|\hat{\vartheta}_{l}
    + \widetilde{C}_{4} \gamma \operatorname{sg}( m)\frac{ k\!\cdot\! m}{| m|}\hat{\vartheta}_{l}.
\end{equation}
Also, $\widetilde{C}_{2}$, $\widetilde{C}_{3}$, and $\widetilde{C}_{4}$
denote some generic constants independent of the summation.
For $p,q\in\mathscr{P}$, we define
\begin{equation}
    C_{p}(\mathcal{U},V^q)\stackrel{\mathrm{def}}{=}
    \sum_{\substack{ k+ l=m\\k=nq, m=n^\prime p\\m|p|=k|q|}}
    \sigma^{\operatorname{sg}(p)}_{klm}V^{\operatorname{sg}(q)}_{ k}H^{\operatorname{sg}(p)}_{ m},
\end{equation}
which is the contribution of $V^q$ on $V^p$.

At this point, we introduce the following functions defined for $z\in\mathbb{T}$
\begin{equation}
    \begin{aligned}
        {\vv}^p(z)&\stackrel{\mathrm{def}}{=}
        \mathrm{i}\widetilde{C}_1\sum_{n\in\mathbb{Z}^{\ast}}
        V^{\operatorname{sg}(p)}_{np}\mathrm{e}^{\mathrm{i}nz},\\
        c_p(\mathcal{U},{\vv}^q)(z)&\stackrel{\mathrm{def}}{=}
        \sum_{\substack{ k+ l=m\\k=nq, m=n^\prime p\\m|p|=k|q|}}
        \sigma^{\operatorname{sg}(p)}_{klm}\mathrm{e}^{\mathrm{i}n^\prime z}.
    \end{aligned}
\end{equation}
As pointed out by \cite{Masmoudi_Poincare}, we have for all $s\in\mathbb{R}$
\begin{equation}\label{energy_equivalence_2}
    \begin{aligned}
        2\pi\mathrm{i}\widetilde{C}_1
        \|V^p\|_{H^s(\mathbb{T}^{N}_{a})}
        = |p|^{s}\|{\vv}^p\|_{H^s(\mathbb{T})},\\
        2\pi\mathrm{i}\widetilde{C}_1
        \|V^p\|_{B^s(\mathbb{T}^{N}_{a})}
        = |p|^{s}\|{\vv}^p\|_{B^s(\mathbb{T})},
    \end{aligned}
\end{equation}
which equate the energys of the original problem
to those of the one dwelling in $\mathbb{T}$.
More importantly,
\begin{proposition}[\cite{Masmoudi_Poincare}]
    Solving the system \eqref{LS} for $V=\sum_{p\in\mathscr{P}}$
    $V^{p}$ is equivalent to solving the system of
    infinite coupled viscous Burgers (ICVB) equations
    \begin{equation}\label{ICVB}\tag{ICVB}
        \forall p\in\mathscr{P},\quad
        \partial_t {\vv}^p + |p|\partial_z ({\vv}^p)^2
        - \bar{\mu}|p|^2\partial_z^2 {\vv}^p
        + \sum_{q\in\mathscr{P}}c_p(\mathcal{U},{\vv}^q) = 0.
    \end{equation}
\end{proposition}
As a consequence, once the a priori bounds are proved,
approximation arguments (see \cite{Masmoudi_Poincare}) leads to the existence of \eqref{ICVB}.

Now we turn to the proof of the a priori bound.
The following lemmas may find useful in later estimates,
we state the results and postpone their proof to the end of this section.
The first lemma tells us that the resonance operators
have some cancellation properties which render them “harmless"
in certain energy estimates
\begin{lemma}\label{cancellation_identities}
    For $s\geqslant 0$, we have
    \begin{equation}
        \begin{aligned}
            &\langle\mathcal{Q}_{2r}(\mathcal{U},V),V\rangle_{H^s} = 0,\\
            &\langle\Delta_j \mathcal{Q}_{2r}(\mathcal{U},V),\Delta_j V\rangle_{\mathbb{H}} = 0,\\
            &\langle\mathcal{Q}_{3r}(V,V),V\rangle_{\mathbb{H}} = 0.
        \end{aligned}
    \end{equation}
\end{lemma}
The second lemma concerns some useful estimates of $\mathcal{Q}_{3r}$
\begin{lemma}\label{estimates_for_Q^{3r}}
    \begin{equation}
        \begin{aligned}
            &\langle\mathcal{Q}_{3r}(V,W),W\rangle_{\mathbb{H}}
            \lesssim \|W\|_{L^2}\|W\|_{B^{\frac{1}{2}}}\|V\|_{H^1},\\
            &\|\mathcal{Q}_{3r}(V,V)\|_{B^{s-1}}\lesssim \|V\|_{B^{\frac{1}{2}}}\|V\|_{B^s}.
        \end{aligned}
    \end{equation}
    \begin{proof}
        The proof of these two estimates relies on the observation \eqref{primely_reduced_Q3r},
        whereby estimates relating to $\mathcal{Q}_{3r}(V,W)$ is reducable
        to those of functions defined on the $1$-dimensional torus $\mathbb{T}$.
        We refer to Lemma 9.1 of \cite{Danchin_AJM} for the details.
    \end{proof}
\end{lemma}

\begin{proof}[Proof of Theorem \ref{apriori_limit_system}]
    We only prove a priori bounds here.
    In light of Lemma \ref{cancellation_identities}, a standard $L^2$ estimate yields
    \begin{equation}
        \frac{1}{2}\frac{\mathrm{d}}{\mathrm{d}t}\|V(t)\|^2_{L^2} + \bar{\mu}\|\nabla V(t)\|^2_{L^2} = 0.
    \end{equation}
    Then a time integration proves the first estimate.
    To obtain estimates in Besov spaces, we apply $\Delta_j$ to (LS)
    and then perform $L^2$ estimate.
    Also by Lemma \ref{cancellation_identities} we get
    \begin{equation}
        \frac{1}{2}\frac{\mathrm{d}}{\mathrm{d}t}\|\Delta_j V\|^2_{L^2} + \bar{\mu}\|\nabla \Delta_j V\|^2_{L^2}
        + \langle\Delta_j \mathcal{Q}_{3r}(V,V),\Delta_j V\rangle_{\mathbb{H}} = 0.
    \end{equation}
    A time integration and a multiplication by $2^{j(s-1)}$ gives
    \begin{equation}
        \begin{aligned}
            2^{j(s-1)}&\|\Delta_j V(t)\|_{L^2} + c\bar{\mu}2^{j(s+1)}\|\Delta_j V\|_{L^1_t(L^2)}\\
            \leqslant& 2^{j(s-1)}\|\Delta_j V_{in}\|_{L^2} + 2^{j(s-1)}\|\Delta_j \mathcal{Q}_{3r}(V,V)\|_{L^1_t(L^2)}.
        \end{aligned}
    \end{equation}
    Summing up in $j$ yields
    \begin{equation}
        \|V(t)\|_{B^{s-1}} + c\bar{\mu}\|V\|_{L^1_t(B^{s+1})}\leqslant\|V_{in}\|_{B^{s-1}} + \int_{0}^{t}\|\mathcal{Q}_{3r}(V,V)\|_{B^{s-1}}\mathrm{d}\tau.
    \end{equation}
    By Lemma \ref{estimates_for_Q^{3r}}, we have
    \begin{equation}
        \|\mathcal{Q}_{3r}(V,V)\|_{B^{s-1}}\lesssim \|V\|_{H^1}\|V\|_{B^s},
    \end{equation}
    whence an interpolation and Young's inequality conclude
    \begin{equation}
        \begin{aligned}
            \|V(t)\|&_{B^{s-1}} + c\bar{\mu}\|V\|_{L^1_t(B^{s+1})}\\
            \leqslant&\|V_{in}\|_{B^{s-1}} +
            \frac{c\bar{\mu}}{2}\|V\|_{L^1_t(B^{s+1})} + \int_{0}^{t} \frac{C}{\bar{\mu}}\|V(\tau)\|^2_{H^1}\|V(\tau)\|_{B^{s-1}}\mathrm{d}\tau.
        \end{aligned}
    \end{equation}
    Gronwall lemma then leads to
    \begin{equation}
        \|V(t)\|_{B^{s-1}} + c\bar{\mu}\|V\|_{L^1_t(B^{s+1})}
        \leqslant \mathrm{e}^{\frac{C}{\bar{\mu}}\int_{0}^{t}\|V(\tau)\|^2_{H^1}\mathrm{d}\tau}\|V_{in}\|_{B^{s-1}},
    \end{equation}
    which combined with the $L^2$-estimate
    \eqref{l2_estimate_for_limit_system} gives the Besov energy estimate.
    Existence of solutions in $V\in F^{s}_{T}$ then lies on standard arguments.

    We prove uniqueness of the solutions in $C([0,T];L^2)\cap L^2(0,T;H^1)$. Remark that this property is not known for the incompressible Naiver-Stokes equation in dimension $N\geqslant3$.
    Let $V_1$, $V_2$ be two solutions of \eqref{LS} in $C([0,T];L^2)\cap L^2(0,T;H^1)$. Define $V^\prime \stackrel{\mathrm{def}}{=} V_2-V_1$, then
    \begin{equation}
        \partial_t V^\prime  - \bar{\mu}\Delta V^\prime  + \mathcal{Q}_{2r}(\mathcal{U},V^\prime )
        = - \mathcal{Q}_{3r}(V_1+V_2,V^\prime ).
    \end{equation}
    It satisfies the estimate
    \begin{equation}
        \frac{1}{2}\frac{\mathrm{d}}{\mathrm{d}t}\|V^\prime \|^2_{L^2} + \bar{\mu}\|\nabla V^\prime \|^2_{L^2}
        = -\langle\mathcal{Q}_{3r}(V_1+V_2,V^\prime ),V^\prime \rangle_{\mathbb{H}}.
    \end{equation}
    According to Lemma \ref{estimates_for_Q^{3r}} and the embedding $H^1\hookrightarrow B^{\frac{1}{2}}$, we have
    \begin{equation}
        |\langle\mathcal{Q}_{3r}(V_1+V_2,V^\prime ),V^\prime \rangle_{\mathbb{H}}|
        \leqslant\frac{C}{\bar{\mu}}(\|V_1+V_2\|_{H^1})^2\|V^\prime \|^2_{L^2} + \bar{\mu}\|\nabla V^\prime \|^2_{L^2}.
    \end{equation}
    Then a use of Gronwall lemma entails $V^\prime \equiv 0$.
\end{proof}

\begin{proof}[Proof of Lemma \ref{cancellation_identities}]
    Performing an $H^s$ inner product between $\mathcal{Q}_{2r}(\mathcal{U},V)$
    and $V$, we have
    \begin{equation}
        \langle\mathcal{Q}_{2r}(\mathcal{U},V),V\rangle_{H^s}
            =\sum_{(\gamma, m)}\sum_{\substack{ k+ l=m\\\alpha \operatorname{sg}( k)=\gamma \operatorname{sg}( m)\\| k|=| m|}}
            |m|^{2s}V^{\alpha}_{ k}\overline{V}^{\gamma}_{ m}\sigma^{\gamma}_{klm},
    \end{equation}
    Noticing the symmetry that $H^{\alpha}_{ k}=\bar{H}^{\alpha}_{- k}$, $V^{\alpha}_{ k}=\overline{V}^{\alpha}_{- k}$,
    we exchange $\alpha$ and $\gamma$, and change $ k$ to $- m$ and change $ m$ to $- k$.
    Then $\sigma^{\gamma}_{klm}$ changes to $-\sigma^{\gamma}_{klm}$,
    because the relation $\hat{\omega}_{ l}\cdot m=\hat{\omega}_{ l}\cdot k$
    is guaranteed by the divergence-free property of $\omega$.
    Note that $ l= m- k$ is invariant under this index change.
    Therefore,
    \begin{equation}
        \begin{aligned}
            \langle\mathcal{Q}_{2r}(\mathcal{U},V),V\rangle_{H^s}
            &= -\sum_{(\gamma, m)}\sum_{\substack{ k+ l=m\\\alpha \operatorname{sg}( k)=\gamma \operatorname{sg}( m)\\| k|=| m|}}
            |m|^{2s}V^{\gamma}_{-m}\overline{V}^{\alpha}_{-k}\sigma^{\gamma}_{klm},\\
            &= -\sum_{(\gamma, m)}\sum_{\substack{ k+ l=m\\\alpha \operatorname{sg}( k)=\gamma \operatorname{sg}( m)\\| k|=| m|}}
            |m|^{2s}V^{\alpha}_{ k}\overline{V}^{\gamma}_{ m}\sigma^{\gamma}_{klm},
        \end{aligned}
    \end{equation}
    whence
    \begin{equation}
        \langle\mathcal{Q}_{2r}(\mathcal{U},V),V\rangle_{H^s}
        = -\langle\mathcal{Q}_{2r}(\mathcal{U},V),V\rangle_{H^s}.
    \end{equation}
    Therefore $\langle\mathcal{Q}_{2r}(\mathcal{U},V),V\rangle_{H^s}=0$.
    Similar arguments give
    \begin{equation}
        \langle\Delta_j \mathcal{Q}_{2r}(\mathcal{U},V),\Delta_j V\rangle_{\mathbb{H}} = 0.
    \end{equation}

    Next, a $\mathbb{H}$-inner product between
    $\mathcal{Q}_{3r}(V,V)$ and $V$ gives
    \begin{equation}
        \begin{aligned}
            \langle\mathcal{Q}_{3r}(V,V),V\rangle_{\mathbb{H}} =
            \sum_{(\gamma, m)}
            \sum_{\substack{ k+ l= m\\ \operatorname{sg}( k)| k| + \operatorname{sg}( l)| l| = \operatorname{sg}( m)| m|}}
            \chi^{\gamma}_{m}
            V^{\gamma}_{ k}V^{\gamma}_{ l}\overline{V}^{\gamma}_{ m},
        \end{aligned}
    \end{equation}
    Notice that if $(V^{\gamma}_{ k},V^{\gamma}_{ l},V^{\gamma}_{ m})$ is a resonant triplet, then it's also true for
    $(V^{\gamma}_{ k},V^{\gamma}_{- m},V^{\gamma}_{- l})$ and $(V^{\gamma}_{- m},V^{\gamma}_{ l},V^{\gamma}_{- k})$.
    Hence
    \begin{equation}
        \begin{aligned}
            \langle\mathcal{Q}_{3r}(V,V),V\rangle_{\mathbb{H}} =& \frac{1}{3}\times\sum_{(\gamma, m)}
            \sum_{\substack{ k+ l= m\\ \operatorname{sg}( k)| k| + \operatorname{sg}( l)| l| =\\ \operatorname{sg}( m)| m|}}
            (\chi^{\gamma}_{m}+\chi^{\gamma}_{-k}+\chi^{\gamma}_{-l})
            V^{\gamma}_{ k}V^{\gamma}_{ l}\overline{V}^{\gamma}_{ m}
            =& 0.
        \end{aligned}
    \end{equation}

\end{proof}
\section{Appendix}
In this section we give a brief introduction to Bony's paradifferential calculus(see work by \cite{Bony}), and we prove some estimates that we need to use.
The product of u and v is defined by
\begin{equation}
    T_u v \stackrel{\text { def }}{=} \sum_{q \in \mathbb{Z}} S_{q-1} u \Delta_q v.
\end{equation}
Then we have following formal decomposition:
\begin{equation}
    \begin{aligned}
        &uv = T_u v + T_v u + R(u,v) \text{ with }\\
        &R(u,v) \stackrel{\text { def }}{=} \sum_{q \in \mathbb{Z} } \Delta_q u \sum_{|q^{\prime} -q| \leqslant 1 } \Delta_{q^{\prime}} V .
    \end{aligned}
\end{equation}
We have following continuity result for the paraproduct. We refer to \cite{Danchin_ARMA} for the proof of it.
\begin{lemma}For all $s_1,\  s_2,\  t_1,\  t_2$ such that $s_1 \leqslant N/2$ and $s_2 \leqslant N/2$, the following estimate holds:
    \begin{equation}
        \left\|T_u v\right\|_{\widetilde{B}^{s_1+t_1-\frac{N}{2}, s_2+t_2-\frac{N}{2}}} \lesssim\|\u\|_{\widetilde{B}^{s_1, s_2}}\|v\|_{\widetilde{B}^{t_1, t_2}} .\\
    \end{equation}
 If  $\min \left(s_1+t_1, s_2+t_2\right)>0$, then
 \begin{equation}
    \|R(u, v)\|_{\widetilde{B}^{s_1+t_1-\frac{N}{2}, s_2+t_2-\frac{N}{2}}} \lesssim\|\u\|_{\widetilde{B}^{s_1, s_2}}\|v\|_{\widetilde{B}^{t_1, t_2}} .
 \end{equation}
If $u \in L^{\infty}$,
\begin{equation}
    \left\|T_u v\right\|_{\widetilde{B}^{t_1, t_2}} \lesssim\|\u\|_{L^{\infty}}\|v\|_{\widetilde{B}^{t_1, t_2}},
\end{equation}
and, if $\min \left(t_1, t_2\right)>0$, then
\begin{equation}
    \|R(u, v)\|_{\widetilde{B}^{t_1, t_2}} \lesssim\|\u\|_{L^{\infty}}\|v\|_{\widetilde{B}^{t_1, t_2}} .
\end{equation}
\end{lemma}
Remark that Lemma \ref{product_estimate} can be easily deduced from this Lemma.
We now state the estimates for the  convection terms which are needed in the a priori estimates.

\begin{lemma} \label{convection_estimates}
 Let $F$ be an homogeneous smooth function of degree $m \in \mathbb{R}$. Suppose that $-N / 2<s_1, t_1, s_2, t_2 \leqslant 1+N / 2$. The following two estimates hold:
\begin{equation}
    \begin{aligned}
    &\begin{aligned}
        & \left<\left(F(D) \Delta_q(v \cdot \nabla a) ,\  F(D) \Delta_q a\right)\right> \\
        & \quad \lesssim c_q 2^{-q\left(\widetilde{\psi}^{s_1, s_2}(q)-m\right)}
        \left(\widetilde{\phi}^{s_1,s_2}(q) \right)^{-1}
        \|v\|_{B^{\frac{N}{2}+1}}\|a\|_{\widetilde{B}^{s_1, s_2}}\left\|F(D) \Delta_q a\right\|_{L^2}, \\
    \end{aligned}\\
    &\begin{aligned}
        & \left<\left(F(D) \Delta_q(v \cdot \nabla a) ,\  \Delta_q b\right)+\left(\Delta_q(v \cdot \nabla b) \mid F(D) \Delta_q a\right)\right> \\
        & \lesssim c_q\|v\|_{B^{\frac{N}{2}+1}}
         \times\left(2^{-q \widetilde{\psi}^{t_1, t_2}(q)}
         \left(\widetilde{\phi}^{t_1,t_2}(q)\right)^{-1}
         \left\|F(D) \Delta_q a\right\|_{L^2}\|b\|_{\widetilde{B}^{t_1, t_2}}\right. \\
        & \left.+2^{-q\left(\widetilde{\psi}^{s_1, s_2}(q)-m\right)}
        \left(\widetilde{\phi}^{s_1,s_2}(q)\right)^{-1}
        \|a\|_{\widetilde{B}^{s_1, s_2}}\left\|\Delta_q b\right\|_{L^2}\right),
    \end{aligned}
\end{aligned}
\end{equation}
where the function $\widetilde{\psi}^{s, t}$ and $\widetilde{\phi}^{s,t}$ are  defined as in Lemma \ref{Localising} and ${c_q} \in \ell_1$
\end{lemma}
Proof of this Lemma can be found in \cite{Danchin_ARMA} with minor modification.

\begin{lemma} \label{Qe_estimate}
    The following estimates hold for $0 < s < \frac{N}{2}+1$:
    \begin{equation} \label{Qs_estimate}
        \left\| \mathcal{Q} (W^{1},W^{2}) \right\|_{H^{s-2}}
        \lesssim
        \min\left(\left\| W^{1} \right\|_{H^{s-1}}  \left\| W^{2} \right\|_{H^{\frac{N}{2}}} \cap L^{\infty} ,\ \left\| W^{1} \right\|_{H^{s-1}} \left\| W^{2}  \right\|_{H^{\frac{N}{2}} \cap L^{\infty}} \right),
    \end{equation}
    \begin{equation} \label{e2r_estimate}
        \left\| \mathcal{Q}^{\varepsilon}_{2r}(\mathcal{U} , B) \right\|_{H^{s-2}}
        \lesssim
        \min\left(\left\|\mathcal{U}\right\|_{H^{s-1}} \left\|B\right\|_{B^{\frac{N}{2}}} ,\ \left\|B\right\|_{H^{s-1}} \left\| \mathcal{U} \right\|_{B^{\frac{N}{2}}} \right),
    \end{equation}
    \begin{equation} \label{e3r_estimate}
        \left\| \mathcal{Q}^{\varepsilon}_{3r}(A , B) \right\|_{H^{s-2}}
        \lesssim
        \min\left(\left\|A\right\|_{H^{s-1}} \left\|B\right\|_{B^{\frac{N}{2}}} ,\ \left\|B\right\|_{H^{s-1}} \left\| A \right\|_{B^{\frac{N}{2}}} \right),
    \end{equation}
    where $\mathcal{U} \in \text{ker} \mathcal{A}$ with $\widebar{\mathcal{U}}=0$, and $A,\ B \in \text{ker}^{\perp} \mathcal{A}$.
\end{lemma}
\begin{proof}
    For \eqref{Qs_estimate} use the expression of $\mathcal{Q}$
    For simplicity of description, we denote $\rho^{1}$ the first component of $W^{1}$, $\u^{1}$ the $(2,\cdots,d+1)$ components, and $\theta^{1}$ the $d+2$ component of $W^{1}$.
    Similarly, $(\rho^{2},\ \u^{2},\ \theta^{2}) $ the components of $W^{2}$.
    For the first component of $\mathcal{Q}\left(W^{1} , W^{2}\right) $
    \begin{equation*}
        \begin{aligned}
            \left\| \operatorname{div}( \rho^1 \u^2 ) + \operatorname{div}(\rho^2 \u^1) \right\|_{H^{s-2}}
            \lesssim \left\|\rho^1 \u^2\right\|_{H^{s-1}} + \left\|\rho^2 \u^1\right\|_{H^{s-1}},
        \end{aligned}
    \end{equation*}
    then Lemma \ref{product_estimate} yields that
    \begin{equation*}
        \begin{aligned}
            \left\|\rho^1 \u^2\right\|_{H^{s-1}} + \left\|\rho^2 \u^1\right\|_{H^{s-1}}
            &\lesssim \min\left(\left\|\rho^1\right\|_{H^{s-1}}\left\| \u^2\right\|_{H^{\frac{N}{2}} \cap L^{\infty}} ,\left\|\rho^1\right\|_{H^{\frac{N}{2}} \cap L^{\infty}}\left\| \u^2\right\|_{H^{s-1}} \right) \\
            &+\min\left(\left\|\rho^2\right\|_{H^{s-1}}\left\| \u^1\right\|_{H^{\frac{N}{2}}\cap L^{\infty}} ,\left\|\rho^2\right\|_{H^{\frac{N}{2}} \cap L^{\infty}}\left\| \u^1\right\|_{H^{s-1}} \right) .
        \end{aligned}
    \end{equation*}
    For the terms in the second and third components of $\mathcal{Q}\left(W^{1},\ W^{2}\right) $, we have
    \begin{equation*}
        \begin{aligned}
        \left\| \u^1 \cdot \nabla \u^2 \right\|_{H^{s-2}} +  \left\| \u^2 \cdot \nabla \u^1 \right\|_{H^{s-2}}
        &\lesssim
        \min\left(\left\|\nabla \u^2\right\|_{H^{s-2}}\left\| \u^1\right\|_{H^{\frac{N}{2}} \cap L^{\infty}} ,\left\|\nabla \u^2\right\|_{H^{\frac{N}{2}-1}}\left\| \u^1\right\|_{H^{s-1}} \right)\\
        &+\min\left(\left\|\u^2\right\|_{H^{s-1}}\left\|\nabla \u^1\right\|_{H^{\frac{N}{2}-1}} ,\left\|\u^2\right\|_{H^{\frac{N}{2}} \cap \infty}\left\|\nabla \u^1\right\|_{H^{s-2}} \right) .
    \end{aligned}
    \end{equation*}
    Same methods apply to $\u^1 \cdot \nabla \theta^2+\u^2 \cdot \nabla \theta^1$, $\rho^1 \nabla \theta^2 + \rho^2 \nabla \theta^1$, $\rho^1 \nabla \rho^2 + \rho^2 \nabla \rho^1$ and $\theta^1 \nabla \theta^2 +\theta^2 \nabla \theta^1$.
    Gathering together all estimates above we have
    \begin{equation*}
        \lVert \mathcal{Q}\left(W^{1} ,\ W^{2}\right) \rVert_{H^{s-2}}
        \lesssim
        \min\left(\left\|W^{1}\right\|_{H^{s-1}} \left\| W^{2} \right\|_{H^{\frac{N}{2}} \cap L^{\infty}} ,\ \left\| W^{2} \right\|_{H^{s-1}} \left\| W^{1} \right\|_{H^{\frac{N}{2}} \cap L^{\infty}} \right).
    \end{equation*}

    For \eqref{e2r_estimate}, by definition of $\mathcal{Q}^{\varepsilon}_{2r}$ 
    \begin{equation}
        \begin{aligned}
            \left\| \mathcal{Q}^{\varepsilon}_{2r}(\mathcal{U} , B) \right\|_{H^{s-2}}
            \lesssim \left\| \mathcal{Q}\left(\mathcal{U} ,\ \mathrm{e}^{-\frac{t}{\varepsilon}\mathcal{A}} B\right) \right\|_{H^{s-2}}.
        \end{aligned}
    \end{equation}
    We then use \eqref{Qs_estimate} to get
    \begin{equation*}
        \left\| \mathcal{Q}\left(\mathcal{U} ,\ \mathrm{e}^{-\frac{t}{\varepsilon}\mathcal{A}} B\right) \right\|_{H^{s-2}}
        \lesssim
        \min\left(\left\|\mathcal{U}\right\|_{H^{s-1}} \left\| \mathrm{e}^{-\frac{t}{\varepsilon}\mathcal{A}}B \right\|_{H^{\frac{N}{2}} \cap L^{\infty}} ,\ \left\| \mathrm{e}^{-\frac{t}{\varepsilon}\mathcal{A}} B \right\|_{H^{s-1}} \left\| \mathcal{U} \right\|_{H^{\frac{N}{2}} \cap L^{\infty}} \right).
    \end{equation*}
    Using Lemma \ref{embedding} and Lemma \ref{norm_conserve} 
    \begin{equation*}
        \left\| \mathrm{e}^{-\frac{t}{\varepsilon}\mathcal{A}} B \right\|_{H^{s-1}} = \left\| B \right\|_{H^{s-1}},\ \left\| \mathrm{e}^{-\frac{t}{\varepsilon}\mathcal{A}} B \right\|_{H^{\frac{N}{2}} \cap L^{\infty}} \leqslant \left\|\mathrm{e}^{-\frac{t}{\varepsilon}\mathcal{A}} B \right\|_{B^{\frac{N}{2}}}
        =\left\| B \right\|_{B^{\frac{N}{2}}},
    \end{equation*}
    we finally get \eqref{e2r_estimate}.

    The proof for \eqref{e3r_estimate} is quite similar noticing that
    \begin{equation*}
        \begin{aligned}
            \left\| \mathcal{Q}^{\varepsilon}_{3r}(A , B) \right\|_{H^{s-2}}
            \lesssim \left\| \mathcal{Q}\left(\mathrm{e}^{-\frac{t}{\varepsilon}\mathcal{A}} A ,\  \mathrm{e}^{-\frac{t}{\varepsilon}\mathcal{A}} B\right) \right\|_{H^{s-2}}.
        \end{aligned}
    \end{equation*}
\end{proof}
We shall also use following lemma for the estimate of nonlinear heat equation, which can be found in \cite{Danchin_AJM}.
\begin{lemma}\label{structure_lemma}
    Let $s\in\mathbb{R}$, $\sigma>0$, $0 < T \leq +\infty$ and $\mathscr{Q}$ be a time dependent nonlinear operator on $H^{s-1}$ valued in $H^{s-2}$ with estimate
    \begin{equation}
        \|\mathscr{Q}(V)(t)\|_{H^{s-2}}\leqslant\mathcal{K}(t)\|V\|_{H^{s-1}},
    \end{equation}
    where $\mathcal{K}\in L^2(0,T)$.
    Let $V\in \widetilde{L}^{\infty}(0,T;H^{s-1})\cap L^2(0,T;H^{s})$ be a solution of the system
    \begin{equation}
        \left\{
            \begin{array}{l}
                \partial_t V + \mathscr{Q}(V) - \sigma  \Delta V  = F,\\
                V|_{t=0}=V_{in},
            \end{array}
        \right.
    \end{equation}
    where $F\in L^2(0,T;H^{s-2})$ and $V_{in}\in H^{s-1}$.
    Then, there exists a positive constant $c$ such that for all $t\in[0,T]$,
    \begin{equation}
        \begin{aligned}
            \|V\|^2_{\widetilde{L}^\infty_{t}(\dot{H}^{s-1})}
            + c\sigma \|V\|^2_{L^2_t( \dot{H}^{s})}
            \leqslant & \mathrm{e}^{\frac{4}{c\sigma}\int_{0}^{t}\mathcal{K}^2(\tau)\mathrm{d}\tau}\\
            &\times\left(2\|V_{in}\|^2_{H^{s-1}}+\frac{16}{c\sigma}\|F\|^2_{L^2_t(\dot{H}^{s-2})}\right).
        \end{aligned}
    \end{equation}
\end{lemma}

\section*{Acknowledgment}
In the process of this work, N. Jiang discussed several times with Dr. Li Lai and Prof. Baoping Liu of Peking University on the estimates of small divisor problem. Dr. Li Lai introduced Lemma 5.3 to him. Prof. Baoping Liu had several insightful comments on the small divisor  estimates on the 3-wave resonance set. We deeply appreciate their help here. 

 \bibliographystyle{amsplain}
\bibliography{references.bib}
\end{document}